\documentclass[a4paper,english]{article}
\usepackage[T1]{fontenc}
\usepackage[latin9]{inputenc}
\usepackage{verbatim}
\usepackage{amsmath}
\usepackage{amsthm}
\usepackage{amssymb}
\usepackage{esint}

\makeatletter

\providecommand{\tabularnewline}{\\}

\theoremstyle{plain}
\newtheorem{thm}{\protect\theoremname}[section]
  \theoremstyle{definition}
  \newtheorem{example}[thm]{\protect\examplename}
  \theoremstyle{plain}
  \newtheorem{conjecture}[thm]{\protect\conjecturename}
  \theoremstyle{plain}
  \newtheorem{cor}[thm]{\protect\corollaryname}
  \theoremstyle{definition}
  \newtheorem{defn}[thm]{\protect\definitionname}
  \theoremstyle{remark}
  \newtheorem{notation}[thm]{\protect\notationname}
  \theoremstyle{plain}
  \newtheorem{lem}[thm]{\protect\lemmaname}
  \theoremstyle{remark}
  \newtheorem{rem}[thm]{\protect\remarkname}
  \theoremstyle{plain}
  \newtheorem{prop}[thm]{\protect\propositionname}
  \theoremstyle{remark}
  \newtheorem{claim}[thm]{\protect\claimname}
\newcommand{\lyxaddress}[1]{
\par {\raggedright #1
\vspace{1.4em}
\noindent\par}
}

\DeclareMathOperator{\id}{id}
\DeclareMathOperator{\diam}{diam}
\DeclareMathOperator{\supp}{supp}

\DeclareMathOperator{\sdim}{s-dim}
\DeclareMathOperator{\eigen}{eigen}
\DeclareMathOperator{\stab}{stab}

\DeclareMathOperator{\sat}{sat}
\DeclareMathOperator{\lip}{Lip}
\DeclareMathOperator{\vol}{vol}

\DeclareMathOperator{\bdim}{bdim}
\DeclareMathOperator{\pdim}{dim_P}

\DeclareMathOperator*{\spn}{span}

\DeclareMathOperator{\rank}{rank}
\DeclareMathOperator{\Gal}{Gal}

\newcommand{\hide}[1]{}
\newcommand{\conv}{\mbox{\LARGE{$.$}}}

\date{}

\makeatletter
\def\blfootnote{\xdef\@thefnmark{}\@footnotetext}
\makeatother

\makeatother

\usepackage{babel}
  \providecommand{\claimname}{Claim}
  \providecommand{\conjecturename}{Conjecture}
  \providecommand{\corollaryname}{Corollary}
  \providecommand{\definitionname}{Definition}
  \providecommand{\examplename}{Example}
  \providecommand{\lemmaname}{Lemma}
  \providecommand{\notationname}{Notation}
  \providecommand{\propositionname}{Proposition}
  \providecommand{\remarkname}{Remark}
\providecommand{\theoremname}{Theorem}

\begin{document}

\title{On self-similar sets with overlaps and inverse theorems for entropy
in $\mathbb{R}^{d}$}

\author{Michael Hochman}
\maketitle
\begin{abstract}
We\blfootnote{Supported by ERC grant 306494, partially supported by ISF grant 1409/11}
study self-similar sets and measures on $\mathbb{R}^{d}$. Assuming
that the defining iterated function system $\Phi$ does not preserve
a proper affine subspace, we show that one of the following holds:
(1) the dimension is equal to the trivial bound (the minimum of $d$
and the similarity dimension $s$); (2) for all large $n$ there are
$n$-fold compositions of maps from $\Phi$ which are super-exponentially
close in $n$; (3) there is a non-trivial linear subspace of $\mathbb{R}^{d}$
that is preserved by the linearization of $\Phi$ and whose translates
typically meet the set or measure in full dimension. In particular,
when the linearization of $\Phi$ acts irreducibly on $\mathbb{R}^{d}$,
either the dimension is equal to $\min\{s,d\}$ or there are super-exponentially
close $n$-fold compositions. We give a number of applications to
algebraic systems, parametrized systems, and to some classical examples.

The main ingredient in the proof is an inverse theorem for the entropy
growth of convolutions of measures on $\mathbb{R}^{d}$, and the growth
of entropy for the convolution of a measure on the orthogonal group
with a measure on $\mathbb{R}^{d}$. More generally, this part of
the paper applies to smooth actions of Lie groups on manifolds.
\end{abstract}
\tableofcontents{}

\section{\label{sec:Introduction}Introduction}

Self similar sets and measures are among the simplest fractal objects:
their defining property is that the whole is made up of finitely many
objects similar to it, i.e. identical to the whole except for scaling,
rotation and translation. When these smaller copies are sufficiently
separated from each other the small-scale structure is relatively
easy to understand and in particular the Hausdorff dimension can be
computed precisely in terms of the defining similitudes. Without separation,
however, things are significantly more complicated, and it is an open
problem to compute the dimension. Many special cases of this problem
have received attention, including the Erd\"{o}s problem on Bernoulli
convolutions, Furstenberg's projection problem for the 1-dimensional
Sierpinski gasket (now settled), the Keane-Smorodinsky the ${0,1,3}$-problem,
and ``fat'' Sierpinski gaskets (for more on these, see below). 

For self-similar sets and measures in $\mathbb{R}$ there is a longstanding
conjecture predicting that the dimension will be ``as large as possible'',
subject to the combinatorial constraints, unless there are exact overlaps,
i.e. unless some of the (iterated) small-scale copies of the original
coincide. In recent work \cite{Hochman2014} we introduced methods
from additive combinatorics to this problem and obtained a partial
result towards the conjecture, showing that if the dimension is ``too
small'' then there are super-exponentially close pairs of small-scale
copies. In particular, for some important classes of self-similar
sets, e.g. those defined by similarities with algebraic coefficients,
this resolves the conjecture. 

In the present paper we treat the general case of self-similar sets
and measures in $\mathbb{R}^{d}$. Easy examples show that in the
higher-dimensional setting the conjecture above is false as stated
(Example \ref{example}). The main new feature of the problem is that
the linear parts of the defining similarities may act reducibly on
$\mathbb{R}^{d}$, and ``excess dimension'' may accumulate on non-trivial
invariant subspaces and produce dimension loss. To correct this we
propose here a modified version of the conjecture that takes this
possibility into account (Conjecture \ref{conjecture}), and prove
a weak version of it (Theorem \ref{thm:main-individual-RD-entropy}),
analogous to the main result of \cite{Hochman2014}. We give various
applications, in particular we show that the modified conjecture holds
when the linear action is irreducible and the coefficients of the
similarities are algebraic.

As in the 1-dimensional case, a central ingredient in the proof is
an inverse theorem about the structure of probability measures on
$\mathbb{R}^{d}$ whose convolutions have essentially the same entropy
as the original (Theorem \ref{thm:inverse-thm-Rd}). In fact, what
we really need is a result of this type for the convolution of a measure
on $\mathbb{R}^{d}$ with a measure on the similarity group, or one
on the isometry group (Theorem \ref{thm:inverse-theorem-for-isometries}).
These results are of independent interest, and provide a versatile
tool for analyzing smooth images of product measures. We take the
opportunity to develop these methods here, in particular stating results
for convolutions in Lie groups and their actions (Theorems \ref{thm:inverse-theorem-for-isometries},
\ref{thm:generalized-inverse-theorem} and the subsequent corollaries).

\subsection{\label{sub:Self-similar-sets}Setup: Self-similar sets and measures}

Let $G$ denote the group of similarities of $\mathbb{R}^{d}$, namely
maps $x\mapsto rUx+a$ for $r\in(0,\infty)$, $a\in\mathbb{R}^{d}$
and $U$ a $d\times d$ orthogonal matrix; we denote the map simply
by $\varphi=rU+a$. In this paper an iterated function system means
a finite family $\Phi=\{\varphi_{i}\}_{i\in\Lambda}\subseteq G$ consisting
of contractions, so $\varphi_{i}=r_{i}U_{i}+a_{i}$ with $0<r_{i}<1$.
A self similar set is the attractor of such a system, defined as the
unique compact set $\emptyset\neq X\subseteq\mathbb{R}$ satisfying
\begin{equation}
X=\bigcup_{i\in\Lambda}\varphi_{i}X.\label{eq:self-similar-set}
\end{equation}
The self-similar measure determined by $\Phi$ and a positive probability
vector $(p_{i})_{i\in\Lambda}$ is the unique Borel probability measure
$\mu$ on $\mathbb{R}^{d}$ satisfying
\[
\mu=\sum_{i=1}^{k}p_{i}\cdot\varphi_{i}\mu.
\]
Here and throughout, $\varphi\mu=\mu\circ\varphi^{-1}$ denotes the
push-forward of $\mu$ by $\varphi$.

It is a classical problem to understand the small-scale structure
of self-similar sets and measures, and especially their dimension.
We shall write $\dim A$ for the Hausdorff dimension of $A$ and define
the dimension of a finite Borel measure $\theta$ by\footnote{This is the lower Hausdorff dimension. Many other notions of dimension
exist but since self-similar measures are exact dimensional \cite{FengHu09},
for them all the major ones coincide. }
\[
\dim\theta=\inf\{\dim E\,:\,\theta(E)>0\}.
\]

The textbook case of self-similar sets and measures occurs when the
images $\varphi_{i}X$ are disjoint, or satisfy some weaker separation
assumption (e.g. the open set condition). Then the dimension can be
computed exactly: $\dim X$ is equal to the similarity dimension\footnote{The similarity dimension depends on the IFS $\Phi$ rather than the
attractor, but we prefer the shorter notation $\sdim X$ in which
$\Phi$ is implicit. The meaning should always be clear from the context.
A similar comment holds for the similarity dimension of measures.} $\sdim X$, i.e. the unique $s\geq0$ solving the equation $\sum|r_{i}|^{s}=1$,
and $\dim\mu$ is equal to the similarity dimension of $\mu$, defined
by 
\[
\sdim\mu=\frac{\sum p_{i}\log p_{i}}{\sum p_{i}\log r_{i}}.
\]

It is when the images $\varphi_{i}X$ have more substantial overlap
that the problem becomes very challenging. The similarity dimension,
and the dimension $d$ of the ambient space $\mathbb{R}^{d}$, still
constitute upper bounds. Thus one always has 
\begin{eqnarray}
\dim X & \leq & \min\{d,\sdim X\}\label{eq:similarity-dimension-bound}\\
\dim\mu & \leq & \min\{d,\sdim\mu\}.\label{eq:similarity-bound-for-measures}
\end{eqnarray}
In general little more is known. In fact, we usually cannot even determine
whether or not equality holds in \eqref{eq:similarity-dimension-bound}
and \eqref{eq:similarity-bound-for-measures}. There is one exception
to this, which arises from combinatorial coincidences of cylinder
sets. For $i=i_{1}\ldots i_{n}\in\Lambda^{n}$ write
\[
\varphi_{i}=\varphi_{i_{1}}\circ\ldots\circ\varphi_{i_{n}}.
\]
One says that exact overlaps occur if there is an $n$ and distinct
$i,j\in\Lambda^{n}$ such that $\varphi_{i}=\varphi_{j}$ (in particular
the images $\varphi_{i}X$ and $\varphi_{j}X$ coincide).\footnote{If $i\in\Lambda^{k}$, $j\in\Lambda^{m}$ and $\varphi_{i}=\varphi_{j}$,
then $i$ cannot be a proper prefix of $j$ and vice versa, because
the maps are all contractions. Thus $ij,ji\in\Lambda^{k+m}$ are distinct,
and $\varphi_{ij}=\varphi_{ji}$. This shows that our definition is
equivalent to one asking for coincidence of compositions of possibly
different lengths. Stated differently, exact overlaps means that the
semigroup generated by the $\varphi_{i}$, $i\in\Lambda$, is not
freely generated by them.} If this occurs then the attractor (or self-similar measure) can be
expressed using an IFS $\Psi$ which is a proper subset of $\{\varphi_{i}\}_{i\in\Lambda^{n}}$,
and a strict inequality in \eqref{eq:similarity-dimension-bound}
and \eqref{eq:similarity-bound-for-measures} may follow from the
trivial bounds \eqref{eq:similarity-dimension-bound} and \eqref{eq:similarity-bound-for-measures}
applied to the IFS $\Psi$.

\subsection{\label{sub:Main-results}Main results}

Define the distance between similarities $\psi=rU+a$ and $\psi'=r'U'+a'$
by\footnote{In \cite{Hochman2014} we used the stronger metric in which the term
$|\log r-\log r'|$ is replaced by the discrete distance $\delta_{r,r'}$.
One could do the same here but the metric above is better suited in
some of the generalizations presented in Section \ref{thm:inverse-theorem-for-isometries}
and is good enough for our applications, so we restrict ourselves
to it. } 
\begin{equation}
d(\psi,\psi')=|\log r-\log r'|+\left\Vert U-U'\right\Vert +\left\Vert a-a'\right\Vert .\label{eq:47}
\end{equation}
Here $\left\Vert \cdot\right\Vert $ denotes the Euclidean or operator
norm as appropriate. Given an IFS $\Phi=\{\varphi_{i}\}_{i\in\Lambda}$,
let 
\begin{equation}
\Delta_{n}=\min\{d(\varphi_{i},\varphi_{j})\,:\,i,j\in\Lambda^{n}\,,\,i\neq j\}.\label{eq:48}
\end{equation}
Note that exact overlaps occur if and only if $\Delta_{n}=0$ for
all large $n$, and it is easy to see that $\Delta_{n}\rightarrow0$
at least exponentially fast (this is an easy consequence of contraction).
Convergence may or may not be faster than this, but we note that in
some cases there is an exponential lower bound $\Delta_{n}\geq c^{n}>0$.

The main result of \cite{Hochman2014} was a step towards the folklore
conjecture that when $d=1$, the occurrence of exact overlaps is the
only mechanism which can lead to a strict inequality in\emph{ \eqref{eq:similarity-dimension-bound}
}and \emph{\eqref{eq:similarity-bound-for-measures}}. Specifically,
we proved the following \cite[Corollary 1.2]{Hochman2014}:
\begin{thm}
\label{thm:1-dim-case}For a self-similar set $X\subseteq\mathbb{R}$,
if $\dim X<\min\{1,\sdim X\}$ then $\Delta_{n}\rightarrow0$ super-exponentially,
i.e. $-\frac{1}{n}\log\Delta_{n}\rightarrow\infty$. The same conclusion
holds if $\dim\mu<\min\{1,\sdim\mu\}$ for a self-similar measure
$\mu$ on $X$.
\end{thm}
When $d\geq2$, the analogous conjecture and analogous theorem are
both false. A trivial class of counterexamples arise when the maps
in $\Phi$ preserve a non-trivial affine subspace $V<\mathbb{R}^{d}$,
which is equivalent to having $X\subseteq V$. In this case, if $\sdim X>\dim V$,
then the trivial bound gives 
\[
\dim X\leq\min\{\dim V,\sdim X\}=\dim V<\min\{d,\sdim X\},
\]
even though there may be no exact overlaps.

We say that $\Phi$ is affinely irreducible if the only trivial affine
subspaces are simultaneously preserved by all $\varphi_{i}\in\Phi$.
The following example shows that affine irreducibility is also not
enough for an analog of Theorem \ref{thm:1-dim-case} to hold.
\begin{example}
\label{example}Begin with the IFS $\Phi=\{\varphi_{\pm}\}$ on $\mathbb{R}$
given by $\varphi_{\pm}(x)=\lambda^{-1}x\pm1$, where $\lambda=1.6956\ldots$
is the real root of $t^{3}-t^{2}-2=0$. This example, due to Garsia
\cite{Garsia1962}, has the property that $\Delta_{n}\geq c\cdot2^{-n}$,
and the attractor is the interval $[-\frac{\lambda}{\lambda-1},\frac{\lambda}{\lambda-1}]$.
Let $\Phi^{3}=\{\varphi_{i}\}_{i\in\{\pm\}^{3}}$ denote the IFS consisting
of all three-fold compositions of the maps $\varphi_{+},\varphi_{-}$.
Then $\Phi^{3}$ has the same attractor and all the maps in $\Phi^{3}$
contract by the same ratio $\lambda^{-3}<1/2$. Now let $\Psi=\{\varphi_{-}^{3},\varphi_{+}^{3}\}$,
where $\varphi^{3}=\varphi\circ\varphi\circ\varphi$. Then $\Psi$
is an IFS with the same contraction ratio $\lambda^{-3}$ as $\Phi^{3}$
but it satisfies the strong separation condition (its attractor $Y$
is the disjoint union of $\varphi_{+}^{3}Y$ and $\varphi_{-}^{3}Y$),
and hence $\dim Y=\log2/\log\lambda^{3}<1$. Finally, take the product
IFS $\Gamma=\Phi^{3}\times\Psi$, consisting of all maps of the form
$(x,y)\mapsto(\varphi x,\psi y)$ for $\varphi\in\Phi^{3}$, $\psi\in\Psi$.
The attractor $Z$ of $\Gamma$ is just the product $Z=[-\frac{\lambda}{\lambda-1},\frac{\lambda}{\lambda-1}]\times Y$
of the attractors of $\Phi^{3}$ and $\Psi$, and its dimension is
$1+\log2/\log\lambda$. We can compute the similarity dimension of
$Z$ using $\lambda^{2}>2$ and $\lambda^{3}-\lambda^{2}-2=0$: 
\[
\sdim Z=\frac{\log|\Gamma|}{\log\lambda^{3}}=\frac{\log16}{\log\lambda^{3}}=\frac{\log16}{\log(2+\lambda^{2})}<2.
\]
We therefore have (using $\lambda<2$): 
\[
\dim Z=1+\frac{\log2}{\log\lambda}<\frac{\log16}{\log\lambda^{3}}=\min\{2,\sdim Z\}.
\]
On the other hand, since both $\Phi^{3}$ and $\Psi$ have exponential
lower bounds on the distance between cylinders, there is also an exponential
lower bound for $\Gamma$. Thus, the example shows that a strict inequality
in \eqref{eq:self-similar-set} with neither exact overlaps or even
super-exponential concentration of cylinders.
\end{example}
Two things stand out about this example. First, the foliation of $\mathbb{R}^{2}$
by horizontal lines is preserved by all maps in $\Gamma$, and, second,
the excess similarity dimension is being ``absorbed'' in the intersection
of the attractor of $\Gamma$ with these lines. Indeed, in these intersections
we are seeing essentially the 1-dimensional IFS $\Phi$, and we are
not getting all of the potential dimension out of it, since its similarity
dimension is $>1$ but attractor is ``trapped'' in a line. We do,
however, have the maximal possible dimension for the intersection
of $Z$ with those horizontal lines that intersect it. 

For an IFS $\Phi=\{\varphi_{i}\}_{i\in\Lambda}$ on $\mathbb{R}^{d}$,
we say that a linear subspace $V<\mathbb{R}^{d}$ is $D\Phi$-invariant
if it is invariant under the orthogonal parts (i.e. differentials)
$U_{i}=D\varphi_{i}$ of $\varphi_{i}\in\Phi$, and nontrivial if
$0<\dim V<d$. If every $D\Phi$-invariant subspace is trivial then
$\Phi$ is said to be \emph{linearly irreducible}. The discussion
above suggests the following:
\begin{conjecture}
\label{conjecture}Let $X\subseteq\mathbb{R}^{d}$ be the attractor
of an affinely irreducible IFS $\Phi\subseteq G$. Then one of the
following must hold:
\begin{enumerate}
\item [(i)] $\dim X=\min\{d,\sdim X\}$.
\item [(ii)] There are exact overlaps.
\item [(iii)] There is a non-trivial $D\Phi$-invariant linear subspace
$V\leq\mathbb{R}^{d}$ and $x\in X$ such that 
\[
\dim(X\cap(V+x))=\dim V.
\]

\end{enumerate}
\end{conjecture}
One might even conjecture a stronger form of (iii), e.g. that the
set of points $x$ in question is of full dimension in $X$, or is
large in some other sense.

The main result of this paper, Theorem \ref{thm:main-individual-RD-entropy-1},
confirms a weakened version of Conjecture \ref{conjecture}:
\begin{thm}
\label{thm:main-sets-Rd}Let $X\subseteq\mathbb{R}^{d}$ be the attractor
of an affinely irreducible IFS $\Phi\subseteq G$. Then one of the
following must hold:
\begin{enumerate}
\item [(i')] $\dim X=\min\{d,\sdim X\}$.
\item [(ii')] $\Delta_{n}\rightarrow0$ super-exponentially.
\item [(iii')] There exists a non-trivial $D\Phi$-invariant linear subspace
$V\leq\mathbb{R}^{d}$ and $x\in X$ such that
\[
\dim(X\cap(V+x))=\dim V.
\]

\end{enumerate}
\end{thm}
The alternatives are not exclusive (all three may hold simultaneously).

The theorem follows, as in the one-dimensional case, from a more precise
statement about the entropy of the measure at small scales. We require
some notation. The level-$n$ dyadic partition $\mathcal{D}_{n}$
of $\mathbb{R}$ is the partition into intervals $[k/2^{n},(k+1)/2^{n})$,
$k\in\mathbb{Z}$. The level-$n$ dyadic partition of $\mathbb{R}^{d}$
is given by 
\[
\mathcal{D}_{n}^{d}=\{I_{1}\times\ldots\times I_{d}\,:\,I_{i}\in\mathcal{D}_{n}\}.
\]
We omit the superscript $d$ when it is clear from the context. 

For a probability measure $\nu$ and partitions $\mathcal{E},\mathcal{F}$
of the underlying probability space we write $H(\nu,\mathcal{E})=-\sum_{E\in\mathcal{E}}\nu(E)\log\nu(E)$
and $H(\nu,\mathcal{E}|\mathcal{F})=H(\nu,\mathcal{E\lor\mathcal{F}})-H(\nu,\mathcal{F})$
for the entropy and conditional entropy of $\nu$ with respect to
$\mathcal{E}$ (conditioned on $\mathcal{F}$, respectively). Here
$\mathcal{E}\lor\mathcal{F}$ is the common refinement of the partitions
$\mathcal{E},\mathcal{F}$. We also write $H(\nu)$ for the entropy
of an atomic measure $\nu$ with respect to the partition into points. 

It is convenient to parametrize $G$ as a subset of $\mathbb{R}\times M_{d}(\mathbb{R})\times\mathbb{R}^{d}$,
with $(t,U,a)$ corresponding to $2^{-t}U+a\in G$. Then the level-$n$
dyadic partition $\mathcal{D}_{n}^{G}$ of $G$ is defined as the
partition induced from the corresponding level-$n$ partition of $\mathbb{R}\times M_{d}(\mathbb{R})\times\mathbb{R}^{d}\cong\mathbb{R}^{1+d^{2}+d}$.
We also introduce the partitions $\mathcal{E}_{n}^{G}$ of $G$ induced
by the dyadic partition according to the translation part of the similarities,
which in the parametrization $G\subseteq\mathbb{R}\times M_{d}(\mathbb{R})\times\mathbb{R}^{d}$
is 
\[
\mathcal{E}_{n}^{G}=\{(\mathbb{R}\times M_{d}(\mathbb{R})\times D)\cap G\,:\,D\in\mathcal{D}_{n}^{d}\}.
\]
Note that $\mathcal{D}_{n}^{G}$ refines $\mathcal{E}_{n}^{G}$. 

Given a self-similar measure $\mu=\sum_{i\in\Lambda}p_{i}\cdot\varphi_{i}\mu$
and assuming all $p_{i}>0$, let 
\[
\nu^{(n)}=\sum_{i\in\Lambda^{n}}p_{i}\cdot\delta_{\varphi_{i}}.
\]
This is a probability measure on $G$, but if we fix $\widetilde{x}$
in the attractor of $\Phi$ then the push-forward of $\nu^{(n)}$
via $g\mapsto g\widetilde{x}$ is the natural ``$n$-th generation''
approximation of $\mu$, given by\footnote{It is also common to approximated $\mu$ ``at scale $\rho$'' by
putting the appropriate mass on the points $\varphi_{i_{1}\ldots i_{m}}(x)$,
where $i_{1}\ldots i_{m}\in\Lambda^{*}$ are the sequences of minimal
length such that $\varphi$$_{i_{1}\ldots i_{m}}$ contracts by at
least $\rho$. We could use this approximation instead of $\nu^{(n)}$,
but this would lead to messier notation and have little advantage.}
\[
\widetilde{\nu}^{(n)}=\sum_{i\in\Lambda^{n}}p_{i}\cdot\delta_{\varphi_{i}(\widetilde{x})}
\]
(this measure depends on the choice of $\widetilde{x}$ but this is
of little consequence). Let 
\[
r=\prod_{i\in\Lambda}r_{i}^{p_{i}}
\]
denote the (geometric) average contraction and for $n\in\mathbb{N}$
let 
\[
n'=[n/\log(1/r)],
\]
so that $2^{-n'}\sim r^{n}$.

Now, it is not hard to show $|H(\nu^{(n)},\mathcal{E}_{n'}^{G})-H(\widetilde{\nu}^{(n)},\mathcal{D}_{n'})|=O(1)$
(in fact if we take $\widetilde{x}=0$ then two entropies are identical),
and since it is easily seen that $\frac{1}{n'}H(\widetilde{\nu}^{(n)},\mathcal{D}_{n'})\rightarrow\dim\mu$,
one concludes 
\[
\lim_{n\rightarrow\infty}\frac{1}{n'}H(\nu^{(n)},\mathcal{E}_{n'}^{G})=\dim\mu.
\]
Observe that when there are no exact overlaps, $\nu^{(n)}$ consists
of $|\Phi|^{n}$ atoms whose masses are all the products $p_{i_{1}}\cdot\ldots\cdot p_{i_{n}}$,
and hence $H(\nu^{(n)})=n\cdot(-\sum p_{i}\log p_{i})$. Thus for
fixed $n$, 
\[
\frac{1}{n'}H(\nu^{(n)},\mathcal{D}_{k})\rightarrow\frac{-\sum p_{i}\log p_{i}}{\log r}=\sdim\mu\qquad\mbox{as }k\rightarrow\infty,
\]
and if there is a strict inequality in \eqref{eq:similarity-bound-for-measures}
we would have 
\begin{eqnarray*}
\frac{1}{n'}H(\nu^{(n)},\mathcal{D}_{k}^{G}|\mathcal{E}_{n'}^{G}) & = & \frac{1}{n'}H(\nu^{(n)},\mathcal{E}_{k}^{G})-\frac{1}{n'}H(\nu^{(n)},\mathcal{D}_{k}^{G})\\
 & \rightarrow & \sdim\mu-\dim\mu\qquad\qquad\qquad\mbox{as }k\rightarrow\infty\\
 & > & 0
\end{eqnarray*}
Therefore it is possible to choose $k=k(n)$ such that the ``excess''
$\frac{1}{n'}H(\nu^{(n)},\mathcal{D}_{k(n)}^{G}|\mathcal{E}_{n'}^{G})$
remains bounded away from $0$ as $n\rightarrow\infty$. It is natural
to ask at what rate this excess entropy emerges, that is, how fast
$k(n)$ must grow for this to hold. The following theorem shows that
it must grow at least super-linearly.
\begin{thm}
\label{thm:main-individual-RD-entropy}Let $\mu\in\mathcal{P}(\mathbb{R}^{d})$
be a self-similar measure for an affinely irreducible IFS $\Phi$.
Then one of the following must hold:
\begin{enumerate}
\item [(i'')] $\dim\mu=\min\{d,\sdim\mu\}$.
\item [(ii'')] $\lim_{n\rightarrow\infty}\frac{1}{n'}H(\nu^{(n)},\mathcal{D}_{qn}^{G}|\mathcal{E}_{n'}^{G})=0$
for all $q>1$.
\item [(iii'')] There is a non-trivial $D\Phi$-invariant linear subspace
$V\leq\mathbb{R}^{d}$ such that for $\mu$-a.e. $x$, the conditional
measure $\mu_{V+x}$ on $V+x$ satisfies $\dim\mu_{V+x}=\dim V$.
\end{enumerate}
\end{thm}
In fact, the second or third alternatives must hold irrespective of
the validity of (i''). The usefulness of the theorem, however, lies
in the fact that if (i'') fails and (ii'') holds then $\Delta_{n}\rightarrow0$
super-exponentially. 

Theorem \ref{thm:main-individual-RD-entropy} and Theorem \ref{thm:main-sets-Rd}
are usually applied by ruling out (iii') or (iii''), and then working
out the implications for the dimension. One trivial way to rule it
out is to just assume it:
\begin{cor}
If $\Phi$ is a linearly irreducible IFS $\Phi$, then its attractor
$X$ satisfies (i') or (ii'), and every self-similar measure $\mu$
for $\Phi$ satisfies (i'') or (ii'').
\end{cor}
As there are no non-trivial linear subspaces of $\mathbb{R}$, every
IFS acts linearly irreducibly, and we have recovered the main results
of \cite{Hochman2014} (Theorem \ref{thm:1-dim-case} above).

We say that $rU+a\in G$ is algebraic if $r$ and all the coordinates
of $U$ and $a$ are algebraic numbers over $\mathbb{Q}$, and we
say that an IFS $\Phi\subseteq G$ is algebraic if all of its elements
are. If $\Phi$ is an algebraic IFS without exact overlaps, and we
take $\widetilde{x}=0$, then for each $n$,  $\Delta_{n}$ is a polynomial
in the algebraic parameters defining the maps on $\Phi$ and has degree
$n$ and height at most exponential in $n$. This implies an exponential
lower bound $\Delta_{n}\geq c^{n}$; this is a well known fact but
we include a proof in Section \ref{sub:Applications-Rd-proofs}. Thus
we have ruled out (ii') ad (ii''), and obtained the following:
\begin{cor}
\label{cor:irr-algebraic-IFS}Let $\Phi$ be an algebraic IFS acting
linearly irreducibly on $\mathbb{R}^{d}$ and without exact overlaps.
Then $\dim\mu=\min\{d,\sdim\mu\}$ for every fully supported self-similar
measure $\mu$ of $\Phi$, and $\dim X=\min\{d,\sdim X\}$.
\end{cor}
Our arguments are purely Euclidean and do not utilize any non-elementary
properties of the orthogonal or similarity groups. However, the nature
of these groups depends crucially on the dimension $d$. For $d\leq2$
the orthogonal group of $\mathbb{R}^{d}$ is abelian (and the similarity
group is solvable). In particular, the set $\mathcal{U}_{n}=\{U_{i}\}_{i\in\Lambda^{n}}$
of the orthogonal parts of $\varphi_{i}$, $i\in\Lambda^{n}$, is
of polynomial size in $n$, and does not contribute to the entropy
$H(\nu^{(n)},\mathcal{D}_{qn}^{G}|\mathcal{E}_{n'}^{G})$ (for the
same reason, the contraction ratios do not contribute asymptotically
to the entropy). For $d\geq3$ the orthogonal group is a virtually
simple Lie group with strong expansion properties, and typically $|\mathcal{U}_{n}|$
is exponential in $n$. Our methods do not make use of any special
properties of the orthogonal group, but concurrently and independently
with our work, Lindenstrauss and Varj\'{u} utilized the work of Bourgain
and Gamburd \cite{BourgainGamburd2008} and of de Saxce \cite{deSaxce2014}
on spectral gap of random walks on the orthogonal group to prove the
following result. 
\begin{thm}
[Lindenstrauss-Varj\'{u}, \cite{LindenstraussVarju2012}] Let $U_{1},\ldots,U_{k}\in SO(d)$
and $p=(p_{1},\ldots,p_{k})$ a probability vector. Suppose that the
operator $f\mapsto\sum_{i=1}^{k}p_{i}f\circ U_{i}$ on $L^{2}(SO(d))$
has a spectral gap. Then there is a number $\widetilde{r}<1$ such
that for every choice $\widetilde{r}<r_{1},\ldots,r_{k}<1$, and for
any $a_{1},\ldots,a_{k}\in\mathbb{R}^{d}$, the self similar measure
with weights $p$ for the IFS $\{r_{i}U_{i}+a_{i}\}_{i=1}^{k}$ is
absolutely continuous with respect to Lebesgue measure on $\mathbb{R}^{d}$.
\end{thm}
The spectral gap hypothesis can currently be verified when the entries
are algebraic and $U_{i}$ generate a dense subgroup of $O(d)$, but
is conjectured to hold much more generally.

Compare this theorem to Corollary \ref{cor:irr-algebraic-IFS}: The
former ensures absolute continuity (which is a stronger property than
full dimension), but only when the contraction of the IFS is uniformly
close enough to $1$, while the latter ensures that the dimension
is $d$ as soon as there is no dimension obstruction (i.e. as soon
as $\sdim\mu\geq d$), but does not give absolute continuity. It is
probable that absolute continuity holds under the same assumptions
but this remains open.

There are other cases in which possibility (iii') of Theorem \ref{thm:main-sets-Rd}
or (iii'') of Theorem \ref{thm:main-individual-RD-entropy}, can be
ruled out. A trivial case is when the attractor $X$ of $\Phi$ satisfies
$\dim X<k$, and all $D\Phi$-invariant subspaces have dimension $\geq k$.
Another case is when $\Phi$ consists of homotheties (i.e. the orthogonal
parts $U_{i}$ of the contractions are identities), and for every
line $\ell$ in $\mathbb{R}^{d}$ we have 
\[
\sum_{i\,:\,\varphi_{i}(X)\cap\ell\neq\emptyset}r_{i}<1
\]
Then elementary covering considerations show that $\dim(X\cap\ell)<1$
for every line $\ell\subseteq\mathbb{R}^{d}$, and consequently (iii')
(and hence (iii'')) fails for every subspace $V$. Similarly, if $\Phi$
consists of homotheties and $\mu=\sum p_{i}\cdot\varphi_{i}\mu$ is
a self-similar measure such that for every line $\ell$,
\[
\frac{\sum_{i\,:\,\varphi_{i}(X)\cap\ell\neq\emptyset}p_{i}\log p_{i}}{\sum_{i\,:\,\varphi_{i}(X)\cap\ell\neq\emptyset}p_{i}\log r_{i}}<1
\]
then one can deduce that $\dim\mu_{\ell+x}<1$ for $\mu$-a.e. $x$,
which by Marstrand's slice theorem rules out (iii''). Another alternative
is to show that the linear images onto $(d-1)$-planes have dimension
greater than $\dim\mu-1$, in which case Dimension conservation \cite{Furstenberg08}
implies that in every dimension, the conditional measure on a.e. line
has dimension $<1$. 

Unfortunately such arguments do not always apply, and we know of no
general method to exclude (iii') and (iii''). See Theorem \ref{thm:fat-Sierpinski}
and the discussion surrounding it.

\subsection{\label{sub:Applications-Rd}Parametric families}

Suppose that $I$ is a set of parameters and that for $t\in I$ we
are given an IFS $\Phi_{t}=\{\varphi_{i,t}\}$, where $\varphi_{i,t}(x)=r_{i}(t)U_{i}(t)x+a_{i}(t)$
for  functions $r_{i},U_{i},a_{i}$ defined on $I$. For $i,j\in\Lambda^{n}$
let 
\[
\Delta_{i,j}(t)=\varphi_{i,t}(0)-\varphi_{j,t}(0).
\]
Then $\left\Vert \Delta_{i,j}(t)\right\Vert $ is the third term in
the definition \eqref{eq:47} of $d(\varphi_{i,t},\varphi_{j,t})$,
and hence, writing $\Delta_{n}(t)$ for the quantity defined as in
\eqref{eq:48} for the system $\Phi_{t}$, we have 
\[
\min\{\left\Vert \Delta_{i,j}(t)\right\Vert \,:\,i,j\in\Lambda^{n}\mbox{ distinct}\}\leq\Delta_{n}(t).
\]
This gives the following formal consequence of Theorem \ref{thm:main-individual-RD-entropy}:
\begin{thm}
\label{thm:description-of-exceptional-params-Rd}Let $\{\Phi_{t}\}_{t\in I}$
be a parametric family of IFSs on $\mathbb{R}^{d}$. Let $E\subseteq I$
be the set 
\[
E=\bigcap_{\varepsilon>0}\left(\bigcup_{N=1}^{\infty}\,\bigcap_{n>N}\left(\bigcup_{i,j\in\Lambda^{n}}\Delta_{i,j}^{-1}((-\varepsilon^{n},\varepsilon^{n})^{d})\right)\right),
\]
and let $F\subseteq I$ be the set of parameters $t$ for which $\Phi_{t}$
is linearly reducible. Then for $t\in I\setminus(E\cup F)$, every
self-similar measure $\mu$ for $\Phi_{t}$ satisfies $\dim\mu=\min\{d,\sdim\mu\}$
and similarly for the attractor of $\Phi_{t}$.
\end{thm}
The main case of interest is when $I\subseteq\mathbb{R}^{m}$. Then,
under rather mild assumptions, the set $E$ of (potential) exceptions
can be shown to be quite small. For $i,j\in\Lambda^{\mathbb{N}}$
let
\[
\Delta_{i,j}(t)=\lim_{n\rightarrow\infty}\Delta_{i_{1}\ldots i_{n},j_{1}\ldots j_{n}}(t).
\]

\begin{thm}
\label{thm:main-parametric-Rd}Let $I\subseteq\mathbb{R}^{m}$ be
connected and compact and let $\{\Phi_{t}\}_{t\in I}$ be a parametrized
family of IFSs for which the associated functions $r_{i}(\cdot)$,
$U_{i}(\cdot)$ and $a_{i}(\cdot)$ are real-analytic on a neighborhood
of $I$. Suppose that 
\[
\forall i,j\in\Lambda^{\mathbb{N}}\;\left(i\neq j\;\implies\;\Delta_{i,j}\not\equiv0\right).
\]
Then the set $E$ of the previous theorem has Hausdorff and packing
dimension $\leq m-1.$ In particular if $\Phi_{t}$ is linearly irreducible
for all $t\in I$, then outside a set of parameters $t$ of dimension
$\leq m-1$ (and in particular for Lebesgue-a.e. parameter), the attractor
and self-similar measures of $\Phi_{t}$ have the expected dimension
(i.e. equality holds in equations \eqref{eq:similarity-dimension-bound}
and \eqref{eq:similarity-bound-for-measures}).
\end{thm}
The condition $\Delta_{i,j}\not\equiv0$ rules out trivial cases.
For instance the theorem cannot be expected to apply when $\Phi_{t}=\Phi$
does not depend on $t$ and the system $\Phi$ has exact overlaps,
in which case there are indeed distinct $i,j\in\Lambda^{\mathbb{N}}$
with $\Delta_{i,j}\equiv0$.

If $I\subseteq\mathbb{R}^{m}$ and the IFS is in $\mathbb{R}^{d}$,
and $m\geq d$, then we expect that for each $i,j\in\Lambda^{\mathbb{N}}$
there typically will be a sub-manifold $I_{i,j}\subseteq I$ of dimension
$m-d$ on which $\Delta_{i,j}=0$ for $i\in I_{i,j}$. Thus, the dimension
bound on $E$ that one expects is $m-d$ rather than the bound $m-1$
appearing in the theorem above. However, the hypothesis $\Delta_{i,j}\not\equiv0$
in itself is certainly not enough to guarantee this bound. To see
this, begin with any a $1$-parameter family $\{\Phi_{u}\}_{u\in[0,1]}$
of linearly irreducible IFSs in $\mathbb{R}^{2}$, and define a two-parameter
family by $\Phi_{(s,t)}=\Phi_{(s-t)^{2}}$, $(s,t)\in[0,1]^{2}$.
One might expect, by the logic above, that $\dim E=m-d=0$. But, evidently,
on the 1-dimensional subspace $V=\{s=t\}$ we have $\Phi_{(s,t)}=\Phi_{0},$
and if the attractor of $\Phi_{0}$ happens to satisfy \eqref{eq:similarity-dimension-bound}
with a strict inequality, then $\dim E\geq1\neq0=m-d$.

It is natural to suggest that, assuming linear irreducibility of the
IFS, the ``correct'' bound for $E$ is
\begin{equation}
\dim E\leq m-\sup\left\{ \dim\Delta_{i,j}^{-1}(0)\,:\,i,j\in\Lambda^{\mathbb{N}}\,,\,i\neq j\right\} .\label{eq:105}
\end{equation}
For $m=d=1$, the bound proved in \cite{Hochman2014} coincides with
this one. The difficulty in higher dimension is that the zero sets
of real-analytic functions, and the behavior of the functions near
them, are not so well understood (for real-analytic functions on the
line things are simple: the zero set consists of isolated points,
away from which the function grows polynomially in a well-understood
manner). It seems likely that having effective bounds on the constants
in \L{}ojasiewicz's inequality \cite{Lojasiewicz1962} might advance
the matter but this seems to a difficult question in itself. What
we prove here is that the bound \eqref{eq:105} holds if one makes
an assumption analogous to the classical transversality assumption.
\begin{thm}
\label{thm:main-parametric-Rd-2}Let $I\subseteq\mathbb{R}^{m}$ be
compact and let $\{\Phi_{t}\}_{t\in I}$ be a parametrized family
of IFSs for which the associated functions $r_{i}(\cdot)$, $U_{i}(\cdot)$
and $a_{i}(\cdot)$ are real-analytic on a neighborhood of $I$. Suppose
that there exists an $r\in\mathbb{N}$ such that for every distinct
pair $i,j\in\Lambda^{n}$ and $t\in I$, 
\[
\Delta_{i,j}(t)=0\qquad\implies\qquad\rank\left(D\Delta_{i,j}(t)\right)\geq r.
\]
Then the set $E$ of Theorem \ref{thm:description-of-exceptional-params-Rd}
has Hausdorff and packing dimension $\leq m-r.$ 
\end{thm}
As noted above, it is likely that there is room for improvement in
these results.

\subsection{Applications}

We demonstrate the use of Theorems \ref{thm:main-parametric-Rd} and
\ref{thm:main-parametric-Rd-2} for families of self-similar measures
in which one varies the translations, contractions, or the IFS. Proofs
are given in Section \ref{sub:Applications-Rd-proofs}.

Let $X_{\Phi}\subseteq\mathbb{R}^{d}$ denote the attractor of an
IFS $\Phi$. 
\begin{thm}
\label{thm:generic-IFS}For a finite set $\Lambda$ and $d\in\mathbb{N}$
let $IFS_{\Lambda}\subseteq G(d)^{\Lambda}$ denote the set $|\Lambda|$-tuples
of contracting similarities , which we identify with the set of IFSs
indexed by $\Lambda$. Then 
\[
\dim\{\Phi\in IFS_{\Lambda}\,:\,\dim X_{\Phi}<\min\{d,\sdim_{\Phi}X_{\Phi}\}\}\leq\dim IFS_{\Lambda}-1.
\]
In particular, $\dim X_{\Phi}=\dim\{1,\sdim X_{\Phi}\}$ for a.e.
IFS $\Phi\in IFS_{\Lambda}$.
\end{thm}
If one fixes the linear parts of the similarity maps and varies the
translation part, one obtains a version of results by Simon and Solomyak
\cite{SimonSolomyak2002}:
\begin{thm}
\label{thm:generic-translation-IFS}Let $\{U_{i}\}_{i\in\Lambda}$
be orthogonal maps acting irreducibly on $\mathbb{R}^{d}$ and fix
$0<r_{i}<1$, $i\in\Lambda$, satisfying the condition
\[
i\neq j\qquad\implies\qquad r_{i}+r_{j}<1.
\]
Then there is a subset $A\subseteq(\mathbb{R}^{d})^{\Lambda}$ with
$\dim(\mathbb{R}^{d})^{\Lambda}\setminus A\leq d|\Lambda|-d$, and
such that for $a\in A$ the attractor of $\Phi=\{r_{i}U_{i}+a_{i}\}_{i\in\Lambda}$
satisfies $\dim X_{\Phi}=\dim\{1,\sdim X_{\Phi}\}$. In particular
this is true for a.e. $a\in(\mathbb{R}^{d})^{\Lambda}$. 
\end{thm}
The condition on the contraction ratios plays a similar role in \cite[Theorem 2.1(c)]{SimonSolomyak2002}
and the forthcoming book \cite{SimonSolomyak2014}, where it is used
in conjunction with the transversality method. It is needed to control
the rank of $D\Delta_{i,j}$, which in our setting is required in
order to apply Theorem \ref{thm:main-parametric-Rd-2}. It is not
clear to what extent the restriction on the contractions in necessary,
but without the irreducibility condition it certainly is, as follows
from \cite[Proposition 3.3]{SimonSolomyak2002}.

Another variant of these results concerns projections of self-similar
measures defined by homotheties. This is a variant of Marstrand's
theorem and Furstenberg's projection problem \cite{Kenyon97,Hochman2014}:
\begin{thm}
\label{thm:generic-projection-IFS}Let $X\in\mathcal{P}(\mathbb{R}^{d})$
be a self-similar set defined by an IFS consisting of homotheties
and satisfying strong separation. Let $k<d$ and let $\Pi_{d,k}$
denote the set of orthogonal projections from $\mathbb{R}^{d}$ to
$k$-dimensional subspaces. Then 
\[
\dim\{\pi\in\Pi_{d,k}\,:\,\dim\pi X=\min\{k,\dim X\}\}\leq\dim\Pi_{d,k}-k.
\]

\end{thm}
A particularly interesting family are the Bernoulli convolutions with
nonuniform contraction. Namely, for $0<\beta,\gamma<1$ let $\lambda_{\beta,\gamma}$
denote the self-similar measure of maximal dimension for the IFS $\{x\mapsto\beta x,x\mapsto\gamma x+1\}$.
Let $S\subseteq(0,1)^{2}$ be the set of $(\beta,\gamma)$ for which
$\sdim\lambda_{\beta,\gamma}>1$; it is expected that $\lambda_{\beta,\gamma}$
is absolutely continuous for a.e. $(\beta,\gamma)\in S$, but this
has been established only in certain restricted parameter ranges,
e.g. \cite{Neunhauserer2001}. 
\begin{thm}
\label{thm:nonuniform-BC}$\dim\lambda_{\beta,\gamma}=\min\{1,\sdim\lambda_{\beta,\gamma}\}$
outside a set of parameters $(\beta,\gamma)\in(0,1)^{2}$ of Hausdorff
(and packing) dimension $1$. In particular this holds for Lebesgue-a.e.
pair $(\beta,\gamma)\in(0,1)^{2}$.
\end{thm}
Finally, our results can be applied to a higher-dimensional analogs
of the Bernoulli convolutions problem, namely the ``fat Sierpinski
gasket'', first studied by Simon and Solomyak \cite{SimonSolomyak2002}.
For $\lambda\in(0,1)$ consider the system of contractions $\{\varphi_{i}\}_{i=a,b,c}$
where $a,b,c$ are the vertices of an equilateral triangle in $\mathbb{R}^{2}$
and $\varphi_{u}(x)=\lambda x+u$. The classical Sierpinski gasket
arises from the choice $\lambda=1/2$, and in general when $0\leq\lambda\leq1/2$
the open set condition is satisfied and the dimension of the attractor
$S_{\lambda}$ is equal to the similarity dimension. When $\lambda>2/3$
the attractor has non-empty interior, and this remains true for $\lambda\geq\lambda_{*}$,
where $\lambda_{*}\approx0.6478$ is the real root of $x^{3}-x^{2}+x=1/2$;
see Broomhead-Montaldi-Sidorov \cite{BroomheadMontaldiSidorov2004}.
For $1/2<\lambda\leq\lambda_{*}$, however, the the dimension is known
only for certain special algebraic parameters and for Lebesgue-typical
$\lambda$ in a certain sub-range, and similarly for absolute continuity
of the appropriate self-similar measures. See Jordan \cite{Jordan2005}
and Jordan-Pollicott \cite{JordanPollicott2006}. 
\begin{thm}
\label{thm:fat-Sierpinski}$\dim S_{\lambda}=\min\{2,\sdim S_{\lambda}\}$
for $\lambda\in(0,1)$ outside a set of Hausdorff (and packing) dimension
$0$.
\end{thm}
The last result is an immediate consequence of Theorem \ref{thm:main-parametric-Rd}
using the fact that $S_{\lambda}$ can be written also as the attractor
of a linearly irreducible IFS (the one given above is reducible).
The possibility of such a presentation of $S_{\lambda}$ comes from
its rotational symmetries. Interestingly, our method do not give comparable
results even for very slight variants of $S_{\lambda}$, e.g. the
fat Sierpinski gaskets studied in  \cite{JordanPollicott2006}.

\subsection{\label{sub:Organization}Organization and notation}

A key ingredient in our argument is played by on the growth of entropy
of measures under convolution. This subject is developed in the next
three sections: Section \ref{sec:Additive-combinatorics} introduces
the statements and basic definitions, Section \ref{sec:Entropy-concentration-uniformity-saturation}
contains preliminaries on entropy, saturation, concentration and convolutions,
and Section \ref{sec:Convolutions} proves the main results on convolutions.
In Section \ref{sec:convolutions-with-isometries} we extend the results
to convolutions of a measure on $\mathbb{R}^{d}$ with a measure on
the isometry group. Finally, in Section \ref{sec:Self-similar-sets-and-measures-Rd}
we state and prove our main theorem on self-similar sets and measures
and their applications.

Some notation: $\mathbb{N}=\{1,2,3,\ldots\}$. All logarithms are
to base $2$. $\mathcal{P}(X)$ is the space of probability measures
on $X$, endowed with the weak-{*} topology if appropriate. We follow
standard ``big $O$'' notation: $O_{\alpha}(f(n))$ is an unspecified
function bounded in absolute value by $C\cdot f(n)$ for some constant
$C=C(\alpha)$ depending on $\alpha$. Similarly $o(1)$ is a quantity
tending to $0$ as the relevant parameter $\rightarrow\infty$. We
implicitly suppress all dependence of constants on the dimension $d$
of $\mathbb{R}^{d}$. Thus $O(1)$ sometimes means $O_{d}(1)$. We
sometimes write $-O(\cdot)$ instead of $+O(\cdot)$ to indicate that
the error may be negative but formally the two notations are equivalent.

The statement ``for all $s$ and $t>t(s),\ldots$'' should be understood
as saying ``there exists a function $t(\cdot)$ such that for all
$s$ and $t>t(s),\ldots$''. The function $t(\cdot)$ will change
between contexts, when we want a persistent name we will designate
the function as $t_{1}(\cdot)$, $t_{2}(\cdot)$, $t_{*}(\cdot)$,
etc. 

For the reader's convenience we summarize our main notation in the
table below.

\noindent \bigskip{}
\begin{tabular}{ll}
\hline 
$d$  & Dimension of the ambient Euclidean space.\tabularnewline
$B_{r}(x)$  & The open Euclidean ball of radius $r$ around $x$\tabularnewline
$\left\Vert x\right\Vert ,$$\left\Vert A\right\Vert $ & Euclidean norm of $x\in\mathbb{R}^{d}$, operator norm of $A\in M_{d}(\mathbb{R})$\tabularnewline
$\dim$  & Hausdorff dimension of sets and measures\tabularnewline
$\Phi=\{\varphi_{i}\}_{i\in\Lambda}$ & Iterated Function system, Section \ref{sub:Self-similar-sets}\tabularnewline
$X$ & Attractor of $\Phi$. Usually assume $0\in X\subseteq[0,1)$, Section
\ref{sub:Self-similar-sets}\tabularnewline
$\mu$ & Self-similar measure (usually), Section \ref{sub:Self-similar-sets}\tabularnewline
$\varphi_{i_{1}\ldots i_{n}}$, $p_{i_{1}\ldots i_{n}}$ & $\varphi_{i_{1}}\circ\varphi_{i_{2}}\circ\ldots\varphi_{i_{n}}$ and
$p_{i_{1}}\cdot p_{i_{2}}\cdot\ldots\cdot p_{i_{n}}$\tabularnewline
$\nu^{(n)}$ & $\sum_{i\in\Lambda^{n}}p_{i}\cdot\delta_{\varphi_{i}(0)}$, the $n$-th
approximation of $\mu$\tabularnewline
$\mathcal{D}_{n}^{k}$ & $n$-th level dyadic partition of $\mathbb{R}^{k}$  ($k=d$ by default);
Section \ref{sub:Main-results}\tabularnewline
$\mathcal{D}_{n}^{G}$ & Dyadic partition of $G\subseteq\mathbb{R}^{+}\times M_{d}(\mathbb{R})\times\mathbb{R}^{d}$,
Section \ref{sub:Main-results}\tabularnewline
$\mathcal{E}_{n}^{G}$ & Dyadic partition of $G$ by translation part, Section \ref{sub:Main-results}\tabularnewline
$\mathcal{P}(X)$ & Space of probability measures on $X$.\tabularnewline
$\mu_{x,n},\mu^{x,n}$ & Component measures (raw and re-scaled), Section \ref{sub:Component-measures}\tabularnewline
$S_{t}$ & Scaling map: $S_{t}(x)=2^{t}x$\tabularnewline
$\tau_{z}$ & Translation map: $\tau_{s}(x)=x+s$\tabularnewline
$\mathbb{P}_{i\in I}$, $\mathbb{E}_{i\in I}$ & Distribution and expectation over components, Section \ref{sub:Component-measures}\tabularnewline
$H(\mu,\mathcal{B})$ & Shannon entropy, Section \ref{sub:Preliminaries-on-entropy}\tabularnewline
$H(\mu,\mathcal{B}|\mathcal{C})$ & Conditional entropy, Section \ref{sub:Preliminaries-on-entropy}\tabularnewline
$H_{m}(\mu)$ & $\frac{1}{m}H(\mu,\mathcal{D}_{m})$, Section \ref{sub:Preliminaries-on-entropy}\tabularnewline
$G,G_{0}$ & The groups of similarities and isometries, respectively.\tabularnewline
$\pi_{V}$ & Orthogonal projection to $V$\tabularnewline
$V^{(\varepsilon)}$ & $\varepsilon$-neighborhood of $V$\tabularnewline
$d(U,V)$ & Distance between linear subspaces of $U,V\leq\mathbb{R}^{d}$, Section
\ref{sub:Geometry-of-subspaces}\tabularnewline
$\sqsubseteq$ & Subset relation restricted to unit ball, Section \ref{sub:Geometry-of-subspaces}\tabularnewline
$\angle(U,V)$ & (Modified) angle between linear subspaces, Section \ref{sub:Geometry-of-subspaces}\tabularnewline
$\mu*\eta$ & Convolution of probability measure on $\mathbb{R}^{d}$.\tabularnewline
$\nu\conv x$, $\nu\conv\mu$ & Action/convolution of $\nu\in\mathcal{P}(G)$ on $x\in\mathbb{R}^{d}$,
$\mu\in\mathcal{P}(\mathbb{R}^{d})$.\tabularnewline
$m(\mu)$, $\Sigma(\mu)$ & Mean and covariance matrix of measure $\mu\in\mathcal{P}(\mathbb{R}^{d})$,
Section \ref{sub:Covariance-matrices}\tabularnewline
$\lambda_{i}(\mu),\lambda_{i}(\Sigma)$ & Eigenvalues of measure or covariance matrix, Section \ref{sub:Covariance-matrices}\tabularnewline
$\eigen_{1\ldots r}(\Sigma)$ & Span of top $r$ eigenvectors of $\Sigma$ (for measure, $\Sigma=\Sigma(\mu)$),
Section \ref{sub:Covariance-matrices}\tabularnewline
$\sat(\eta,\varepsilon,n,m)$ & Set of $(V,\varepsilon,m)$-saturated subspaces at level $n$, Section
\ref{sub:saturation-and-self-similar-measures}\tabularnewline
\hline 
\end{tabular}

\subsection*{Acknowledgment}

I am grateful to Pablo Shmerkin and Boris Solomyak for their many
helpful comments, and to Ariel Rapaport for his contribution to the
argument in Section \ref{sub:Entropy-and-dimension}. Part of this
work was done during a visit to Microsoft Research in Redmond, Washington,
and I would like to thank Yuval Peres and the members of the theory
group for their hospitality.

\section{\label{sec:Additive-combinatorics}An inverse theorem for the entropy
of convolutions}

\subsection{Entropy and additive combinatorics}

A subject of independent interest and central to our work is an analysis
of the growth of the entropy of measures under convolution, either
with other measures or with measures on the group of isometries (or
similarities). This topic will occupy us for a large part of the paper. 

We begin with a discussion of convolutions on Euclidean space, leaving
generalizations to later. It is convenient to introduce the normalized
scale-$n$ entropy 
\[
H_{n}(\mu)=\frac{1}{n}H(\mu,\mathcal{D}_{n}).
\]
This normalization makes $H_{n}(\mu)$ a finite-scale surrogate for
the dimension of $\mu$. In particular, for $\mu\in\mathcal{P}([0,1)^{d})$
we have 
\[
0\leq H_{n}(\mu)\leq d,
\]
with equality holding for all $n$ if and only if $\mu$ is Lebesgue
measure on $[0,1)^{d}$, and in general for measures $\mu$ of bounded
support,
\[
0\leq H_{n}(\mu)\leq d+O(\frac{1}{n}),
\]
where the constant depends logarithmically on the diameter of the
support. 

Our aim is to obtain structural information about measures $\mu,\nu$
for which $\mu*\nu$ is small in the sense that 
\begin{equation}
H_{n}(\mu*\nu)\leq H_{n}(\mu)+\delta,\label{eq:mean-entropy-growth}
\end{equation}
where $\delta>0$ is small but fixed, and $n$ is large. This problem
is a relative of classical ones in additive combinatorics concerning
the structure of sets $A,B$ whose sumset $A+B=\{a+b\,:\,a\in A\,,\,b\in B\}$
is appropriately small. The general principle is that when the sum
is small, the sets should have some algebraic structure. Results to
this effect are known as inverse theorems. For example the Freiman-Rusza
theorem asserts that if $|A+B|\leq C|A|$ then $A,B$ are close, in
a manner depending on $C$, to generalized arithmetic progressions\footnote{A generalized arithmetic progression is an injective affine image
of a box in a higher-dimensional lattice.} (the converse is immediate). See e.g \cite{TaoVu2006}. 

The entropy of a discrete measure corresponds to the logarithm of
the cardinality of a set, and convolution is the analog for measures
of the sumset operation. Thus the analog of the condition $|A+A|\leq C|A|$
is 
\begin{equation}
H_{n}(\mu*\mu)\leq H_{n}(\mu)+O(\frac{1}{n}).\label{eq:O-of-1-entropy-growth}
\end{equation}
An entropy version of Freiman's theorem was recently proved by Tao
\cite{Tao2010}, who showed that if $\mu$ satisfies \eqref{eq:O-of-1-entropy-growth}
then it is close, in an appropriate sense, to a uniform measures on
a (generalized) arithmetic progression. 

The condition \eqref{eq:mean-entropy-growth}, however, is significantly
weaker than \eqref{eq:O-of-1-entropy-growth} even when $\nu=\mu$,
and it is harder to draw conclusions from it about the global structure
of $\mu$. Consider the following example. Start with an arithmetic
progression of length $n_{1}$ and gap $\varepsilon_{1}$, and put
the uniform measure on it. Now split each atom $x$ into an arithmetic
progression of length $n_{2}$ and gap $\varepsilon_{2}<\varepsilon_{1}/n_{2}$,
starting at $x$ (so the entire gap fits in the space between $x$
and the next atom). Repeat this procedure $N$ times with parameters
$n_{i},\varepsilon_{i}$, and call the resulting measure $\mu$. Let
$k$ be such that $\varepsilon_{N}$ is of order $2^{-k}$. It is
not hard to verify that we can have $H_{k}(\mu)=1/2$ but $|H_{k}(\mu)-H_{k}(\mu*\mu)|$
arbitrarily small. This example is actually the uniform measure on
a (generalized) arithmetic progression, as predicted by Freiman-type
theorems, but as we allow the rank $N$ to grow, the entropy growth
can be made arbitrarily small. Furthermore, if one conditions $\mu$
on an exponentially small subset of its support one gets another example
with the similar properties that is quite far from a generalized arithmetic
progression. 

Our main contribution to this matter is Theorem \ref{thm:inverse-thm-Rd}
below, which shows that constructions like the one above are, in a
certain statistical sense, the only way that \eqref{eq:mean-entropy-growth}
can occur. We note that there is a substantial existing literature
on the growth condition $|A+B|\leq|A|^{1+\delta}$, which is the sumset
analog of \eqref{eq:mean-entropy-growth}. Such a condition appears
in the sum-product theorems of Bourgain-Katz-Tao \cite{BourgainKatzTao2004}
and in the work of Katz-Tao \cite{KatzTao2001}, and in the Euclidean
setting more explicitly in Bourgain's work on the Erd\H{o}s-Volkmann
conjecture \cite{Bourgain2003} and Marstrand-like projection theorems
\cite{Bourgain2010}. However we have not found a result in the literature
that meets our needs and, in any event, we believe that the formulation
given here will find further applications.

\subsection{\label{sub:Concentration-and-saturation}Concentration and saturation
on subspaces}

We begin by discussing global properties of measures that lead to
the inequality in \eqref{eq:mean-entropy-growth}, and formulate discrete
analogs of them.

For a linear subspace $V\leq\mathbb{R}^{d}$ we say that a measure
$\mu$ is absolutely continuous on a translate $V'$ of $V$ if it
is absolutely continuous with respect to the $\dim V$-dimensional
volume (Hausdorff measure) $\lambda_{V'}$ on $V'$. Suppose that
$\mu\in\mathcal{P}(\mathbb{R}^{d})$ is compactly supported on a translate
$V_{1}$ of $V$, and is absolutely continuous there. Then the Lebesgue
differentiation theorem implies that $\mu(B_{r}(x))=c_{x}\cdot(r^{\dim V}+o(1))$
as $r\rightarrow0$, and it follows that 
\begin{equation}
H_{n}(\mu)=\dim V-o(1)\qquad\mbox{as }n\rightarrow\infty.\label{eq:33}
\end{equation}
If $\nu\in\mathcal{P}(\mathbb{R}^{d})$ is compactly supported on
another translate $V_{2}$ of $V$, then $\nu*\mu$ is supported on
$V_{3}=V_{1}+V_{2}$, which is a translate of $V$, and is absolutely
continuous there. Thus it also satisfies \eqref{eq:33}, and consequently
$H_{n}(\mu*\nu)=H_{n}(\mu)+o(1)$: i.e., at small scales there is
negligible entropy growth, and \eqref{eq:mean-entropy-growth} is
satisfied.

More generally, let $W=V^{\perp}$ be the orthogonal complement of
$V$ and write $\pi_{W}$ for the orthogonal projection to $W$. Suppose
$\mu\in\mathcal{P}(\mathbb{R}^{d})$ is compactly supported and its
conditional measures on the translates of $V$ are absolutely continuous,
that is, $\mu=\int\mu_{w}d\theta(w)$ where $\theta=\pi_{W}\mu$ and
$\mu_{w}$ is $\theta$-a.s. supported and absolutely continuous on
$\pi_{W}^{-1}(w)=V+w$. Then instead of \eqref{eq:33}, one can show
that 
\begin{equation}
H_{n}(\mu)=H_{n}(\pi_{W}(\mu))+\dim V-o(1)\qquad\mbox{as }n\rightarrow\infty,\label{eq:34}
\end{equation}
and, if $\nu$ is compactly supported on a translate of $V$, then
$\mu*\nu$ again has absolutely continuous conditional measures on
translates of $V$, and it projects to a translate of $\theta$, so
it satisfies the same relation \eqref{eq:34}. Again, we have $H_{n}(\nu*\mu)=H_{n}(\mu)+o(1)$,
and \eqref{eq:mean-entropy-growth} is satisfied.

This discussion motivates the following finite-scale analogs. For
$A\subseteq\mathbb{R}^{d}$ and $\varepsilon>0$ denote the $\varepsilon$-neighborhood
of $A$ by 
\[
A^{(\varepsilon)}=\{x\in\mathbb{R}^{d}\,:\,d(x,A)<\varepsilon\}.
\]

\begin{defn}
\label{def:concentration}Let $V\leq\mathbb{R}^{d}$ be a linear subspace
and $\varepsilon>0$. A measure $\mu\in\mathcal{P}(\mathbb{R}^{d})$
is $(V,\varepsilon)$-concentrated if there is a translate $W$ of
$V$ such that $\mu(W^{(\varepsilon)})\geq1-\varepsilon$. 
\end{defn}
Note that $(V,\varepsilon)$-concentration does not imply that the
measure is supported near $V$ itself, only near a translate of it.
Next, discretizing \eqref{eq:33} we have
\begin{defn}
\label{def:uniform}Let $V\leq\mathbb{R}^{d}$ be a linear subspace,
$\varepsilon>0$ and $m\in\mathbb{N}$. A measure $\mu\in\mathcal{P}(\mathbb{R}^{d})$
is $(V,\varepsilon)$-uniform at scale $m$, or $(V,\varepsilon,m)$-uniform,
if it is $(V,2^{-m})$-concentrated and $H_{m}(\mu)>\dim V-\varepsilon$.
\end{defn}
Finally, discretizing \eqref{eq:34}, we have:
\begin{defn}
\label{def:saturation}Let $V\leq\mathbb{R}^{d}$ be a linear subspace,
$W=V^{\perp}$ its orthogonal complement, and $\varepsilon>0$. A
probability measure $\mu\in\mathcal{P}(\mathbb{R}^{d})$ is $(V,\varepsilon)$-saturated
at scale $m$, or $(V,\varepsilon,m)$-saturated, if
\[
H_{m}(\mu)\geq H_{m}(\pi_{W}\mu)+\dim V-\varepsilon.
\]

\end{defn}
There are obvious relations between the notions above: being nearly
uniform implies saturation, and saturation implies being essentially
a convex combination of nearly uniform measures. Furthermore, as one
would expect from the discussion above, if we convolve a measure which
is highly concentrated on a subspace with another measure which is
uniform or saturated on that subspace at some scale, there will be
little entropy growth at that scale. For precise statements see Sections
\ref{sub:More-about-concentration-uniformity-continuity} and \ref{sub:saturation-and-concentration-2}.

\subsection{\label{sub:Component-measures}Component measures }

Let $\mathcal{D}_{n}(x)\in\mathcal{D}_{n}$ denote the unique level-$n$
dyadic cell containing the point $x\in\mathbb{R}^{d}$. For $D\in\mathcal{D}_{n}$
let $T_{D}:\mathbb{R}^{d}\rightarrow\mathbb{R}^{d}$ be the unique
homothety mapping $D$ to $[0,1)^{d}$. Recall that if $\mu\in\mathcal{P}(\mathbb{R})$
then $T_{D}\mu$ is the push-forward of $\mu$ through $T_{D}$ .
\begin{defn}
For $\mu\in\mathcal{P}(\mathbb{R}^{d})$ and a dyadic cell $D$ with
$\mu(D)>0$, the (raw) $D$-component of $\mu$ is
\[
\mu_{D}=\frac{1}{\mu(D)}\mu|_{D},
\]
and the (rescaled) $D$-component is 
\[
\mu^{D}=\frac{1}{\mu(D)}T_{D}(\mu|_{D}).
\]
For $x\in\mathbb{R}^{d}$ with $\mu(\mathcal{D}_{n}(x))>0$ we write
\begin{eqnarray*}
\mu_{x,n} & = & \mu_{\mathcal{D}_{n}(x)}\\
\mu^{x,n} & = & \mu^{\mathcal{D}_{n}(x)}.
\end{eqnarray*}
These measures, as $x$ ranges over all possible values for which
$\mu(\mathcal{D}_{n}(x))>0$, are called the level-$n$ components
of $\mu$.
\end{defn}
Our results on the multi-scale structure of $\mu\in\mathbb{R}^{d}$
are stated in terms of the behavior of random components of $\mu$,
defined as follows.\footnote{Definition \ref{def:component-distribution} is motivated by Furstenberg's
notion of a CP-distribution \cite{Furstenberg70,Furstenberg08,HochmanShmerkin2012},
which arise as limits as $N\rightarrow\infty$ of the distribution
of components of level $1,\ldots,N$. These limits have a useful dynamical
interpretation but in our finitary setting we do not require this
technology.}
\begin{defn}
\label{def:component-distribution}Let $\mu\in\mathcal{P}(\mathbb{R}^{d})$. 
\begin{enumerate}
\item A random level-$n$ component, raw or rescaled, is the random measure
$\mu_{D}$ or $\mu^{D}$, respectively, obtained by choosing $D\in\mathcal{D}_{n}$
with probability $\mu(D)$; equivalently, this is the random measure
$\mu_{x,n}$ or $\mu^{x,n}$, respectively, with $x$ chosen according
to $\mu$.
\item For a finite set $I\subseteq\mathbb{N}$, a random level-$I$ component,
raw or rescaled, is chosen by first choosing $n\in I$ uniformly,
and then (conditionally independently on the choice of $n$) choosing
a raw or rescaled level-$n$ component.
\end{enumerate}
\end{defn}
\begin{notation}
When the symbols $\mu^{x,i}$ and $\mu_{x,i}$ appear inside an expression
$\mathbb{P}\left(\ldots\right)$ or $\mathbb{E}\left(\ldots\right)$,
they will always denote random variables drawn according to the component
distributions defined above. The range of $i$ will be specified as
needed. When dealing with components of several measures $\mu,\nu$,
we assume all choices of components are independent unless otherwise
stated. 
\end{notation}
The definition is best understood with some examples. For $\mathcal{A},\mathcal{B}\subseteq\mathcal{P}([0,1]^{d})$,
and writing $1_{\mathcal{A}}$ for the indicator function of $\mathcal{A}$,
we have
\begin{eqnarray*}
\mathbb{P}_{i=n}\left(\mu^{x,i}\in\mathcal{A}\right) & = & \int1_{\mathcal{A}}(\mu^{x,n})\,d\mu(x)\\
\mathbb{P}_{0\leq i\leq n}\left(\mu^{x,i}\in\mathcal{A}\right) & = & \frac{1}{n+1}\sum_{i=0}^{n}\int1_{\mathcal{A}}(\mu^{x,i})\,d\mu(x)\\
\mathbb{P}_{i=n}\left(\mu^{x,i}\in\mathcal{A}\,,\,\nu^{y,i}\in\mathcal{B}\right) & = & \int\int1_{\mathcal{A}}(\mu^{x,n})\cdot1_{\mathcal{B}}(\nu^{y,n})\,d\mu(x)\,d\nu(y).
\end{eqnarray*}
This notation implicitly defines $x,i$ as random variables. Thus
if $\mathcal{A}_{0},\mathcal{A}_{1},\ldots\subseteq\mathcal{P}([0,1]^{d})$
and $D\subseteq[0,1]^{d}$ we could write 
\[
\mathbb{P}_{0\leq i\leq n}\left(\mu^{x,i}\in\mathcal{A}_{i}\mbox{ and }x\in D\right)=\frac{1}{n+1}\sum_{i=0}^{n}\mu\left(x\,:\,\mu^{x,i}\in\mathcal{A}_{i}\mbox{ and }x\in D\right).
\]
Similarly, for $f:\mathcal{P}([0,1)^{d})\rightarrow\mathbb{R}$ and
$I\subseteq\mathbb{N}$,
\[
\mathbb{E}_{i\in I}\left(f(\mu^{x,i})\right)=\frac{1}{|I|}\sum_{i\in I}\int f(\mu^{x,i})\,d\mu(x).
\]
We use similar expectation notation to average a sequence $a_{n},\ldots,a_{n+k}\in\mathbb{R}$:
\[
\mathbb{E}_{n\leq i\leq n+k}\left(a_{i}\right)=\frac{1}{k+1}\sum_{i=n}^{n+k}a_{i}.
\]
We note in particular one trivial identity that will be used repeatedly
later on:
\begin{equation}
\mu=\mathbb{E}_{i=n}\left(\mu_{x,i}\right).\label{eq:components-average-to-the-whole}
\end{equation}

Component distributions have the convenient property that they are
almost invariant under repeated sampling, i.e. choosing components
of components. More precisely, for $\mu\in\mathcal{P}(\mathbb{R}^{d})$
and $m,n\in\mathbb{N}$, let $\mathbb{P}_{n}^{\mu}$ denote the distribution
of components $\mu^{x,i}$, $0\leq i\leq n$, as defined above; and
let $\mathbb{Q}_{n,m}^{\mu}$ denote the distribution on components
obtained by first choosing a random component $\mu^{x,i}$, $0\leq1\leq n$,
as above, and then, conditionally on $\theta=\mu^{x,i}$, choosing
a component $\theta^{y,j}$, $i\leq j\leq i+m$ with the usual distribution
(note that $\theta^{y,j}=\mu^{y,j}$ is indeed a component of $\mu$). 
\begin{lem}
\label{lem:distribution-of-components-of-components}Given $\mu\in\mathcal{P}(\mathbb{R}^{d})$
and $m,n\in\mathbb{N}$, the total variation distance between $\mathbb{P}_{n}^{\mu}$
and $\mathbb{Q}_{n,m}^{\mu}$ satisfies 
\[
\left\Vert \mathbb{P}_{n}^{\mu}-\mathbb{Q}_{n,m}^{\mu}\right\Vert =O(\frac{m}{n})
\]
In particular if $\mathcal{A},\mathcal{B}\subseteq\mathcal{P}([0,1)^{d})$
and $\varepsilon,\delta>0$ are such that
\begin{eqnarray}
\mathbb{P}_{0\leq i\leq n}(\mu^{x,y}\in\mathcal{A}) & > & 1-\varepsilon\nonumber \\
\mathbb{P}_{i\leq j\leq i+m}(\theta^{y,i}\in\mathcal{B}) & > & 1-\delta\qquad\mbox{ for every }\theta\in\mathcal{A}\label{eq:106}
\end{eqnarray}
Then
\[
\mathbb{P}_{0\leq i\leq n}(\mu^{x,i}\in\mathcal{B})>1-\varepsilon-\delta-O(\frac{m}{n})
\]
\end{lem}
\begin{proof}
Observe that both $\mathbb{P}_{n}^{\mu}$ and $\mathbb{Q}_{n,m}^{\mu}$
produce a component $\mu_{z,k}$ by choosing $z$ according to $\mu$,
and independently choosing a level $k\in\mathbb{N}$. The difference
is that $\mathbb{P}_{n}^{\mu}$ chooses $k$ uniformly in the range
$0,\ldots,n$, whereas for $\mathbb{Q}_{n,m}^{\mu}$, an elementary
calculation shows that with probability $1-O(m/n)$ is choses $k$
uniformly in the range $m,m+1,\ldots,n$, and with probability $O(m/n)$
it is chooses $k\in\{0,1,,\ldots,m-1\}\cup\{n+1,\ldots,n+m\}$ (one
can easily determine the distribution in this case but it is not relevant
here). This gives the first statement.

For the second statement, what we want to show is that $\mathbb{P}_{n}^{\mu}(\mathcal{B})>1-\varepsilon-\delta-O(m/n)$.
This will follow from the first statement if we show that $\mathbb{Q}_{n,m}^{\mu}(\mathcal{B})>1-\varepsilon-\delta$.
Let $\theta=\mu^{x,i}$ and $\theta^{y,j}$ be as in the previous
paragraph, so $\theta^{y,j}$ is distributed according to $\mathbb{Q}_{n,m}^{\mu}$
. By the law of total probability and our hypotheses, 
\begin{eqnarray*}
\mathbb{Q}_{n,m}^{\mu}(\mathcal{B}) & = & \mathbb{P}(\theta^{y,j}\in\mathcal{B})\\
 & \geq & \mathbb{P}(\theta^{y,j}\in\mathcal{B}|\mu^{x,i}\in\mathcal{A})\cdot\mathbb{P}(\mu^{x,i}\in\mathcal{A})\\
 & > & (1-\delta)(1-\varepsilon)
\end{eqnarray*}
and the claim follows.
\end{proof}
Similar statements hold for raw components and components of measures
on the similarity group. We omit the proofs, which are the same.

\subsection{\label{sub:inverse-theorem}An inverse theorem for convolutions on
$\mathbb{R}^{d}$}

Our main result on entropy growth is that the global obstructions
described at the beginning of Section \ref{sub:Concentration-and-saturation}
are the only local obstructions.
\begin{thm}
\label{thm:inverse-thm-Rd}For every $R,\varepsilon>0$ and $m\in\mathbb{N}$
there is a $\delta=\delta(\varepsilon,R,m)>0$ such that for every
$n>n(\varepsilon,R,\delta,m)$, the following holds: if $\mu,\nu\in\mathcal{P}([-R,R]^{d})$
and 
\[
H_{n}(\mu*\nu)<H_{n}(\mu)+\delta,
\]
then there exists a sequence $V_{0},\ldots,V_{n}\leq\mathbb{R}^{d}$
of subspaces such that 
\begin{eqnarray}
\mathbb{P}_{_{0\leq i\leq n}}\left(\begin{array}{c}
\mu^{x,i}\mbox{ is }(V_{i},\varepsilon,m)\mbox{-saturated and}\\
\nu^{y,i}\mbox{ is }(V_{i},\varepsilon)\mbox{-concentrated}
\end{array}\right) & > & 1-\varepsilon.\label{eq:inverse-theorem-conclusion}
\end{eqnarray}

\end{thm}
The proof of the theorem is given in Section \ref{sub:Proof-of-inverse-theorem}. 
\begin{rem}
{}
\begin{enumerate}
\item The dependence of $\delta$ on $\varepsilon,m$ is effective, but
the bounds we obtain are certainly far from optimal, and we do not
pursue this topic. Also note that the theorem is not a characterization
(this is already the case in dimension 1, see discussion after \cite[Theorem 2.7]{Hochman2014}).
\item We have assumed that $\mu,\nu\in\mathcal{P}([-R,R]^{d}])$ but the
theorem can be extended to measures with unbounded support having
finite entropy by an approximation argument, see also \cite[Section 5.5]{Hochman2014}. 
\item An application of Markov's inequality shows that (up to replacing
$\varepsilon$ by $\sqrt{\varepsilon}$) equation \eqref{eq:inverse-theorem-conclusion}
is equivalent to
\begin{eqnarray}
\mathbb{P}_{_{0\leq i\leq n}}\left(\mu^{x,i}\mbox{ is }(V_{i},\varepsilon,m)\mbox{-saturated and}\right) & > & 1-\varepsilon\label{eq:inverse-thm-1st-alternative}\\
\mathbb{P}_{_{0\leq i\leq n}}\left(\nu^{y,i}\mbox{ is }(V_{i},\varepsilon)\mbox{-concentrated}\right) & > & 1-\varepsilon.\label{eq:inverse-theorem-2nd-alternative}
\end{eqnarray}

\item There is no assumption in the theorem on the entropy of $\nu$, but
if $H_{n}(\nu)$ is sufficiently close to $0$ the conclusion will
automatically hold with $V_{i}=\{0\}$ (indeed, a small value of $H_{n}(\nu)$
implies that with high probability $\nu^{y,i}$ will be highly concentrated
on $\{0\}$, so \eqref{eq:inverse-thm-1st-alternative} holds, and
\eqref{eq:inverse-theorem-2nd-alternative} is automatic, every measure
is $(\{0\},\varepsilon,m\})$-saturated).
\item The version of Theorem \ref{thm:inverse-thm-Rd} given in \cite{Hochman2014}
for the case $d=1$ had a somewhat different, but equivalent, appearance.
The statement there was that for small enough $\delta>0$, if $H_{n}(\mu*\nu)\leq H_{n}(\mu)+\delta$,
then there exist disjoint sets $I,J\subseteq\{0,\ldots,n\}$ with
$|I\cup J|>(1-\varepsilon)n$ such that \eqref{eq:inverse-thm-1st-alternative}
holds for $V_{i}=\mathbb{R}$ when the expectation is conditioned
on $i\in I$, and \eqref{eq:inverse-theorem-2nd-alternative} holds
for $V_{i}=\{0\}$ when the expectation is conditioned on $i\in J$.
Indeed, if such $I,J\subseteq\{0,\ldots,n\}$ are given, observe that
by setting $V_{i}=\mathbb{R}$ for $i\in I$ and $V_{i}=\{0\}$ for
$i\in J$, and defining $V_{i}$ arbitrarily on the at most $\varepsilon n$
remaining $i$, equations \eqref{eq:inverse-thm-1st-alternative}
and \eqref{eq:inverse-theorem-2nd-alternative} will hold for slightly
larger $\varepsilon$ also without conditioning on $I,J$, because
every measure is $(\mathbb{R},\varepsilon)$-concentrated and $(\{0\},\varepsilon,m)$-saturated.
Thus the version in \cite{Hochman2014} implies the $d=1$ case of
Theorem \ref{thm:inverse-thm-Rd}. Conversely, assuming subspaces
$V_{i}$ as in Theorem \ref{thm:inverse-thm-Rd}, we recover the version
from \cite{Hochman2014} by setting $I=\{i\,:\,V_{i}=\mathbb{R}\}$
and $J=\{j\,:\,V_{j}=\{0\}\}$ and adjusting $\varepsilon$.
\end{enumerate}
\end{rem}
Specializing to self-convolutions and using some of the basic relations
between saturation, concentration and uniformity, one deduces a multi-scale
Freiman-type result:
\begin{thm}
For every $\varepsilon>0$ and $m\in\mathbb{N}$, there is a $\delta=\delta(\varepsilon,m)>0$
such that for every $n>n(\varepsilon,\delta,m)$ and every $\mu\in\mathcal{P}([0,1)^{d})$,
if 
\[
H_{n}(\mu*\mu)<H_{n}(\mu)+\delta,
\]
then there exists a sequence $V_{0},\ldots,V_{n}<\mathbb{R}^{d}$
such that 
\begin{eqnarray*}
\mathbb{P}_{_{0\leq i\leq n}}\left(\mu^{x,i}\mbox{ is }(V_{i},\varepsilon,m)\mbox{-uniform}\right) & > & 1-\varepsilon.
\end{eqnarray*}

\end{thm}

\subsection{\label{sub:inverse-theorem-isometries}An inverse theorem for isometries
acting on $\mathbb{R}^{d}$}

Recall that $G=G(d)$ denotes the group of similarities of $\mathbb{R}^{d}$.
For $g=rU+a$ we write $r_{g}=r,U_{g}=U$ and $a_{g}=a$. The dyadic
partitions $\mathcal{D}{}_{n}^{G}$ and $\mathcal{E}_{n}^{G}$ of
$G$ were defined in Section \ref{sub:Main-results} using the identification
of $G$ with a subset of $\mathbb{R}\times M_{d}(\mathbb{R})\times\mathbb{R}^{d}$.
For $\nu\in\mathcal{P}(G)$ and for $g\in G$, $n\in\mathbb{N}$,
we define the raw component $\nu_{g,n}$ in terms of the partition
$\mathcal{D}_{n}^{G}$,
\[
\nu_{g,n}=c\cdot\nu|_{\mathcal{D}_{n}^{G}(g)},
\]
where $c$ is a normalizing constant. We adopt the same notation and
conventions for these components as laid out in Section \ref{sub:Component-measures}.

It is not natural in this context to define ``rescaled'' components.
When we need to rescale we shall do so explicitly using the maps $S_{t}\in G$,
\[
S_{t}x=2^{t}x.
\]

For $\nu\in\mathcal{P}(G)$ and $\mu\in\mathcal{P}(\mathbb{R}^{d})$
we write $\nu\conv\mu$ for the push-forward of $\nu\times\mu$ via
$(\varphi,x)\mapsto\varphi(x)$, and similarly for $x\in\mathbb{R}^{d}$
write $\nu\conv x$ for the push-forward of $\nu$ by $g\mapsto gx$.
Our aim is to understand when the entropy of $\nu\conv\mu$ is large
relative to the entropy of $\mu$, for $\nu\in\mathcal{P}(G)$ and
$\mu\in\mathcal{P}(\mathbb{R}^{d})$. 

While our methods are able to treat this setting, it is more transparent
if we assume that $\nu$ is supported on the isometry group $G_{0}<G$,
and we shall mostly restrict our attention to this case. 

The statement we would like to make is that, if $\nu\in\mathcal{P}(G_{0})$
and $\mu\in\mathcal{P}(\mathbb{R}^{d})$, and if $\nu$ is of large
entropy, then $\mu\conv\nu$ will have substantially more entropy
than $\mu$, at small enough scales, unless certain specific obstructions
occur. In the present setting the obvious global obstruction is that
$\mu$ may be close to uniform on an orbit of a subgroup $H<G_{0}$,
and $\nu$ supported on $H$ or a left coset of $H$. However, locally,
this situation is not very different from the one we have already
seen, and it is more natural to study the concentration of $\mu$
on affine subspaces, as in the Euclidean case. This is because the
orbit of a point $x\in\mathbb{R}^{d}$ under a closed subgroup $H<G_{0}$
is a finite union of smooth manifolds, and at small scales these look
like affine subspaces of $\mathbb{R}^{d}$ (essentially, the tangent
hyperplanes of the manifolds). Thus we continue to state our results
in terms of the concentration on subspaces of (the components of)
$\mu$ and (the components of) the image of $\nu$ under the action.

Even so, there are several complications related to the phenomenon
above. The first is demonstrated by the following example. Let $d=2$,
let $\mu$ be the uniform measure on the circle $\{x\in\mathbb{R}^{2}\,:\,\left\Vert x\right\Vert _{2}=1\}$,
and let $\nu$ be the uniform measure on the group of rotations about
the origin. Then $\nu\conv\mu=\mu$, so there is no entropy growth.
In this case, as predicted in the previous paragraph, the components
$\mu^{x,n}$ become saturated on lines when $n$ is large, but the
line varies according to the point $x$ (the distribution of directions
for $x\sim\mu$ is of course uniform). In contrast, recall from Theorem
\ref{thm:inverse-thm-Rd} that, for convolutions of measures on $\mathbb{R}^{d}$,
at each scale there was a single subspace on which, with high probability,
all components of $\mu$ at a given level became saturated, irrespective
of their spatial positions.

Another complication is the possibility that at small scales $\mu$
indeed becomes saturated, and $\nu$ concentrated, on subspaces, but
that these subspaces are trivial. In the Euclidean setting such an
occurrence was possible only if $\nu$ had nearly vanishing entropy,
since if $H_{n}(\nu)$ is substantial then the components of $\nu$
cannot with high probability be highly concentrated on points. In
the current setting, however, this cannot be ruled out. To see this
let $\mu=\delta_{0}$ and let $\nu$ be normalized Haar measure on
the orthogonal group $O(d)=\stab_{G_{0}}(0)$. Then $\nu\conv\mu=\mu$,
so there is no entropy growth, and $\nu$ has large entropy at all
scales, but the components of $\mu$ are not saturated on any non-trivial
subspace. Thus the theorem above applies, but $V_{i}=\{0\}$. This
type of situation can be avoided, however, if no part of the measure
$\mu$ is close to a proper affine subspace. To make this quantitative
we introduce the following definition:
\begin{defn}
\label{def:epsilon-sigma-continuity}$\mu\in\mathcal{P}(\mathbb{R}^{d})$
is $(\varepsilon,\sigma)$-non-affine if $\mu(V^{(\sigma)})<\varepsilon$
for every proper affine subspace $V\leq\mathbb{R}^{d}$. 
\end{defn}
We can now state the inverse theorem. Informally, it says that if
$\nu\conv\mu$ does not have substantially more entropy than $\mu$,
then, to most components of $\mu$ and $\nu$ at a moderately small
scale, we can associate a subspace (depending on the components in
question) such that the sub-components of the components typically
become concentrated or saturated on this subspace. Furthermore, these
subspaces will frequently be non-trivial if $\mu$ is not too close
to being supported on a proper affine subspace of $\mathbb{R}^{d}$.
Here is the precise formulation:
\begin{thm}
\label{thm:inverse-theorem-for-isometries}For every $\varepsilon>0$,
$R>0$ and $m\in\mathbb{N}$, there exists $\delta=\delta(\varepsilon,R,m)>0$
such that for every $k>k(\varepsilon,R,m)$ and every $n>n(\varepsilon,R,m,k)$,
the following holds. For every $\nu\in\mathcal{P}(G_{0})$ and $\mu\in\mathcal{P}([-R,R]^{d})$
that are supported on balls of radius $R$, either 
\[
H_{n}(\nu\conv\mu)>H_{n}(\mu)+\delta,
\]
or else, to every pair of level-$k$ components $\widetilde{\nu}$
of $\nu$ and $\widetilde{\mu}$ of $\mu$ we can assign a sequence
of subspaces $V_{i}=V_{i}(\widetilde{\nu},\widetilde{\mu})<\mathbb{R}^{d}$,
$0\leq i\leq n$, such that with probability at least $1-\varepsilon$
over the choice of $\widetilde{\mu},\widetilde{\nu}$, 
\[
\mathbb{P}_{0\leq i\leq n}\left(\begin{array}{c}
\widetilde{\mu}^{x,i}\mbox{ is }(V_{i},\varepsilon,m)\mbox{-saturated and }\\
S_{i}U_{g}^{-1}(\widetilde{\nu}_{g,i}\conv x)\mbox{ is }(V_{i},\varepsilon)\mbox{-concentrated}
\end{array}\right)>1-\varepsilon
\]
If in addition $\mu$ is $((\varepsilon/5d)^{2(d+1)},\sigma)$-non-affine
for some $\sigma>0$, and the relation among parameters takes $\sigma$
into account, then for those $\widetilde{\nu},\widetilde{\mu}$ in
the set of good components above, then for those $\widetilde{\nu},\widetilde{\mu}$
in the set of good components above, 
\[
\frac{1}{n+1}\sum_{i=0}^{n}\dim V_{i}>\frac{1}{d+1}H_{n}(\widetilde{\nu})-\varepsilon,
\]
and 
\begin{equation}
\mathbb{E}_{i=k}\left(\frac{1}{n+1}\sum_{j=0}^{n}\dim V_{j}(\nu_{g,i},\mu_{x,i})\right)>\frac{1}{d+1}H(\nu)-\varepsilon\label{eq:46}
\end{equation}
\end{thm}
\begin{rem}
{}
\begin{enumerate}
\item Given $\varepsilon$, the assumption that $\mu$ is $((\varepsilon/5d)^{2(d+1)},\sigma)$-non-affine
is global, and imposes no restriction on the structure of $\mu$ below
at scales smaller than $O(\varepsilon^{2(d+1)})$). Indeed, if $\mu\in\mathcal{P}(\mathbb{R}^{d})$
does not give mass to any affine subspace, then for any $\tau$ it
is $(\tau,\sigma)$-non-concentrated for some $\sigma>0$. Thus, if
we fix $\mu$ in advance, then for every $\varepsilon,m$ the conclusion
of the theorem holds automatically for suitable parameters $\delta,k,n$,
and all $\nu\in\mathcal{P}(G)$.
\item The average in \eqref{eq:46} is over all pairs of components $\nu_{g,k},\mu_{x,k}$,
not only those for which the first part of the conclusion holds. But
the total mass of the exceptional components is at most $\varepsilon$,
and $\dim V_{i}\leq d$, so the exceptional components contribute
$O(\varepsilon)$ to the average, which is of the same order as the
error term. Thus we get an equivalent statement if in \eqref{eq:46}
we average only over only the ``good'' components from the first
part of the theorem.
\end{enumerate}
\end{rem}
The proof of the theorem is based on a linearization argument which
allows us to apply Theorem \ref{thm:inverse-thm-Rd} from the Euclidean
setting. See Section \ref{sub:proof-of-inverse-theorems-for-G}.

\subsection{Generalizations}

It is possible to apply our methods also to convolutions in Lie groups,
actions of Lie groups on manifolds, and more general settings. Let
$I\subseteq\mathbb{R}^{d_{1}}$ and $J\subseteq\mathbb{R}^{d_{2}}$
be closed balls and $f:I\times J\rightarrow\mathbb{R}^{d}$ a $C^{1}$
map. For $z=(x,y)\in I\times J$ we can write the differential $Df(z):\mathbb{R}^{d_{1}+d_{2}}\rightarrow\mathbb{R}^{d}$
in matrix form, as
\[
Df(z)=[A_{z},B_{z}]:\mathbb{R}^{d_{1}+d_{2}}\rightarrow\mathbb{R}^{d},
\]
where $A_{z}\in M_{d\times d_{1}}$ and $B_{z}\in M_{d\times d_{2}}$. 
\begin{thm}
\label{thm:generalized-inverse-theorem}Let $f:I\times J\rightarrow\mathbb{R}^{d}$
be as above. For every $\varepsilon>0$ and $m\in\mathbb{N}$ there
exists $\delta=\delta(f,\varepsilon,m)>0$ such that for every $k>k(f,\varepsilon,m)$
and every $n>n(f,\varepsilon,m,k)$, the following holds. Let $\nu\in\mathcal{P}(I)$
and $\mu\in\mathcal{P}(J)$. Then either
\begin{equation}
H_{n}(f(\mu\times\nu))>\int H_{n}(f(\mu\times\delta_{y}))\,d\nu(y)+\delta\label{eq:entropy-growth-in-generalized-inv-thm}
\end{equation}
or else, for independently chosen level-$k$ components $\widetilde{\mu},\widetilde{\nu}$
of $\mu,\nu$, respectively, with probability at least $1-\varepsilon$
there are subspaces $V_{0},\ldots,V_{n}<\mathbb{R}^{d}$ such that
\[
\mathbb{P}_{0\leq i\leq n}\left(\begin{array}{c}
A_{x,y}\widetilde{\mu}^{x,i}\mbox{ is }(V_{i},\varepsilon,m)\mbox{-saturated and }\\
B_{x,y}\widetilde{\nu}^{y,i}\mbox{ is }(V_{i},\varepsilon)\mbox{-concentrated}
\end{array}\right)>1-\varepsilon
\]
and
\[
\frac{1}{n+1}\sum_{i=0}^{n}\dim V_{i}>c\int H_{n}(f(\delta_{x}\times\nu))\,d\widetilde{\mu}(x).
\]

\end{thm}
Note that since $I\times J$ is compact the norms of $A_{x,y}$ and
$B_{x,y}$ are bounded over $(x,y)\in I\times J$, and since $\varepsilon$
may be small and $m$ large relative to these norms, we have not bothered
to re-scale the measures $A_{x,y}\widetilde{\mu}^{x,i},B_{x,y}\widetilde{\nu}^{y,i}$
to compensate for their contraction/expansion (the distortion caused
by these matrices is also one reason for the dependence of the parameters
on $f$, the other being the speed of linear approximation). The proof
is given in Section \ref{sub:Generalizations}. 

We note two important special cases.
\begin{cor}
Let $G<GL_{d}(\mathbb{R})\subseteq\mathbb{R}^{d^{2}}$ be a matrix
group acting by left multiplication on $\mathbb{R}^{d}$. Let $\nu\in\mathcal{P}(G)$
and $\mu\in\mathcal{P}(\mathbb{R}^{d})$ be measures of bounded support.
Then for every $\varepsilon>0$ and $m\in\mathbb{N}$ there is a $\delta=\delta(\nu,\mu,\varepsilon,m)>0$,
such that for $k>k(\nu,\mu,\varepsilon,m,\delta)$ and $n>n(\nu,\mu,\varepsilon,m,\delta,k)$,
either
\[
H_{n}(\nu\conv\mu)>H_{n}(\mu)+\delta,
\]
or else, for independently chosen level-$k$ components $\widetilde{\mu},\widetilde{\nu}$
of $\mu,\nu$, respectively, with probability at least $1-\varepsilon$
there are subspaces $V_{0},\ldots,V_{n}<\mathbb{R}^{d}$ such that
\[
\mathbb{P}_{0\leq i\leq n}\left(\begin{array}{c}
y\conv\widetilde{\mu}^{x,i}\mbox{ is }(V_{i},\varepsilon,m)\mbox{-saturated and }\\
\widetilde{\nu}^{y,i}\conv x)\mbox{ is }(V_{i},\varepsilon)\mbox{-concentrated}
\end{array}\right)>1-\varepsilon
\]
(The dependence of $\delta,k,\nu$ on the measures depends only on
their support and is uniform on compact sets). 

If in addition $\mu$ is $((\varepsilon/3d)^{d+1},\sigma)$-non-affine
for some $\sigma>0$, then for $\delta,k,n$ which also depend on
$\sigma$, we also have 
\[
\frac{1}{n+1}\sum_{i=0}^{n}\dim V_{i}>c\cdot H_{n}(\widetilde{\nu})-\varepsilon.
\]
for a constant $c$ depending only on $d,\sigma$ and the support
of $\nu$.
\end{cor}

\begin{cor}
Let $G<GL_{d}(\mathbb{R})\subseteq\mathbb{R}^{d^{2}}$ be a matrix
group and $\mu,\nu\in\mathcal{P}(G)$ measures of bounded support.
Then for every $\varepsilon>0$ and $m\in\mathbb{N}$ there is a $\delta>0$
such that for every large enough $k$ and all suitably large enough
$n$, either
\[
H_{n}(\mu*\nu)>H_{n}(\mu)+\delta,
\]
or else, for an independently chosen pair of raw level-$k$ components
$\widetilde{\mu},\widetilde{\nu}$ of $\mu,\nu$, respectively, with
probability $>1-\varepsilon$, there are subspaces $V_{0},\ldots,V_{n}<\mathbb{R}^{d^{2}}$
such that
\[
\mathbb{P}_{0\leq i\leq n}\left(\begin{array}{c}
y*\widetilde{\mu}^{x,i}\mbox{ is }(V_{i},\varepsilon,m)\mbox{-saturated and }\\
\widetilde{\nu}^{y,i}*x\mbox{ is }(V_{i},\varepsilon)\mbox{-concentrated}
\end{array}\right)>1-\varepsilon
\]
and
\[
\frac{1}{n+1}\sum_{i=0}^{n}\dim V_{i}>c\cdot H_{n}(\nu)-\varepsilon.
\]
for a suitable constant $c$.
\end{cor}
Both corollaries follow from the the previous theorem by taking $f(x,y)=yx$
to be the appropriate action map; for the first corollary an additional
argument is needed to produce the constant $c$. The dependence of
the paramerets on the measures is only through their supports: If
we fix a large ball in advance and assume the measures are supported
on it, then the parameters depend only on the ball, not the measures. 
\begin{rem}
{\ }
\begin{enumerate}
\item It is important to note the order of quantifiers in the theorem and
corollaries: In the theorem all parameters depend on the function
$f$, and in the corollaries the function $f$ is the action map restricted
to the (compact) product of the supports of $\nu$ and $\mu$, which
are fixed before the other parameters. The reason this works is that
once the functions is fixed and the measures are fixed, and compactly
supported, the speed with which the function $f$ approaches its linearzation
is uniform, hence, at small enough scales, we are essentially dealing
with linear convolutions rather than a non-linear image. 
\item In some applications the order of quantifiers above is not sufficient
and it is necessary to obtain statements that are uniform over many
functions or independent of the support of the measures. Then a more
quantitative analysis is needed. Such an example can be found in \cite{HochmanSolomyak2016}.
\item One can formulate the corollaries in abstract Lie groups using partitions
introduced from local coordinates, or using general theorem on the
existence of similar partitions in doubling metric spaces, see e.g.
\cite{Kaenmaki2012}.
\item When dealing with more general group actions one would also like to
relax the condition that the measures be compactly supported. But
in doing so one must take into account how various properties of the
action affect the dependence between parameters in the theorem. For
example they are sensitive to the speed at which the action approaches
its linearization (which may differ from point to point), how well
the an element of the group is determined by its action on $k$-tuples,
and how sensitive the latter procedure is to changes in the $k$-tuple.
It turns out that the cleanest approach is to choose a left-invariant
Riemmanian metric on the group and dyadic partition adapted to it.
For a detailed development of this approach in one example we refer
the reader to \cite{HochmanSolomyak2016}.
\end{enumerate}
\end{rem}

\section{\label{sec:Entropy-concentration-uniformity-saturation}Entropy,
concentration, uniformity and saturation}

This section presents without proof some standard results about entropy,
followed by a more detailed analysis of concentration, saturation
and uniformity.

\subsection{\label{sub:Preliminaries-on-entropy}Preliminaries on entropy}

The Shannon entropy of a probability measure $\mu$ with respect to
a countable partition $\mathcal{E}$ is given by
\[
H(\mu,\mathcal{E})=-\sum_{E\in\mathcal{E}}\mu(E)\log\mu(E),
\]
where the logarithm is in base $2$ and $0\log0=0$. The conditional
entropy with respect to a countable partition $\mathcal{F}$ is
\[
H(\mu,\mathcal{E}|\mathcal{F})=\sum_{F\in\mathcal{F}}\mu(F)\cdot H(\mu_{F},\mathcal{E}),
\]
where $\mu_{F}=\frac{1}{\mu(F)}\mu|_{F}$ is the conditional measure
on $F$. For a discrete probability measure $\mu$ we write $H(\mu)$
for the entropy with respect to the partition into points, and for
a probability vector $\alpha=(\alpha_{1},\ldots,\alpha_{k})$ we write
\[
H(\alpha)=-\sum\alpha_{i}\log\alpha_{i}.
\]
and for $0<\varepsilon<1$ we abbreviate 
\[
H(\varepsilon)=H((\varepsilon,1-\varepsilon))
\]
Note that if $0<\varepsilon<1/2$ then $H(\varepsilon)=O(\varepsilon\log(1/\varepsilon))$.

We collect here some standard properties of entropy.
\begin{lem}
\label{lem:entropy-combinatorial-properties}Let $\mu,\nu$ be probability
measures on a common space, $\mathcal{E},\mathcal{F}$ partitions
of the underlying space and $\alpha\in[0,1]$.
\begin{enumerate}
\item \label{enu:entropy-positivity}$H(\mu,\mathcal{E})\geq0$, with equality
if and only if $\mu$ is supported on a single atom of $\mathcal{E}$. 
\item \label{enu:entropy-combinatorial-bound}If $\mu$ is supported on
$k$ atoms of $\mathcal{E}$ then $H(\mu,\mathcal{E})\leq\log k$,
with equality if and only if each of these atoms has mass $1/k$.
\item \label{enu:entropy-refining-partitions}If $\mathcal{F}$ refines
$\mathcal{E}$ (i.e. $\forall\;F\in\mathcal{F}\;\exists E\in\mathcal{E}\,s.t.\,F\subseteq E$)
then $H(\mu,\mathcal{F})\geq H(\mu,\mathcal{E})$.
\item \label{enu:entropy-conditional-formula}If $\mathcal{E}\lor\mathcal{F}=\{E\cap F\,:\,E\in\mathcal{E}\,,\,F\in\mathcal{F}\}$
denotes the join of $\mathcal{E}$ and $\mathcal{F}$, then 
\[
H(\mu,\mathcal{E}\lor\mathcal{F})=H(\mu,\mathcal{F})+H(\mu,\mathcal{E}|\mathcal{F}),
\]
in particular
\[
H(\mu,\mathcal{E}\lor\mathcal{F})\leq H(\mu,\mathcal{E})+H(\mu,\mathcal{F}).
\]

\item \label{enu:entropy-concavity}$H(\cdot,\mathcal{E})$ and $H(\cdot,\mathcal{E}|\mathcal{F})$
are concave.
\item \label{enu:entropy-almost-convexity}$H(\cdot,\mathcal{E})$ obeys
the ``convexity'' bound 
\[
H(\sum\alpha_{i}\mu_{i},\mathcal{E})\leq\sum\alpha_{i}H(\mu_{i},\mathcal{E})+H(\alpha).
\]
and similarly after conditioning on $\mathcal{F}$.
\end{enumerate}
\end{lem}
In particular, we note that for $\mu\in\mathcal{P}([0,1)^{d})$ we
have the bounds $H(\mu,\mathcal{D}_{m})\leq md$ and $H(\mu,\mathcal{D}_{n+m}|\mathcal{D}_{n})\leq md$.

Although the function $(\mu,m)\mapsto H(\mu,\mathcal{D}_{m})$ is
not continuous in the weak-{*} topology on measures, the following
estimates provide usable substitutes.
\begin{lem}
\label{lem:entropy-weak-continuity-properties}Let $\mu,\nu\in\mathcal{P}(\mathbb{R}^{d})$,
let $\mathcal{E},\mathcal{F}$ be partitions of $\mathbb{R}^{d}$,
and $m,m'\in\mathbb{N}$. 
\begin{enumerate}
\item \label{enu:entropy-approximation} Given a compact $K\subseteq\mathbb{R}^{d}$
and $\mu\in\mathcal{P}(K)$, there is a neighborhood $U\subseteq\mathcal{P}(K)$
of $\mu$ such that $|H(\nu,\mathcal{D}_{m})-H(\mu,\mathcal{D}_{m})|=O_{d}(1)$
for $\nu\in U$.
\item \label{enu:entropy-combinatorial-distortion} If each $E\in\mathcal{E}$
intersects at most $k$ elements of $\mathcal{F}$ and vice versa,
then $|H(\mu,\mathcal{E})-H(\mu,\mathcal{F})|=O(\log k)$.
\item \label{enu:entropy-geometric-distortion} If $f,g:\mathbb{R}^{d}\rightarrow\mathbb{R}^{k}$
and $\left\Vert f(x)-g(x)\right\Vert \leq C2^{-m}$ for $x\in\mathbb{R}^{d}$
then $|H(f\mu,\mathcal{D}_{m})-H(g\mu,\mathcal{D}_{m})|\leq O_{C,k}(1)$.
\item \label{enu:entropy-translation} If $\nu(\cdot)=\mu(\cdot+x_{0})$
then $\left|H(\mu,\mathcal{D}_{m})-H(\nu,\mathcal{D}_{m})\right|=O_{d}(1)$.
\item \label{enu:entropy-change-of-scale} If $|m'-m|\leq C$, then $\left|H(\mu,\mathcal{D}_{m})-H(\mu,\mathcal{D}_{m'})\right|\leq O_{C,d}(1)$.
\end{enumerate}
\end{lem}
We will use some easy corollaries of Lemma \ref{lem:entropy-combinatorial-properties}
\eqref{enu:entropy-concavity} and \eqref{enu:entropy-almost-convexity}. 
\begin{lem}
\label{lem:entropy-total-variation-perturbation}Let $\mu,\nu\in\mathcal{P}([-r,r]^{d})$,
let $\delta>0$, and let $\theta=(1-\delta)\mu+\delta\nu$. Then for
partitions $\mathcal{A},\mathcal{B}$ of $\mathbb{R}^{d}$ we have
\begin{eqnarray*}
|H(\theta,\mathcal{A})-H(\mu,\mathcal{A})| & \leq & H(\delta)+\delta|H(\mu,\mathcal{A})-H(\nu,\mathcal{A})|,\\
|H(\theta,\mathcal{A}|\mathcal{B})-H(\mu,\mathcal{A}|\mathcal{B})| & \leq & H(\delta)+\delta|H(\mu,\mathcal{A}|\mathcal{B})-H(\nu,\mathcal{A}|\mathcal{B})|.
\end{eqnarray*}

\end{lem}
Recall that the total variation distance between $\mu,\nu\in\mathcal{P}(\mathbb{R}^{d})$
is
\[
\left\Vert \mu-\nu\right\Vert =\sup_{A}|\mu(A)-\nu(A)|,
\]
where the supremum is over Borel sets $A$. This is a complete metric
on $\mathcal{P}(\mathbb{R}^{d})$. It follows from standard measure
theory%
{} that given $\mu,\nu$ there are probability measures $\tau,\mu',\nu'$
such that $\mu=(1-\delta)\tau+\delta\mu'$ and $\nu=(1-\delta)\tau+\delta\nu'$,
where $\delta=\frac{1}{2}\left\Vert \mu-\nu\right\Vert $. Combining
this with Lemma \ref{lem:entropy-combinatorial-properties} \eqref{enu:entropy-concavity}
and \eqref{enu:entropy-almost-convexity}, we have
\begin{lem}
\label{lem:entropy-total-variation-continuity}If $\mathcal{A},\mathcal{B}$
are partitions of $\mathbb{R}^{d}$, and if $\mu,\nu\in\mathcal{P}(\mathbb{R}^{d})$
are supported on at most $k$ atoms of each partition and $\left\Vert \mu-\nu\right\Vert <\varepsilon$,
then 
\begin{eqnarray*}
|H(\mu,\mathcal{A})-H(\nu,\mathcal{A})| & < & 2k\varepsilon+2H(\frac{1}{2}\varepsilon).\\
|H(\mu,\mathcal{A}|\mathcal{B})-H(\nu,\mathcal{A}|\mathcal{B})| & < & 2k\varepsilon+2H(\frac{1}{2}\varepsilon).
\end{eqnarray*}
In particular, if $\mu,\nu\in\mathcal{P}([0,1)^{d})$, then
\[
|H_{m}(\mu)-H_{m}(\nu)|<2d\varepsilon+\frac{2H(\frac{1}{2}\varepsilon)}{m}.
\]

\end{lem}

\subsection{\label{sub:Global-entropy-from-local-entropy}Global entropy from
local entropy}

Recall from Section \ref{sub:Component-measures} the definition of
the raw and re-scaled components $\mu_{x,n}$, $\mu^{x,n}$, and note
that 
\[
H(\mu^{x,n},\mathcal{D}_{m})=H(\mu_{x,n},\mathcal{D}_{n+m}).
\]
Also, 
\begin{eqnarray*}
\mathbb{E}_{i=n}\left(H_{m}(\mu^{x,i})\right) & = & \int\frac{1}{m}H(\mu^{x,n},\mathcal{D}_{m})\,d\mu(x)\\
 & = & \frac{1}{m}\int H(\mu_{x,n},\mathcal{D}_{n+m})\,d\mu(x)\\
 & = & \frac{1}{m}H(\mu,\mathcal{D}_{n+m}\,|\,\mathcal{D}_{n}).
\end{eqnarray*}
The following basic lemmas enable us to get bounds on the scale-$n$
entropy of a measure, or a convolution of measures, in terms of the
average scale-$m$ entropy of their components or convolution of their
components, when $m\ll n$.
\begin{lem}
\label{lem:entropy-local-to-global}For $r\geq1$, $\mu\in\mathcal{P}([-r,r]^{d})$
and integers $m<n$,
\begin{eqnarray*}
H_{n}(\mu) & = & \mathbb{E}_{0\leq i\leq n}\left(H_{m}(\mu^{x,i})\right)+O(\frac{m+\log r}{n}).
\end{eqnarray*}

\end{lem}

\begin{lem}
\label{lem:entropy-of-convolutions-via-component-convolutions}For
$r>0$, $\mu,\nu\in\mathcal{P}([-r,r]^{d})$ and integers $m<n$,
\begin{eqnarray*}
H_{n}(\mu*\nu) & \geq & \mathbb{E}_{0\leq i\leq n}\left(\frac{1}{m}H(\mu_{x,i}*\nu_{y,i},\mathcal{D}_{i+m}|\mathcal{D}_{i})\right)+O(\frac{m+\log r}{n})\\
 & \geq & \mathbb{E}_{0\leq i\leq n}\left(H_{m}(\mu^{x,i}*\nu^{y,i})\right)+O(\frac{1}{m}+\frac{m+\log r}{n}).
\end{eqnarray*}

\end{lem}
For proofs see \cite[Section 3.2]{Hochman2014}, or the proof of the
following variant, which is essentially the same as the Euclidean
case.
\begin{lem}
\label{lem:local-to-global-on-G}Let $\nu\in\mathcal{P}(G_{0})$ and
$\mu\in\mathcal{P}(\mathbb{R}^{d})$ be supported on sets of diameter
$r$. Then for $m<n$,
\[
H_{n}(\nu\conv\mu)\geq\mathbb{E}_{0\leq i\leq n}\left(H_{i,m}(\nu_{g,i}\conv\mu)\right)-O(\frac{1}{m}+\frac{m+\log r}{n}).
\]
\end{lem}
\begin{proof}
We can assume that $n=n_{0}m$, since replacing $n$ by the closest
multiple of $m$ results in a change of $O(m/n)$ to $H_{n}(\nu\conv\mu)$,
which is absorbed in the error term. Let us also introduce a parameter
$0\leq k<m$. Then
\begin{eqnarray*}
H_{n}(\nu\conv\mu) & = & \frac{1}{n}H(\nu\conv\mu,\mathcal{D}_{n})\\
 & = & \frac{1}{n}H(\nu\conv\mu,\mathcal{D}_{k+n})+O(\frac{k}{n})\\
 & = & \frac{1}{n}H(\nu\conv\mu,\mathcal{D}_{k})+\frac{1}{n}H(\nu\conv\mu,\mathcal{D}_{k+n}|\mathcal{D}_{k})+O(\frac{m}{n})
\end{eqnarray*}
Since $\nu$ is supported on a set of diameter $O(1)$ and $\mu$
on a set of diameter $O(r)$, also $\nu\conv\mu$ is supported on
a set of diameter $O(r)$, so the trivial entropy bound gives
\[
\frac{1}{n}H(\nu\conv\mu,\mathcal{D}_{k})=O(\frac{\log r+m}{n})
\]
We next evaluate $\frac{1}{n}H(\nu\conv\mu,\mathcal{D}_{k+n}|\mathcal{D}_{k})$.
Recalling our assumption $n=n_{0}m$ and the definition of conditional
entropy, we have 
\[
\frac{1}{n}H(\nu\conv\mu,\mathcal{D}_{k+n}|\mathcal{D}_{k})=\frac{1}{n}\sum_{j=0}^{n_{0}-1}H(\nu\conv\mu,\mathcal{D}_{k+(j+1)m}|\mathcal{D}_{k+jm})
\]
For each $j$ we have the identities $\nu=\mathbb{E}_{i=j}(\nu_{g,i})$
and $\mu=\mathbb{E}_{i=j}(\mu_{x,i})$, which implies $\nu\conv\mu=\mathbb{E}_{i=j}(\nu_{g,i}\conv\mu)$.
By concavity of entropy, we get
\begin{eqnarray*}
\frac{1}{n}\sum_{j=1}^{n_{0}}H(\nu\conv\mu,\mathcal{D}_{k+jm}|\mathcal{D}_{k+(j-1)m}) & = & \frac{1}{n}\sum_{j=0}^{n_{0-1}}H(\mathbb{E}_{i=k+jm}(\nu_{g,i}\conv\mu),\mathcal{D}_{k+(j+1)m}|\mathcal{D}_{k+jm})\\
 & \geq & \frac{1}{n}\sum_{j=1}^{n_{0}}\mathbb{E}_{i=k+jm}\left(H(\nu_{g,i}\conv\mu,\mathcal{D}_{k+(j+1)m}|\mathcal{D}_{k+jm})\right)\\
 & = & \frac{1}{n}\sum_{j=1}^{n_{0}}\mathbb{E}_{i=k+jm}\left(H(\nu_{g,i}\conv\mu,\mathcal{D}_{k+(j+1)m})-H(\nu_{g,i}\conv\mu,\mathcal{D}_{k+jm})\right)
\end{eqnarray*}
Since $\nu_{g,i}\conv\mu$ is supported on a set of diameter $O(2^{-i})$,
for $i=k+jm$ we have $H(\nu_{g,i}\conv\mu\mathcal{D}_{k+jm})=O(1)$.
Thus the total sum of error terms in the sum above is $O(n_{0})$,
which upon dividing by $n$ is $O(n_{0}/n)=O(1/m)$. The discussion
so far shows that 
\begin{eqnarray*}
\frac{1}{n}H(\nu\conv\mu,\mathcal{D}_{n}) & \geq & \frac{1}{n}\sum_{j=1}^{n_{0}}\mathbb{E}_{i=k+jm}\left(H(\nu_{g,i}\conv\mu,\mathcal{D}_{k+(j+1)m})\right)-O(\frac{1}{m}+\frac{m+\log r}{n})
\end{eqnarray*}
Averaging now over $k=0,\ldots,m$ gives
\begin{eqnarray*}
\frac{1}{n}H(\nu\conv\mu,\mathcal{D}_{n}) & = & \frac{1}{m}\sum_{k=0}^{m-1}\frac{1}{n}H(\nu\conv\mu,\mathcal{D}_{k+n})-O(\frac{m}{n})\\
 & \geq & \frac{1}{m}\sum_{k=0}^{m-1}\frac{1}{n}\sum_{j=1}^{n_{0}}\mathbb{E}_{i=k+jm}\left(H(\nu_{g,i}\conv\mu,\mathcal{D}_{k+(j+1)m})\right)-O(\frac{1}{n}+\frac{m+\log r}{n})\\
 & = & \frac{1}{n}\sum_{j=1}^{n}\frac{1}{m}\mathbb{E}_{i=j}\left(H(\nu_{g,i}\conv\mu,\mathcal{D}_{i+m})\right)-O(\frac{1}{n}+\frac{m+\log r}{n})\\
 & = & \mathbb{E}_{1\leq i\leq n}\left(\frac{1}{m}H(\nu_{g,i}\conv\mu,\mathcal{D}_{i+m})\right)-O(\frac{1}{n}+\frac{n+\log r+k}{n})
\end{eqnarray*}
as claimed. 
\end{proof}
We also need the following variant of Lemma \ref{lem:local-to-global-on-G}:
\begin{lem}
\label{lem:skipping-k-scales}Let $\nu\in\mathcal{P}(G_{0})$ and
$\mu\in\mathcal{P}(\mathbb{R}^{d})$ be supported on balls of diameter
$r$. Then for every $k,n\in\mathbb{N}$,
\[
H_{n}(\nu\conv\mu)\geq\mathbb{E}_{i=k}(H_{n}(\nu_{g,i}\conv\mu))+O_{R,k}(\frac{1}{n})
\]
and in particular
\[
\mathbb{E}_{i=k}(H_{n}(\nu_{g,i}\conv\mu)-H_{n}(\mu_{x,i}))\leq H_{n}(\nu\conv\mu)-H_{n}(\mu)+O_{R,k}(\frac{1}{n})
\]
\end{lem}
\begin{proof}
Since $\mu,\nu$ are supported on balls of radius $R$, so is $\nu\conv\mu$,
so the scale-$k$ entropies of all these measures is $O_{R,k}(1)$.
It follows that 
\[
H_{n}(\nu\conv\mu)=\frac{1}{n}H(\nu\conv\mu,\mathcal{D}_{n}|\mathcal{D}_{k})+O_{R,k}(\frac{1}{n})
\]
By concavity of conditional entropy,
\begin{eqnarray*}
\frac{1}{n}H(\nu\conv\mu,\mathcal{D}_{n}|\mathcal{D}_{k}) & = & \frac{1}{n}H(\mathbb{E}_{i=k}(\nu_{g,i}\conv\mu),\mathcal{D}_{n}|\mathcal{D}_{k})\\
 & \geq & \mathbb{E}_{i=k}(\frac{1}{n}H(\nu_{g,i}\conv\mu,\mathcal{D}_{n}|\mathcal{D}_{k}))
\end{eqnarray*}
But $\nu_{g,i}\conv\mu$ is supported on a set of diameter $O(2^{-i})$,
so (taking $i=k$),
\begin{eqnarray*}
\frac{1}{n}H(\nu_{g,k}\conv\mu,\mathcal{D}_{n}|\mathcal{D}_{k}) & = & \frac{1}{n}H(\nu_{g,k}\conv\mu,\mathcal{D}_{n})+O(\frac{1}{n})\\
 & = & H_{n}(\nu_{g,k}\conv\mu)+O(\frac{1}{n})
\end{eqnarray*}
Combining the last three equations gives the first claim. For the
second claim, note that we have
\begin{eqnarray*}
H_{n}(\mu) & = & \frac{1}{n}H(\mu,\mathcal{D}_{n}|\mathcal{D}_{k})+O_{R,k}(\frac{1}{n})\\
 & = & \frac{1}{n}\mathbb{E}_{i=k}(H(\mu_{x,i},\mathcal{D}_{n}))+O_{R,k}(\frac{1}{n})\\
 & = & \mathbb{E}_{i=k}(H_{n}(\mu_{x,i}))+O_{R,k}(\frac{1}{n})
\end{eqnarray*}
(the first inequality again because $\mu$ is supported on a set of
diameter $O(R)$). Subtracting this expression for $H_{n}(\mu)$ from
the previous one for $H_{n}(\nu\conv\mu)$ gives the claim.
\end{proof}

\subsection{\label{sub:More-about-concentration-uniformity-continuity}First
lemmas on concentration, uniformity, saturation }

We consider here some basic connections between uniform, concentrated
and saturated measures. We make the statements as general as possible,
but in some cases, especially when dealing with uniform measures,
it is necessary to assume that the support of the measures is bounded,
and the constants in the error terms may depend on the diameter of
the support. Since we are interested in the asymptotics as $m\rightarrow\infty$
we rarely make the dependence explicit, but it can be read off from
the proofs. 

Given partitions $\mathcal{E}$ and $\mathcal{F}$ of sets $X,Y$,
respectively, write 
\[
\mathcal{E}\otimes\mathcal{F}=\{E\times F\,:\,E\in\mathcal{E}\,,\,F\in\mathcal{F}\}
\]
for the product partition of $X\times Y$. We will also often identify
$\mathcal{E}$ with the partition $\mathcal{E}\otimes\{Y\}$ of $X\times Y$,
and similarly $\mathcal{F}$ with the partition $\{X\}\otimes\mathcal{F}$
of $X\times Y$.

For a linear subspace $V\leq\mathbb{R}^{d}$ we write $\mathcal{D}_{n}^{V}$
for the level-$n$ dyadic partition on $V$ with respect to some fixed
(but arbitrary) orthogonal coordinate system in $V$, which we usually
do not specify.

Let $V\leq\mathbb{R}^{d}$ be a linear subspace and $W=V^{\perp}$,
and let $\mathcal{D}'_{m}=\mathcal{D}_{m}^{V}\otimes\mathcal{D}_{m}^{W}$
denote the product partition of $\mathbb{R}^{d}\cong V\times W$.
Each element of $\mathcal{D}_{m}$ intersects at most $O(1)$ elements
of $\mathcal{D}'_{m}$, and vice versa, so by Lemma \ref{lem:entropy-weak-continuity-properties}
\eqref{enu:entropy-combinatorial-distortion}, 
\[
|H(\mu,\mathcal{D}_{m})-H(\mu,\mathcal{D}'_{m})|=O(1).
\]
The same is true for the induced partitions on $W$, so, writing $\pi_{W}$
for the orthogonal projection to $W$,
\[
|H(\pi_{W}\mu,\mathcal{D}_{m})-H(\pi_{W}\mu,\mathcal{D}'_{m})|=O(1)
\]
and also
\[
|H(\pi_{W}\mu,\mathcal{D}_{m})-H(\pi_{W}\mu,\mathcal{D}_{m}^{W})|=O(1).
\]
Recall that we identify $\mathcal{D}_{m}^{V},\mathcal{D}_{m}^{W}$
with the partitions $\pi_{V}^{-1}\mathcal{D}_{m}^{V}$, $\pi_{W}^{-1}\mathcal{D}_{m}^{W}$
of $\mathbb{R}^{d}$, respectively. With this identification we have
$\mathcal{D}'_{m}=\mathcal{D}_{m}^{V}\lor\mathcal{D}_{m}^{W}$, and
\[
H(\pi_{W}\mu,\mathcal{D}'_{m})=H(\mu,\mathcal{D}_{m}^{W}).
\]
From this discussion we have the following immediate consequence:
\begin{lem}
\label{lem:saturation-under-coordinate-change}With the above notation,
a measure $\mu\in\mathcal{P}(\mathbb{R}^{d})$ is $(V,\varepsilon+O(1/m),m)$-saturated
if and only if 
\[
\frac{1}{m}H(\mu,\mathcal{D}_{m}^{V}|\mathcal{D}_{m}^{V^{\perp}})\geq\dim V-(\varepsilon+O(\frac{1}{m})).
\]

\end{lem}
From similar considerations we have
\begin{lem}
\label{lem:saturation-transformed-by-similarity}If $\mu\in\mathcal{P}(\mathbb{R}^{d})$
is $(V,\varepsilon,m)$-saturated and $g=2^{t}U+a\in G$ is a similarity,
then $g\mu$ is $(UV,\varepsilon+O(|t|/m),m)$-saturated; and similarly
for uniformity. 
\end{lem}
One way to get saturated measures is from uniform measures:
\begin{lem}
\label{lem:uniform-measures-are-saturated}If $\mu\in\mathcal{P}([-r,r]^{d})$
is $(V,\varepsilon,m)$-uniform then it is $(V,O_{r}(\varepsilon+1/m),m)$-saturated.\end{lem}
\begin{proof}
By uniformity, we can write $\mu=(1-\varepsilon)\mu'+\varepsilon\mu''$
, where $\mu'$ is supported on the $2^{-m}$-neighborhood of a translate
of $V$. By concavity of conditional entropy, 
\begin{eqnarray*}
H(\mu,\mathcal{D}_{m}|\mathcal{D}_{m}^{V^{\perp}}) & \geq & (1-\varepsilon)H(\mu',\mathcal{D}_{m}|\mathcal{D}_{m}^{V^{\perp}})\\
 & \geq & H(\mu',\mathcal{D}_{m}|\mathcal{D}_{m}^{V^{\perp}})-\varepsilon H(\mu',\mathcal{D}_{m}).
\end{eqnarray*}
Since $\mu$, and hence $\mu'$, is supported on at most $O(r^{d}\cdot2^{m})$
atoms of $\mathcal{D}_{m}$, we have $H(\mu',\mathcal{D}_{m})=O(m\log r)$,
and the inequality above becomes 
\[
H(\mu,\mathcal{D}_{m}|\mathcal{D}_{m}^{V^{\perp}})\geq H(\mu',\mathcal{D}_{m}|\mathcal{D}_{m}^{V^{\perp}})-\varepsilon O(m\log r).
\]
Since $\mu'$ is supported on a $2^{-m}$-neighborhood of a translate
of $V$, it is supported on $O(1)$ atoms of $\mathcal{D}_{m}^{V^{\perp}}$,
so $H(\mu',\mathcal{D}_{m}^{V^{\perp}})=O(1)$, hence 
\begin{eqnarray*}
H(\mu',\mathcal{D}_{m}|\mathcal{D}_{m}^{V^{\perp}}) & \geq & H(\mu',\mathcal{D}_{m})-H(\mu',\mathcal{D}_{m}^{V^{\perp}})\\
 & \geq & H(\mu',\mathcal{D}_{m})-O(1).
\end{eqnarray*}
Finally, by Lemma \ref{lem:entropy-total-variation-perturbation}
applied to $\mu=(1-\varepsilon)\mu+\varepsilon\mu''$, and using the
bound $O(r^{d}2^{m})$ on the number of $\mathcal{D}_{m}$-atoms supporting
$\mu',\mu''$ and uniformity of $\mu$, 
\begin{eqnarray*}
H(\mu',\mathcal{D}_{m}) & > & H(\mu)-\varepsilon(m+O(\log r))-H(\varepsilon)\\
 & > & m\dim V-\varepsilon(m+\log r)-H(\varepsilon)
\end{eqnarray*}
Putting it all together, using $H(\varepsilon)\leq1$ and dividing
by $m$ gives the claim.
\end{proof}
Another way to get saturated measures is to take convex combinations
of saturated measures:
\begin{lem}
\label{lem:saturation-passes-to-convex-combinations}A convex combination
of $(V,\varepsilon,m)$-saturated measures on $\mathbb{R}^{d}$ is
$(V,\varepsilon+O(1/m),m)$-saturated.\end{lem}
\begin{proof}
Immediate from Lemma \ref{lem:saturation-under-coordinate-change}
and concavity of conditional entropy (Lemma \ref{lem:entropy-combinatorial-properties}
\eqref{enu:entropy-concavity}).
\end{proof}
Combining the two lemmas above gives the following:
\begin{cor}
\label{cor:saturation-from-uniform-measures}A convex combination
of $(V,\varepsilon,m)$-uniform measures on $[-r,r]^{d}$ is $(V,O_{r}(\varepsilon+1/m),m)$-saturated.
\end{cor}
Saturation is also stable under small perturbations in the total variation
metric:
\begin{lem}
\label{lem:saturation--of-total-variation-perturbations}Let $\mu,\nu\in\mathcal{P}([-r,r]^{d})$.
If $\mu$ is $(V,\varepsilon,m)$-saturated and $\left\Vert \mu-\nu\right\Vert <\delta$
then $\nu$ is $(V,\varepsilon',m)$-saturated for $\varepsilon'=\varepsilon+O(\delta\log r+1/m)$.\end{lem}
\begin{proof}
Take $\mathcal{\mathcal{A}}=\mathcal{D}_{m}^{V}\lor\mathcal{D}_{m}^{V^{\perp}}$
and $\mathcal{\mathcal{B}}=\mathcal{D}_{m}^{V^{\perp}}$ in Lemma
\ref{lem:entropy-total-variation-continuity}, and use Lemma \ref{lem:saturation-under-coordinate-change}.
\end{proof}
Finally, we shall need an entropy bound for concentrated measures.
\begin{lem}
\label{lem:entropy-of-V-conventrated-measures}If $\mu\in\mathcal{P}([-r,r]^{d})$
is $(V,2^{-m})$-concentrated then $H_{m}(\mu)\leq\dim V+O_{r}(\frac{\log m}{m})$.\end{lem}
\begin{proof}
Write $\mu=(1-2^{-m})\mu_{1}+2^{-m}\mu_{2}$ where $\mu_{1}\in\mathcal{P}(W^{(2^{-m})})$
for some translate $W$ of $V$ and $\mu_{2}\in\mathcal{P}([-r,r]^{d})$.
Since $H_{m}(\mu_{i})=O_{r}(1)$ for $i=1,2$, by Lemma \ref{lem:entropy-total-variation-perturbation}
it suffices for us to show that $H_{m}(\mu_{1})\leq\dim V+O_{r}(1/m)$.
This again follows from the fact that $W^{(2^{-m})}\cap[-r,r]^{d}$
intersects $O(r^{d}2^{m})$ atoms of $\mathcal{D}_{m}$ and the trivial
entropy bound.
\end{proof}

\subsection{Concentration and saturation of components}

In this section all measures are supported on $[0,1)^{d}$.
\begin{lem}
\label{lem:saturation-passes-to-components}If $\mu\in\mathcal{P}([0,1)^{d})$
is $(V,\varepsilon,n)$-saturated, then for every $1\leq m<n$, 
\[
\mathbb{P}_{0\leq i\leq n}\left(\mu^{x,i}\mbox{ is }(V,\varepsilon',m)\mbox{-saturated}\right)>1-\varepsilon',
\]
where \textup{$\varepsilon'=\sqrt{d\varepsilon+O(\frac{m}{n})}$.}\end{lem}
\begin{proof}
Without loss of generality, we may assume that $\mathcal{D}_{n}=\mathcal{D}_{n}^{V}\lor\mathcal{D}_{n}^{W}$
where $W=V^{\perp}$ (Lemma \ref{lem:saturation-under-coordinate-change}).
By the fact that $\mu$ is $(V,\varepsilon,n)$-saturated and by Lemma
\ref{lem:entropy-local-to-global}, we have
\begin{multline*}
\dim V+H_{n}(\pi_{W}\mu)-\varepsilon\leq\\
\begin{aligned}\leq & H_{n}(\mu)\\
= & \mathbb{E}_{0\leq i\leq n}\left(H_{m}(\mu^{x,i})\right)+O(\frac{m}{n})\\
= & \mathbb{E}_{0\leq i\leq n}\left(H_{m}(\mu^{x,i},\mathcal{D}_{m}^{W})\right)+\mathbb{E}_{0\leq i\leq n}\left(\frac{1}{m}H\left(\mu^{x,i},\mathcal{D}_{m}|\mathcal{D}_{m}^{W}\right)\right)+O(\frac{m}{n})\\
= & \mathbb{E}_{0\leq i\leq n}\left(H_{m}(\pi_{W}(\mu^{x,i}))\right)+\mathbb{E}_{0\leq i\leq n}\left(\frac{1}{m}H\left(\mu^{x,i},\mathcal{D}_{m}|\mathcal{D}_{m}^{W}\right)\right)+O(\frac{m}{n}).
\end{aligned}
\end{multline*}
Since $(\pi_{W}\mu)_{y,i}$ is the convex combination (using the natural
weights) of $\pi_{W}(\mu_{D})$ over those $D\in\mathcal{D}_{i}$
with $D\cap\pi_{W}^{-1}(y)\neq\emptyset$ (recall that we are assuming
$\mathcal{D}_{n}=\mathcal{D}_{n}^{V}\lor\mathcal{D}_{n}^{W}$), concavity
of entropy implies 
\begin{eqnarray*}
H_{n}(\pi_{W}\mu) & = & \mathbb{E}_{0\leq i\leq n}\left(H_{m}((\pi_{W}\mu)^{y,i})\right)+O(\frac{m}{n})\\
 & \geq & \mathbb{E}_{0\leq i\leq n}\left(H_{m}(\pi_{W}(\mu^{x,i})\right)+O(\frac{m}{n}).
\end{eqnarray*}
Combining these we have
\[
\mathbb{E}_{0\leq i\leq n}\left(\frac{1}{m}H\left(\mu^{x,i},\mathcal{D}_{m}|\mathcal{D}_{m}^{W}\right)\right)\geq\dim V-(\varepsilon+O(\frac{m}{n})).
\]
But we also have the trivial bound $\frac{1}{m}H(\mu^{x,i},\mathcal{D}_{m}|\mathcal{D}_{m}^{W})\leq\dim V\leq d$.
Combining this and the last inequality, the lemma follows by Markov's
inequality.
\end{proof}
The analogous statement for concentration is valid at individual scales
(rather than for typical scales between $0$ and $n$, as above):
\begin{lem}
\label{lem:concentrated-measures-have-concentrated-components}If
$\mu\in\mathcal{P}([0,1)^{d})$ is $(V,\varepsilon)$-concentrated
and $1\leq m\leq\log(1/\varepsilon)$, then
\[
\mathbb{P}_{i=m}\left(\mu^{x,i}\mbox{ is }(V,\sqrt{2^{m}\varepsilon})\mbox{-concentrated}\right)>1-\sqrt{2^{-m}\varepsilon}.
\]
\end{lem}
\begin{proof}
Let $W=V+v_{0}$ be such that $\mu(W^{(\varepsilon)})>1-\varepsilon$.
For a dyadic cube $D$ write $T_{D}$ for the surjective homothety
$D\rightarrow[0,1)^{d}$ and let $W^{D}=T_{D}(W)$. Clearly, for any
$D\in\mathcal{D}_{m}$ we have $T_{D}(W^{(\varepsilon)})=(W^{D})^{(2^{m}\varepsilon)}$.
Take $\delta=\sqrt{2^{m}\varepsilon}\leq1$ and let $\mathcal{E}\subseteq\mathcal{D}_{m}$
denote the family of cells $D$ such that 
\[
\mu_{D}(D\setminus W^{(\varepsilon)})=\mu^{D}([0,1]^{d}\setminus(W^{D})^{(2^{m}\varepsilon)})>\delta.
\]
It follows that 
\[
\varepsilon\geq\mu([0,1]^{d}\setminus W^{(\varepsilon)})\geq\sum_{D\in\mathcal{E}}\mu(D\setminus W^{(\varepsilon)})>\delta\cdot\mu(\cup\mathcal{E}),
\]
so $\mu(\cup\mathcal{E})<\varepsilon/\delta=\sqrt{2^{-m}\varepsilon}$.
Hence $\mu(\cup(\mathcal{D}_{m}\setminus\mathcal{E}))>1-\sqrt{2^{-m}\varepsilon}$,
and the conclusion follows.
\end{proof}
We often will want to change the scale at which measures are saturated.
Clearly if $\delta<\varepsilon$ and $\mu$ is $(V,\delta)$-concentrated,
then it is also $(V,\varepsilon)$-concentrated. However for $\delta<\varepsilon$
and $k>m$ it is in general not true that if $\mu$ is $(V,\delta,k)$-saturated
then $\mu$ is also $(V,\varepsilon,m)$-saturated (though of course
it certainly is $(V,\varepsilon,k)$-saturated). The issue is that
the first few scales do not greatly affect the entropy at a fine scale.
In order to allow such change of parameters we will pass to components,
using the lemmas above. We will also need a simple covering argument
for intervals of $\mathbb{Z}$:
\begin{lem}
\label{lem:covering-by-intervals}Let $I\subseteq\{0,\ldots,n\}$
and $m\in\mathbb{N}$ be given. Then there is a subset $I'\subseteq I$
such that $I\subseteq I'+[0,m]$ and $[i,i+m]\cap[j,j+m]=\emptyset$
for distinct $i,j\in I'$.\end{lem}
\begin{proof}
Define $I'$ inductively. Begin with $I'=\emptyset$ and, at each
successive stage, if $I\setminus\bigcup_{i\in I'}[i,i+m]\neq\emptyset$
then add its least element to $I'$. Stop when $I\subseteq\bigcup_{i\in I'}[i,i+m]$. \end{proof}
\begin{prop}
\label{prop:concentration-and-saturation-change-of-parameters}For
every $\varepsilon>0$ and $m\in\mathbb{N}$, if $k>k(\varepsilon,m)$
and $0<\delta<\delta(\varepsilon,m,k)$, then for all large enough
$n>n(\varepsilon,m,k,\delta)$, the following holds. Let $\nu,\mu\in\mathcal{P}(\mathbb{R}^{d})$
and let $V_{0},V_{1},\ldots,V_{n}\leq\mathbb{R}^{d}$ be linear subspaces
such that 
\begin{equation}
\mathbb{P}_{0\leq i\leq n}\left(\begin{array}{c}
\mu^{x,i}\mbox{ is }(V_{i},\delta,k)\mbox{-saturated and}\\
\nu^{y,i}\mbox{ is }(V_{i},\delta)\mbox{-concentrated}
\end{array}\right)>1-\delta.\label{eq:100a}
\end{equation}
Then there are linear subspaces $V'_{0},\ldots,V'_{n}\leq\mathbb{R}^{d}$
such that 
\begin{equation}
\mathbb{P}_{0\leq i\leq n}\left(\begin{array}{c}
\mu^{x,i}\mbox{ is }(V'_{i},\varepsilon,m)\mbox{-saturated and}\\
\nu^{y,i}\mbox{ is }(V'_{i},\varepsilon)\mbox{-concentrated}
\end{array}\right)>1-\varepsilon.\label{eq:100b}
\end{equation}
Furthermore if $V_{i}=V$ is independent of $i$ then we can take
$V'_{i}=V$.\end{prop}
\begin{proof}
Fix $\delta,k$ and suppose that \eqref{eq:100a} holds for some $n$.
Let $I\subseteq\{0,\ldots,n\}$ denote the set of indices $u$ such
that
\[
\mathbb{P}_{i=u}\left(\begin{array}{c}
\mu^{x,i}\mbox{ is }(V_{i},\delta,k)\mbox{-saturated and}\\
\nu^{y,i}\mbox{ is }(V_{i},\delta)\mbox{-concentrated}
\end{array}\right)>1-\sqrt{\delta}.
\]
By Markov's inequality, 
\[
|I|\geq(1-\sqrt{\delta})(n+1)
\]
Let $I'\subseteq I$ be chosen as in the previous lemma with parameter
$k$, so $I\subseteq I'+[0,k]$ and $[i,i+k]\cap[j,j+k]=\emptyset$
for distinct $i,j\in I'$. If $j=i+u$ for some $i\in I'$ and $0\leq u\leq k$,
define $V'_{j}=V_{i}$. Define $V'_{j}$ arbitrarily for other $j$.
Note that when $V_{i}=V$ is independent of $i$ then also $V'_{j}=V$
for $j$ as above, in which case we can set $V'_{i}=V$ for all $i$
and satisfy the last assertion in the statement.

To see that this choice works (assuming the parameters satisfy the
proper relations), note that for any pair of components $\theta=\mu^{x,i},\eta=\nu^{y,i}$
such that $\theta$ is $(V_{i},\delta,k)$-saturated and $\eta$ is
$(V_{i},\delta)$-concentrated, we have by Lemmas \ref{lem:saturation-passes-to-components}
and \ref{lem:concentrated-measures-have-concentrated-components}
that 
\begin{eqnarray*}
\mathbb{P}_{i\leq j\leq i+k}(\theta^{w,j}\mbox{ is }(V'_{j},\sqrt{d\delta+O(\frac{m}{k})},m)\mbox{-saturated}) & > & 1-\sqrt{d\delta+O(\frac{m}{k})}\\
\mathbb{P}_{i\leq j\leq i+k}(\eta^{z,j}\mbox{ is }(V'_{j},\sqrt{2^{k}\delta})\mbox{-concentrated}) & > & 1-\sqrt{2^{-k}\delta}.
\end{eqnarray*}
so
\[
\mathbb{P}_{i\leq j\leq i+k}\left(\begin{array}{c}
\theta^{w,j}\mbox{ is }(V'_{j},\sqrt{d\delta+O(\frac{m}{k})},m)\mbox{-saturated and}\\
\eta^{z,j}\mbox{ is }(V'_{j},\sqrt{2^{k}\delta})\mbox{-concentrated}
\end{array}\right)>1-O(\sqrt{\delta+\frac{m}{k}}).
\]
Write $U=\bigcup_{i\in I'}[i,i+k]$. The union is disjoint by assumption,
so the bounds above combine to give
\[
\mathbb{P}_{i\in U}\left(\begin{array}{c}
\theta^{w,i}\mbox{ is }(V'_{j},\sqrt{d\delta+O(\frac{m}{k})},m)\mbox{-saturated and}\\
\eta^{zij}\mbox{ is }(V'_{j},\sqrt{2^{k}\delta})\mbox{-concentrated}
\end{array}\right)>1-O(\sqrt{\delta+\frac{m}{k}}).
\]
Let $V=U\cap[0,n]$. Then we have the trivial inequalities 
\begin{eqnarray*}
\mathbb{P}_{i\in V}(\ldots) & \geq & \mathbb{P}_{i\in U}(\ldots)-\frac{|U\setminus V|}{|U|}\\
\mathbb{P}_{0\leq i\leq n}(\ldots) & \geq & \frac{|V|}{n+1}\mathbb{P}_{i\in U}(\ldots).
\end{eqnarray*}
Since $I\subseteq U\subseteq[0,n+k]$ to we have $|U\setminus V|\leq k$
and $|U|\geq(1-\sqrt{\delta})(n+1)$, so combining the identities
above with the previous inequality we get 
\[
\mathbb{P}_{0\leq i\leq n}\left(\begin{array}{c}
\theta^{w,i}\mbox{ is }(V'_{j},\sqrt{d\delta+O(\frac{m}{k})},m)\mbox{-saturated and}\\
\eta^{zij}\mbox{ is }(V'_{j},\sqrt{2^{k}\delta})\mbox{-concentrated}
\end{array}\right)>1-O(\sqrt{\delta+\frac{m}{k}})-O(\frac{k}{n}).
\]
Thus if $k$ is large enough relative to $\varepsilon,m$; $\delta$
is small enough relative to $\varepsilon,k$; and $n$ is large enough
relative to $\varepsilon,k$, we obtain \eqref{eq:100b}.
\end{proof}
We remark that the use of Lemma \ref{lem:covering-by-intervals} and
Markov's inequality in the proof is rather crude, and one might want
to use Lemma \ref{lem:distribution-of-components-of-components} instead.
This would have shown that one can associate to most components a
subspace on which it is suitably concentrated and saturated, but the
subspaces would generally depend on the component, and not just on
the level it belongs to. The argument above gives the desired uniformity
across each level.

\subsection{The space of subspaces}

Let $B_{r}(x)$ denote the open Euclidean ball of radius $r$ around
$x\in\mathbb{R}^{d}$, and, as before, for $A\subseteq\mathbb{R}^{d}$
let $A^{(\varepsilon)}=\{x\in\mathbb{R}^{d}\,:\,d(x,A)<\varepsilon\}$.
Define a metric on the space of linear subspaces $V,W\leq\mathbb{R}^{d}$
by 
\begin{equation}
d(V,W)=\inf\{\varepsilon>0\,:\,V\cap B_{1}(0)\subseteq W^{(\varepsilon)}\mbox{ and }W\cap B_{1}(0)\subseteq V^{(\varepsilon)}\}\label{eq:metric-on-subspaces}
\end{equation}
This is just the Hausdorff metric on the intersections of $V,W$ with
the closed unit ball, so the induced topology on the space of linear
subspaces of $\mathbb{R}^{d}$ is compact (note that it decomposes
into $d+1$ connected components, corresponding to the dimensions
of the subspaces). It is also the same as the distance $\left\Vert \pi_{V}-\pi_{W}\right\Vert $,
where $\left\Vert \cdot\right\Vert $ denotes the operator norm and
$\pi_{V},\pi_{W}$ the orthogonal projections.

It will be convenient to write 
\[
A\sqsubseteq A'\qquad\mbox{if}\qquad A\cap B_{1}(0)\subseteq A'.
\]
This a transitive, reflexive relation. In this notation, the distance
between subspaces $V,W\leq\mathbb{R}^{d}$ defined above is 
\[
d(V_{1},V_{2})=\inf\{\varepsilon>0\,:\,V_{1}\sqsubseteq V_{2}^{(\varepsilon)}\mbox{ and }V_{2}\sqsubseteq V_{1}^{(\varepsilon)}\}.
\]
Define the ``angle'' between subspaces $V_{1},V_{2}$ by $\angle(V_{1},V_{2})=0$
if $V_{1}\subseteq V_{2}$ or $V_{2}\subseteq V_{1}$; otherwise set
$W=V_{1}\cap V_{2}$ and 
\[
\angle(V_{1},V_{2})=\inf\{\left\Vert v_{1}-v_{2}\right\Vert \,:\,v_{1}\in V_{1}\cap W^{\perp}\;,\;v_{2}\in V_{2}\cap W^{\perp}\;,\;\left\Vert v_{1}\right\Vert =\left\Vert v_{2}\right\Vert =1\}.
\]
This is not the usual notion of angle, but it agrees with the usual
definition up to a multiplicative constant, and is more convenient
to work with.

The following properties are elementary and we omit their proof.
\begin{lem}
\label{lem:subspace-lemmas}Let $V,W\leq\mathbb{R}^{d}$ be linear
subspaces and $\varepsilon>0$.
\begin{enumerate}
\item $d(V,W)\leq1$ with equality if and only if $V\cap W^{\perp}\not=\{0\}$
or $W\cap V^{\perp}\neq\{0\}$. In particular if $\dim W>\dim V$
then $W\not\sqsubseteq V^{(1)}$ and $d(V,W)=1$.
\item If $0<\varepsilon<1$ and $V\sqsubseteq W^{(\varepsilon)}$ then $\pi_{W}:V\rightarrow W$
is injective, $\dim V\leq\dim W$, and if $\dim V=\dim W$ then $W\sqsubseteq V^{(\varepsilon)}$
and $d(V,W)\leq\varepsilon$.
\item $\angle(V,W)\leq\sqrt{2}\cdot d(V,W)$.
\item \label{enu:subspace-lemmas-4}If $V\not\sqsubseteq W^{(\varepsilon)}$
then there exists a vector $v\in V$ with $\angle(\mathbb{R}v,W)\geq\varepsilon$. 
\end{enumerate}
\end{lem}
We collect some elementary implications for concentration, uniformity
and saturation:
\begin{lem}
\label{lem:concentration-saturation-uniformity-under-change-of-subspace}Let
$\mu\in\mathcal{P}([0,1)^{d})$ and $V,W\leq\mathbb{R}^{d}$.
\begin{enumerate}
\item \label{enu:concentration-passes-to-nearby-subspaces}If $\mu$ is
$(V,\varepsilon)$ concentrated and $d(W,V)<\delta$, then $\mu$
is $(W,\varepsilon+\sqrt{d}\delta)$-concentrated.
\item \label{enu:uniformity-passes-to-nearby-subspaces}If $\mu$ is $(V,\varepsilon,m+1)$-uniform
and $d(W,V)<\frac{1}{\sqrt{d}}2^{-(m+1)}$, then $\mu$ is $(W,\varepsilon,m)$-uniform.
\item \label{enu:saturation-passes-to-nearby-subspaces}If $\mu$ is $(V,\varepsilon,m)$-saturated
and $d(W,V)<2^{-m}$, then $\mu$ is $(W,\varepsilon+O(1/m),m)$-saturated.
\item \label{enu:saturation-subspaces}If $\mu$ is $(V,\varepsilon,m)$-saturated
and $W\leq V$ is a subspace then $\mu$ is $(W,\varepsilon+O(1/m),m)$-saturated.
\item \label{enu:saturation-direct-sum}If $\mu\in\mathcal{P}([0,1)^{d})$
is both $(V_{1},\varepsilon,m)$, and $(V_{2},\varepsilon,m)$-saturated,
and $\angle(V_{1},V_{2})>\delta>0$, then $\mu$ is $(V_{1}+V_{2},\varepsilon',m)$-saturated,
where $\varepsilon'=2\varepsilon+O(\frac{1}{m}\log(\frac{1}{\delta}))$.
\end{enumerate}
\end{lem}
\begin{proof}
If $d(W,V)<\delta$ then $V\cap B_{1}(0)\subseteq W^{(\delta)}\cap B_{1}(0)$,
so $V\cap B_{\sqrt{d}}(0)\subseteq W^{(\sqrt{d}\delta)}\cap B_{\sqrt{d}}(0)$.
It follows that if $(V+v)\cap[0,1)^{d}\neq\emptyset$ then $V+v\cap[0,1)^{d}\subseteq(W+v)^{(\sqrt{d}\cdot\delta)}$
(we use the fact that the diameter of $[0,1)^{d}$ is $\sqrt{d}$),
so $(V^{(\varepsilon)}+v)\cap[0,1)^{d}\subseteq(W+v)^{(\varepsilon+\sqrt{d}\delta)}$.
The first claim follows.

For (2), observe that if $d(W,V)<2^{-(m+1)}$ and $\mu$ is $(V,2^{-(m+1)}/\sqrt{d})$-concentrated,
then by the first claim, $\mu$ is $(W,2^{-m})$-concentrated. Since
by assumption $H_{m}(\mu,\mathcal{D}_{m})>\dim V-\varepsilon$, and
$d(V,W)<2^{-(m+1)}$ implies $\dim W=\dim V$, we have shown that
$\mu$ is $(V,\varepsilon,m)$-uniform.

For (3), note that $d(V,W)<2^{-m}$ implies that $\left\Vert \pi_{V^{\perp}}-\pi_{W^{\perp}}\right\Vert <2^{-m}$,
so $|H(\mu,\mathcal{D}_{m}|\mathcal{D}_{m}^{V^{\perp}})-H(\mu,\mathcal{D}_{m}|\mathcal{D}_{m}^{W^{\perp}})|=O(1)$,
and the claim follows.

For (4), we may assume $W\neq V$. Let $W'<V$ denote the orthogonal
complement of $W$ in $V$ and write $\mathbb{R}^{d}$ as the orthogonal
direct sum $W\oplus W'\oplus V^{\perp}$. Without loss of generality
we may assume $\mathcal{D}_{m}=\mathcal{D}_{m}^{W}\lor\mathcal{D}_{m}^{W'}\lor\mathcal{D}_{m}^{V^{\perp}}$
(Lemma \ref{lem:saturation-under-coordinate-change}); by doing so
we implicitly increased $\varepsilon$ by $O(1/m)$. Since $\mu$
is $(V,\varepsilon,m)$-saturated,
\[
\frac{1}{m}H(\mu,\mathcal{D}_{m}|\mathcal{D}_{m}^{V^{\perp}})\geq\dim V-\varepsilon.
\]
Since $\mathcal{D}_{m}$ refines $\mathcal{D}_{m}^{W^{\perp}}=\mathcal{D}_{m}^{W'\oplus V^{\perp}}$
which in turn refines $\mathcal{D}_{m}^{V^{\perp}}$, we have 
\begin{eqnarray*}
H(\mu,\mathcal{D}_{m}|D_{m}^{V^{\perp}}) & = & H(\mu,\mathcal{D}_{m}\lor\mathcal{D}_{m}^{W'}|\mathcal{D}_{m}^{V^{\perp}})\\
 & = & H(\mu,\mathcal{D}_{m}^{W'}|\mathcal{D}_{m}^{V^{\perp}})+H(\mu,\mathcal{D}_{m}|\mathcal{D}_{m}^{W^{\perp}}).
\end{eqnarray*}
Inserting this into the inequality above gives 
\[
\frac{1}{m}H(\mu,\mathcal{D}_{m}|\mathcal{D}_{m}^{W^{\perp}})\geq\dim V-H(\mu,\mathcal{D}_{m}^{W'}|\mathcal{D}_{m}^{V^{\perp}})-\varepsilon.
\]
Since $\frac{1}{m}H(\mu,\mathcal{D}_{m}^{W'}|\mathcal{D}_{m}^{V^{\perp}})\leq\dim W'+O(1/m)=\dim V-\dim W+O(1/m)$,
this is precisely $(W,\varepsilon+O(1/m),m)$-saturation of $\mu$.

We turn to (5). Let $V'_{2}=V_{2}\cap(V_{1}\cap V_{2})^{\perp}$,
so that $V'_{2}<V_{2}$, $V_{1}\cap V'_{2}=\{0\}$, $\angle(V_{1},V'_{2})=\angle(V_{1},V_{2})>\delta$
and $V_{1}+V'_{2}=V_{1}+V_{2}$. By (4) we can replace $V_{2}$ by
$V'_{2}$ at the cost of increasing $\varepsilon$ by $O(1/m)$. Thus,
we may assume from the start that $V_{1}\cap V_{2}=\{0\}$. 

Write $V=V_{1}\oplus V_{2}$ (this is an algebraic, not an orthogonal,
sum) and $W=V^{\perp}$. We can assume without loss of generality
that $\mathcal{D}_{m}=\mathcal{D}_{m}^{V}\lor\mathcal{D}_{m}^{W}$.
Also let $\mathcal{E}_{m}=\mathcal{D}_{m}^{V_{1}}\lor\mathcal{D}_{m}^{V_{2}}\lor\mathcal{D}_{m}^{W}$
be the partition corresponding to the direct sum $\mathbb{R}^{d}=V_{1}\oplus V_{2}\oplus W$. 

By Lemma \ref{lem:saturation-under-coordinate-change}, we must show
that 
\[
H_{m}(\mu,\mathcal{D}_{m}^{V}|\mathcal{D}_{m}^{W})\geq\dim V-2\varepsilon-O(\frac{\log(1/\delta)}{m}).
\]
Because of the assumption $\angle(V_{1},V_{2})>\delta$, the partitions
of $\mathcal{D}_{m}^{V_{1}}\lor\mathcal{D}_{m}^{V_{2}}$ and $\mathcal{D}_{m}^{V}$
of $V$, and also the corresponding partitions of $\mathbb{R}^{d}$,
have the property that each atom of one intersects $O(1/\delta)$
atoms of the other. Thus 
\[
\left|H_{m}(\mu,\mathcal{D}_{m}^{V_{1}}\lor\mathcal{D}_{m}^{V_{2}}|\mathcal{D}_{m}^{W})-H(\mu,\mathcal{D}_{m}^{V}|\mathcal{D}_{m}^{W})\right|=O(\log(1/\delta)),
\]
so it is sufficient for us to prove that
\begin{equation}
H_{m}(\mu,\mathcal{D}_{m}^{V_{1}}\lor\mathcal{D}_{m}^{V_{2}}|\mathcal{D}_{m}^{W})\geq\dim V-2\varepsilon-O(\frac{\log(1/\delta)}{m}).\label{eq:74}
\end{equation}
Now, 
\begin{equation}
\frac{1}{m}H(\mu,\mathcal{D}_{m}^{V_{1}}\lor\mathcal{D}_{m}^{V_{2}}|\mathcal{D}_{m}^{W})=\frac{1}{m}H(\mu,\mathcal{D}_{m}^{V_{1}}|\mathcal{D}_{m}^{W})+\frac{1}{m}H(\mu,\mathcal{D}_{m}^{V_{2}}|\mathcal{D}_{m}^{V_{1}}\lor\mathcal{D}_{m}^{W}).\label{eq:75}
\end{equation}
Since $W\subseteq V_{1}^{\perp}$, we can assume that the partition
$\mathcal{D}_{m}^{V_{1}^{\perp}}$ refines $\mathcal{D}_{m}^{W}$.
Using the fact that $\mu$ is $(V_{1},\varepsilon,m_{1})$-saturated,
we get a bound for the first term on the right hand side of the above
identity:
\[
\frac{1}{m}H(\mu,\mathcal{D}_{m}^{V_{1}}|\mathcal{D}_{m}^{W})\geq\frac{1}{m}H(\mu,\mathcal{D}_{m}^{V_{1}}|\mathcal{D}_{m}^{V_{1}^{\perp}})\geq\dim V_{1}-\varepsilon.
\]
As for the second term, again using the fact that each atom of $\mathcal{D}_{m}^{V_{1}}\lor\mathcal{D}_{m}^{V_{2}}\lor\mathcal{D}_{m}^{W}$
intersects $O(1/\delta)$ atoms of $\mathcal{D}_{m}^{V_{2}}\lor\mathcal{D}_{m}^{V_{2}^{\perp}}$
and vice versa, and similarly for $\mathcal{D}_{m}^{V_{1}}\lor\mathcal{D}_{m}^{W}$
and $\mathcal{D}_{m}^{V_{2}^{\perp}}$, we have 
\begin{eqnarray*}
\frac{1}{m}H(\mu,\mathcal{D}_{m}^{V_{2}}|\mathcal{D}_{m}^{V_{1}}\lor\mathcal{D}_{m}^{W}) & = & \frac{1}{m}H(\mu,\mathcal{D}_{m}^{V_{2}}|\mathcal{D}^{V_{2}^{\perp}})-O(\frac{\log(1/\delta)}{m})\\
 & \geq & \dim V_{2}-\varepsilon-O(\frac{\log(1/\delta)}{m}).
\end{eqnarray*}
Combining the last two inequalities and \eqref{eq:75} gives the desired
inequality \eqref{eq:74}
\end{proof}

\subsection{\label{sub:Geometry-of-subspaces}Geometry of thickened subspaces }

In this section we develop some methods for understanding unions and
intersections of thickened subspaces. We require some elementary linear
algebra estimates.
\begin{lem}
Let $v_{1},\ldots,v_{k}\in\mathbb{R}^{d}$ with $\left\Vert v_{i}\right\Vert \leq1$,
and suppose that 
\[
d(v_{i},\spn\{v_{1},\ldots,v_{i-1}\})>\delta\qquad\mbox{ for all }1\leq i\leq k.
\]
Then for any $v=\sum t_{i}v_{i}$ we have $\left\Vert (t_{1},\ldots,t_{k})\right\Vert \leq\sqrt{k}\cdot2^{k}\left\Vert v\right\Vert /\delta^{k}$.\end{lem}
\begin{proof}
We first claim for every $1\leq i\leq k$ that
\[
|t_{i}|\leq(1+\frac{1}{\delta})^{k-i+1}\left\Vert v\right\Vert 
\]
This we show by induction on $k$. For $k=1$ it is trivial. In general
set $V_{i}=\spn\{v_{1},\ldots,v_{i}\}$ and $W_{i}=V_{i}^{\perp}$.
By hypothesis, $\left\Vert \pi_{W_{i-1}}(v_{i})\right\Vert >\delta$
for all $1\leq i\leq k$. Thus 
\[
\left\Vert v\right\Vert \geq\left\Vert \pi_{W_{k-1}}(v)\right\Vert =|t_{k}|\cdot\left\Vert \pi_{W_{k-1}}(v_{k})\right\Vert >t_{k}\delta
\]
so $|t_{k}|<\left\Vert v\right\Vert /\delta$, and the claim holds
for $i=k$. Now, 
\begin{eqnarray*}
\left\Vert \sum_{i=1}^{k-1}t_{i}v_{i}\right\Vert  & = & \left\Vert \pi_{V_{k-1}}(v)-t_{k}\pi_{V_{k-1}}(v_{k})\right\Vert \\
 & \leq & \left\Vert v\right\Vert +|t_{k}|\\
 & \leq & (1+\frac{1}{\delta})\left\Vert v\right\Vert 
\end{eqnarray*}
Thus, by the induction hypothesis, for $1\leq i\leq k-1$,
\begin{eqnarray*}
|t_{i}| & \leq & (1+\frac{1}{\delta})^{(k-1)-i+1}\left\Vert \sum_{i=1}^{k-1}t_{i}v_{i}\right\Vert \\
 & \leq & (1+\frac{1}{\delta})^{k-i+1}\left\Vert v\right\Vert 
\end{eqnarray*}
as claimed. It remains to note that
\[
\left\Vert t\right\Vert \leq\sqrt{k}\left\Vert t\right\Vert _{\infty}\leq\sqrt{k}(1+\frac{1}{\delta})^{k}\left\Vert v\right\Vert 
\]
The claim follows (note that $\delta<1$).
\end{proof}
It will be convenient to introduce notation for the constant 
\[
p_{k}=k\cdot2^{k}
\]

\begin{cor}
\label{cor:span-is-close-to-V-if-basis-is-close-to-V}Let $V\leq\mathbb{R}^{d}$
be a subspace and $v_{1},\ldots,v_{k}\in V^{(\varepsilon)}$ with
$\left\Vert v_{i}\right\Vert \leq1$. If $d(v_{i},\spn\{v_{1},\ldots,v_{i-1}\})>\delta$
for all $1\leq i\leq k$, then $\spn\{v_{1},\ldots,v_{k}\}\sqsubseteq V^{(p_{k}\varepsilon/\delta^{k})}$.\end{cor}
\begin{proof}
Write $W=\spn\{v_{i}\}$ and let $w=\sum t_{i}v_{i}\in W$ be a unit
vector. Write $t=(t_{1},\ldots,t_{k})$. Then by the last lemma, $\left\Vert t\right\Vert _{2}\leq2^{k}\sqrt{k}/\delta^{k}$.
Thus
\begin{eqnarray*}
d(w,V) & \leq & \sum|t_{i}|\cdot d(v_{i},V)\\
 & < & \varepsilon\cdot\left\Vert t\right\Vert _{1}\\
 & \leq & \varepsilon\cdot\sqrt{k}\cdot\left\Vert t\right\Vert _{2}\\
 & \leq & p_{k}\cdot\frac{\varepsilon}{\delta^{k}},
\end{eqnarray*}
where we used the hypothesis and the general inequality $\left\Vert u\right\Vert _{1}\leq\sqrt{k}\left\Vert u\right\Vert _{2}$.\end{proof}
\begin{cor}
Suppose that $E,V\leq\mathbb{R}^{d}$ are subspaces such that $E\sqsubseteq V^{(\varepsilon)}$,
and $e\in V^{(\varepsilon)}\cap B_{1}(0)$ is such that $d(e,E)>\delta>0$.
Then $E'=E\oplus\mathbb{R}e$ satisfies $E'\sqsubseteq V^{(8\varepsilon/\delta^{2})}$.\end{cor}
\begin{proof}
Every vector in $E'$ belongs to a subspace of the form $\mathbb{R}e'\oplus\mathbb{R}e$
for some $e'\in E$, so it is enough to show $\mathbb{R}e'\oplus\mathbb{R}e\sqsubseteq V^{(8\varepsilon/\delta^{2})}$.
But the pair $e',e$ satisfies the assumptions of the previous corollary
with $k=2$. Since $p_{2}=8$, the claim follows.\end{proof}
\begin{cor}
\label{cor:sintersection-of-thick-subsapces}Suppose that $E,V,W\leq\mathbb{R}^{d}$
are subspaces such that $E\sqsubseteq V^{(\varepsilon)}\cap W^{(\varepsilon)}$,
and $e\in(V^{(\varepsilon)}\cap W^{(\varepsilon)})\cap B_{1}(0)$
is such that $d(e,E)>\delta>0$. Let $E'=E\oplus\mathbb{R}e$. Then
$E'\sqsubseteq V^{(8\varepsilon/\delta^{2})}\cap W^{(8\varepsilon/\delta^{2})}$.\end{cor}
\begin{proof}
Immediate from the lemma. 
\end{proof}
Proposition \ref{prop:minimal-engulfing-subspace} below takes a family
$\mathcal{W}$ of subspaces and finds an essentially minimal subspace
that almost-contains all $W\in\mathcal{W}$. The basic step in the
proof is to do this for two subspaces, and this is given by the next
corollary.
\begin{cor}
\label{cor:intersection-contains-subspace}Given $\varepsilon>0$
let $\varepsilon_{k}=4\varepsilon^{1/3^{k}}$. Then for any $V,W\leq\mathbb{R}^{d}$,
there is a $0\leq k\leq d$ and a $k$-dimensional subspace $E\sqsubseteq V^{(\varepsilon_{k})}\cap W^{(\varepsilon_{k})}$
such that $V^{(\varepsilon_{k})}\cap W^{(\varepsilon_{k})}\sqsubseteq E^{(\varepsilon_{k+1})}$.\end{cor}
\begin{proof}
Let $E$ be a subspace of maximal dimension satisfying $E\sqsubseteq V^{(\varepsilon_{\dim E})}\cap W^{(\varepsilon_{\dim E})}$
(such subspaces exist, e.g. $\{0\}$). Let $k=\dim E$. If $V^{(\varepsilon_{k})}\cap W^{(\varepsilon_{k})}\not\sqsubseteq E^{(\varepsilon_{k+1})}$
then by the previous corollary we can replace $E$ by $E'=E+\mathbb{R}e$
for some $e\in\partial B_{1}(0)\cap(V^{(\varepsilon_{k})}\cap W^{(\varepsilon_{k})}\setminus E^{(\varepsilon_{k+1})})$
and $E'$ will satisfy 
\[
E'\sqsubseteq V^{(8\varepsilon_{k}/\varepsilon_{k+1}^{2})}\cap W^{(8\varepsilon_{k}/\varepsilon_{k+1}^{2})}\subseteq V^{(\varepsilon_{k+1})}\cap W^{(\varepsilon_{k+1})}
\]
where we have used 
\[
\frac{8\varepsilon_{k}}{\varepsilon_{k+1}^{2}}=\frac{8\cdot4\varepsilon^{1/3^{k}}}{4^{2}\varepsilon^{2/3^{k+1}}}=2\cdot\varepsilon^{1/3^{k+1}}=\frac{1}{2}\varepsilon_{k+1}
\]
But $\dim E'=\dim E+1$, which contradicts the maximality of $E$. \end{proof}
\begin{prop}
\label{prop:minimal-engulfing-subspace}Let $\varepsilon>0$ and $\varepsilon_{k}=4\varepsilon^{1/3^{k}}$.
Then for any family $\mathcal{W}$ of subspaces of $\mathbb{R}^{d}$,
there is a subspace $V\leq\mathbb{R}^{d}$ such that $W\sqsubseteq V^{(\varepsilon_{d})}$
for all $W\in\mathcal{W}$, and if $\widetilde{V}$ is a subspace
such that $W\sqsubseteq\widetilde{V}^{(\varepsilon)}$ for all $W\in\mathcal{W}$,
then $V\sqsubseteq\widetilde{V}^{(\varepsilon_{d})}$. \end{prop}
\begin{proof}
We may assume that $\varepsilon_{d}<1$ since otherwise the statement
is trivial (any subspace $V$ will do). Let $V$ be a subspace of
minimal dimension such that $W\sqsubseteq V^{(\varepsilon_{d-\dim V})}$
for all $W\in\mathcal{W}$ (such subspaces exist, e.g. $V=\mathbb{R}^{d}$).
Write $k=d-\dim V$. We can assume $k<d$ since the case $k=d$ corresponds
to $V=\{0\}$, and then the conclusion is trivial. 

We claim that $V$ is the desired subspace. First, $\varepsilon_{k}\leq\varepsilon_{d}$,
so we have $W\sqsubseteq V^{(\varepsilon_{k})}\sqsubseteq V^{(\varepsilon_{d})}$
for all $W\in\mathcal{W}$, which is the first property.

For the second property of $V$, suppose that there is a subspace
$\widetilde{V}\leq\mathbb{R}^{d}$ such that $W\sqsubseteq\widetilde{V}^{(\varepsilon)}\sqsubseteq\widetilde{V}^{(\varepsilon_{k})}$
for $W\in\mathcal{W}$, but such that $V\not\sqsubseteq\widetilde{V}^{(\varepsilon_{d})}$.
Let $E$ be a subspace of maximal dimension satisfying  $E\sqsubseteq V^{(\varepsilon_{k+1})}\cap\widetilde{V}^{(\varepsilon_{k+1})}$.
Clearly $\dim E\leq\dim V$ (since $E\sqsubseteq V^{(\varepsilon_{k+1})}$
and $\varepsilon_{k+1}\leq\varepsilon_{d}<1$), and we cannot have
$\dim E=\dim V$ because then we would have $V\sqsubseteq E^{(\varepsilon_{k})}\sqsubseteq\widetilde{V}^{(\varepsilon_{k}+\varepsilon_{k+1})}\sqsubseteq\widetilde{V}^{(\varepsilon_{d})}$,
contrary to assumption. So $\dim E<\dim V$. Thus, by the definition
of $V$, there exists a $W\in\mathcal{W}$ with $W\not\sqsubseteq E^{(\varepsilon_{k+1})}$.
Choose a vector $e\in(B_{1}(0)\cap W)\setminus E^{(\varepsilon_{k+1})}$,
so that $d(e,E)\geq\varepsilon_{k+1}$, and note that since $W\sqsubseteq V^{(\varepsilon_{k})}\cap\widetilde{V}^{(\varepsilon_{k})}$
we also have $e\in V^{(\varepsilon_{k})}\cap\widetilde{V}^{(\varepsilon_{k})}$.
Thus, by Corollary \ref{cor:sintersection-of-thick-subsapces} (with
$\varepsilon_{k}$ and $\varepsilon_{k+1}$ in the role of $\varepsilon,\delta$),
the subspace $E'=E\oplus\mathbb{R}e$ satisfies $E'\sqsubseteq V^{(\varepsilon_{k+1})}\cap\widetilde{V}^{(\varepsilon_{k+1})}$.
But $\dim E'>\dim E$, which contradicts the definition of $E$. We
conclude that $V\subseteq\widetilde{V}^{(\varepsilon_{d})}$, as desired.
\end{proof}
We note that the proof actually shows $V\sqsubseteq W^{(\varepsilon_{d-\dim V})}$
for all $W\in\mathcal{W}$ and that any $\widetilde{V}$ with this
property satisfies $V\sqsubseteq\widetilde{V}^{(\varepsilon_{d-\dim V+1})}$.

From the last proposition we can derive a dual version: for any family
$\mathcal{W}$ of subspaces and any $\varepsilon>0$, there is a subspace
$V$ such that $V\sqsubseteq W^{(\varepsilon_{d})}$ for all $W\in\mathcal{W}$
and any other subspace $\widetilde{V}$ with this property satisfies
$\widetilde{V}\sqsubseteq V^{(\varepsilon_{d})}$. To see this, observe
that $U_{1}\sqsubseteq U_{2}^{(\varepsilon)}$ if and only if $U_{2}^{\perp}\sqsubseteq(U_{2}^{\perp})^{(\varepsilon)}$,
and apply the previous proposition to $\mathcal{W}^{\perp}=\{W^{\perp}\,:\,W\in\mathcal{W}\}$.
However, in a later application we will want to present the subspace
$V$ as an intersection of a small number of (neighborhoods of) subspaces
from $\mathcal{W}$. This is provided for in the following proposition.
\begin{prop}
\label{prop:maximal-common-subspace}Let $\varepsilon>0$ and $\delta=8^{d-1}\varepsilon^{1/3^{d^{2}}}$.
Then for any family $\mathcal{W}$ of subspaces of $\mathbb{R}^{d}$,
there is a subspace $V\leq\mathbb{R}^{d}$ such that $V\sqsubseteq W^{(\delta)}$
for every $W\in\mathcal{W}$, and subspaces $W_{1},\ldots,W_{k}\in\mathcal{W}$
with $k\leq d-\dim V$ such that $\bigcap_{i=1}^{k}W_{i}^{(\varepsilon)}\sqsubseteq V^{(\delta)}$.
In particular, if $V'$ is any other subspace satisfying $V'\sqsubseteq W^{(\varepsilon)}$
for every $W\in\mathcal{W}$, then $V'\sqsubseteq V^{(\delta)}$. 

Furthermore, if we are given an increasing sequence $\mathcal{W}^{1}\subseteq\mathcal{W}^{2}\subseteq\ldots$
with each $\mathcal{W}^{i}$ a family of subspaces of $\mathbb{R}^{d}$,
then we can assign $V^{i}$ to $\mathcal{W}^{i}$ as above in such
a way that $V^{i+1}\sqsubseteq(V^{i})^{(\delta)}$.\end{prop}
\begin{proof}
Fix $\varepsilon,\delta,\mathcal{W}$ as in the statement. We shall
recursively choose finite sequences of subspaces $W_{1},W_{2},\ldots\in\mathcal{W}$
and $V_{0},V_{1},\ldots\leq\mathbb{R}^{d}$, and of real numbers $\delta_{0},\delta_{1},\ldots>0$,
such that $\bigcap_{j=1}^{i}W_{j}^{(\varepsilon)}\sqsubseteq V_{i}^{(\delta_{i})}$.

Begin with $V_{0}=\mathbb{R}^{d}$ and $\delta_{0}=\varepsilon$.
Now for $j\geq1$ suppose we have defined $V_{i},W_{i},\delta_{i}$
for $i<j$. Let $\delta_{j}^{*}=8(\delta_{j-1})^{1/3^{d}}$. If $V_{j-1}\sqsubseteq W^{(\delta_{j}^{*})}$
for all $W\in\mathcal{W}$, we terminate the construction. Otherwise,
choose $W_{j}\in\mathcal{W}$ such that $V_{j-1}\not\sqsubseteq W_{j}^{(\delta_{j}^{*})}$.
Apply Corollary \ref{cor:intersection-contains-subspace} to the subspaces
$V_{j-1}$, $W_{j}$ with the parameter $\delta_{j-1}$. We obtain
a subspace $V_{j}\leq\mathbb{R}^{d}$ and real numbers $\delta_{j-1}\leq\delta'_{j}\leq\delta_{j}\leq4(\delta_{j-1})^{1/3^{d}}$
satisfying 
\begin{equation}
V_{j}\sqsubseteq V_{j-1}^{(\delta'_{j})}\cap W_{j}^{(\delta'_{j})}\label{eq:120}
\end{equation}
and 
\[
V_{j-1}^{(\delta'_{j})}\cap W_{j}^{(\delta'_{j})}\sqsubseteq V_{j}^{(\delta_{j})}
\]
(in the notation of the corollary, $\delta'_{j}=\varepsilon_{k}$
and $\delta_{j}=\varepsilon_{k+1}$, but if $k=d$ we can take $\delta'_{j}=\delta_{j}=\varepsilon_{d}$).
Since $\varepsilon\leq\delta_{j-1}\leq\delta'_{j}$ and, by the induction
hypothesis, $\bigcap_{i=0}^{j-1}W_{j}^{(\varepsilon)}\sqsubseteq V_{j-1}^{(\delta_{j-1})}$,
the last equation implies that $\bigcap_{i=0}^{j}W_{j}^{(\varepsilon)}\sqsubseteq V_{j}^{(\delta_{j})}$,
and the conditions of the construction are satisfied.

We now claim that $\dim V_{j}<\dim V_{j-1}$ as long as they are defined.
Indeed, suppose the construction completed the $j$-th step of the
construction without terminating, so $V_{j-1}\not\sqsubseteq W_{j}^{(\delta_{j}^{*})}$.
In particular this means that $\delta_{j}\leq\delta_{j}^{*}<1$. Now,
we know that $V_{j}\sqsubseteq V_{j-1}^{(\delta_{j})}$, which together
with $\delta_{j}<1$ implies $\dim V_{j}\leq\dim V_{j-1}$. Suppose
that equality held. Then, again using $\delta_{j}<1$, we would have
the reverse containment $V_{j-1}\sqsubseteq V_{j}^{(\delta_{j})}$.
This, together with $V_{j}\sqsubseteq W_{j}^{(\delta'_{j})}$ and
$\delta'_{j}\leq\delta_{j}$, implies $V_{j-1}\sqsubseteq W_{j}^{(2\delta_{j})}$.
Since $2\delta_{j}\leq\delta_{j}^{*}$, this contradicts the assumption
$V_{j-1}\not\sqsubseteq W_{j}^{(\delta_{j}^{*})}$, so we must have
$\dim V_{j}<\dim V_{j-1}$.

Since $\dim V_{j}$ is strictly decreasing, the procedure terminates
after completing some $k\leq d$ iterations, which in our numbering
means it completed step $k-1$ and terminated at step $k$. This means
that $V_{k-1}\sqsubseteq W{}^{(\delta_{k}^{*})}$ for all $W\in\mathcal{W}$
and $\bigcap_{i=0}^{k-1}W_{i}^{(\varepsilon)}\sqsubseteq V_{k}^{(\delta_{k-1})}$.
Observe that 
\[
\delta_{k-1}\leq\delta_{k-1}^{*}<8^{k-1}(\delta_{0})^{1/(k-1)d}\leq\delta
\]
(since $\delta_{0}=\varepsilon$). Hence for $V=V_{k}$ we have $W\sqsubseteq V^{(\delta)}$
for all $W\in\mathcal{W}$ and $\bigcup_{i=1}^{k}W_{i}^{(\varepsilon)}\sqsubseteq V^{(\delta)}$,
as desired.

The statement about $V'$ is immediate from the first statement of
the lemma.

Finally, for the last part, we note that in the construction we may
first exhaust the subspaces in $\mathcal{W}^{1}$, obtaining $V^{1}$,
then move on to those in $\mathcal{W}^{2}$ obtaining possibly a different
$V^{2}$, etc. The containment relation follows from \eqref{eq:120}.
\end{proof}

\subsection{\label{sub:saturation-and-concentration-2}Minimally concentrated
and maximally saturated subspaces }

Our goal in this section is to identify, given a measure and associated
parameters, a subspace $V$ on which it is in a sense most concentrated,
and one on which it is most saturated, relative to the parameters.
By this we mean that if $\widetilde{V}$ is another subspace for which
the measure is concentrated or saturated, relative to comparable parameters,
then $\widetilde{V}$ is, respectively, essentially contained in,
or essentially contains, $V$.

The existence of a ``minimal'' subspace on which a given measure
concentrates is proved by a variation on the argument in Proposition
\ref{prop:minimal-engulfing-subspace}:
\begin{prop}
\label{prop:concentration-subspace-of-a-measure}Let $\varepsilon>0$
and $\varepsilon_{k}=4\varepsilon^{1/3^{k}}$, and assume that $\varepsilon_{d}<1/2$.
Then for any $\eta\in\mathcal{P}([0,1)^{d})$, there is a subspace
$V\leq\mathbb{R}^{d}$ such that $\eta$ is $(V,\sqrt{d}\cdot\varepsilon_{d})$-concentrated,
and if $W$ is any subspace such that $\eta$ is $(W,\varepsilon)$-concentrated,
then $V\sqsubseteq W^{(\varepsilon_{d})}$.\end{prop}
\begin{proof}
We can assume $\varepsilon_{d}<1$. Choose a subspace $V\leq\mathbb{R}^{d}$
of minimal dimension such that $\eta$ is $(V,\sqrt{d}\cdot\varepsilon_{d-\dim V})$-concentrated
(the family of such subspaces is non-empty, e.g. $V=\mathbb{R}^{d}$).
Write $k=d-\dim V$ note that we can assume $k<d$ since otherwise
$V=\{0\}$ and the claim is trivial.

We claim that $V$ is the desired subspace. Suppose that $\eta$ is
$(W,\varepsilon)$-concentrated (and hence $(W,\varepsilon_{k})$-concentrated)
but that $V\not\sqsubseteq W^{(\varepsilon_{d})}$. Let $E\sqsubseteq V^{(\varepsilon_{k+1})}\cap W^{(\varepsilon_{k+1})}$
be a subspace of maximal dimension. Then $\dim E<\dim V$ so by the
definition of $V$ the measure $\eta$ is not $(E,\sqrt{d}\cdot\varepsilon_{k+1})$-concentrated.
Now, consider translates of $V^{(\varepsilon_{k})}+v$ and $W^{(\varepsilon_{k})}+w$
which cover all but $\varepsilon_{k}$ and $\varepsilon$ of the mass
of $\eta$, respectively. Choose $u\in(V^{(\varepsilon_{k})}+v)\cap(W^{(\varepsilon_{k})}+w)$
(the intersection is non-empty because it has $\eta$-mass at least
$1-2\varepsilon_{k}>0$), and observe that $V^{(\varepsilon_{k})}+v\subseteq V^{(2\varepsilon_{k})}+u$
and $W^{(\varepsilon_{k})}+w\subseteq W^{(2\varepsilon_{k})}+u$.
Hence 
\[
\eta\left([0,1]^{d}\cap(V^{(2\varepsilon_{k})}\cap W^{(2\varepsilon_{k})}+u\right)>1-2\varepsilon_{k}.
\]
Now consider the translate $[0,1]^{d}\cap(E^{(\sqrt{d}\varepsilon_{k})}+u)$.
It covers at most $1-\varepsilon_{k+1}$ of the mass of $\eta$, which,
since since $2\varepsilon_{k}<\varepsilon_{k+1}$, is less than the
mass of the previous intersection. Thus, translating back to the origin
and scaling by $1/\sqrt{d}$ (so that $[0,1]^{d}+u$ is mapped into
the unit ball), we find that there exists a point $e\in(B_{1}(0)\cap V^{(2\varepsilon_{k}/\sqrt{d})}\cap W^{(2\varepsilon_{k}/\sqrt{d})})\setminus E^{(\varepsilon_{k})}$.
By Corollary \ref{cor:sintersection-of-thick-subsapces} the subspace
$E'=E+\mathbb{R}e$ satisfies $E'\sqsubseteq V^{(\varepsilon_{k+1})}\cap W^{(\varepsilon_{k+1})}$
(we have used that $8\cdot2\varepsilon_{k}/(\sqrt{d}\varepsilon_{k+1}^{2})<\varepsilon_{k+1}$).
But $\dim E'>\dim E$, contradicting the choice of $E$. We conclude
therefore $V\sqsubseteq W^{(\varepsilon_{d})}$, as desired.
\end{proof}
We turn to the analog of Proposition \ref{prop:concentration-subspace-of-a-measure},
which provides a ``maximal'' subspace on which a given measure is
saturated to a certain degree. The argument is again similar to the
measureless case.
\begin{prop}
\label{prop:saturation-subspace-of-a-measure}Given $m\in\mathbb{N}$
and $\theta\in\mathcal{P}([0,1)^{d}d)$, there is a subspace $V\leq\mathbb{R}^{d}$
such that $\theta$ is $(V,O(\frac{\log m}{m}),m)$-saturated, and
if $W$ is any subspace such that $\theta$ is $(W,\frac{1}{m},m)$-saturated,
then $W\sqsubseteq V^{(O((\log m)/m))}$.\end{prop}
\begin{proof}
Write $\delta_{k}=C2^{k}k\log(m)/m$ where $C>1$ is large enough
to serve as the implicit a constant in the big-$O$ expressions we
invoke below. Note that $\delta_{k}<\delta_{k+1}$ and $\delta_{d}=O_{d}(\frac{\log m}{m})$.
Let $V$ be a subspace of maximal dimension such that $\theta$ is
$(V,\delta_{\dim V},m)$-saturated (such subspaces exist, e.g. $V=\{0\}$).
Write $k=\dim V$ and suppose $\theta$ is $(W,1/m,m)$-saturated
for some $W$. If $W\not\sqsubseteq V^{(\delta_{d})}$ then certainly
$W\not\sqsubseteq V^{(\delta_{k})}$, so by Lemma \ref{lem:subspace-lemmas}\eqref{enu:subspace-lemmas-4},
there is a subspace $W'\subseteq W$ with $\angle(V,W')>\delta_{k}$.
By Lemma \ref{lem:concentration-saturation-uniformity-under-change-of-subspace}
\eqref{enu:saturation-passes-to-nearby-subspaces} $\theta$ is $(W',(1+C)/m,m)$-saturated,
and since $(1+C)/m<\delta_{k}$ it is $(W',\delta_{k},m)$-saturated.
By Lemma \ref{lem:concentration-saturation-uniformity-under-change-of-subspace}
\eqref{enu:saturation-direct-sum}, $\theta$ is $(V+W',2\delta_{k}+\frac{C}{m}\log(\frac{1}{\delta_{k}}),m)$-saturated.
Since $2\delta_{k}+\frac{C}{m}\log(\frac{1}{\delta_{k}})<\delta_{k+1}$
the measure $\theta$ is $(V+W',\delta_{k+1},m)$-saturated. Since
$V'=V+W'$ has dimension at least $1+k$ and $\theta$ is $(V',\delta_{\dim V'},m)$-saturated,
this contradicts the definition of $V$.
\end{proof}

\subsection{Measures with uniformly concentrated components}

When a measure has the property that at each level the components
are with high probability concentrated on a subspace, one may expect
the subspace to vary slowly between levels. This is the content of
the following proposition, which may be applied to the conclusion
of Theorem \ref{thm:inverse-thm-Rd}, but is also needed in the theorem's
proof.
\begin{prop}
\label{prop:uniformizing-concentration-subspaces}Let $0<\varepsilon<1$
 and set $\delta=3\cdot8^{d-1}\varepsilon^{1/(4\cdot3^{d^{2}})}$.
Let $\eta\in\mathcal{P}([0,1)^{d})$ and $n\in\mathbb{N}$, and suppose
that for every $n\leq k\leq n+\frac{1}{2}\log(1/\varepsilon)$ there
is given a linear subspace $W_{k}\leq\mathbb{R}^{d}$ satisfying
\begin{equation}
\mathbb{P}_{i=k}(\eta^{x,i}\mbox{ is }(W_{k},\varepsilon)\mbox{-concentrated})>1-\varepsilon.\label{eq:59a}
\end{equation}
Then there are subspaces $V_{k}\leq W_{k}$ such that for $n\leq k\leq\frac{1}{2}\log(1/\varepsilon)$,
\begin{equation}
\mathbb{P}_{i=k}(\eta^{x,i}\mbox{ is }(V_{k},\delta)\mbox{-concentrated})>1-2d\sqrt{\varepsilon},\label{eq:59b}
\end{equation}
and $V_{j}\sqsubseteq V_{i}^{(\delta)}$ for all $n\leq i\leq j\leq n+\frac{1}{2}\log(1/\varepsilon)$.\end{prop}
\begin{proof}
Write $N=[\frac{1}{2}\log(1/\varepsilon)]$. For each $n\leq i\leq n+N$
set $\mathcal{W}^{i}=\{W_{j}\,:\,n\leq j\leq i\}$ and apply Proposition
\ref{prop:maximal-common-subspace} the with parameter $\varepsilon^{1/4}$
to obtain a subspace $V_{i}$ satisfying $V_{i}\sqsubseteq W_{j}^{(\delta/3)}$
and $V_{i}\sqsubseteq V_{j}^{(\delta/3)}$ for $n\leq j\leq i$, and
$r(i)\leq d$ subspaces $W_{i,1},\ldots,W_{i,r(i)}\in\mathcal{W}_{i}$
such that $\bigcap_{j=1}^{r(i)}W_{i,j}^{(\varepsilon^{1/4})}\sqsubseteq V_{i}^{(\delta/3)}$. 

Now, given $i$ and $1\leq j\leq r(i)$, there is by definition a
$n\leq k=k(i,j)\leq i$ such that $W_{i,j}=W_{k(i,j)}$. For every
component $\theta=\eta^{x,k}$ in the event in \eqref{eq:59a}, we
can apply Lemma \ref{lem:concentrated-measures-have-concentrated-components}
(using $i-k(i,j)\leq\frac{1}{2}\log(1/\varepsilon)$) to get 
\[
\mathbb{P}_{u=i}(\theta^{x,u}\mbox{ is }(V_{k},\varepsilon^{1/4})\mbox{-concentrated})>1-\sqrt{\varepsilon}.
\]
Thus by \eqref{eq:59a}, 
\[
\mathbb{P}_{u=i}(\eta^{x,u}\mbox{ is }(V_{k},\varepsilon^{1/4})\mbox{-concentrated})>1-2\sqrt{\varepsilon}.
\]
Hence,
\begin{eqnarray*}
\mathbb{P}_{u=i}(\eta^{x,u}\mbox{ is }(V_{k(i,j)},\varepsilon^{1/4})\mbox{-concentrated for all }1\leq j\leq r(i)) & > & 1-2r(i)\sqrt{\varepsilon}\\
 & \geq & 1-2d\sqrt{\varepsilon}.
\end{eqnarray*}
Finally, if $\theta=\eta^{x,i}$ is in the event above then, using
$\bigcap_{j=1}^{r(i)}W_{i,j}^{(\varepsilon^{1/4})}=\bigcap_{j=1}^{r(i)}W_{k(i,j)}^{(\varepsilon^{1/4})}\sqsubseteq V_{i}^{(\delta/3)}$
we have 
\begin{eqnarray*}
\theta(V_{i}^{(\delta/3)}) & \geq & 1-\sum_{j=1}^{r(i)}(1-\theta(W_{k(i,j)}^{(\varepsilon^{1/4})}))\\
 & \geq & 1-r(i)\cdot\varepsilon^{1/4}.
\end{eqnarray*}
Since $r(i)\leq d$ and $d\varepsilon^{1/4}\leq\delta/3$, this means
that $\theta$ is $(V_{i},\delta/3)$-concentrated. Since this is
true for components $\theta=\eta^{x,i}$ with probability $>1-2d\sqrt{\varepsilon}$,
we have established \eqref{eq:59b}, in fact with $\delta/3$ instead
of $\delta$.

Finally, we show that we can assume $V_{k}\leq W_{k}$. If $\varepsilon$
is so large that $\delta\geq1$ there is nothing to prove since we
can take $V_{k}=W_{k}$ from the start, so assume $\delta<1$. From
this and the relation $V_{i}\sqsubseteq W_{i}^{(\delta/3)}$ it follows
that $\pi_{W_{i}}$ is injective on$V_{i}$ and satisfies $d(V_{i},\pi_{W_{i}}V_{i})\leq\delta/3$.
Thus $V_{i}^{(\delta/3)}\sqsubseteq(\pi_{W_{i}}V_{i})^{(\delta)}$,
so if a measure $\theta$ is $(V_{i},\delta/3)$-concentrated, it
is also $(\pi_{W_{i}}V_{i},\delta)$-concentrated. It follows that
if we replace $V_{i}$ by $\pi_{W_{i}}V_{i}$, we still will have
\eqref{eq:59b}, as desired. Also, since $V_{j}\sqsubseteq V_{i}^{(\delta/3)}$
for $n\leq i<j$ before the modification, and each subspace moves
by at most $\delta/3$, after the change we have $V_{j}\sqsubseteq V_{i}^{(\delta)}$
for $n\leq i<j$, as desired.\end{proof}
\begin{cor}
\label{cor:uniformizing-concentration-subspaces}For every $\ell\in\mathbb{N}$
and $0<\varepsilon<1$ the following holds with $\delta=3\cdot8^{d-1}\varepsilon^{1/(4\cdot3^{d^{2}})}$.
Let $\eta\in\mathcal{P}([0,1]^{d})$ and $N>\frac{1}{2}\log(1/\varepsilon)$,
and suppose that for each $0\leq q\leq N$ there is given a subspace
$W_{q}\leq\mathbb{R}^{d}$ such that 
\[
\mathbb{P}_{i=q}(\eta^{x,i}\mbox{ is }(W_{q},\varepsilon)\mbox{-concentrated})>1-\varepsilon.
\]
Then there are subspaces $V_{q}\leq W_{q}$ such that
\[
\mathbb{P}_{i=q}(\eta^{x,i}\mbox{ is }(V_{q},\delta)\mbox{-concentrated})>1-2d\sqrt{\varepsilon},
\]
and 
\[
\frac{1}{N+1}\#\left\{ 0\leq i\leq N\,:\,d(V_{i},V_{i-\ell})\leq\delta\right\} \geq1-\frac{2(d+1)\ell}{\log(1/\varepsilon)}.
\]
(Note that the conclusion is of interest only when $\ell$ is small
compared to $\log(1/\varepsilon)$).\end{cor}
\begin{proof}
We may assume that $\delta<1$, otherwise the statement is trivial.

Let $m=[\frac{1}{2}\log(1/\varepsilon)]$. For each $k<[N/m]$ write
$I_{k}=\{mk,mk+1,\ldots,m(k+1)-1\}$ and for $k=[N/m]$ write $I_{k}=\{m[N/m],\ldots,N\}$.
For each $k\leq[N/m]$, apply the previous proposition with $n=km$
to find subspaces $V_{q}\leq W_{q}$, $q\in I_{k}$, such that $V_{j}\sqsubseteq V_{i}^{(\delta)}$
for all $i<j$ in $I_{k}$. This defines $V_{q}$ for all $0\leq q\leq N$.

Fix $k$. If $i<j$ are in $I_{k}$ then $V_{j}\sqsubseteq V_{i}^{(\delta)}$
(since $\delta<1$), hence $\dim V_{j}\leq\dim V_{i}$, and if $\dim V_{i}=\dim V_{j}$
then $d(V_{i},V_{j})\leq\delta$ (since $V_{j}\sqsubseteq V_{i}^{(\delta)}$).
Let $i_{0}=mk$ and let $i_{u+1}\in I_{k}$ denote the least index
such that $\dim V_{i_{u+1}}<\dim V_{i_{u}}$. There are at most $d$
such indices. It follows from the above that if $j+\ell\in I_{k}$
and $d(V_{j},V_{j-\ell})\geq\delta$ then $i_{u}\leq j<i_{u}+\ell$
for some $u$. There are at most $(d+1)\ell$ such indices $j$, so
\[
\#\left\{ i\,:\,i+\ell\in I_{k}\mbox{and }d(V_{i},V_{i-\ell})\geq\delta\right\} \leq(d+1)\ell.
\]
As the sets $I_{0},\ldots,I_{[N/m]}$ are disjoint and cover $\{0,\ldots,N+1\}$,
the bound above applies to each of them, so
\[
\#\left\{ 0\leq i<N\,:\,d(V_{i},V_{i-\ell})\geq\delta\right\} \leq(\frac{N}{m}+1)(d+1)\ell.
\]
Dividing by $N+1$ and using $N>\frac{1}{2}\log(1/\varepsilon)$ gives
the desired bound.
\end{proof}

\section{\label{sec:Convolutions} The inverse theorem in $\mathbb{R}^{d}$}

Our goal in this section is to prove Theorem \ref{thm:inverse-thm-Rd}.

\subsection{Elementary properties of convolutions}

We begin with the obvious.
\begin{lem}
\label{lem:entropy-monotonicity-under-convolution}For $m\in\mathbb{N}$
and $\mu,\nu\in\mathcal{P}(\mathbb{R}^{d})$ ,
\[
H_{m}(\mu*\nu)\geq H_{m}(\mu)-O(\frac{1}{m}).
\]
Also, if $\mu$ is $(V,\varepsilon,m)$-saturated then $\mu*\nu$
is $(V,\varepsilon',m)$-saturated, where $\varepsilon'=\varepsilon+O(1/m)$.\end{lem}
\begin{proof}
Notice that $\mu*\delta_{y}(A)=\mu(A-y)$, so that $H(\mu*\delta_{y},\mathcal{D}_{m})=H(\mu,\mathcal{D}_{m}+y)$,
where $\mathcal{D}_{m}+y=\{[a+y,b+y)\,:\,[a,b)\in\mathcal{D}_{m}\}$.
Thus by Lemma \ref{lem:entropy-weak-continuity-properties} \eqref{enu:entropy-translation},
we have $H_{m}(\mu*\delta_{y})\geq H_{m}(\mu)-O(\frac{1}{m})$. Since
$\mu*\nu=\int\mu*\delta_{y}d\nu(y)$, concavity of entropy implies
$H_{m}(\mu*\nu)\geq H_{m}(\mu)-O(\frac{1}{m})$. The second part follows
using the same relation and Lemma \ref{lem:saturation-passes-to-convex-combinations}.\end{proof}
\begin{cor}
\label{cor:entropy-growth-from-non-saturation}Let $\mu,\nu\in\mathcal{P}(\mathbb{R}^{d})$,
$m\in\mathbb{N}$, and let $V\leq\mathbb{R}^{d}$ be a linear subspace.
Suppose that $\mu$ is not $(V,2\varepsilon,m)$-saturated, and that
$\nu$ is $(V,\varepsilon,m)$-saturated. Then
\[
H_{m}(\mu*\nu)>H_{m}(\mu)+\varepsilon',
\]
where $\varepsilon'=\varepsilon-O(1/m)$.\end{cor}
\begin{proof}
Write $W=V^{\perp}$. By the previous lemma (with the roles of $\mu,\nu$
reversed),
\[
H_{m}(\mu*\nu)\geq H_{m}(\pi_{W}(\mu*\nu))+\dim V-\left(\varepsilon+O(1/m)\right).
\]
Since $\pi_{W}$ is linear, $\pi_{W}(\mu*\nu)=\pi_{W}\mu*\pi_{W}\nu$,
by the previous lemma $H_{m}(\pi_{W}(\mu*\nu))\geq H_{m}(\pi_{W}\mu)-O(1/m)$.
Inserting this in the last inequality and using the assumption that
$H_{m}(\mu)\leq H_{m}(\pi_{W}\mu)+\dim V-2\varepsilon$, and absorbing
another $O(1/m)$ into the error term, we have
\begin{eqnarray*}
H_{m}(\mu*\nu) & \geq & H_{m}(\pi_{W}\mu)+\dim V-\left(\varepsilon+O(1/m)\right)\\
 & \geq & H_{m}(\mu)+\left(\varepsilon-O(1/m)\right),
\end{eqnarray*}
as claimed.\end{proof}
\begin{lem}
\label{lem:convolution-of-concenrated-subspaces}Let $\mu\in\mathcal{P}(\mathbb{R}^{d})$
be $(V,\varepsilon)$-concentrated, $0<\varepsilon<1$. Then $\mu^{*k}$
is $(V,(1-\varepsilon k))$-concentrated for all $k\in\mathbb{N}$
with $\varepsilon k<1$. \end{lem}
\begin{proof}
Let $\mu=\varepsilon\mu_{1}+(1-\varepsilon)\mu_{2}$ with $\mu_{1}$,$\mu_{2}$
probability measures and $\mu_{1}$ supported on a translate of $V^{(\varepsilon)}$.
Then we can write $\mu^{\times k}=(1-\varepsilon)^{k}\mu_{1}^{\times k}+(1-(1-\varepsilon)^{k})\nu_{k}$
for some probability measure $\nu_{k}$, so, writing $\pi_{k}(x_{1}\ldots x_{k})=\sum_{i=1}^{k}x_{i}$,
we have 
\[
\mu^{*k}=\pi_{k}\mu^{\times k}=(1-\varepsilon)^{k}\pi_{k}\mu_{1}^{\times k}+(1-(1-\varepsilon)^{k})\pi_{k}\nu_{k}
\]
Since $\mu_{1}$ is supported on a translate of $V^{(\varepsilon)}$,
the measure $\mu_{1}^{*k}=\pi_{k}\mu^{\times k}$ is supported on
a translate of $\sum_{i=1}^{k}V^{(\varepsilon)}=V^{(\varepsilon k)}$.
So the splitting of $\mu^{*k}\,$ above shows that $(1-\varepsilon)^{k}$
of the mass of $\mu^{*k}$ is supported on an $\varepsilon k$-neighborhood
of a translate of $V$. Since $(1-\varepsilon)^{k}\geq1-\varepsilon k$,
the claim follows.
\end{proof}

\subsection{\label{sub:Covariance-matrices}\label{sub:Normal-measures-and-Berry-Esseen}Mean,
covariance and concentration}

A rough but convenient way to describe the distribution of a measure
is via its mean and covariance matrix. In this section we develop
some basic properties of these objects and their relation to concentration.

By a covariance matrix we shall mean a $d\times d$ real symmetric
matrix with non-negative eigenvalues (we do not require them to be
positive). We denote the eigenvalues of such a matrix $\Sigma$ by
\[
\lambda_{1}(\Sigma)\geq\lambda_{2}(\Sigma)\geq\ldots\geq\lambda_{d}(\Sigma).
\]
set $\lambda_{k}=0$ for $k>d$, preserving monotonicity. Define $\eigen_{1\ldots r}(\Sigma)$
to be the span in $\mathbb{R}^{d}$ of the eigenvectors corresponding
to eigenvalues $\geq\lambda_{r}(\Sigma)$. Note that if $\lambda_{r}(\Sigma)=\lambda_{r+1}(\Sigma)$
then $\dim(\eigen_{1\ldots r}(\Sigma))>r$. 

It is advantageous to think of a covariance matrix as the positive
semi-definite bi-linear form which it determines. The correspondence
between these objects is not one-to-one: The matrix determines the
form but the form determines the matrix only given the standard basis.
Nevertheless, given the inner product, the form determines the eigenvalues
and eigenspaces, and we are primarily interested in these; since the
inner product is always fixed in our discussion, we will not lose
much by thinking in terms of linear forms, and use the same notation
for both. One advantage of this approach is that a bi-linear form
can be restricted to a linear subspace, giving another bi-linear form,
which is positive semi-definite if the original one was. 
\begin{lem}
\label{lem:eigenvalues-of-restricted-form}Let $\Sigma$ be a positive
semidefinite form on $\mathbb{R}^{d}$ and $U\leq\mathbb{R}^{d}$
as subspace. Suppose that $u_{1},\ldots,u_{k}$ is an orthonormal
basis for $U$ and that $\Sigma(u_{i},u_{i})<\varepsilon$ for $i=1,\ldots,d$.
Then $\lambda_{1}(\Sigma|_{U\times U})\leq d\varepsilon$.\end{lem}
\begin{proof}
Let $u=\sum a_{i}u_{i}\in U$ be a unit vector, write $a=(a_{1},\ldots,a_{d-r+1})$,
so that $\left\Vert a\right\Vert _{2}=1$. Using Cauchy-Schwartz inequality
for the ``semi-inner product'' $\left\langle v,w\right\rangle \mapsto v^{T}\Sigma w$
(which may not be positive, but satisfies the requirements for the
weak inequality), and again to get $\sum|a_{i}|\leq\sqrt{d}\left\Vert a\right\Vert _{2}$,
we have
\begin{eqnarray*}
u^{T}\Sigma_{\mu}u & = & \sum_{i,j}a_{i}a_{j}u_{i}^{T}\Sigma_{\mu}u_{j}\\
 & \leq & \sum_{i,j}a_{i}a_{j}\sqrt{(u_{i}^{T}\Sigma_{\mu}u_{i})(u_{j}^{T}\Sigma_{\mu}u_{j})}\\
 & < & \sum_{i,j}a_{i}a_{j}\varepsilon\\
 & \leq & \varepsilon\cdot(\sum|a_{i}|)(\sum|a_{j}|)\\
 & \leq & \varepsilon\cdot d\cdot\left\Vert a\right\Vert _{2}^{2}\\
 & = & d\cdot\varepsilon
\end{eqnarray*}
This proves the claim.
\end{proof}
For $\mu\in\mathcal{P}(\mathbb{R}^{d})$, the mean\emph{ }of $\mu$
is 
\[
m=m(\mu)=\int x\,d\mu(x),
\]
and the covariance\emph{  }matrix of $\mu$ is 
\[
\Sigma(\mu)=\int(x-m)(x-m)^{T}\,d\mu(x).
\]
In this case we abbreviate 
\[
\lambda_{i}(\mu)=\lambda_{i}(\Sigma(\mu)),
\]
and similarly $\eigen_{1\ldots r}(\mu)$. We note that scaling a measure
by $r$ results in multiplying its covariance matrix by $r^{2}$,
an operation which does not affect the eigenvalues or eigenspaces.
\begin{lem}
\label{lem:restricted-form-as-form-of-projected-measure}Let $\mu\in\mathcal{P}(\mathbb{R}^{d})$,
write $\Sigma=\Sigma(\mu)$, and let $U\leq\mathbb{R}^{d}$ be a linear
subspace. Then $\Sigma|_{U}=\Sigma(\pi_{U}\mu)$ (the equality is
of bi-linear forms on $U$).\end{lem}
\begin{proof}
Write $m=m(\mu)$ and $m_{U}=\pi_{U}m=m(\pi_{U}\mu)$ (the last equality
is immediate). For vectors $u,v\in U$, we now have
\begin{eqnarray*}
u^{T}\Sigma v & = & u^{T}\left(\int(x-m)(x-m)^{T}d\mu(x)\right)v\\
 & = & \int u^{T}(x-m)(x-m)^{T}vd\mu(x)\\
 & = & \int\left\langle u,x-m\right\rangle \left\langle v,x-m\right\rangle d\mu(x)\\
 & = & \int\left\langle u,\pi_{U}(x-m)\right\rangle \left\langle v,\pi_{U}(x-m)\right\rangle d\mu(x)\\
 & = & \int\left\langle u,x-m_{U}\right\rangle \left\langle v,x-m_{U}\right\rangle d\pi_{U}\mu(x)\\
 & = & u^{T}\Sigma(\pi_{U}\mu)v.
\end{eqnarray*}
This proves the claim.
\end{proof}
A measure $\mu$ is supported on an $r$-dimensional affine subspace
of $\mathbb{R}^{d}$ if and only if $\lambda_{i}(\mu)=0$ for $i>r$,
in which case it is supported on a translate of $\eigen_{1\ldots r}\mu$.
We will use a quantitative version of this fact: 
\begin{lem}
\label{lem:concentration-from-covariance-matrix}Let $\mu\in\mathcal{P}(\mathbb{R}^{d})$
and write $\lambda_{i}=\lambda_{i}(\mu)$ and $V_{r}=\eigen_{1\ldots r}(\mu)$.
\begin{enumerate}
\item $\mu$ is $(V_{r},O(\lambda_{r+1}^{1/3}))$-concentrated. 
\item \label{enu:concentration-from-covariance-matrix-concentrated-case}If
$\mu\in\mathcal{P}([0,1])$ is $(V,\varepsilon)$-concentrated for
some $r$-dimensional subspace $V$ and $\varepsilon>0$, then $\lambda_{r+1}=O(\varepsilon)$
and $\mu$ is $(V_{r},O(\varepsilon^{1/3}))$-concentrated.
\item \label{enu:concentration-from-covariance-matrix-random-case}Let $\mu=\mu_{\omega}\in\mathcal{P}([0,1]^{d})$
be a random measure defined on some probability space $(\Omega,\mathcal{F},\mathbb{P})$.
Set $A=\mathbb{E}(\Sigma(\mu))$. If $\lambda_{r+1}(A)<\varepsilon$,
then, writing $V=\eigen_{1,\ldots,r}(A)$,
\[
\mathbb{P}\left(\mu\mbox{ is }(V,O(\varepsilon^{1/6}))\mbox{-concentrated}\right)>1-O(\sqrt{\varepsilon}).
\]

\end{enumerate}
\end{lem}
\begin{proof}
Let $\xi$ be an $\mathbb{R}^{d}$-valued random variable distributed
according to $\mu$ and let $m=m(\mu)$ and $\Sigma=\Sigma(\mu)$.
Identifying column vectors with $d\times1$ matrices and scalars with
$1\times1$ matrices, for $u\in\mathbb{R}^{d}$ we have
\begin{eqnarray*}
\mathbb{E}\left(\left\langle u,\xi-m\right\rangle ^{2}\right) & = & \mathbb{E}\left(u^{T}(\xi-m)(\xi-m)^{T}u\right)\\
 & = & u^{T}\mathbb{E}\left((\xi-m)(\xi-m)^{T}\right)u\\
 & = & u^{T}\Sigma u.
\end{eqnarray*}
For a subspace $W$, let $\eta_{W}$ denote the rotation-invariant
probability measure on the unit sphere in $W$. Then there exists
a constant $c=c(r)$ such that 
\[
d(\xi,m+V_{r})^{2}=c\cdot\int\left\langle u,\xi-m\right\rangle ^{2}\,d\eta_{V_{r}^{\perp}}(u).
\]
Therefore,
\begin{eqnarray*}
\mathbb{E}\left(d(\xi,m+V_{r})^{2}\right) & = & \mathbb{E}\left(c\cdot\int\left\langle u,\xi-m\right\rangle ^{2}\,d\eta_{V_{r}^{\perp}}(u).\right)\\
 & = & c\cdot\int\mathbb{E}\left(\left\langle u,\xi-m\right\rangle ^{2}\right)\,d\eta_{V_{r}^{\perp}}(u).\\
 & \leq & c\cdot\lambda_{r+1}.
\end{eqnarray*}
because $u^{T}\Sigma u\leq\lambda_{r+1}$ for every unit vector in
$V_{r}^{\perp}$. Now (1) follows from Markov's inequality.

For (2), fix $V$ as in the statement. Since $r+1+\dim V^{\perp}>d$,
we must have $\dim\left(\eigen_{1,\ldots,r+1}\cap V^{\perp}\right)\geq1$.
Fix a unit vector $w\in\eigen_{1,\ldots,r+1}\cap V^{\perp}$. Then
\begin{eqnarray}
\mathbb{E}\left(d(\xi,m+V)^{2}\right) & = & \mathbb{E}\left(\sup_{u\in V^{\perp}}\frac{\left\langle u,\xi-m\right\rangle ^{2}}{\left\Vert u\right\Vert ^{2}}\right)\nonumber \\
 & \geq & \sup_{u\in V^{\perp}}\mathbb{E}\left(\frac{\left\langle u,\xi-m\right\rangle ^{2}}{\left\Vert u\right\Vert ^{2}}\right)\nonumber \\
 & \geq & \mathbb{E}\left(\left\langle w,\xi-m\right\rangle ^{2}\right)\nonumber \\
 & = & w^{T}\Sigma(\mu)w\nonumber \\
 & \geq & \lambda_{r+1}.\label{eq:51}
\end{eqnarray}
On the other hand, since $\mu\in\mathcal{P}([0,1]^{d})$ we have $\left\Vert \xi\right\Vert \leq\sqrt{d}$
$\mu$-a.s., hence, writing $\delta=\varepsilon(1+2\sqrt{d})$, 
\begin{eqnarray}
\mathbb{E}\left(d(\xi,m+V)^{2}\right) & \leq & \delta^{2}\mathbb{P}(\xi\in(m+V)^{(\delta)})+d\cdot\mathbb{P}(\xi\in[0,1]^{d}\setminus(m+V)^{(\delta)})\nonumber \\
 & \leq & \delta^{2}+d\cdot\mathbb{P}(\xi\in[0,1]^{d}\setminus(m+V)^{(\delta)}).\label{eq:50}
\end{eqnarray}
Finally, since $\mu$ is $(V,\varepsilon)$-concentrated, there is
a translate $U$ of $V$ such that $\mu(U^{(\varepsilon)})>1-\varepsilon$.
Hence
\begin{eqnarray*}
m & = & \mathbb{E}(\xi)\\
 & = & \mu(U^{(\varepsilon)})\mathbb{E}(\xi|\xi\in U^{(\varepsilon)})+(1-\mu(U^{(\varepsilon)})\mathbb{E}(\xi|\xi\in\mathbb{R}^{d}\setminus U^{(\varepsilon)}).
\end{eqnarray*}
Since $U^{(\varepsilon)}$ is convex, $\mathbb{E}(\xi|\xi\in U^{(\varepsilon)})\in U^{(\varepsilon)}$.
Also, since $\left\Vert \xi\right\Vert \leq\sqrt{d}$, both expectations
on the right hand side of the last equation have magnitude at most
$\sqrt{d}$. Thus 
\[
d(m,U^{(\varepsilon)})\leq\left\Vert m-\mathbb{E}(\xi|U^{(\varepsilon)})\right\Vert \leq2\varepsilon\sqrt{d}.
\]
Therefore $U^{(\varepsilon)}\subseteq m+V^{(\varepsilon+2\varepsilon\sqrt{d})}=m+V^{(\delta)}$,
and consequently
\[
\mathbb{P}(\xi\in[0,1]^{d}\setminus(m+V)^{(\delta)})\leq\mathbb{E}(\xi\notin U^{(\varepsilon)})<\varepsilon.
\]
Combined with \eqref{eq:51} and \eqref{eq:50} this proves the first
part of (2), the second part now follows from (1).

We turn to (3). Let $U=(\eigen_{1\ldots r}A)^{\perp}\leq\eigen_{r+1,\ldots,d}A$.
Also for brevity write $\Sigma_{\mu}=\Sigma(\mu)$. For any unit vector
$u\in U$, we have 
\[
\varepsilon>u^{T}Au=\mathbb{E}(u\Sigma_{\mu}u^{T})
\]
Since $u^{T}\Sigma_{\mu}u\geq0$, by Markov's inequality,
\[
\mathbb{P}(u^{T}\Sigma_{\mu}u>\sqrt{\varepsilon})<\sqrt{\varepsilon}
\]
Now fix an orthonormal basis $u_{1}\ldots u_{\ell}$ of $U$ (so $\ell\leq d-(r+1)$).
By the last inequality,
\[
\mathbb{P}(u_{i}^{T}\Sigma_{\mu}u_{i}\leq\sqrt{\varepsilon}\mbox{ for all }i=1,\ldots,d)\geq1-d\sqrt{\varepsilon}
\]
By the Lemma \ref{lem:eigenvalues-of-restricted-form}, the condition
in the event above implies that $\lambda_{1}(\Sigma_{\mu}|_{U\times U})\leq d\sqrt{\varepsilon}$,
where $\Sigma_{\mu}|_{U\times U}$ is the restriction of the quadratic
for $\Sigma_{\mu}$ to $U\times U$, and by Lemma \ref{lem:restricted-form-as-form-of-projected-measure},
$\Sigma_{\mu}|_{U\times U}=\Sigma(\pi_{U}\mu)$ (as linear forms on
$U$). Combined with the previous probability estimate we get
\begin{equation}
\mathbb{P}\left(\lambda_{1}(\Sigma_{\pi_{U}\mu})\leq d\sqrt{\varepsilon}\right)\geq1-d\sqrt{\varepsilon}\label{eq:21}
\end{equation}
By the first part of this lemma, for $\mu$ in the event in \eqref{eq:21},
$\pi_{U}\mu$ is $(\{0\},O(\varepsilon^{1/6}))$-concentrated, and
this is the same as saying that $\mu$ is $(V,O(\varepsilon^{1/6}))$-concentrated,
as claimed.
\end{proof}
Recall the definition of the distance between linear subspaces \eqref{eq:metric-on-subspaces}.
We shall use the following basic fact, which we state without proof.
\begin{lem}
\label{lem:continuity-of-covariance}The maps $\Sigma\mapsto\lambda_{i}(\Sigma)$
are continuous on the set of positive semi-definite matrices. Furthermore,
given $\tau>\sigma>0$ and $1\leq r\leq d$, the map $\Sigma\mapsto\eigen_{1\ldots r}\Sigma$
is continuous on the compact space of positive semi-definite matrices
$\Sigma$ satisfying $\lambda_{r}(\Sigma)\geq\tau$ and $\lambda_{r+1}(\Sigma)\leq\sigma$.
\end{lem}

\subsection{Gaussian measures and the Berry-Esseen-Rotar estimate}

The standard $d$-dimensional Gaussian measure $\gamma=\gamma_{d}$
is given by $\gamma(A)=\int_{A}\varphi(x)dx$, where $\varphi=\varphi_{d}$
is $\varphi(x)=(2\pi)^{d/2}\exp(-\frac{1}{2}\left\Vert x\right\Vert ^{2})$.
The mean and covariance are $0$ and $I$ (the $d\times d$ identity
matrix), respectively. Given a $d\times d$ covariance matrix $\Sigma$
and $m\in\mathbb{R}^{d}$, write $\Sigma=BB^{T}$. The\emph{ }Gaussian
measure with mean $m\in\mathbb{R}^{d}$ and covariance $\Sigma$ is
the push-forward of $\gamma$ by the map $x\mapsto Bx+m$ and is denoted
$N(m,\Sigma)$. When $\Sigma$ is non-singular its density with respect
to Lebesgue is 
\[
f(x)=\frac{1}{\sqrt{(2\pi)^{d}\det\Sigma}}\exp(-\frac{1}{2}(x-m)^{T}\Sigma^{-1}(x-m)).
\]
When $\Sigma$ is singular and $r$ is such that $\lambda_{r}(\Sigma)>0$,
$\lambda_{r+1}(\Sigma)=0$, one obtains a similar formula for the
density on the affine space $V=\eigen_{1\ldots r}(\Sigma)+m$ with
respect to the $r$-dimensional Hausdorff measure on $V$. In particular,
if $\mu=N(m,\Sigma)$ and $\nu$ is the push-forward of $\mu$ through
the map $x\mapsto rx$, then $\nu=N(rm,r^{2}\Sigma)$. 

If $\mu_{1},\ldots,\mu_{k}$ are measures then $\mu=\mu_{1}*\ldots*\mu_{k}$
has mean $m(\mu)=\sum_{i=1}^{k}m(\mu_{i})$ and covariance $\Sigma(\mu)=\sum_{i=1}^{k}\Sigma(\mu_{i})$.
If $\mu_{i}=N(m_{i},\Sigma_{i})$ then $\mu_{1}*\ldots*\mu_{k}=N(\sum m_{i},\sum\Sigma_{i})$. 

The central limit theorem asserts that, for $\mu_{1},\mu_{2},\ldots\in\mathcal{P}(\mathbb{R}^{d})$
which are not too concentrated on subspaces, the convolutions $\mu_{1}*\ldots*\mu_{k}$
can be re-scaled so that the resulting measure is close to a Gaussian
measure. The Berry-Esseen estimate and its variants quantify the rate
of this convergence. The following multi-dimensional variant is due
to Rotar \cite{Rotar1970}. 
\begin{thm}
\label{thm:Berry-Esseen-Rotar}Let $\mu_{1},\ldots,\mu_{k}$ be probability
measures on $\mathbb{R}^{d}$ with finite third moments $\rho_{i}=\int\left\Vert x\right\Vert ^{3}\,d\mu_{i}(x)$.
Let $\mu=\mu_{1}*\ldots*\mu_{k}$ and let $\gamma$ be the Gaussian
measure with the same mean and covariance matrix as $\mu$. Then for
any convex Borel set $D\subseteq\mathbb{R}^{d}$,
\[
|\mu(D)-\gamma(D)|\leq C_{1}\cdot\frac{\sum_{i=1}^{k}\rho_{i}}{\lambda_{d}(\mu)^{3/2}},
\]
where $C_{1}=C_{1}(d)$. In particular, if $\rho_{i}\leq C$ and $\lambda_{d}(\mu_{i})\geq c$
for constants $c,C>0$ then 
\[
|\mu(D)-\gamma(D)|=O_{c,C}(k^{-1/2}).
\]

\end{thm}

\subsection{\label{sub:Estimating-modulus-of-continuity}Multi-scale analysis
of repeated self-convolutions}

If $\mu\in\mathcal{P}(\mathbb{R}^{d})$ is supported on a subspace
$V\leq\mathbb{R}^{d}$ but not on a smaller subspace, and if the support
is bounded, then $\mu^{*k}$ becomes increasingly smooth as a measure
on $V$, in the sense that, by the central limit theorem, it converges
(after suitable re-scaling) to a Gaussian on $V$. In this section
we prove a localized version of this statement which applies with
high probability to the components of the measure. Specifically, for
$\mu\in\mathcal{P}(\mathbb{R}^{d})$ of bounded support, for every
$\delta>0$ and integer scale $m$, there are subspaces $V_{0},V_{1},\ldots$
such that typical level-$i$ components of $\mu$ are $(\delta,2^{-m})$-concentrated
on $V_{i}$, and, when $k$ is large, a typical level-$i$ component
of $\mu^{*k}$ is $(V_{i},\delta,m)$-saturated. 

For a linear subspace $V\leq\mathbb{R}^{d}$ let $\pi_{V}$ denote
the orthogonal projection $\mathbb{R}^{d}\rightarrow V$. Recall our
convention that $\lambda_{i}=0$ for $i>d$, and in what follows define
$\lambda_{0}(\Sigma)=d$, so that when $\mu\in\mathcal{P}([0,1)^{d})$
and $\Sigma=\Sigma(\mu)$ the sequence $(\lambda_{i}(\mu))_{i=0}^{\infty}$
is monotone. 
\begin{prop}
\label{pro:entropy-convolution-estimate}Let $\sigma>0$, $\delta>0$,
$R>0$ and\footnote{In the one-dimensional case in \cite{Hochman2014} there wan no requirement
that $m$ be large. The reason this is necessary in the multi-dimensional
case is that, even when $\mathcal{\mu\in P}([0,1]^{d})$ is Lebesgue
measure on $V\cap[0,1]^{d}$ for an affine subspace $V$, we do not
generally have $H_{m}(\mu)=\dim V$, but rather only $H_{m}(\mu)=\dim V-o(1)$.
One can change coordinates so that if $\mathcal{D}_{m}$ is defined
in the new coordinates, $H_{m}(\mu)=\dim V$, but the coordinate change
itself incurs an $O(1/m)$ loss for $H_{m}(\cdot)$. } $m>m(\delta,R)$. Then there exists an integer $p=p_{0}(\sigma,\delta,R,m)$
such that for all $k\geq k_{0}(\sigma,\delta,R,m)$ and all $0\leq\rho<\rho_{0}(\sigma,\delta,R,m,k)$,
the following holds:

Let $\mu_{1},\ldots,\mu_{k}\in\mathcal{P}([-R,R]^{d})$, let $\mu=\mu_{1}*\ldots*\mu_{k}$
and $V=\eigen_{1,\ldots,r}\mu$ for some $0\leq r\leq d$, and suppose
that $\lambda_{r}(\mu)\geq\sigma k$ and $\lambda_{r+1}(\mu)\leq\rho$.
Then
\begin{equation}
\mathbb{P}_{i=p-[\log\sqrt{k}]}\left(\mu^{x,i}\mbox{ is }(V,\delta,m)\mbox{-uniform}\right)>1-\delta.\label{eq:3}
\end{equation}
\end{prop}
\begin{rem}
Instead of $\lambda_{r+1}(\mu)<\rho$ we could require $\mu$ to be
$(V,\rho)$-concentrated. This would give a formally equivalent statement
(using Lemma \ref{lem:concentration-from-covariance-matrix} \eqref{enu:concentration-from-covariance-matrix-concentrated-case}).\end{rem}
\begin{proof}
It is a general fact that, for an absolutely continuous probability
measure $\gamma$, for $\gamma$-a.e. $x$, as $p\rightarrow\infty$
the components $\gamma^{x,p}$ converge weak-{*} to Lebesgue measure
on $[0,1]^{d}$, and in particular 
\begin{equation}
\mathbb{E}_{i=p}(H_{m}(\gamma^{x,i}))\rightarrow d\qquad\mbox{as }p\rightarrow\infty\label{eq:local-entropy-of-abs-cont-measures}
\end{equation}
(this is a consequence of the martingale convergence theorem). There
is no guaranteed rate of convergence, but if $\gamma$ has a continuous
density function $f$, then convergence holds at every $x$ for which
$f(x)>0$, and the rate depends only on $f(x)$ and on the modulus
of continuity of $f$ at $x$. In particular, when $f\in C^{1}$ has
a smooth density $f$, the convergence rate at $x$ is controlled
by $f(x)$ and the bounds on $\left\Vert \nabla f(x)\right\Vert $
near $x$. Thus, for any compact family $\mathcal{E}\subseteq M_{d}(\mathbb{R})$
of non-singular co-variance matrices and any compact $K\subseteq\mathbb{R}^{d}$,
convergence in \eqref{eq:local-entropy-of-abs-cont-measures} is uniform
as $\gamma$ ranges over the Gaussians $\gamma$ with mean 0 and co-variance
matrix $\Sigma\in\mathcal{E}$, and $x$ ranges over $K$. Furthermore,
given $\mathcal{E}$ we can choose a compact $K_{2}\subseteq\mathbb{R}^{d}$
so that it has arbitrarily large mass uniformly for such $\gamma$.
Summarizing, given $0<\sigma,\delta<1$, there is a $p=p_{0}(\sigma,\delta,m)$
such that, for any Gaussian $\gamma$ with $\sigma\leq\lambda_{d}(\gamma)\leq1/\sigma$,
\begin{equation}
\mathbb{P}_{i=p}\left(H_{m}(\gamma^{x,i})>d-\delta\right)>1-\delta.\label{eq:2}
\end{equation}
In addition, by Lemma \ref{lem:entropy-weak-continuity-properties}
\eqref{enu:entropy-approximation}, there is a weakly open neighborhood
$\mathcal{U}_{\delta}\subseteq\mathcal{P}(\mathbb{R}^{d})$ of these
Gaussians such that the inequality continues to be valid for all $\gamma\in\mathcal{U}_{\delta}$.

Next, let $\mu=\mu_{1}*\ldots*\mu_{k}$ be as in the statement of
the proposition and first assume $r=d$, so $\lambda_{d}(\mu)\geq\sigma k$.
The third moments of the $\mu_{i}$ are bounded by $O_{R}(1)$, because
$\mu_{i}\in\mathcal{P}([-R,R]^{d})$. Thus by Theorem \ref{thm:Berry-Esseen-Rotar},
if $k$ is large enough in a manner depending only on $\delta$, the
scaling $\mu'$ of $\mu$ given by $\mu'(A)=\mu(2^{[\log\sqrt{k}]}\cdot A)$
will belong to $U_{\delta}$. Thus we obtain \eqref{eq:2} for $\mu'$.
Scaling everything back by a factor of $2^{[\log\sqrt{k}]}$ we obtain
\eqref{eq:3}.

Now consider the case that $\Sigma(\mu)$ is singular, i.e. $\lambda_{d}(\mu)=0$.
Fix $r,\rho$ and $V=\eigen_{1\ldots r}\mu$ as in the statement of
the proposition, and let $\pi=\pi_{V}$ denote the orthogonal projection
to $V$. Then the argument in the last paragraph applies in $V$ to
the measure $\pi\mu=\pi\mu_{1}*\ldots*\pi\mu_{k}$ and ensures that
\[
\mathbb{P}_{i=p-[\log\sqrt{k}]}\left(H_{m}(\pi\mu^{x,i})>r-\delta/2-O(1/m)\right)>1-\delta/2.
\]
The $O(1/m)$ term arises because we have transferred the entropy
bound from the dyadic partition $\mathcal{D}_{m}^{V}$ on $V$ to
the dyadic partition $\mathcal{D}_{m}$ of $\mathbb{R}^{d}$. But,
as we are assuming that $m$ is large relative to $\delta$, we can
absorb this term in $\delta$ and assume that
\begin{equation}
\mathbb{P}_{i=p-[\log\sqrt{k}]}\left(H_{m}(\pi_{V}\mu^{x,i})>r-\delta/2\right)>1-\delta/2.\label{eq:20}
\end{equation}
Now, the hypothesis $\lambda_{r+1}(\mu)\leq\rho$ means that $\mu$
is $(V,\sqrt{\rho})$-concentrated (Lemma \ref{lem:concentration-from-covariance-matrix})
and so for $\rho$ small enough in a manner depending on the other
parameters, a $(1-\delta/2)$-fraction of the components $\mu^{x,p-[\log\sqrt{k}]}$
are $(V,2^{-m})$-concentrated (Lemma \eqref{lem:convolution-of-concenrated-subspaces}).
In fact by taking $\rho$ small we can ensure an arbitrarily high
degree of concentration. Furthermore, if enough of the mass of such
a component $\mu^{x,p-[\log\sqrt{k}]}$ is concentrated on a small
enough neighborhood of $V$, then on this neighborhood $\pi$ will
be close enough to the identity map (in the supremum norm on continuous
self-maps of $[0,1]^{d}$) that Lemma \ref{lem:entropy-weak-continuity-properties}
\eqref{enu:entropy-geometric-distortion} will imply $|H_{m}(\mu^{x,p-[\log\sqrt{k}]})-H_{m}(\pi_{V}\mu^{x,p-[\log\sqrt{k}]})|<\delta/2$.
Combined with \eqref{eq:20}, we obtain \eqref{eq:3}.
\end{proof}
We now specialize to convolutions of a single measure.
\begin{prop}
\label{prop:saturation-of-components-of-convolution}Let $\sigma,\delta>0$
and $m>m(\delta)$. Then there exists $p=p_{1}(\sigma,\delta,m)$
such that for sufficiently large $k\geq k_{1}(\sigma,\delta,m)$ and
sufficiently small $0<\rho\leq\rho_{1}(\sigma,\delta,m,k)$, the following
holds. 

Let $\mu\in\mathcal{P}(\mathbb{R}^{d})$, fix an integer $i_{0}\geq0$,
and write 
\[
A=\mathbb{E}_{i=i_{0}}(\Sigma(\mu^{x,i})).
\]
If $\lambda_{r}(A)>\sigma$ and $\lambda_{r+1}(A)<\rho$ for some
$1\leq r\leq d$, then, setting $V=\eigen_{1\ldots r}(A)$, $\nu=\mu^{*k}$
and $j_{0}=i_{0}-[\log\sqrt{k}]+p$, we have 
\[
\mathbb{P}_{j=j_{0}}\left(\nu^{x,j}\mbox{ is }(V,\delta,m)\mbox{-saturated}\right)>1-\delta.
\]
\end{prop}
\begin{proof}
Fix $\sigma$, $\delta$, $m$, $k$, $\rho$, $\mu$, $A$ , $V$,
$i_{0}$ as in the statement, we will see that if the stated relationships
hold and $p$ is defined as in the statement, then the conclusion
holds.

Let $\widetilde{\mu}$ denote the $k$-fold self-product $\widetilde{\mu}=\mu\times\ldots\times\mu$
and $\pi:(\mathbb{R}^{d})^{k}\rightarrow\mathbb{R}^{d}$ the map 
\[
\pi(x_{1},\ldots,x_{k})=\sum_{i=1}^{k}x_{i}.
\]
Then $\nu=\pi\widetilde{\mu}$, and, since $\widetilde{\mu}=\mathbb{E}_{i=i_{0}}\left(\widetilde{\mu}_{x,i}\right)$,
we also have by linearity $\nu=\mathbb{E}_{i=i_{0}}\left(\pi(\widetilde{\mu}_{x,i})\right)$.
Thus, by Corollary \ref{cor:saturation-from-uniform-measures} and
an application of Markov's inequality, there is a $\delta_{1}>0$,
depending only on $\delta$ and $d$, such that if $m$ is large enough
as a function of $\delta_{1}$ then the proposition will follow if
we show that with probability $>1-\delta_{1}$ over the choice of
the component $\widetilde{\mu}_{x,i_{0}}$ of $\widetilde{\mu}$,
the measure $\tau=\pi(\widetilde{\mu}_{x,i_{0}})$ satisfies
\[
\mathbb{P}_{j=j_{0}}\left(\tau^{y,j}\mbox{ is }(V,\delta_{1},m)\mbox{-uniform}\right)>1-\delta_{1}.
\]
If we manage to define a random subspace $W=W(\widetilde{\mu}_{x,i_{0}})$
such that
\[
\mathbb{P}_{j=j_{0}}\left(d(W,V)<\frac{1}{\sqrt{d}}2^{-(m+1)}\mbox{ and }\tau^{y,j}\mbox{ is }(W,\delta_{1},m+1)\mbox{-uniform}\right)>1-\delta_{1},
\]
then the previous inequality follows by applying Lemma \ref{lem:concentration-saturation-uniformity-under-change-of-subspace}
to each component $\eta^{y,i}$ in the last event (we use here the
assumption that $m$ is large relative to $\delta_{1}$). We thus
aim to define $W$ such that \eqref{eq:5} holds.

Set $\eta=\pi(\widetilde{\mu}^{x,i_{0}})$ and notice that, with $\tau$
as before, the distribution of the components $\tau^{y,j_{0}}$ is
the same as the distribution of the components of $\eta^{z,j_{0}-i_{0}}$.
Thus what we really aim to prove is that we 
\begin{equation}
\mathbb{P}_{j=j_{0}-i_{0}}\left(d(W,V)<\frac{1}{\sqrt{d}}2^{-(m+1)}\mbox{ and }\eta^{y,j}\mbox{ is }(W,\delta_{1},m+1)\mbox{-uniform}\right)>1-\delta_{1}.\label{eq:5}
\end{equation}

A random component $\widetilde{\mu}^{x,i_{0}}$ is itself a product
measure $\widetilde{\mu}^{x,i_{0}}=\mu^{x_{1},i_{0}}\times\ldots\times\mu^{x_{k},i_{0}}$
(here $x=(x_{1},\ldots,x_{k})$), and the marginal measures $\mu^{x_{j},i_{0}}$
of this product are distributed independently according to the distribution
of the re-scaled components of $\mu$ at level $i_{0}$. Recall that
\begin{eqnarray}
\Sigma(\pi(\mu^{x_{1},i_{0}}\times\ldots\times\mu^{x_{k},i_{0}})) & = & \sum_{j=1}^{k}\Sigma(\mu^{x_{j},i_{0}})\label{eq:covariance-of-component-projection}
\end{eqnarray}
Fixing a parameter $\delta_{2}$ which will depend on $\sigma,\delta_{1}$,
by the weak law of large numbers, if $k$ is large enough in a manner
depending on $\delta_{2}$, then with probability $>1-\delta_{2}$
over the choice of $\widetilde{\mu}_{x,i_{0}}$ we will have\footnote{We use here the fact that we have a uniform bound for the rate of
convergence in the weak law of large numbers for i.i.d. random variables
$X_{1},X_{2},\ldots$. In fact, the rate can be bounded in terms of
the mean and variance of $X_{n}$. Here $X_{n}$ are matrix-valued
(they are distributed like the covariance matrix of the level-$i_{0}$
components of $\mu$), and therefore the mean and variance of the
components of $X_{n}$ can be bounded independently of the measure
$\mu\in\mathcal{P}([0,1)^{d})$.}
\begin{equation}
\left\Vert \frac{1}{k}\Sigma(\pi\widetilde{\mu}^{x,i_{0}})-A\right\Vert <\delta_{2}.\label{eq:4}
\end{equation}
Using Lemma \ref{lem:continuity-of-covariance} and the fact that
$\mu_{r}(A)>\sigma$ and $\lambda_{r+1}(A)<\rho$, and assuming as
we may that $\rho<k/4$, we can choose $\delta_{2}$ in a manner depending
on $\sigma,\delta_{1}$, in such a way that \eqref{eq:4} implies
\begin{eqnarray*}
\lambda_{r}(\pi\widetilde{\mu}^{x,i_{0}}) & > & \frac{k\sigma}{2}\\
\lambda_{r+1}(\pi\widetilde{\mu}^{x,i_{0}}) & < & \frac{k\sigma}{4}
\end{eqnarray*}
and such that, if we write $W_{x,i_{0}}=\eigen_{1\ldots r}\Sigma(\pi\widetilde{\mu}_{x,i_{0}})=\eigen_{1\ldots r}\Sigma(\pi\widetilde{\mu}^{x,i_{0}})$,
then 
\[
d(W_{x,i_{0}},V)<\frac{1}{\sqrt{d}}2^{-(m+1)}
\]
Thus, assuming that $k$ is large enough, we have shown 
\begin{equation}
\mathbb{P}_{i=i_{0}}\left(\lambda_{r}(\pi\widetilde{\mu}^{x,i})>\frac{k\sigma}{2}\mbox{ and }d(W_{x,i},V)<\frac{1}{\sqrt{d}}2^{-(m+1)}\right)>1-\delta_{2}\label{eq:23}
\end{equation}
Next, fix such a $k$. By hypothesis $\lambda_{r+1}(A)=\lambda_{r+1}(\mathbb{E}_{i=i_{0}}(\Sigma(\mu^{x,i})))<\rho$,
so by Lemma \ref{lem:concentration-from-covariance-matrix} \eqref{enu:concentration-from-covariance-matrix-random-case},
\[
\mathbb{P}_{i=i_{0}}\left(\mu^{x,i}\mbox{ is }(V,O(\rho^{1/6}))\mbox{-concentrated}\right)>1-O(\sqrt{\rho})
\]
Using again the fact that $\widetilde{\mu}^{x,i_{0}}$ is a product
of $k$ independent copies of level-$i_{0}$ components of $\mu$,
the last inequality implies
\begin{equation}
\mathbb{P}_{i=i_{0}}\left(\mbox{all marginals of }\widetilde{\mu}^{x,i}\mbox{ are }(V,O(\rho^{1/6})\mbox{-concentrated}\right)>1-O(k\sqrt{\rho})\label{eq:22}
\end{equation}
If $\widetilde{\mu}^{x,i_{0}}$ is in the event above, then all its
marginals are $(V,O(\rho^{1/6}))$-concentrated, so by Lemma \ref{lem:convolution-of-concenrated-subspaces},
$\pi(\widetilde{\mu}^{x,i_{0}})$ is $(V,O(k\rho^{1/6}))$-concentrated.
By Lemma \ref{lem:concentration-from-covariance-matrix} \eqref{enu:concentration-from-covariance-matrix-concentrated-case},
$\lambda_{r+1}(\pi(\widetilde{\mu}^{x,i_{0}}))<O_{k}(\rho^{1/6})$,
and we conclude that
\[
\mathbb{P}_{i=i_{0}}\left(\lambda_{r+1}(\pi(\widetilde{\mu}^{x,i_{0}}))\leq O_{k}(\rho^{1/6})\right)>1-O(k\sqrt{\rho})
\]
Combining this with \eqref{eq:23} and assuming that $\rho$ is sufficiently
small relative to $\delta_{2}$, we have
\[
\mathbb{P}_{i=i_{0}}\left(\begin{aligned}\lambda_{r}(\pi\widetilde{\mu}^{x,i}) & >\;\frac{k\sigma}{2}\\
\lambda_{r+1}(\pi(\widetilde{\mu}^{x,i})) & <\;O_{k}(\rho^{1/6})\\
d(W_{x,i_{0}},V) & <\;\frac{1}{\sqrt{d}}2^{-(m+1)}
\end{aligned}
\right)>1-2\delta_{2}
\]
Let $\widetilde{\mu}^{x,i_{0}}$ belong to the event above. Let us
recall the dependences of the parameters: $\delta$ is given and determines
$\delta_{1}$, then $m$ is large relative to $\delta_{1}$, then
$\delta_{2}$ small depending on $\sigma,\delta_{!}$, then $k$ is
correspondingly large, and $\rho$ correspondingly small. So we can
assume that $k$ is large enough, and $\rho$ small enough, to apply
Proposition \ref{pro:entropy-convolution-estimate} with parameters
$\delta_{1},m+1$ and $\sigma/2$, and conclude that there is a $p=p(\delta_{1},m+1,k)=p(\delta,m,k)$
such that, writing $\eta=\pi\widetilde{\mu}^{x,i_{0}}$, 
\[
\mathbb{P}_{j=j_{0}-i_{0}}\left(\eta^{y,j}\mbox{ is }(W_{x,i_{0}},\delta_{1},m+1)\mbox{-uniform}\right)>1-\delta_{1}.
\]
This and the estimate above on the probability that $d(W_{x,i_{0}},V)<\;\frac{1}{\sqrt{d}}2^{-(m+1)}$
give \eqref{eq:5}, which is what we wanted.\end{proof}
\begin{thm}
\label{thm:saturation-of-repeated-convolutions}Let $\delta>0$ and
$m\in\mathbb{N}$. Then there exists\footnote{In \cite{Hochman2014} the corresponding statement holds for all large
enough $k$. The reason the size of $k$ must be restricted is, roughly,
that if $\mu$ is concentrated extremely near a subspace $V$ then
it will remain so for a reasonable number of convolutions, but too
many convolutions will make it drift away from $V$.} a $0\leq k\leq k_{2}(\delta,m)$ such that for all sufficiently large
$n\geq n_{2}(\delta,m,k)$, the following holds: For any $\mu\in\mathcal{P}(\mathbb{R}^{d})$
there is a sequence $V_{0},\ldots,V_{n}$ of subspaces\footnote{The corresponding theorem in \cite{Hochman2014} is stated differently,
in terms of disjoint subsets $I,J\subseteq\{1,\ldots,n\}$. See remark
after Theorem \ref{thm:inverse-thm-Rd}.} of $\mathbb{R}^{d}$ such that, writing $\nu=\mu^{*k}$,
\[
\mathbb{P}_{0\leq i\leq n}\left(\nu^{x,i}\mbox{ is }(V_{i},\delta,m)\mbox{-saturated}\right)>1-\delta
\]
and
\[
\mathbb{P}_{0\leq i\leq n}\left(\mu^{x,i}\mbox{ is }(V_{i},\delta)\mbox{-concentrated}\right)>1-\delta.
\]
\end{thm}
\begin{proof}
It is a formal consequence of Proposition \ref{prop:concentration-and-saturation-change-of-parameters}
that we may assume that $m$ is large in a manner depending on $\delta$.
We also may assume that $\delta<1/2$. Also, since we are free to
take $R$ large relative to $n$, we can assume that $\mu$ is supported
on $[0,1)^{d}$.

Let $k_{1}(\cdot),$ $p_{1}(\cdot)$, $\rho_{1}(\cdot)$ be as in
Proposition \ref{prop:saturation-of-components-of-convolution}. We
assume, without loss of generality, that these functions are monotone
in each of their arguments. 

Let $c>1$ denote a constant good for all previous big-$O$ bounds.

The proof will depend on a function $\widetilde{\rho}:(0,d]\rightarrow(0,d]$
such that $\widetilde{\rho}(\sigma)$ is small in a manner depending
on $\sigma,\delta,m$. Specifically, we require that $\widetilde{\rho}$
satisfy the following inequalities, where $\exp_{2}(y)=2^{y}$ (for
concreteness, on could define $\widetilde{\rho}(\sigma)$ to be one-half
the minimum of the right-hand sides): 
\begin{eqnarray}
\widetilde{\rho}(\sigma) & < & \sigma,\label{eq:rho-tilde-1}\\
\widetilde{\rho}(\sigma) & < & \;\rho_{1}(\sigma,\delta/2,m,k_{1}(\sigma,\delta/2,m)),\label{eq:rho-tilde-2}\\
\widetilde{\rho}(\sigma) & < & \;\frac{\delta^{12}}{c^{6}(2d)^{12}},\label{eq:rho-tilde-3}\\
\widetilde{\rho}(\sigma) & < & \;\frac{1}{c^{6}}\exp_{2}(-\frac{24(d+1)\cdot([\log\sqrt{k_{1}(\sigma,\delta/2,m)}]-p_{1}(\sigma,\delta/2,m))}{\delta/2}),\label{eq:rho-tilde-4}\\
\widetilde{\rho}(\sigma) & < & \;\frac{1}{c^{4}}\cdot(\frac{\delta}{(\sqrt{d}+1)\cdot3\cdot8^{d-1}})^{24\cdot3^{d^{2}}}.\label{eq:rho-tilde-5}
\end{eqnarray}

As before define $\lambda_{0}(\Sigma_{i})=d$ and $\lambda_{d+1}(\Sigma_{i})=0$.
Fix $n$ and $\mu$, we shall later see how large an $n$ is desirable.
For $0\leq q\leq n$ write
\[
\Sigma_{q}=\mathbb{E}_{i=q}\left(\Sigma(\mu^{x,i})\right).
\]
Define a sequence $\sigma_{0}>\sigma_{1}>\ldots$ by $\sigma_{0}=d$
and $\sigma_{i}=\widetilde{\rho}(\sigma_{i-1})$ (the sequence is
decreasing because of \eqref{eq:rho-tilde-1}). For a covariance matrix
$\Sigma$ and $s\in\mathbb{N}$, set 
\[
N_{s}(\Sigma)=\#\{1\leq j\leq d\;:\;\lambda_{j}(\Sigma)\in(\sigma_{s},\sigma_{s-1}]\}.
\]

\begin{claim}
There is an $s\leq\left\lceil 1+2d/\delta\right\rceil $ satisfying
\[
\mathbb{P}_{0\leq q\leq n}(N_{s}(\Sigma_{q})=0)>1-\frac{\delta}{2}.
\]
\end{claim}
\begin{proof}
Note that $\sum_{r=1}^{\infty}N_{r}(\Sigma_{q})=d$, so 
\begin{equation}
\sum_{s=1}^{\left\lceil 1+2d/\delta\right\rceil }\mathbb{E}_{0\leq q\leq n}(N_{i}(\Sigma_{q}))=\mathbb{E}_{0\leq q\leq n}(\sum_{s=1}^{\left\lceil 1+2d/\delta\right\rceil }N_{i}(\Sigma_{q}))\leq d.\label{eq:expected-covariance-bound}
\end{equation}
Thus there must exist an $s\leq\left\lceil 1+2d/\delta\right\rceil $
such that 
\[
\mathbb{E}_{0\leq q\leq n}(N_{s}(\Sigma_{q}))\leq\frac{d}{\left\lceil 1+2d/\delta\right\rceil }<\frac{\delta}{2}.
\]
Since $N_{i}(\cdot)$ is integer valued, we have 
\[
\mathbb{P}_{0\leq q\leq n}(N_{s}(\Sigma_{q})\geq1)\leq\mathbb{E}_{0\leq q\leq n}(N_{s}(\Sigma_{q})),
\]
so this is the desired $s$.
\end{proof}
Fix an $s$ that satisfies the conclusion of the lemma, write 
\begin{eqnarray*}
\sigma & = & \sigma_{s-1}\\
\rho & = & \sigma_{s}\\
 & = & \widetilde{\rho}(\sigma),
\end{eqnarray*}
and set
\[
k=k_{1}(\sigma,\frac{\delta}{2},m).
\]
Note that $k$ is bounded above by some expression $k_{2}(\delta,m)$
(also depending implicitly on the choice of the function $\widetilde{\rho}$),
as in the statement, since its largest possible value occurs for $s=[1+2d/\delta]$,
and once the function $\widetilde{\rho}$ is fixed, the magnitude
$\sigma$, and hence $k$, is bounded.

Let
\[
I=\{0\leq q\leq n\,:\,N_{s}(\Sigma_{q})=0\}.
\]
By our choice of $s$, 
\[
|I|\geq(1-\frac{\delta}{2})(n+1).
\]
For $q\in I$ let $1\leq r_{q}\leq d$ denote the smallest integer
such that 
\[
\lambda_{r_{q}}(\Sigma_{q})\geq\sigma\qquad\mbox{and}\qquad\lambda_{r_{q}+1}(\Sigma_{q})<\rho,
\]
which exists by definition, and set 
\[
W_{q}=\eigen_{1,\ldots,r_{q}}(\Sigma_{q}).
\]
We define $W_{q}=\mathbb{R}^{d}$ for $q\notin I$. Finally, write
\[
\ell=[\log\sqrt{k}]-p_{1}(\sigma,\frac{\delta}{2},m).
\]

\begin{claim}
For $q\in I$,
\begin{eqnarray}
\mathbb{P}_{i=q}\left(\nu^{x,i-\ell}\mbox{ is }(W_{i},\frac{\delta}{2},m)\mbox{-saturated}\right) & > & 1-\frac{\delta}{2}\label{eq:8}\\
\mathbb{P}_{i=q}\left(\mu^{x,i}\mbox{ is }(W_{i},c\rho^{1/6})\mbox{-concentrated}\right) & > & 1-c\rho^{1/6}.\label{eq:8-prime}
\end{eqnarray}
\end{claim}
\begin{proof}
The first inequality follows from Proposition \ref{prop:saturation-of-components-of-convolution}
and our choice of parameters, specifically the definition of $\ell$
and assumption \eqref{eq:rho-tilde-2}. The second follows from Lemma
\ref{lem:concentration-from-covariance-matrix} (3) applied to the
random component $\mu^{x,i}$, since $W_{q}=\eigen_{1,\ldots,r_{q}},$
$\mathbb{E}_{i=q}\left(\lambda_{r_{q}+1}(\mu^{x,i})\right)<\rho$.
\end{proof}
This is almost what we want, except that in \eqref{eq:8} the level
of the component is shifted by $\ell$ (that is, $\nu^{x,i-\ell}$
appears instead of $\nu^{x,i}$). To correct this we apply Corollary
\ref{cor:uniformizing-concentration-subspaces} to \eqref{eq:8-prime}
with parameter $c\rho^{1/2}$ (we can do this since we are assuming
that $n$ is large relative to $\rho$). Then, writing 
\[
\rho'=3\cdot8^{d-1}c^{1/(4\cdot3^{d^{2}})}\rho^{1/(24\cdot3^{d^{2}})},
\]
we conclude the there are subspaces $W'_{i}\leq W_{i}$ such that
for all $0\leq q\leq n$, 
\begin{eqnarray}
\mathbb{P}_{i=q}\left(\mu^{x,i}\mbox{ is }(W'_{i},\rho')\mbox{-concentrated}\right) & > & 1-2d\sqrt{c\rho^{1/6}}\nonumber \\
 & > & 1-\delta\label{eq:8-prime-2}
\end{eqnarray}
(the last inequality by \eqref{eq:rho-tilde-3}), and
\begin{eqnarray*}
\frac{1}{n+1}\#\left\{ 0\leq q\leq n\,:\,d(W'_{q},W'_{q-\ell})\leq\rho'\right\}  & \geq & 1-\frac{2(d+1)\ell}{\log(1/c\rho^{1/6})}\\
 & > & 1-\frac{\delta}{2}
\end{eqnarray*}
(the last inequality in by assumption \eqref{eq:rho-tilde-4}). Let
\[
J=\{i\in I\,,\,d(W'_{i},W'_{i-\ell})\leq\rho'\}.
\]
Since $\frac{1}{n+1}|I|>1-\delta/2$, the previous equation implies
that
\begin{equation}
\frac{1}{n+1}|J|\geq1-\delta.\label{eq:55}
\end{equation}
Now, for any $\ell\leq q\leq n$, applying \eqref{eq:8-prime-2} to
$q-\ell$ we have
\[
\mathbb{P}_{i=q}\left(\nu^{x,i-\ell}\mbox{ is }(W'_{i-\ell},\rho')\mbox{-concentrated}\right)>1-\delta.
\]
Assuming also $q\in J$, we also have $d(W'_{q-\ell},W'_{q})\leq\rho'$,
so by Lemma \ref{lem:concentration-saturation-uniformity-under-change-of-subspace}
(1) applied to each component $\nu^{x,i-\ell}$ in the event above,
\[
\mathbb{P}_{i=q}\left(\nu^{x,i-\ell}\mbox{ is }(W'_{i},(\sqrt{d}+1)\rho')\mbox{-concentrated}\right)>1-\delta\qquad\mbox{for }q\in J.
\]
Our assumption \eqref{eq:rho-tilde-5} implies that $(\sqrt{d}+1)\rho'<\delta$,
and the last inequality yields 
\begin{equation}
\mathbb{P}_{i=q}\left(\nu^{x,i-\ell}\mbox{ is }(W'_{i},\delta)\mbox{-concentrated}\right)>1-\delta\qquad\mbox{for }q\in J.\label{eq:56}
\end{equation}
On the other hand for $q\in J$ we have $q\in I$ and so \eqref{eq:8}
holds. Since $W'_{q}\leq W_{q}$, by Lemma \ref{lem:concentration-saturation-uniformity-under-change-of-subspace}
(4), we have
\[
\mathbb{P}_{i=q}\left(\mu^{x,i-\ell}\mbox{ is }(W'_{i},\frac{\delta}{2}+O(\frac{1}{m}),m)\mbox{-saturated}\right)>1-\delta\qquad\mbox{for }q\in J.
\]
Since we are assuming $m$ large enough relative to $\delta$, this
implies
\begin{equation}
\mathbb{P}_{i=q}\left(\mu^{x,i-\ell}\mbox{ is }(W'_{i},\delta,m)\mbox{-saturated}\right)>1-\delta\qquad\mbox{for }q\in J.\label{eq:57}
\end{equation}
In conclusion, if we define $V_{i}=W'_{i+\ell}$ and replace $I$
by  $(J-\ell)\cap[0,n]$, then equations \eqref{eq:55}, \eqref{eq:56},
and \eqref{eq:57} give the desired conclusion, assuming that $m,\delta$
have the appropriate relationship to each other and to $\varepsilon$,
and that $n$ and is large enough.
\end{proof}

\subsection{\label{sub:Kaimanovitch-Vershik-Tao-theorem}The Ka\u\i{}manovich-Vershik
lemma}

The second ingredient in our proof of Theorem \ref{thm:inverse-thm-Rd}
is the following entropy analog of the Pl\"{u}nnecke-Rusza inequality:
\begin{lem}
[Ka\u\i{}manovich-Vershik, \cite{KaimanovichVershik1983}]\label{thm:Kaimanovitch-Vershik-Tao}
Let $\Gamma$ be a countable abelian group and let $\mu,\nu\in\mathcal{P}(\Gamma)$
be probability measures with $H(\mu)<\infty$, $H(\nu)<\infty$. Let
\[
\delta_{k}=H(\mu*(\nu^{*(k+1)}))-H(\mu*(\nu^{*k})).
\]
Then $\delta_{k}$ is non-increasing in $k$. In particular, 
\[
H(\mu*(\nu^{*k}))\leq H(\mu)+k\cdot(H(\mu*\nu)-H(\nu)).
\]

\end{lem}
This lemma first appears in a study of random walks on groups by Ka\u\i{}manovich
and Vershik \cite{KaimanovichVershik1983}. It was more recently rediscovered
and applied in additive combinatorics by Madiman and co-authors \cite{Madiman2008,MadimanMarcusTetali2012},
and in a weaker form independently by Tao \cite{Tao2010}. For a proof
using our notation see \cite{Hochman2014}. 

For non-discrete measures in $\mathbb{R}^{d}$ we have the following
analog:
\begin{cor}
\label{cor:non-discrete-KVT-theorem}Let $\mu,\nu\in\mathcal{P}(\mathbb{R}^{d})$
with $H_{n}(\mu),H_{n}(\nu)<\infty$. Then
\[
H_{n}(\mu*(\nu^{*k}))\leq H_{n}(\mu)+k\cdot\left(H_{n}(\mu*\nu)-H_{n}(\mu)\right)+O(\frac{k}{n}).
\]

\end{cor}
The error term arises in the same way as in Lemma \ref{lem:entropy-monotonicity-under-convolution}.
For the proof see \cite{Hochman2014} (the passage from $\mathbb{R}$
to $\mathbb{R}^{d}$ requires only notational changes).

\subsection{\label{sub:Proof-of-inverse-theorem}Proof of the inverse theorem}

We now prove Theorem \ref{thm:inverse-thm-Rd}, which we re-state
for convenience. 
\begin{thm}
\label{thm:inverse-theorem-with-large-m}For every $\varepsilon>0$,
$R>0$ and $m\in\mathbb{N}$, there exists $\delta=\delta(\varepsilon,R,m)>0$
such that for all $n>n(\varepsilon,R,m,\delta)$, the following holds:
if $\nu,\mu\in\mathcal{P}([-R,R]^{d})$ and 
\[
H_{n}(\mu*\nu)<H_{n}(\mu)+\delta,
\]
then there exists a sequence $V_{0},\ldots,V_{n}\leq\mathbb{R}^{d}$
of subspaces such that 
\begin{eqnarray*}
\mathbb{P}_{_{0\leq i\leq n}}\left(\begin{array}{c}
\mu^{x,i}\mbox{ is }(V_{i},\varepsilon,m)\mbox{-saturated and}\\
\nu^{x,i}\mbox{ is }(V_{i},\varepsilon)\mbox{-concentrated}
\end{array}\right) & > & 1-\varepsilon
\end{eqnarray*}
\end{thm}
\begin{proof}
Fix $\varepsilon>0$. It is a formal consequence of Proposition \ref{prop:concentration-and-saturation-change-of-parameters}
that it suffices for us to prove the theorem with the assumption that
$m$ is large in a manner depending on $\varepsilon$. We can also
assume that $\varepsilon<1/2$, and that $\varepsilon$ is small with
respect to $d$. Also, as we are free to choose $n$ large relative
to $R$, the distribution on components depends negligibly on dyadic
scales greater than $0$, and the scale-$n$ entropy of $\mu$ and
$\mu*\nu$ differs negligibly from the same entropy conditioned on
$\mathcal{D}_{0}$. Thus, without loss of generality, we can assume
that the measures are supported on $[0,1)^{d}$, and we omit mention
of $R$ from now on. 

Choose $k=k_{2}(\varepsilon,m)$ as in Theorem \ref{thm:saturation-of-repeated-convolutions}.
We shall show that the conclusion holds if $n$ is large relative
to the previous parameters.

Let $\mu,\nu\in\mathcal{P}([0,1)^{d})$. Denote 
\[
\tau=\nu^{*k}
\]
Assuming $n$ is large enough, Theorem \ref{thm:saturation-of-repeated-convolutions}
provides us with subspaces $V_{0},\ldots,V_{n}\subseteq\mathbb{R}^{d}$
such that 
\begin{equation}
\mathbb{P}_{0\leq i\leq n}\left(\nu^{x,i}\mbox{ is }(V_{i},\varepsilon)\mbox{-concentrated}\right)\geq1-\varepsilon,\label{eq:11}
\end{equation}
and
\[
\mathbb{P}_{0\leq i\leq n}\left(\tau^{x,i}\mbox{ is }(V_{i},\varepsilon,m)\mbox{-saturated}\right)>1-\varepsilon.
\]
If it holds that 
\begin{equation}
\mathbb{P}_{0\leq i\leq n}\left(\mu^{x,i}\mbox{ is }(V_{i},2\varepsilon,m)\mbox{-saturated}\right)>1-2\varepsilon\label{eq:12}
\end{equation}
then we are done, since \eqref{eq:11} and \eqref{eq:12} together
are the second alternative of the theorem we want to prove (with a
multiple of $\varepsilon$ instead of $\varepsilon$, but this is
formally equivalent). 

Otherwise, by Lemma \ref{cor:entropy-growth-from-non-saturation}
and the above we have
\begin{multline*}
\mathbb{P}_{0\leq i\leq n}\left(H_{m}(\mu^{x,i}*\tau^{y,i})>H_{m}(\mu^{x,i})+\varepsilon-O(\frac{1}{m})\right)\\
\begin{aligned} & \qquad\geq\quad\mathbb{P}_{0\leq i\leq n}\left(\mu^{x,i}\mbox{ is not }(V_{i},2\varepsilon,m)\mbox{-saturated and }\tau^{y,i}\mbox{ is }(V_{i},\varepsilon,m)\mbox{-saturated}\right)\\
 & \qquad>\quad\mathbb{P}_{0\leq i\leq n}\left(\mu^{x,i}\mbox{ is not }(V_{i},2\varepsilon,m)\mbox{-saturated}\right)\\
 & \qquad\quad\quad-\;\left(1-\mathbb{P}_{0\leq i\leq n}\left(\tau^{y,i}\mbox{ is }(V_{i},\varepsilon,m)\mbox{-saturated}\right)\right)\\
 & \qquad>\quad2\varepsilon-(1-(1-\varepsilon))\\
 & \qquad=\quad\varepsilon.
\end{aligned}
\end{multline*}
Let $\delta'(\mu^{x,i},\tau^{y,i})=H_{m}(\mu^{x,i}*\tau^{y,i})-H_{m}(\mu^{x,i})$.
By the previous calculation, with probability at least $\varepsilon$
we have $\delta'\geq\varepsilon-O(1/m)$, and by Lemma \ref{lem:entropy-monotonicity-under-convolution}
we always have $\delta'\geq-O(1/m)$. Thus
\[
\mathbb{E}_{0\leq i\leq n}\left(\delta'(\mu^{x,i},\tau^{y,i})\right)\geq\varepsilon^{2}-O(\frac{1}{m})
\]
Thus, by Lemmas \ref{lem:entropy-local-to-global} and \ref{lem:entropy-of-convolutions-via-component-convolutions},
\begin{eqnarray*}
H_{n}(\mu*\tau) & > & \mathbb{E}_{0\leq i<n}\left(H_{m}(\mu^{x,i}*\tau^{y,i})\right)-O(\frac{m}{n})\\
 & \geq & \mathbb{E}_{0\leq i\leq n}\left(H_{m}(\mu^{x,i})+\delta'(\mu^{x,i},\tau^{y,i})\right)-O(\frac{m}{n})\\
 & \geq & \mathbb{E}_{0\leq i\leq n}\left(H_{m}(\mu^{x,i})\right)+\varepsilon^{2}-O(\frac{1}{m}+\frac{m}{n})\\
 & = & H_{n}(\mu)+\varepsilon^{2}-O(\frac{1}{m}+\frac{m}{n}).
\end{eqnarray*}
So, assuming that $m$ is large and $n$ larger still, all in a manner
depending on $\varepsilon,d$, we have 
\[
H_{n}(\mu*\tau)>H_{n}(\mu)+\frac{\varepsilon^{2}}{2}.
\]
On the other hand, by Corollary \ref{cor:non-discrete-KVT-theorem}
above,
\[
H_{n}(\mu*\tau)\leq H_{n}(\mu)+k\cdot\left(H_{n}(\mu*\nu)-H_{n}(\mu)\right)+O(\frac{k}{n}).
\]
Assuming that $n$ is large enough in a manner depending on $d$,
$\varepsilon$ and $k$, this and the previous inequality give 
\[
H_{n}(\mu*\nu)\geq H_{n}(\mu)+\frac{\varepsilon^{2}}{3k}.
\]
This completes the proof of Theorem \ref{thm:inverse-thm-Rd} with
$\delta=\varepsilon^{2}/3k$.
\end{proof}

\section{\label{sec:convolutions-with-isometries}Inverse theorem for the
action of the isometry group on $\mathbb{R}^{d}$}

In this section we prove the inverse theorems for convolutions $\nu\conv\mu$
for $\mu\in\mathcal{P}(\mathbb{R}^{d})$ and $\nu\in\mathcal{P}(G_{0})$,
where $G_{0}$ is the group of isometries of $\mathbb{R}^{d}$. Our
strategy is to linearize the action $G_{0}\times\mathbb{R}^{d}\rightarrow\mathbb{R}^{d}$
and apply the Euclidean inverse theorem. 

Recall that the elements of $G_{0}$ are denoted $g=U+a$, with $U$
an orthogonal matrix, $a\in\mathbb{R}^{d}$, and $gx=Ux+a$. Given
$g\in G_{0}$ we denote the associated matrix and vector by $U_{g}$
and $a_{g}$. Also recall that $S_{t}$ is the scalar map $S_{t}(x)=2^{t}x$,
and introduce the translation map 
\[
\tau_{s}(x)=x+s.
\]

\subsection{Concentration and saturation on random subspaces}

This section contains additional technical results on concentration
and saturation of components of a measures. Our first goal is to show
that if two measures $\eta,\theta\in\mathcal{P}(\mathbb{R}^{d})$
are such that with high probability pairs of components $\eta^{x,i}$,
$\theta^{y,i}$ are highly concentrated and saturated, respectively,
on a subspace $V=V^{(i,x,y)}$, then we can assume that $V^{(i,x,y)}$
is essentially independent of $x,y$. 
\begin{prop}
\label{prop:derandomizing-subspaces}For every $\varepsilon>0$ and
$m\in\mathbb{N}$, there are $\varepsilon'=\varepsilon'(\varepsilon,m)\rightarrow0$
and $m'=m'(\varepsilon,m)\rightarrow\infty$ as $\varepsilon\rightarrow0$
and $m\rightarrow\infty$, such that the following holds. Suppose
that $\theta,\eta$ are independent $\mathcal{P}(\mathbb{R}^{d})$-valued
random variables defined on a probability space $(\Omega,\mathcal{F},\mathbb{P})$
and that $V=V(\theta,\eta)\leq\mathbb{R}^{d}$ is a linear subspace
determined by the random measures $\theta,\eta$ (hence $V$ is random).
If 
\begin{equation}
\mathbb{P}\left(\theta\mbox{ is }(V,\varepsilon',m')\mbox{-saturated and }\eta\mbox{ is }(V,\varepsilon')\mbox{-concentrated}\right)>1-\varepsilon'.\label{eq:76a}
\end{equation}
Then there is a deterministic subspace $V_{*}$ such that 
\begin{equation}
\mathbb{P}_{0\leq i\leq m'}\left(\theta^{x,i}\mbox{ is }(V_{*},\varepsilon,m)\mbox{-saturated and }\eta^{y,i}\mbox{ is }(V_{*},\varepsilon)\mbox{-concentrated}\right)>1-\varepsilon.\label{eq:76b}
\end{equation}
(the probability in the last equation is over both the measures $\theta,\eta$
and $i,x,y$, independently).\end{prop}
\begin{proof}
It is a formal consequence of Proposition \ref{prop:concentration-and-saturation-change-of-parameters}
that it is enough to prove the statement under the assumption that
$m$ is large relative to $\varepsilon$. 

Fix $\varepsilon$. Assume $m,m'$ large relative to $\varepsilon$,
and $\varepsilon'$ small relative to the other parameters. We shall
show that  \eqref{eq:76a} implies \eqref{eq:76b}. 

Apply Proposition \ref{prop:concentration-subspace-of-a-measure}
to $\eta$ with parameter $\varepsilon'$. We obtain a random subspace
$V_{c}=V_{c}(\eta)$ and a constant $C_{c}\geq1$ (which we will later
assume is large compared to another constant $D$) such that, for
\[
\delta_{c}=C_{c}\cdot(\varepsilon')^{1/3^{d}},
\]
the measure $\eta$ is $(V_{c},\delta_{c})$-concentrated and if $\eta$
is $(W,\varepsilon')$-concentrated then $V_{c}\sqsubseteq W^{(\delta_{c})}$. 

Apply Proposition \ref{prop:saturation-subspace-of-a-measure} to
$\theta$ with parameter $m'$. We obtain a random subspace $V_{s}=V_{s}(\theta)$
such that, for a constant $C_{s}$ and 
\[
\delta_{s}=C_{s}\frac{\log(m')}{m'},
\]
the measure $\theta$ is $(V_{s},\delta_{s},m')$-saturated and if
$\theta$ is $(W,1/m',m')$-saturated, then $W\sqsubseteq V_{s}^{(\delta_{s})}$.
We shall assume that $\varepsilon'<1/m'$. Thus, if $\theta$ is $(W,\varepsilon',m')$-saturated
then $W\sqsubseteq V_{s}^{(\delta_{s})}$.

The random subspace $V$ satisfies \eqref{eq:76a}, so we have
\[
\mathbb{P}\left(V_{c}\sqsubseteq V^{(\delta_{c})}\mbox{ and }V\sqsubseteq V_{s}^{(\delta_{s})}\right)>1-\varepsilon'.
\]
Thus, writing 
\[
\delta=\delta_{c}+\delta_{s},
\]
we have 
\begin{equation}
\mathbb{P}(V_{c}\sqsubseteq V_{s}^{(\delta)})>1-\varepsilon'.\label{eq:73}
\end{equation}

Let $\mathcal{W}=\{W_{1},\ldots,W_{N}\}$ denote a minimal $\delta_{c}$-dense
sequence of subspaces with respect to the metric \eqref{eq:metric-on-subspaces}.
This metric is bi-Lipschitz equivalent to a smooth metric on the compact
manifold of subspaces, so $N\leq D\cdot\delta_{c}^{-[d^{2}/2]}$ for
some universal constant $D>1$ (here $[d^{2}/2]$ is the dimension
of the space \,of subspaces). Let 
\[
\mathcal{W}_{0}=\{W\in\mathcal{W}\,:\,\mathbb{P}(d(V_{c},W)<\delta_{c})>\frac{\delta_{c}}{N}\}.
\]
Apply Proposition \ref{prop:minimal-engulfing-subspace} to $\mathcal{W}_{0}$
with parameter $2\delta$ to obtain the parameter 
\[
\delta'=4\cdot2^{1/3^{d}}\cdot\delta^{1/3^{d}}
\]
and a non-trivial subspace $V_{*}$ such that
\begin{enumerate}
\item [a.] $W\sqsubseteq V_{*}^{(\delta')}$ for all $W\in\mathcal{W}_{0}$,
\item [b.] If $\widetilde{V}_{*}$ is another subspace such that $W\sqsubseteq\widetilde{V}_{*}^{(2\delta)}$
for all $W\in\mathcal{W}_{0}$, then $V_{*}\sqsubseteq\widetilde{V}_{*}^{(2\delta')}$. 
\end{enumerate}
We claim that $V_{*}$ is the desired subspace. Writing $\mathcal{W}_{1}=\mathcal{W}\setminus\mathcal{W}_{0}$,
\begin{eqnarray*}
\mathbb{P}(d(V_{c},W)\geq\delta_{c}\mbox{ for all }W\in\mathcal{W}_{0}) & = & P(V_{c}\notin\bigcup_{W\in\mathcal{W}_{0}}B_{\delta_{c}}(W))\\
 & \leq & \mathbb{P}(V_{c}\in\bigcup_{W\in\mathcal{W}_{1}}B_{\delta_{c}}(W))\\
 & \leq & \sum_{W\in\mathcal{W}_{1}}\mathbb{P}(V_{c}\in B_{\delta_{c}}(W))\\
 & \leq & |\mathcal{W}_{1}|\cdot\frac{\delta_{c}}{N}\\
 & < & \delta_{c},
\end{eqnarray*}
where in the first inequality we used the fact that $\bigcup_{W\in\mathcal{W}_{0}\cup\mathcal{W}_{1}}B_{\delta_{c}}(W)$
covers all subspaces, and in the last line we used $|\mathcal{W}_{1}|\leq|\mathcal{W}|=N$.
Hence
\[
\mathbb{P}(d(V_{c},W)<\delta_{c}\mbox{ for some }W\in\mathcal{W}_{0})>1-\delta_{c}.
\]
Consequently, by property (a) of $V_{*}$ and the fact that $\delta_{c}\leq\delta'$,
\begin{equation}
\mathbb{P}(V_{c}\sqsubseteq V_{*}^{(2\delta')})>1-\delta_{c}.\label{eq:71}
\end{equation}

Since $V_{c}$ is a function of $\eta$ and $V_{s}$ is a function
of $\theta$, and since $\eta,\theta$ are independent, also each
of the pairs $V_{c},\theta$ and $V_{c},V_{s}$ is independent. Therefore,
for a.e. value $\theta_{0}$ of $\theta$,
\begin{equation}
\mathbb{P}(d(V_{c},W)<\delta_{c}|\theta=\theta_{0})=\mathbb{P}(d(V_{c},W)<\delta_{c})>\frac{\delta_{c}}{N}\qquad\mbox{for all }W\in\mathcal{W}_{0}.\label{eq:70}
\end{equation}
Observe that
\begin{eqnarray*}
\frac{\delta_{c}}{N} & \geq & \frac{1}{D}\delta_{c}^{1+[d^{2}/2]}\\
 & = & \frac{C_{c}^{1+[d^{2}/2]}}{D}(\varepsilon')^{(1+[d^{2}/2])/3^{d}}\\
 & > & (\varepsilon')^{1/2},
\end{eqnarray*}
where, to justify the last inequality, we increase the constant $C_{c}$
if necessary to ensure $C_{c}^{1+[d^{2}/2]}/D\geq1$, and note that
$(1+[d^{2}/2])/3^{d}<1/2$. Thus if a fixed measure $\theta_{0}$
satisfies 
\begin{equation}
\mathbb{P}(V_{c}\sqsubseteq V_{s}^{(\delta)}|\theta=\theta_{0})>1-\sqrt{\varepsilon'}\label{eq:49}
\end{equation}
then, by \eqref{eq:70} and \eqref{eq:49}, for all $W\in\mathcal{W}_{0},$
\begin{eqnarray*}
\mathbb{P}(W\sqsubseteq V_{s}^{(2\delta)}|\theta=\theta_{0}) & \geq & \mathbb{P}(d(V_{c},W)<\delta_{c}\mbox{ and }V_{c}\sqsubseteq V_{s}^{(\delta)}|\theta=\theta_{0})\\
 & \geq & \mathbb{P}(d(V_{c}\sqsubseteq V_{s}^{(\delta)}|\theta=\theta_{0})-(1-\mathbb{P}(d(V_{c},W)<\delta_{c}|\theta=\theta_{0}))\\
 & > & (1-\sqrt{\varepsilon'})-(1-\sqrt{\varepsilon'})\\
 & = & 0,
\end{eqnarray*}
Since $V_{s}$ is a function of $\theta$, this says that for $\theta_{0}$
satisfying \eqref{eq:49} we have $W\sqsubseteq V_{s}^{(2\delta)}$
for each $W\in\mathcal{W}_{0}$; consequently, by property (b) of
the definition of $V_{*}$, for such $\theta_{0}$ we have that $V_{*}\sqsubseteq V_{s}^{(2\delta')}$. 

By Markov's inequality and \eqref{eq:73}, the relation \eqref{eq:49}
holds with probability $1-\sqrt{\varepsilon'}$ over the choice of
$\theta_{0}$. Thus we conclude 
\begin{equation}
\mathbb{P}(V_{*}\sqsubseteq V_{s}^{(2\delta')})>1-\sqrt{\varepsilon'}.\label{eq:72}
\end{equation}
Combining \eqref{eq:71} and \eqref{eq:72} and using $\sqrt{\varepsilon'}\leq\delta_{c}$,
we find that
\[
\mathbb{P}(V_{c}\sqsubseteq V_{*}^{(2\delta')}\mbox{ and }V_{*}\sqsubseteq V_{s}^{(2\delta')})>1-2\delta_{c}.
\]

Finally, fix $\eta,\theta$ and associated to $V_{s},V_{c}$ belonging
to this event, we have that $\eta$ is $(V_{c},\delta_{c})$-concentrated.
Therefore by Lemma \ref{lem:concentrated-measures-have-concentrated-components},
\[
\mathbb{P}_{0\leq i\leq m'}\left(\eta^{x,i}\mbox{ is }(V_{c},\sqrt{2^{i}\delta_{c}})\mbox{-concentrated}\right)>1-\sqrt{\delta_{c}}
\]
(here and below the randomness is over $i$, with $\eta,\theta$ fixed).
Since $V_{c}\sqsubseteq V_{*}^{(2\delta')}$, and assuming as we may
that $\varepsilon'$, and hence $\delta_{c}$, is small enough relative
to $\varepsilon,m'$, this implies
\begin{equation}
\mathbb{P}_{0\leq i\leq m'}\left(\eta^{x,i}\mbox{ is }(V_{*},\varepsilon)\mbox{-concentrated}\right)>1-\frac{\varepsilon}{3}.\label{eq:77}
\end{equation}
Similarly, $\theta$ is $(V_{s},\delta_{s},m')$-saturated, so arguing
in the same manner and using  Lemma \ref{lem:saturation-passes-to-components},
\[
\mathbb{P}_{0\leq i\leq m'}\left(\theta^{y,i}\mbox{ is }(V_{s},\sqrt{d\delta_{s}+O(\frac{m}{m'})},m)\mbox{-saturated}\right)>1-O(\sqrt{\delta_{s}+\frac{m}{m'}}).
\]
Assuming $\varepsilon'$ is small enough and $m'$ large enough relative
to $\varepsilon,m$, the constant $\delta'$ can be assumed arbitrarily
small compared to $\varepsilon,m$. Since $V_{*}\sqsubseteq V_{s}^{(2\delta')}$
we have $d(V_{*},\pi_{V_{s}}V_{*})<2\delta'$, so by Lemma \ref{lem:concentration-saturation-uniformity-under-change-of-subspace}
and \eqref{enu:saturation-subspaces}, and assuming the parameters
satisfy the appropriate relationship, we have 
\begin{equation}
\mathbb{P}_{0\leq i\leq m'}\left(\theta^{y,i}\mbox{ is }(V_{*},\varepsilon,m)\mbox{-saturated}\right)>1-\frac{\varepsilon}{3}.\label{eq:78}
\end{equation}
Thus, combining \eqref{eq:77} and \eqref{eq:78} for $\eta,\theta$
in the event in \eqref{eq:72}, we have 
\[
\mathbb{P}_{0\leq i\leq m'}\left(\theta^{y,i}\mbox{ is }(V_{*},\varepsilon,m)\mbox{-saturated and }\eta^{x,i}\mbox{ is }(V_{*},\varepsilon)\mbox{-concentrated}\right)>1-\frac{2}{3}\varepsilon.
\]
Using \eqref{eq:72} and assuming as we may that $\sqrt{\varepsilon'}<\varepsilon/3$,
we obtain \eqref{eq:76b}.\end{proof}
\begin{cor}
\label{cor:derandomizing-subspaces}Let $\varepsilon>0$ and $m\in\mathbb{N}$.
Then there exist $\varepsilon''=\varepsilon''(\varepsilon,m)\rightarrow0$
and $m''=m''(\varepsilon,m)\rightarrow\infty$ as $\varepsilon\rightarrow0$
and $m\rightarrow\infty$, such that for all large enough $n$, the
following holds. Suppose that we are given subspaces $V^{(i,x,y)}$
for $0\leq i\leq n$ and $x,y\in[0,1]^{d}$ and measures $\theta,\eta\in\mathcal{P}([0,1]^{d})$
such that
\[
\mathbb{P}_{0\leq i\leq n}\left(\begin{array}{c}
\theta^{x,i}\mbox{ is }(V^{(i,x,y)},\varepsilon'',m'')\mbox{-saturated and }\\
\eta^{y,i}\mbox{ is }(V^{(i,x,y)},\varepsilon'')\mbox{-concentrated}
\end{array}\right)>1-\varepsilon'',
\]
then there are subspaces $V^{i}\leq\mathbb{R}^{d}$ such that 
\[
\mathbb{P}_{0\leq i\leq n}\left(\begin{array}{c}
\theta^{x,i}\mbox{ is }(V^{i},\varepsilon,m)\mbox{-saturated and }\\
\eta^{y,i}\mbox{ is }(V^{i},\varepsilon)\mbox{-concentrated}
\end{array}\right)>1-\varepsilon.
\]
\end{cor}
\begin{proof}
Apply the previous proposition to obtain $\varepsilon'=\varepsilon'(\frac{1}{2}\varepsilon,m)$
and $m'=m'(\frac{1}{2}\varepsilon,m)$, and set $m''=m'$ and $\varepsilon''=\min\{(\varepsilon')^{2},1/m'\}$.

For $0\leq k\leq n$, let 
\[
p_{k}=\mathbb{P}_{i=k}\left(\begin{array}{c}
\theta^{x,i}\mbox{ is }(V^{(i,x,y)},\varepsilon'',m'')\mbox{-saturated and }\\
\eta^{y,i}\mbox{ is }(V^{(i,x,y)},\varepsilon'')\mbox{-concentrated}
\end{array}\right),
\]
and assume as in the hypothesis that $\frac{1}{n+1}\sum_{k=0}^{n}p_{k}>1-\varepsilon''$.
Let $I\subseteq\{0,\ldots,n\}$ denote the set of $k$ such that $p_{k}>1-\sqrt{\varepsilon''}=1-\varepsilon'$,
so by Markov, $|I|>(1-\varepsilon')(n+1)$.

For $i\in I$, consider the random and independently chosen components
$\theta^{x,i}$, $\eta^{y,i}$ and the subspace $V_{i}=V^{i,x,y}$.
Without loss of generality we may assume that $V^{i,x,y}$ depend
only on $\theta^{x,i}$ and $\eta^{y,i}$, since the only stated property
of $V^{i,x,y}$ involves these measures. From the previous proposition
and our choice of $\varepsilon',m'$, we conclude that there exists
a subspace $V^{i}$ such that 
\[
\mathbb{P}_{i\leq j\leq i+m'}\left(\begin{array}{c}
\theta^{x,i}\mbox{ is }(V^{i},\frac{1}{2}\varepsilon,m)\mbox{-saturated and }\\
\eta^{y,i}\mbox{ is }(V^{i},\frac{1}{2}\varepsilon)\mbox{-concentrated}
\end{array}\right)>1-\frac{\varepsilon}{2}.
\]
The remainder of the argument involves choosing one of these subspaces
$V^{i(j)}$, for every $j\in\bigcup_{u\in I}[u,u+m']$. The details
are identical to the proof of Proposition \ref{prop:concentration-and-saturation-change-of-parameters}.
\end{proof}
We will actually need a more general version of the last corollary,
but, as the proof is identical to the one above, we only give the
statement.
\begin{cor}
\label{cor:derandomizing-subspaces-with-isometries}Let $\varepsilon>0$
and $m\in\mathbb{N}$. Then there exist $\varepsilon''=\varepsilon''(\varepsilon,m)\rightarrow0$
and $m''=m''(\varepsilon,m)\rightarrow\infty$ as $\varepsilon\rightarrow0$,
such that for all large enough $n$, the following holds. Suppose
that $\theta\in\mathcal{P}(G_{0})$, $\eta\in\mathcal{P}(\mathbb{R}^{d})$,
and that for $0\leq i\leq n$, $x\in\supp\eta$ and $g\in\supp\theta$
there are subspaces $V^{(i,x,g)}$ such that
\[
\mathbb{P}_{0\leq i\leq n}\left(\begin{array}{c}
\eta^{x,i}\mbox{ is }(V^{(i,x,y)},\varepsilon'',m'')\mbox{-saturated and }\\
S_{i}U_{g}^{-1}(\theta_{g,i}\conv x)\mbox{ is }(V^{(i,x,y)},\varepsilon'')\mbox{-concentrated}
\end{array}\right)>1-\varepsilon''.
\]
Then there are subspaces $V^{i}\leq\mathbb{R}^{d}$ such that 
\[
\mathbb{P}_{0\leq i\leq n}\left(\begin{array}{c}
\eta^{x,i}\mbox{ is }(V^{i},\varepsilon,m)\mbox{-saturated and }\\
S_{i}U_{g}^{-1}(\theta^{y,i}\conv x)\mbox{ is }(V^{i},\varepsilon)\mbox{-concentrated}
\end{array}\right)>1-\varepsilon.
\]

\end{cor}

\subsection{From concentration of Euclidean components to $G_{0}$-components}

We turn our attention to measures $\eta\in\mathcal{P}(\mathbb{R}^{d})$
of the form $\eta=\theta\conv x$ for some $\theta\in\mathcal{P}(G_{0})$
and $x_{0}\in\mathbb{R}^{d}$. Our goal is to show that the concentration
properties of typical components $\eta^{y,i}$ translates to similar
properties of the ``components'' $\theta_{g,i}\conv x$. The issue
which we must overcome is that $\theta_{g,i}\conv x$ is supported
on $\mathcal{D}_{i}^{G}(g)\conv x$, and this set generally intersects
more than one dyadic cell of $\mathcal{D}_{i}^{d}$. Thus even if
$\eta$ is highly concentrated on a translate of a subspace $W$ on
each of these cells, taken together all one can say is that  $\theta_{g,i}\conv x$
is concentrated on the union of several translates of $W$. 

For a linear subspace $W\leq\mathbb{R}^{d}$ we say that a measure
$\eta\in\mathcal{P}(\mathbb{R}^{d})$ is $(W,\delta)^{m}$-concentrated
if for some $m'\leq m$ there are $m'$ translates $W_{1},\ldots,W_{m'}$
of $W$ such that $\eta(\bigcup_{u=1}^{m'}W_{u}^{(\delta)})\geq1-\delta$.
Thus $(W,\delta)^{1}$-concentration is the same as $(W,\delta)$-concentration.
\begin{lem}
Let $R>0$ , let $\theta\in\mathcal{P}(G_{0})$ and $x\in[-R,R^{d}]$.
Suppose that $\delta>0$, $m\in\mathbb{N}$ and that $\theta\conv x$
is $(W,\delta)^{m}$-concentrated. Then for \textup{$n=[\frac{1}{2}\log(1/\delta)]$
and} $\delta'=O_{R,m}(\frac{\log\log(1/\delta)}{\log(1/\delta)})$
we have 
\[
\mathbb{P}_{0\leq i\leq n}\left(S_{i}(\theta_{g,i}\conv x)\mbox{ is }(W,\delta')\mbox{-concentrated}\right)>1-\delta'.
\]
\end{lem}
\begin{proof}
Although the ``rescaled component'' $\theta^{g,i}$ is not defined,
it will be convenient to introduce the notation 
\[
\theta^{g,i}\conv x=S_{i}(\theta_{g,i}\conv x),
\]
and define the distribution on these ``components'' in the usual
manner. Note that there is a constant $C=C(R)\geq1$ such that $\theta_{g,i}\conv x$
is supported on a set of diameter $\leq C2^{-i}$, and $\theta^{g,i}\conv x$
is supported on a set of diameter $\leq C$.

Let $W_{1},\ldots,W_{m}$ be affine subspaces parallel to $W$ verifying
that $\theta\conv x$ is $(W,\delta)^{m}$-concentrated. We may assume
the $W_{u}$ are distinct. For $u\neq v$ let 
\[
d_{u,v}=d(W_{u},W_{v})=\min\{d(x,y)\,:\,x\in W_{u}\,,\,y\in W_{v}\}.
\]
Notice that for any $1\leq k\leq n$,
\begin{itemize}
\item If $2^{-k}<\frac{1}{4C}d_{u,v}$ for some $u,v$, then for any $g$
the measure $\theta_{g,k}\conv x$ is supported on a set of diameter
at most $C2^{-k}<\frac{1}{4}d_{u,v}$, and on the other hand $\sqrt{\delta}\leq2^{-n}\leq2^{-k}\leq\frac{1}{4}d_{u,v}$,
hence $\theta_{g,k}\conv x$ gives positive mass to at most one the
sets $W_{u}^{(\sqrt{\delta})},W_{v}^{(\sqrt{\delta})}$.
\item Let $I_{g,k}\subseteq\{1,\ldots,m\}$ be the set of indices $u$ such
that $(\theta_{g,k}\conv x)(W_{u}^{(\delta)})>0$. Given $\rho>0$,
if all distinct $u,v\in I_{g,k}$ satisfy $d_{u,v}\leq\rho2^{-k}$,
then there is a translate $W_{g,k}$ of $W$ such that $\bigcup_{u\in I_{g,k}}W_{u}^{(\delta)}\cap\supp(\theta_{g,k}\conv x)\subseteq W_{g,k}^{(\rho2^{-k}+2\delta)}$,
so $\theta^{g,k}\conv x$ is $(W,\rho+2^{k+1}\delta))$-concentrated.
\end{itemize}
Now, for $0\leq k\leq n$ we have identity
\[
\theta\conv x=\mathbb{E}_{i=k}\left(\theta_{g,i}\conv x\right).
\]
Using the hypothesis that $(\theta\conv x)(\bigcup_{u=1}^{m}W_{u}^{(\delta)})>1-\delta$
and Markov's inequality we conclude that
\begin{equation}
\mathbb{P}_{i=k}\left((\theta_{g,i}\conv x)(\bigcup_{u=1}^{m}W_{u}^{(\delta)})>1-\sqrt{\delta}\right)>1-\sqrt{\delta}.\label{eq:111}
\end{equation}
Fix a small parameter $\rho>0$, and suppose that $k$ satisfies 
\begin{equation}
\mbox{For each }1\leq u<v\leq m\mbox{ either }2^{-k}<\frac{1}{4C}d_{u,v}\mbox{ or }d_{u,v}\leq\rho2^{-k},\label{eq:45}
\end{equation}
or, equivalently, that $k$ does not belong to any of the intervals
$J_{u,v}=[\log\frac{\rho}{d_{u,v}},\log\frac{4C}{d_{u,v}})$. Then,
setting 
\[
\sigma=\sigma(\rho)=\max\{\sqrt{\delta},\rho+2^{k+1}\delta\}
\]
the two observations above and \eqref{eq:111} imply
\[
\mathbb{P}_{i=k}\left(\theta^{g,k}\conv x\mbox{ is }(W,\sigma)\mbox{-concentrated}\right)>1-\sqrt{\delta}.
\]
Note that $\delta2^{k}\leq\delta2^{n}\leq\sqrt{\delta}$, so in fact
$\sigma\leq\rho+2\sqrt{\delta}$. 

Next, since the length of $J_{u,v}$ is $\log\frac{4C}{\rho}$ and
there are at most $m(m-1)$ distinct values of $1\leq u,v\leq m$,
the fraction of $0\leq k\leq n$ which satisfy \eqref{eq:45} is at
least $1-m^{2}\log(\frac{4C}{\rho})/n$. Averaging the last equation
over $k=0,\ldots,n$, we conclude that
\begin{eqnarray*}
\mathbb{P}_{0\leq k\leq n}\left(\theta^{g,k}\conv x\mbox{ is }(W,\rho+\sqrt{\delta})\mbox{-concentrated}\right) & > & 1-\sqrt{\delta}-\frac{m^{2}\log(4C/\rho)}{n}\\
 & = & 1-\sqrt{\delta}-O_{m}(\frac{\log(1/\rho)}{\log(1/\delta)}).
\end{eqnarray*}
Choosing $\rho=\frac{1}{\log(1/\delta)}$ gives the desired result.\end{proof}
\begin{prop}
For every $\varepsilon>0$ and $R>0$ there exists $n=n(\varepsilon,R)$
(with $n(\varepsilon,R)\rightarrow\infty$ as $\varepsilon\rightarrow0$)
and an $\delta=\delta(\varepsilon,R)>0$ (with $\delta(\varepsilon,R)\rightarrow0$
as $\varepsilon\rightarrow0$) such that the following holds. Let
$\nu\in\mathcal{P}(G_{0})$ and $x_{0}\in[-R,R]^{d}$, and write $\eta=\nu\conv x_{0}$.
Let $V<\mathbb{R}^{d}$ be a linear subspace and $k\in\mathbb{N}$
such that 
\[
\mathbb{P}_{j=k}\left(\eta^{x,j}\mbox{ is }(V,\delta)\mbox{-concentrated}\right)>1-\delta.
\]
Then 
\[
\mathbb{P}_{k\leq j\leq k+n}\left(S_{j}(\nu_{g,j}\conv x)\mbox{ is }(V,\varepsilon)\mbox{-concentrated}\right)>1-\varepsilon.
\]
\end{prop}
\begin{proof}
Fix $\varepsilon,n,\delta$ for the moment and assume that the hypothesis
holds. Consider the identities 
\[
\mathbb{E}_{j=k}(\eta_{x,j})=\eta=\mathbb{E}_{j=k}(\nu_{g,j}\conv x).
\]
This means that the measures $\nu_{g,k}\conv x$ are ($\nu$-almost-surely
over choice of $g$) absolutely continuous with respect to the weighted
average of the components $\eta_{x,k}$. In fact, since each $\nu_{g,k}\conv x$
is supported on a set that intersects $m=O_{R}(1)$ level-$k$ dyadic
cells, each ``component'' $\nu_{g,k}\conv x$ is absolutely continuous
with respect to the average of these $O(1)$ components $\eta_{x,k}$.
Most of these components are $(V,\delta)$-concentrated, so a Markov-inequality
argument (similar to the one in Lemma \ref{lem:saturation-at-level-n}
below) shows that 
\[
\mathbb{P}_{j=k}\left(\nu_{g,j}\conv x\mbox{ is }(V,2^{-k}\delta')^{m}\mbox{-concentrated}\right)>1-\delta',
\]
where $\delta'\rightarrow0$ as $\delta\rightarrow0$. Equivalently,
\[
\mathbb{P}_{j=k}\left(S_{j}(\nu_{g,j}\conv x)\mbox{ is }(V,\delta')^{m}\mbox{-concentrated}\right)>1-\delta'.
\]
Apply the previous lemma to each component $\theta=\nu_{g,k}$ in
the event above with parameter $\delta'$. Taking $\delta''=O_{R,m}(\log\log(1/\delta')/\log(1/\delta'))$
and $n=[\frac{1}{2}\log(1/\delta')]$, the conclusion is 
\[
\mathbb{P}_{k\leq i\leq k+n}\left(S_{i}(\nu_{g,i}\conv x)\mbox{ is }(V,\delta'')\mbox{-concentrated}\right)>1-\delta'',
\]
and $\delta''$ can be made arbitrarily small by taking $\delta$
small. This is what was claimed.\end{proof}
\begin{prop}
\label{prop:concentration-passes-to-group-components}For every $\delta>0$,
$R>0$ and $m\in\mathbb{N}$, if $m'>m'(\delta,m,R)$, $0<\delta'<\delta'(\delta,m,R)$,
then for all large enough $n$ (depending on previous parameters),
the following holds. Let $\mu\in\mathcal{P}(\mathbb{R}^{d})$, $\nu\in\mathcal{P}(G_{0})$
and $x_{0}\in\mathbb{R}^{d}$, and write $\eta=\nu\conv x_{0}$. Let
$V_{0},V_{1},\ldots,V_{n}\leq\mathbb{R}^{d}$ be linear subspaces,
and suppose that 
\[
\mathbb{P}_{0\leq j\leq n}\left(\begin{array}{c}
\mu^{x,j}\mbox{ is }(V_{j},\delta',m')\mbox{-saturated and}\\
\eta^{y,j}\mbox{ is }(V_{j},\delta')\mbox{-concentrated}
\end{array}\right)>1-\delta'.
\]
Then there are subspaces $V'_{0},V'_{1},\ldots,V'_{n}$ such that
\[
\mathbb{P}_{0\leq j\leq n}\left(\begin{array}{c}
\mu^{x,j}\mbox{ is }(V'_{j},\delta,m)\mbox{-saturated and}\\
S_{j}(\nu_{g,j}\conv x_{0})\mbox{ is }(V'_{j},\delta)\mbox{-concentrated}
\end{array}\right)>1-\delta.
\]
\end{prop}
\begin{proof}
Let $\delta,R,m$ be given. Fix a small auxiliary parameter $\delta_{1}$
which we will specify later and let $m_{1}$ be the number $n(\delta_{1},R)$
from the previous proposition, in particular it can be made arbitrarily
large by making $\delta_{1}$ small. Let $\varepsilon=\varepsilon(\delta_{1})$
as in the previous proposition, and let $\delta_{2}=\frac{1}{2}\varepsilon^{2}$.
Then for all small enough $\delta'$ and all large enough $m'$, the
hypothesis implies, by Proposition \ref{prop:concentration-and-saturation-change-of-parameters},
that there are subspaces $V''_{j}$ such that

\[
\mathbb{P}_{0\leq j\leq n}\left(\begin{array}{c}
\mu^{x,j}\mbox{ is }(V''_{j},\delta_{2},m_{1})\mbox{-saturated and}\\
\eta^{y,j}\mbox{ is }(V''_{j},\delta_{2})\mbox{-concentrated}
\end{array}\right)>1-\delta_{2}.
\]
By Markov's inequality, the set $I\subseteq\{0,\ldots,n\}$ consisting
of $k$ such that 
\[
\mathbb{P}_{j=k}\left(\begin{array}{c}
\mu^{x,j}\mbox{ is }(V''_{j},\delta_{2},m_{1})\mbox{-saturated and}\\
\eta^{y,j}\mbox{ is }(V''_{j},\delta_{2})\mbox{-concentrated}
\end{array}\right)>1-\sqrt{\delta_{2}}
\]
has size $|I|\geq(1-\sqrt{\delta_{2}})(n+1)$. Since $\delta_{2}<\varepsilon^{2}$,
by our choice of $m_{1}$ and $\varepsilon$, for each $k\in I$ we
have 
\[
\mathbb{P}_{k\leq j\leq k+m_{1}}\left(S_{j}(\nu_{g,j}\conv x_{0})\mbox{ is }(V''_{k},\delta_{1})\mbox{-concentrated}\right)>1-\delta_{1}.
\]
Also, applying Lemma \ref{lem:saturation-passes-to-components} to
each $(V''_{k},\delta_{2},m_{1})$-saturated component $\mu^{x,k}$
of $\mu$, we find that 
\[
\mathbb{P}_{k\leq j\leq k+m_{1}}\left(\mu^{x,j}\mbox{ is }(V''_{k},\sqrt{d\sqrt{\delta_{2}}+O(\frac{m}{m_{1}})},m)\mbox{-saturated}\right)>1-\sqrt{d\sqrt{\delta_{2}}+O(\frac{m}{m_{1}})}.
\]
Now, by choosing $\delta_{1}$ small enough we can ensure that $\varepsilon$
and $\delta_{2}$ is small, and $m_{1}$ is large, relative to $\delta,m$.
With suitable choices, one now argues as in the proof of Proposition
\ref{prop:concentration-and-saturation-change-of-parameters} to combine
the last two equations over all $k\in I$ and define $V'_{j}$ with
the desired properties.
\end{proof}

\subsection{Entropy and the $G_{0}$-action on $\mathbb{R}^{d}$}

For $g=U+a$ and $g'=U'+a'$ in $G_{0}$ and $x,x'\in\mathbb{R}^{d}$,
\begin{equation}
g'x'-gx=(U'-U)x+U'(x'-x)+(a'-a).\label{eq:discrete-differentiation-of-G0-action}
\end{equation}
In particular 
\[
\left\Vert gx-g'x'\right\Vert \leq\left\Vert U-U'\right\Vert \left\Vert x\right\Vert +\left\Vert U'\right\Vert \left\Vert x-x'\right\Vert +\left\Vert a-a'\right\Vert ,
\]
so if $g,g'$ are in a common level-$k$ dyadic cell and $x,x'\in[-R,R]^{d}$
are in a common level-$k$ dyadic cell, then $\left\Vert gx-g'x'\right\Vert =O_{R}(2^{-k})$.
In particular if $\nu\in\mathcal{P}(G_{0})$ and $\mu\in\mathcal{P}([-R,R]^{d})$
are both supported on level-$k$ dyadic cells, then $\nu\conv\mu$
is supported on a set of diameter $O_{R}(2^{-k})$. 

For a probability measure $\theta$ on $\mathbb{R}^{d}$ or $G_{0}$
it will be convenient in this section to write
\[
H_{i,n}(\theta)=\frac{1}{n}H(\theta,\mathcal{D}_{i+n}).
\]
(This differs from $H_{i+n}(\theta)$ because we normalize by $1/n$
instead of $1/(i+n)$). In particular $H_{n}(\theta)=H_{0,n}(\theta)$.
By the previous paragraph, if $\theta\in\mathcal{P}(\mathbb{R}^{d})$
and $\nu\in\mathcal{P}(G_{0})$ are supported on level-$i$ dyadic
cells then $\nu\conv\theta$ is supported on $O(1)$ level-$i$ dyadic
cells, so
\[
H_{i,n}(\nu\conv\theta)=\frac{1}{n}H(\nu\conv\theta,\mathcal{D}_{i+n}|\mathcal{D}_{i})+O(\frac{1}{n}).
\]
Also observe that for $\theta$ as above, $H_{i,n}(\theta)=H_{n}(S_{i}\theta)+O(1/n)$
(Lemma \ref{lem:entropy-combinatorial-properties} \eqref{enu:entropy-change-of-scale}).

We now address the issue, described in Section \ref{sub:inverse-theorem-isometries},
of pairs  $\nu\in\mathcal{P}(G_{0})$ and $x\in\mathbb{R}^{d}$ such
that $\nu$ has substantial entropy but $\nu\conv x$ does not (e.g.
because $\nu$ is supported close to $\stab_{G_{0}}(x)$. 
\begin{defn}
For $\sigma>0$ we say that $x_{1},\ldots,x_{d+1}\in\mathbb{R}^{d}$
are $\sigma$-independent if each $x_{i}$ is at distance at least
$\sigma$ from the affine subspace spanned by the others. 
\end{defn}
The action of an element $g\in G$ is determined by its action on
any $(d+1)$-tuple of affinely independent vectors in $\mathbb{R}^{d}$,
in particular of any $\sigma$-independent $(d+1)$-tuple.
\begin{prop}
\label{prop:entropy-passes-to-G-orbits}For every $\varepsilon,\sigma,R>0$,
and $k\in\mathbb{Z}$ and $m\in\mathbb{N}$, the following holds.
For every $\sigma$-independent sequence $x_{1},\ldots,x_{d+1}\in[-R,R]^{d}$
and every $\nu\in\mathcal{P}(G_{0})$ that is supported on a level-$k$
dyadic cell, if 
\[
H_{k,m}(\nu\conv x_{i})<\varepsilon\quad\mbox{ for all }i=1,\ldots,d+1,
\]
then 
\[
H_{k,m}(\nu)<(d+1)\varepsilon+O_{\sigma,R}(\frac{1}{m}).
\]
\end{prop}
\begin{proof}
Since $\nu$ is supported on a level-$k$ dyadic cell, each $\nu\conv x_{i}$
is supported on $O(1)$ level-$k$ dyadic cells, and therefore
\[
H(\nu\conv x_{i},\mathcal{D}_{k})=O(1)
\]
Thus the hypothesis is $\frac{1}{m}H(\nu\conv x_{i},\mathcal{D}_{k+m})<\varepsilon$
for $i=1,\ldots,d+1$, and it is enough to prove that $\frac{1}{m}H(\nu,\mathcal{D}_{k+m}^{G})<(d+1)\varepsilon+O_{\sigma,R}(1/m)$. 

Define the map $f:G_{0}\rightarrow(\mathbb{R}^{d})^{(d+1)}$ by $g\rightarrow(gx_{1},\ldots,gx_{d+1})$.
Then $f$ is a diffeomorphism and one may easily verify that $f$
is uniformly bi-Lipschitz with its image,\footnote{This fact depends of course on the metric with which we endowed $G_{0}$.
In general when applying this type of argument to a non-compact group
this is one point where the choice of metric must be carefully considered.} with Lipschitz constants of $f$ and $f^{-1}$ depending only on
$\sigma$ and $R$. Thus (e.g. by Lemma \ref{lem:entropy-weak-continuity-properties}
\eqref{enu:entropy-combinatorial-distortion} applied to $f^{-1}\mathcal{D}_{k+m}^{d(d+1)}$
and $\mathcal{D}_{k+m}^{G}$), 
\[
|\frac{1}{m}H(f\nu,\mathcal{D}_{k+m}^{d(d+1)})-\frac{1}{m}H(\nu,\mathcal{D}_{k+m}^{G})|=O_{\sigma,R}(\frac{1}{m}).
\]
Let $\pi_{i}:(\mathbb{R}^{d})^{d+1}\rightarrow\mathbb{R}^{d}$ denote
the projection to the $i$-th copy of $\mathbb{R}^{d}$. Then $\nu\conv x_{i}=\pi_{i}(f\nu)$.
Therefore, if $\frac{1}{m}H(\nu\conv x_{i},\mathcal{D}_{k+m})<\varepsilon$
for all $i=1,\ldots,d+1$, then $\frac{1}{m}H(\pi_{i}f\nu,\mathcal{D}_{k+m})<\varepsilon$
for all $i=1,\ldots,d+1$, and so $\frac{1}{m}H(f\nu,\mathcal{D}_{k+m}^{d(d+1)})\leq(d+1)\varepsilon$
(because $\mathcal{D}_{k+m}^{d(d+1)}=\bigvee\pi_{i}^{-1}\mathcal{D}_{k+m}^{d}$,
and using Lemma \ref{lem:entropy-combinatorial-properties} \eqref{enu:entropy-conditional-formula}).
The claim follows.
\end{proof}
Recall the definition of $(\varepsilon,\sigma)$-non-affine measures,
Definition \ref{def:epsilon-sigma-continuity}. 
\begin{lem}
If $\mu\in\mathcal{P}(\mathbb{R}^{d})$ is $(\varepsilon,\sigma)$-non-affine
and $A\subseteq\mathbb{R}^{d}$ is a Borel set with $\mu(A)>((d+1)\varepsilon)^{1/(d+1)}$,
then there exists a $\sigma$-independent sequence $x_{1},\ldots,x_{d+1}\in A$.\end{lem}
\begin{proof}
Let $X_{1},\ldots,X_{d+1}$ be independent $\mathbb{R}^{d}$-valued
random variables, each distributed according to $\mu$. Let $V_{i}$
be the (random) affine subspace spanned by the $d$ vectors $\{X_{j}\}_{j\neq i}$.
For each $i$ the vector $X_{i}$ is independent of $V_{i}$, and
$X_{i}$ is distributed according to $\mu$, so, since $\mu$ is $(\varepsilon,\sigma)$-non-affine,
\[
\mathbb{P}(X_{i}\notin V_{i}^{(\sigma)})=\mu(\mathbb{R}^{d}\setminus V_{i}^{(\sigma)})>1-\varepsilon.
\]
This implies
\[
\mathbb{P}(X_{i}\notin V_{i}^{(\sigma)}\mbox{ for all }i=1,\ldots,d+1)>1-(d+1)\varepsilon.
\]
Therefore, if $\mu(A)>((d+1)\varepsilon)^{1/(d+1)}$, 
\begin{multline*}
\mathbb{P}(X_{i}\notin V_{i}^{(\sigma)}\mbox{ and }X_{i}\in A\mbox{ for all }i=1,\ldots,d+1)\\
\begin{aligned}\geq & \;\mathbb{P}(X_{i}\in A\mbox{ for all }i)-(d+1)\varepsilon\\
\geq & \;\mu(A)^{d+1}-(d+1)\\
> & \;0.
\end{aligned}
\end{multline*}
Any realization $X_{1},\ldots,X_{d+1}$ from the event above is $\sigma$-independent.\end{proof}
\begin{cor}
\label{cor:entropy-passes-to-G-orbits}Let $k\in\mathbb{Z}$ and let
$\nu\in\mathcal{P}(G_{0})$ be supported on a level-$k$ dyadic cell.
Then for every $\varepsilon,\sigma,R>0$, every $(\varepsilon,\sigma)$-non-affine
measure $\mu\in\mathcal{P}([-R,R]^{d})$, and for every $m\in\mathbb{N}$,
\[
\mu\left(x\in\mathbb{R}^{d}\,:\,H_{k,m}(\nu\conv x)>\frac{1}{d+1}H_{k,m}(\nu)-O_{\sigma,R}(\frac{1}{m})\right)>1-((d+1)\varepsilon)^{1/(d+1)}.
\]
\end{cor}
\begin{proof}
Let $c=c(\sigma,R)$ denote the constant in the error term of Proposition
\ref{prop:entropy-passes-to-G-orbits}. Let $A=\{x\in\mathbb{R}^{d}\,:\,H_{k,m}(\nu\conv x)\leq\frac{1}{d+1}H_{k,m}(\nu)-\frac{c}{m}\}$,
we claim that $\mu(A)\leq((d+1)\varepsilon)^{1/(d+1)}$. Otherwise,
by the previous lemma, there is an $(\varepsilon,\sigma)$-non-affine
tuple $x_{1},\ldots,x_{d+1}\in A$. By the Proposition \ref{prop:entropy-passes-to-G-orbits}
applied to $x_{1},\ldots,x_{d+1}$ and using the definition of $A$
we have 
\[
H_{k,m}(\nu)<(d+1)(\frac{1}{d+1}H_{k,m}(\nu)-\frac{c}{m})+\frac{c}{m}<H_{k,m}(\nu),
\]
which is a contradiction.
\end{proof}

\subsection{\label{sub:Components-and-entropy-on-G}Linearization of the $G_{0}$-action}

Next we utilize the differentiability of the action of $G_{0}\times\mathbb{R}^{d}\rightarrow\mathbb{R}^{d}$,
which implies that at small scales, a convolution $\nu\conv\mu$ of
$\nu\in\mathcal{P}(G)$ and $\mu\in\mathcal{P}(\mathbb{R}^{d})$ can
be well approximated by a Euclidean convolution. Since it is easy
to give an elementary argument, we do so.

Let $g_{0}=U_{0}+a_{0}$ and $g=U+a$ be elements of $G_{0}$ and
$x_{0},x\in\mathbb{R}^{d}$. Then we have the identity
\begin{equation}
g\conv x=g\conv x_{0}+U_{0}(x-x_{0})+(U-U_{0})(x-x_{0})\label{eq:discrete-differentiation-identity}
\end{equation}
Assuming further that $g,g_{0}$ belong to a common level-$k$ dyadic
cell in $G_{0}$ and $x,x_{0}$ belong to a common level-$k$ dyadic
cell in $\mathbb{R}^{d}$, we have $\left\Vert U-U_{0}\right\Vert =\left\Vert x-x_{0}\right\Vert =O(2^{-k})$,
so 
\begin{eqnarray}
g\conv x & = & g\conv x_{0}+U_{0}(x-x_{0})+O(2^{-2k}).\label{eq:discrete-differentiation-of-G0-action-2}
\end{eqnarray}
Recall that $\tau_{z}(y)=y+z$ is the translation map. It follows
from the above that if $\nu\in\mathcal{P}(G_{0})$ and $\mu\in\mathcal{P}(\mathbb{R}^{d})$
are supported on the level-$k$ dyadic cells containing $g_{0},x_{0}$
respectively, then for $f\in\lip(\mathbb{R}^{d})$ we have 
\begin{eqnarray*}
\int f\,d(\nu\conv\mu) & = & \int\int f(g\conv x)\,d\nu(g)\,d\mu(x)\\
 & = & \int\int f(g\conv x_{0}+U_{0}(x-x_{0}))\,d\nu(g)\,d\mu(x)+O(2^{-2k}\cdot\left\Vert f\right\Vert _{\lip})\\
 & = & \int f(y)\,d((\nu\conv x_{0})*(U_{0}\tau_{-x_{0}}\mu))(y)+O(2^{-2k}\cdot\left\Vert f\right\Vert _{\lip}).
\end{eqnarray*}
Let $\nu'$, $\mu'$ and $\theta$ be the measures obtained from $\nu\conv x$,
$U_{0}\tau_{x_{0}}\mu$ and $\nu\conv\mu$, respectively, by scaling
them by a factor of $2^{k}$ and translating them so that they are
supported on a closed ball $B$ of radius $O(1)$ at the origin. Define
a metric on $\mathcal{P}(B)$ by 
\[
d(\alpha,\beta)=\sup_{\left\Vert f\right\Vert _{\lip}=1}|\int fd\alpha-\int fd\beta|.
\]
It is well known that $d(\cdot,\cdot)$ is compatible with the weak-{*}
topology on $\mathcal{P}(B)$ (see e.g. \cite[Chapter 14]{Mattila95}),
and the calculation above implies that $d(\nu'*\mu',\theta)=O(2^{-k})$.
Thus when $k$ is large, by Lemma \ref{lem:entropy-weak-continuity-properties}
\eqref{enu:entropy-approximation}, $|H_{m}(\nu'*\mu')-H_{m}(\theta)|=O(1/m)$.
Restating this in terms of the original measure, we have shown:
\begin{lem}
\label{lem:linearization-of-entropy-of-G-convolution}For every $m\in\mathbb{N}$
and $k>k(m)$ the following holds. If $\mu\in\mathcal{P}(\mathbb{R}^{d})$
and $\nu\in\mathcal{P}(G_{0})$ are supported on level-$k$ dyadic
cubes, and $x_{0}\in\supp\mu$ and $g_{0}\in\supp\nu$, then 
\begin{eqnarray*}
H_{k,m}(\nu\conv\mu) & = & H_{k,m}\left((\nu\conv x_{0})*U_{0}\mu\right)+O(\frac{1}{m})\\
 & = & H_{k,m}\left((U_{0}^{-1}(\nu\conv x_{0}))*\mu\right)+O(\frac{1}{m})
\end{eqnarray*}

\end{lem}
We omitted the translation in the statement because it commutes with
convolution, and does not affect entropy more than the error term.
The second line follows from the first by applying $U_{0}^{-1}$ to
the convolution.

Reasoning similarly, let $g,g_{0}\in G_{0}$ belong to a common level-$\ell$
dyadic cube $D$, and $x,x_{0}\in\mathbb{R}^{d}$ belong to a common
level-$k$ dyadic cell. Then, using \eqref{eq:discrete-differentiation-identity},
and the fact that $\left\Vert (U-U_{0})(y-x)\right\Vert =O(2^{-\ell-k})$,
we have
\begin{eqnarray*}
gx & = & gx_{0}+U_{0}(x_{0}-x)+O(2^{-\ell-k})
\end{eqnarray*}
Thus, for $\nu\in\mathcal{P}(D)$ and $f\in\lip(\mathbb{R}^{d})$,
we have
\begin{eqnarray*}
\int f\,d(\nu\conv x) & = & \int f(gx),d\nu(g)\\
 & = & \int f(gx_{0}+U_{0}(x_{0}-x)+O(2^{-\ell-k})\,d\nu(g)\\
 & = & \int f(gx_{0}+U_{0}(x_{0}-x))\,d\nu(g)+O(2^{-\ell-k}\left\Vert f\right\Vert _{\lip}\\
 & = & \int f(y)\,d(\tau_{U_{0}(x_{0}-x)}(\nu\conv x))(y)+O(2^{-\ell-k}\left\Vert f\right\Vert _{\lip}).
\end{eqnarray*}
Now, $\nu\conv x$ and $\nu\conv x_{0}$ are measures supported on
sets of diameter $O_{\left\Vert x\right\Vert ,\left\Vert x_{0}\right\Vert }(2^{-\ell})$
(since $\nu$ is supported on a level-$\ell$ dyadic cell), so re-scaling
by $2^{\ell}$ turns them into ``macroscopic'' measures. The equation
above says that after this the resulting measures are, up to a translation,
$2^{-k}$-close in the weak sense. Therefore,
\begin{lem}
\label{lem:linearization-of-G-orbits}For every $\varepsilon>0$ and
$k>k(\varepsilon)$, if $\nu\in\mathcal{P}(G_{0})$ is supported on
a level-$k$ dyadic cube, $x,y\in\mathbb{R}^{d}$ are in the same
level-$k$ dyadic cube, and $V\leq\mathbb{R}^{d}$ is a linear subspace
then
\[
S_{k}(\nu\conv x)\mbox{ is }(V,\varepsilon)\mbox{-concentrated}\qquad\implies\qquad S_{k}(\nu\conv y)\mbox{ is }(V,2\varepsilon)\mbox{-concentrated}.
\]

\end{lem}

\subsection{\label{sub:proof-of-inverse-theorems-for-G}Proof of the inverse
theorem}

We first prove a version of the inverse theorem \ref{thm:inverse-theorem-for-isometries}
which assumes that $\nu,\mu$ are supported on small dyadic cubes.
These cubes are introduced to ensure that the supports of the measures
are small enough for the linearization machinery to kick in, and the
proof focuses on this aspect of the argument. After the proof we explain
how to get the stronger version, in which the measures have larger
support, and give bounds on the dimensions of the subspaces produced
by the theorem.
\begin{thm}
\label{thm:inverse-theorem-isometries-prelim} For every $\varepsilon>0$,
$R>0$ and $m\in\mathbb{N}$ there is a $\delta=\delta(\varepsilon,R,m)>0$,
such that, for all $k>k(\varepsilon,R,m,\delta)$ and all $n>n(\varepsilon,R,m,\delta,k)$,
the following holds: If $\nu\in\mathcal{P}(G_{0})$ and $\mu\in\mathcal{P}([-R,R]^{d})$
are supported on level-$k$ dyadic cells, then either 
\[
H_{n}(\nu\conv\mu)>H_{n}(\mu)+\delta,
\]
or there is a sequence $V_{k},\ldots,V_{n}$ of subspaces of $\mathbb{R}^{d}$
such that
\[
\mathbb{P}_{0\leq i\leq n}\left(\mu^{x,i}\mbox{ is }(V_{i},\varepsilon,m)\mbox{-saturated}\right)>1-\varepsilon
\]
and for all $x\in\supp\mu$,
\[
\mathbb{P}_{0\leq i\leq n}\left(S_{i}(\nu_{g,i}\conv x)\mbox{ is }(U_{g}V_{i},\varepsilon)\mbox{-concentrated}\right)>1-\varepsilon.
\]
\end{thm}
\begin{rem}
{}\end{rem}
\begin{enumerate}
\item Since $k$ is assumed large relative to $\varepsilon$, by Lemma \ref{lem:linearization-of-G-orbits}
the last condition holds for all $x\in\supp\mu$ if and only if it
holds for some $x\in\supp\mu$, up to a change of a factor of $2$
in the degree of concentration.
\item The measures $\mu$, $\nu$ and $\mu\conv\nu$ are supported on sets
of diameter $O_{R}(2^{-k})$, so when measuring their scale-$n$ entropy
it might seem more natural to rescale them by $O_{R}(2^{k})$. However,
the statement of the theorem is formally unchanged if we do so, since
we are taking $n$ large relative to $k$, and the average entropy
over $n$ scales is negligibly affected by the first $k$ scales. \end{enumerate}
\begin{proof}
Let $\varepsilon,R$ and $m$ be given. 
\begin{enumerate}
\item [i.]Apply Corollary \ref{cor:derandomizing-subspaces-with-isometries}
with parameters $\varepsilon$ and $m$ to obtain parameters $\varepsilon'$
and $m'$ and $n'$.
\item [ii.] Apply Proposition \ref{prop:concentration-passes-to-group-components}
with parameter $\frac{1}{4}\varepsilon'$, $R$ and $m'$ to obtain
parameters $\varepsilon''$ and $m''$.
\item [iii.] Apply the Euclidean inverse theorem (Theorem \ref{thm:inverse-thm-Rd})
with parameters $\varepsilon'',R,m''$, obtaining $\delta'$ and $n''$.
We are free to assume that $\delta'$ is arbitrarily small in a manner
depending on the previous parameters, and that $n''$ is large with
respect to previous parameters. In particular we assume $n''$ is
large relative to 
\[
\delta=(\delta'/2)^{2}.
\]

\item [iv.] Choose $k$ large enough that the conclusions of Lemma \ref{lem:linearization-of-G-orbits}
hold for parameter $\varepsilon''$ and Lemma \ref{lem:linearization-of-entropy-of-G-convolution}
holds for parameter $n'$ (instead of $m$ there). We also assume
that for any $D\in\mathcal{D}_{k}^{G_{0}}$ and $g,h\in D$, the difference
$\left\Vert U_{g}-U_{g'}\right\Vert $ is small enough that if $\theta\in\mathcal{P}([0,1]^{d})$
is a $(U_{g}V,\frac{1}{4}\varepsilon')$-concentrated measure then
it is also $(U_{h}V,\varepsilon')$-concentrated. 
\item [v.]Let $n$ be very large in a manner depending on all previous
parameters. 
\end{enumerate}
Now let $\nu\in\mathcal{P}(G),\mu\in\mathcal{P}(\mathbb{R}^{d})$
be supported on level-$k$ dyadic cells, and suppose that 
\begin{equation}
H_{n}(\nu\conv\mu)\leq H_{n}(\mu)+\delta.\label{eq:101}
\end{equation}
By Lemmas \ref{lem:entropy-local-to-global} and \ref{lem:local-to-global-on-G},
assuming $n$ is large compared to $n''$, \eqref{eq:101} implies
\[
\mathbb{E}_{0\leq i\leq n}\left(H_{i,n''}(\nu_{g,i}\conv\mu)-H_{i,n''}(\mu_{x,i})\right)<2\delta.
\]
By Lemma \ref{lem:linearization-of-entropy-of-G-convolution}, our
choice of $k$ and the fact that $n''$ are large in a manner depending
on $\delta$, 
\[
\mathbb{E}_{0\leq i\leq n}\left(H_{i,n''}((U_{g}^{-1}(\nu_{g,i}\conv x))*\mu_{x,i})-H_{i,n''}(\mu_{x,i})\right)<3\delta.
\]
Since $n''$ is large enough relative to $\delta$, the difference
inside the expectation is essentially non-negative, the is, larger
than $-\delta$ (Lemma \ref{lem:entropy-monotonicity-under-convolution}).
Since $\delta'=\sqrt{4\delta}$, by Markov's inequality we conclude
that
\[
\mathbb{P}_{0\leq i\leq n}\left(H_{i,n''}((U_{g}^{-1}(\nu_{g,i}\conv x))*\mu_{x,i}))\leq H_{i,n''}(\mu_{x,i})+\delta'\right)>1-\delta'.
\]

Fix $g,x$ such that $\nu_{g,i}$ and $\mu_{x,i}$ are in the event
above. Write $\eta=U_{g}^{-1}(\nu_{g,i}\conv x)$ and $\theta=\mu_{x,i}$.
Since
\[
H_{i,n''}(\eta*\theta)\leq H_{i,n''}(\theta)+\delta',
\]
and $\eta$ is supported on a set of diameter $O(R\cdot2^{-i})$,
we can, after implicitly re-scaling by $2^{i}$, apply the Euclidean
inverse theorem (Theorem \ref{thm:inverse-thm-Rd}) and conclude,
by our choice of the parameters $n'',\delta'$, that there are subspaces
$V_{j}=V_{j}^{(i,g,x)}$ for $i\leq j\leq i+n''$, such that
\[
\mathbb{P}_{i\leq j\leq i+n''}\left(\begin{array}{c}
\theta^{y,j}\mbox{ is }(V_{j}^{(i,g,x)},\varepsilon'',m')\mbox{-saturated and}\\
\eta^{z,j}\mbox{ is }(V_{j}^{(i,g,x)},\varepsilon'')\mbox{-concentrated}
\end{array}\right)>1-\varepsilon''.
\]
Since we can assume $n''>n'$, by Proposition \ref{prop:concentration-passes-to-group-components}
and our choice of parameters, writing $\tau=\nu_{g,i}$, 
\[
\mathbb{P}_{i\leq j\leq i+n''}\left(\begin{array}{c}
\theta^{y,j}\mbox{ is }(V_{j}^{(i,g,x)},\frac{1}{2}\varepsilon',m')\mbox{-saturated and}\\
S_{-j}U_{g}^{-1}(\tau_{h,j}\conv x)\mbox{ is }(V_{j}^{(i,g,x)},\frac{1}{4}\varepsilon')\mbox{-concentrated}
\end{array}\right)>1-\frac{1}{2}\varepsilon'
\]
(in the last equation, $g,x$ are fixed, and the randomness is over
$y,h$ and $j$). Recalling that $\mu,\nu$ are supported on level-$k$
dyadic cells and the definition of $k$, we can apply Lemma \ref{lem:linearization-of-G-orbits}
in the event above to replace $\tau_{h,j}\conv x$ by $\tau_{h,j}\conv y$.
As a result the degree of concentration degrades from $\varepsilon'/4$
to $\varepsilon'/2$. Then, since $h,g$ are in the same level $j$
(and hence level-$k$) component, we can exchange $U_{g}$ with $U_{h}$
in the event above with another $\varepsilon'/4$ degradation of the
concentration. After these adjustments we have
\[
\mathbb{P}_{i\leq j\leq i+n''}\left(\begin{array}{c}
\theta^{y,j}\mbox{ is }(V_{j}^{(i,g,x)},\varepsilon',m')\mbox{-saturated and}\\
S_{-j}U_{h}^{-1}(\tau_{h,j}\conv y)\mbox{ is }(V_{j}^{(i,g,x)},\varepsilon')\mbox{-concentrated}
\end{array}\right)>1-\frac{1}{2}\varepsilon'.
\]

So far we have seen that with high probability (at least $1-\delta'$)
over choice of components $\theta=\mu_{x,i}$ and $\tau=\nu_{g,i}$,
we can associate subspaces $V_{j}^{(i,g,x)}$ to a large fraction
(at least $1-\varepsilon'/2$) of the components of $\theta,\tau$
at levels $i,\ldots,i+n''$. These components are also components
of $\nu,\mu$, but each component of $\nu,\mu$ may arise in several
ways as a component of a components. So we have not associated a subspace
to (most) components of $\nu,\mu$, but rather to (most) components
of $\nu,\mu$ we have associated several subspaces. To correct this
we invoke Lemma \ref{lem:distribution-of-components-of-components},
letting us select subspaces $V^{(i,g,x)}$ (no longer depending on
$j$) such that 
\[
\mathbb{P}_{0\leq i\leq n}\left(\begin{array}{c}
\mu^{x,i}\mbox{ is }(V^{(i,g,x)},\varepsilon',m')\mbox{-saturated and}\\
S_{-i}U_{g}^{-1}\nu_{g,i}\conv x\mbox{ is }(V^{(i,g,x)},\varepsilon')\mbox{-concentrated}
\end{array}\right)>1-\frac{1}{2}\varepsilon'-\delta'-O(\frac{n''}{n}).
\]
The right hand side is $>1-\varepsilon'$ assuming as we may that
$\delta'$ is small compared to $\varepsilon'$ and $n$ large relative
to $n''$. Applying Corollary \ref{cor:derandomizing-subspaces-with-isometries},
and by our choice of $\varepsilon'$, there are subspaces $V^{i}$,
independent of $g,x$, such that
\[
\mathbb{P}_{0\leq i\leq n}\left(\begin{array}{c}
\mu^{x,i}\mbox{ is }(V^{i},\varepsilon,m)\mbox{-saturated and}\\
S_{i}U_{g}^{-1}\nu_{g,i}\conv x\mbox{ is }(V^{i},\varepsilon)\mbox{-concentrated}
\end{array}\right)>1-\varepsilon.
\]
This implies the statement.
\end{proof}
We now prove Theorem \ref{thm:inverse-theorem-for-isometries}, which
we repeat for convenience:
\begin{thm}
For every $\varepsilon>0$, $R>0$ and $m\in\mathbb{N}$, there exists
$\delta=\delta(\varepsilon,R,m)>0$ such that for every $k>k(\varepsilon,R,m)$
and every $n>n(\varepsilon,R,m,k)$, the following holds. For every
$\nu\in\mathcal{P}(G_{0})$ and $\mu\in\mathcal{P}([-R,R]^{d})$ that
are supported on balls of radius $R$, either 
\[
H_{n}(\nu\conv\mu)>H_{n}(\mu)+\delta,
\]
or else, to every pair of level-$k$ components $\widetilde{\nu}$
of $\nu$ and $\widetilde{\mu}$ of $\mu$ we can assign a sequence
of subspaces $V_{i}=V_{i}(\widetilde{\nu},\widetilde{\mu})<\mathbb{R}^{d}$,
$0\leq i\leq n$, such that with probability at least $1-\varepsilon$
over the choice of $\widetilde{\mu},\widetilde{\nu}$, 
\[
\mathbb{P}_{0\leq i\leq n}\left(\begin{array}{c}
\widetilde{\mu}^{x,i}\mbox{ is }(V_{i},\varepsilon,m)\mbox{-saturated and }\\
S_{i}U_{g}^{-1}(\widetilde{\nu}_{g,i}\conv x)\mbox{ is }(V_{i},\varepsilon)\mbox{-concentrated}
\end{array}\right)>1-\varepsilon
\]
If in addition $\mu$ is $((\varepsilon/5d)^{2(d+1)},\sigma)$-non-affine
for some $\sigma>0$, and the relation among parameters takes $\sigma$
into account, then for those $\widetilde{\nu},\widetilde{\mu}$ in
the set of good components above, 
\begin{equation}
\frac{1}{n+1}\sum_{i=0}^{n}\dim V_{i}>\frac{1}{d+1}H_{n}(\widetilde{\nu})-\varepsilon\label{eq:112}
\end{equation}
and 
\begin{equation}
\mathbb{E}_{i=k}\left(\frac{1}{n+1}\sum_{j=0}^{n}\dim V_{j}(\nu_{g,i},\mu_{x,i})\right)>\frac{1}{d+1}H(\nu)-\varepsilon\label{eq:114}
\end{equation}
\end{thm}
\begin{proof}
Fix $\varepsilon,R,m$ (the error terms below depend on them but we
suppress it in the notation). Let also $\delta,k,n$ be parameters
whose relations we will specify later, and suppose that 
\[
H_{n}(\nu\conv\mu)<H_{n}(\mu)+\delta
\]
By Lemma \ref{lem:skipping-k-scales},
\[
\mathbb{E}_{i=k}\left(H_{n}(\nu_{g,i}\conv\mu)-H_{n}(\mu_{x,i})\right)<2\delta+O(\frac{1}{n})<3\delta.
\]
By Markov's inequality, assuming $n$ large enough,
\begin{equation}
\mathbb{P}_{i=k}\left(H_{n}(\nu_{g,i}\conv\mu)-H_{n}(\mu_{x,i})<\sqrt{3\delta}\right)>1-\sqrt{3\delta}\label{eq:108}
\end{equation}
Assuming as we may that $\sqrt{3\delta}<\varepsilon$, the last probability
is at least $1-\varepsilon$. 

Now fix a pair of components $\widetilde{\nu},\widetilde{\mu}$ from
the event in \eqref{eq:108}. Assuming that $\sqrt{3\delta}$ is small
relative to $\varepsilon,R,m$ and that $k,n$ are large enough, we
can apply the previous theorem to $\widetilde{\nu}$ and $\widetilde{\mu}$
(which are by definition supported on level-$k$ components) and obtain
corresponding subspaces $V_{i}(\widetilde{\nu},\widetilde{\mu})$,
$k\leq i\leq n$. This proves the first part of the present theorem.

For the second part (bounding the dimensions of the subspaces), suppose
that $\mu$ is $(\sigma',\sigma)$-non-concentrated. Fix an auxiliary
parameter $\varepsilon'$ depending in a manner we shall later determine
on $\varepsilon,\sigma$ and $R$, and run first part using $\varepsilon'$
instead of $\varepsilon$, obtaining associated $\delta,k,n$ etc.,
and a set of level-$k$ components $\widetilde{\nu},\widetilde{\mu}$
of probability at least $1-\varepsilon'$ to which are associated
subspaces $V_{i}(\widetilde{\mu},\widetilde{\nu})$ with the desired
properties w.r.t. $\varepsilon'$. Define $V_{i}(\widetilde{\nu},\widetilde{\mu})=\mathbb{R}^{d}$
for any pair of level-$k$ components $\widetilde{\nu},\widetilde{\mu}$
for which is was not yet defined (i.e. pairs that are not in the event
in \eqref{eq:108}). For $i\geq k$ and components $\nu_{g,i}$ and
$\mu_{x,i}$ set 
\[
V(\nu_{g,i},\mu_{x,i})=V_{i}(\nu_{g,k},\mu_{g,k}).
\]
This is well defined because a level-$i$ component for $i\geq k$
determines uniquely the level-$k$ component it belongs to (on the
other hand we are abusing notation slightly since, strictly speaking,
$\nu_{g,i},\mu_{x,i}$ do not determine $g,x,i$; but as they are
written explicitly, no confusion should occur). 

Observe now that, by the first part of the proof, 
\[
\mathbb{P}_{0\leq i\leq n}\left(S_{i}U_{g}^{-1}\nu_{g,i}\conv x\mbox{ is }(V(\nu_{g,i},\mu_{x,i}),\varepsilon)\mbox{-concentrated}\right)\,d\mu(x)>1-\varepsilon'.
\]
Indeed, if we write $\widetilde{\nu},\widetilde{\mu}$ for the level-$k$
components to which $\nu_{g,i},\mu_{x,i}$ belong, respectively, then
conditioned on $\widetilde{\nu},\widetilde{\mu}$ belonging to the
event in \eqref{eq:108}, the probability of the event above is at
least $1-\varepsilon'$; while conditioned on the complementary event,
the probability is $1$, since then $V(\widetilde{\nu},\widetilde{\mu})=\mathbb{R}^{d}$.
Thus the unconditional probability above is at least $1-\varepsilon'$.

Set
\[
\ell=[\log(1/\varepsilon')]
\]
By Lemma \ref{lem:entropy-of-V-conventrated-measures}, the previous
inequality gives 
\begin{eqnarray*}
\mathbb{P}_{0\leq i\leq n}\left(H_{\ell}(S_{i}U_{g}^{-1}\nu_{g,i}\conv x)<\dim V(\nu_{g,i},\mu_{x,i})+O(\frac{\log\ell}{\ell})\right) & > & 1-\varepsilon'.
\end{eqnarray*}
Since by Lemma \ref{lem:entropy-weak-continuity-properties} \eqref{enu:entropy-combinatorial-distortion},
\[
\left|H_{\ell}(S_{i}U_{g}^{-1}\nu_{g,i}\conv x)-H_{i,\ell}(\nu_{g,i}\conv x)\right|=O(\frac{1}{\ell}),
\]
we obtain
\begin{eqnarray}
\mathbb{P}_{0\leq i\leq n}\left(H_{i,\ell}(\nu_{g,i}\conv x)<\dim V(\nu_{g,i},\mu_{x,i})+O(\frac{\log\ell}{\ell})\right) & > & 1-2\varepsilon'-O(\frac{1}{\ell})\nonumber \\
 & = & 1-O(\frac{1}{\ell}).\label{eq:107}
\end{eqnarray}

We now use the assumption that $\mu$ is $((\varepsilon/5d)^{2(d+1)},\sigma)$-non-affine.
By Corollary \ref{cor:entropy-passes-to-G-orbits}, for every component
$\nu_{g,i}$ of $\nu$, 
\begin{eqnarray*}
\mu\left(x\in\mathbb{R}^{d}\,:\,H_{i,\ell}(\nu_{g,i}\conv x)>\frac{1}{d+1}H_{i,\ell}(\nu_{g,i})-O_{\sigma,R}(\frac{1}{\ell})\right) & > & 1-\frac{1}{5}\varepsilon^{2},
\end{eqnarray*}
Choosing the component $\nu_{g,i}$, $k\leq i\leq n$, at random,
and then $x$ independently according to $\mu$, we conclude that
$H_{i,\ell}(\nu_{g,i}\conv x)>\frac{1}{d+1}H_{i,\ell}(\nu_{g,i})-O_{\sigma,R}(\frac{1}{\ell})$
with probability at least $1-\varepsilon^{2}/5$. Therefore, combined
with \eqref{eq:107}, we have
\begin{eqnarray*}
\mathbb{P}_{0\leq i\leq n}\left(\begin{array}{c}
H_{i,\ell}(\nu_{g,i}\conv x)<\dim V(\nu_{g,i},\mu_{x,i})+O(\frac{\log\ell}{\ell})\\
\mbox{and }H_{i,\ell}(\nu_{g,i}\conv x)>\frac{1}{d+1}H_{i,\ell}(\nu_{g,i})-O_{\sigma,R}(\frac{1}{\ell})
\end{array}\right) & > & 1-O(\frac{1}{\ell})-\frac{1}{5}\varepsilon^{2}
\end{eqnarray*}
Recalling that $\ell=\log(1/\varepsilon')$, by $\varepsilon'$ small
we can assume the error term does not exceed $\varepsilon^{2}/4$.
We obtain 
\begin{eqnarray*}
\mathbb{P}_{0\leq i\leq n}\left(\frac{1}{d+1}H_{i,\ell}(\nu_{g,i})<\dim V(\nu_{g,i},\mu_{x,i})+O_{\sigma,R}(\frac{\log\ell}{\ell})\right) & > & 1-\frac{1}{4}\varepsilon^{2}
\end{eqnarray*}
By Markov's inequality, there is a set $A\subseteq G_{0}\times\mathbb{R}^{d}$
with $\nu\times\mu(A)>1-\varepsilon/2$, such that for every $(g_{0},x_{0})\in A$,
setting $\widetilde{\nu}=\nu_{g_{0},k}$ and $\widetilde{\mu}=\mu_{x_{0},k}$,
\begin{eqnarray*}
\mathbb{P}_{0\leq i\leq n}\left(\frac{1}{d+1}H_{i,\ell}(\widetilde{\nu}_{g,i})<\dim V(\widetilde{\nu}_{g,i},\widetilde{\mu}_{x,i})+O_{\sigma,R}(\frac{\log\ell}{\ell})\right) & > & 1-\frac{1}{2}\varepsilon
\end{eqnarray*}
Outside of the event above we have the trivial bound $\frac{1}{d+1}H_{i,\ell}(\widetilde{\nu}_{g,i})\leq\frac{1}{d+1}d<1$.
On the other hand, by Lemma \ref{lem:entropy-local-to-global} (which
holds also in $G$), 
\[
H_{n}(\widetilde{\nu})=\mathbb{E}_{0\leq i\leq n}(H_{i,\ell}(\widetilde{\nu}_{g,i}))+O(1/\ell+\ell/n)
\]
Finally, since $V(\widetilde{\nu}_{g,i},\widetilde{\mu}_{x,i})=V_{i}(\widetilde{\nu},\widetilde{\mu})$),
the last two equations give 
\begin{eqnarray}
\frac{1}{d+1}H_{n}(\widetilde{\nu}) & = & \mathbb{E}_{0\leq i\leq n}\left(\frac{1}{d+1}H_{i,\ell}(\widetilde{\nu}_{g,i})\right)+O(\frac{1}{\ell}+\frac{\ell}{n})\nonumber \\
 & < & \frac{1}{n+1}\sum_{i=0}^{n}\left(\dim V_{i}(\widetilde{\nu},\widetilde{\mu})+O_{\sigma,R}(\frac{\log\ell}{\ell})\right)+O(\frac{1}{\ell}+\frac{\ell}{n})+\frac{1}{2}\varepsilon\label{eq:113}
\end{eqnarray}
Taking $\varepsilon'$ small (and hence $\ell$ large) relative to
$\varepsilon,R,\sigma$, and $n$ larger, and rearranging, we obtain
\eqref{eq:112}. 

Finally, recall that \eqref{eq:113} holds on a set of pairs of level-$k$
components $\widetilde{\nu},\widetilde{\mu}$ of probability at least
$1-\varepsilon/2$, and recall that 
\[
\mathbb{E}_{i=k}(H_{n}(\nu_{g,i}))=H_{n}(\nu)-O(\frac{k}{n})
\]
The last statement of the theorem, \eqref{eq:114}, follows now by
taking expectation of both sides in \eqref{eq:113} and making $\varepsilon'$
small enough and $n$ large enough.
\end{proof}

\subsection{\label{sub:Generalizations}Generalizations}

To derive Theorem \ref{thm:generalized-inverse-theorem} very few
changes are needed to the convolution case. The $C^{1}$-assumption
of $f$, and the compactness of its domain, easily imply analogs of
Equations \eqref{eq:discrete-differentiation-of-G0-action}, \eqref{eq:discrete-differentiation-of-G0-action-2}
and \eqref{eq:discrete-differentiation-identity} and their consequences
(without quantitative control on the error, but one cannot expect
it in the general setting). In particular, for large enough $k$ and
suitably large $n$, with $\mu\times\nu$-probability at least $1-\delta$
over choice of $(x,y)$ we have 
\[
|H_{n}(f(\mu_{x,k}\times\delta_{y}))-H_{n}(A_{x,y}\mu_{x,k})|<\frac{\delta}{10},
\]
and 
\[
|H_{n}(f(\mu_{x,k}\times\nu_{y,k}))-H_{n}(A_{x,y}\mu_{x,k}*B_{x,y}\nu_{x,k})|<\frac{\delta}{10}
\]
(note that since $n\gg k$, there is no advantage in scaling $1/\left\Vert A\right\Vert $,
$1/\left\Vert B\right\Vert $ by $2^{k}$, as might seem natural). 

By concavity and almost-convexity of entropy (Lemma \ref{lem:entropy-combinatorial-properties}
\eqref{enu:entropy-concavity} and \eqref{enu:entropy-almost-convexity}),
for $n\gg k$ we have
\[
\left|H_{n}(f(\mu\times\nu))-\int H_{n}(f(\mu_{x,k}\times\nu_{y,k}))\,d\mu\times\nu(x,y)\right|<\frac{\delta}{10},
\]
and for every $y$, similarly,
\[
\left|H_{n}(f(\mu\times\delta_{y}))-\int H_{n}(f(\mu_{x,k}\times\delta_{y}))\,d\mu(x)\right|<\frac{\delta}{10}.
\]
Thus the hypothesis \eqref{eq:entropy-growth-in-generalized-inv-thm}
of Theorem \ref{thm:generalized-inverse-theorem} implies that for
any $k$ and $n\gg k$,
\[
\int H_{n}(f(\mu_{x,k}\times\nu_{y,k}))\,d\mu\times\nu(x,y)<\int H_{n}(f(\mu_{x,k}\times\delta_{y}))\,d\mu\times\nu(x,y)+\frac{8}{10}\delta.
\]
By the above, for large $k$ this is
\[
\int H_{n}(A_{x,y}\mu_{x,k}*B_{x,y}\nu_{x,k})\,d\mu\times\nu(x,y)<\int H_{n}(A_{x,y}\mu_{x,k})\,d\mu\times\nu(x,y)+\frac{6}{10}\delta.
\]
Since for large $n$ we essentially have the reverse inequality between
the integrands, we conclude that with high probability at least $1-\delta$
over the components $\widetilde{\mu}=\mu_{x,k}$ and $\widetilde{\nu}=\nu_{y,k}$,
we have
\[
H_{n}(A_{x,y}\widetilde{\mu}*B_{x,y}\widetilde{\nu})<H_{n}(A_{x,y}\widetilde{\mu})+\delta',
\]
where $\delta'$ tends to zero with $\delta$. From here one can apply
the Euclidean inverse theorem to the components $\widetilde{\nu},\widetilde{\mu}$
as we did in the proof of the convolution case, with very few changes
other than notational ones. We omit the details. 

In the special case of actions of matrix groups on $\mathbb{R}^{d}$
or on themselves, one has analogs of Corollary \ref{cor:entropy-passes-to-G-orbits}.
In the first case essentially by the same lemma (using compactness
of the domain of the action function in place of compactness of the
orthogonal group). For a matrix group acting on itself, there are
in fact trivial stabilizers, so there conclusion is automatic.

\section{\label{sec:Self-similar-sets-and-measures-Rd}Self-similar sets and
measures on $\mathbb{R}^{d}$}

The derivation of our main result, Theorem \ref{thm:main-individual-RD-entropy},
from the Theorem \ref{thm:inverse-theorem-for-isometries} (the inverse
theorem for the $G_{0}$-action), follows lines similar to the argument
in \cite{Hochman2014} for $\mathbb{R}$. One new ingredient is the
explicit presence of the isometry group, but this is implicit in the
original argument and the main change is notational. More significant
is the appearance of invariant subspaces in the third alternative
of the theorem. This will require some further analysis, and will
occupy us in the first few subsections.

We remark that our analysis so far, and much of the analysis below,
is of a finitary nature, involving entropies at fine (but finite)
partitions. Certainly we must somewhere connect this to dimension,
specifically to the dimension of the conditional measures of a given
self-similar measure on the family of translates of a subspace (as
in (iii'') of Theorem \ref{thm:main-individual-RD-entropy}). It is
an unfortunate reality that such a connection seems to be available
only when the subspace is invariant under the linearization of the
IFS (see Section \ref{sub:Entropy-and-dimension} below). If such
results were available without invariance, much of the technical work
of the next few sections could be avoided by passing to a limit at
an earlier stage. However, understanding these ``slice'' measures
for general self-similar measures remains an open problem.

\subsection{Almost-invariance and invariance}

We will obtain invariant subspaces from almost invariant ones:
\begin{defn}
A subspace $V\leq\mathbb{R}^{d}$ is $\varepsilon$-invariant under
a subgroup $H<G_{0}$, or $(H,\varepsilon)$-invariant, if $d(hV,V)\leq\varepsilon$
for every $h\in H$. 
\end{defn}
Evidently, $(H,0)$-invariance is $H$-invariance in the usual sense.
Furthermore,
\begin{lem}
\label{lem:invariant-subspace-near-almost-invariant-one}Let $H<G_{0}$
be a closed subgroup. For every $\varepsilon>0$ there is a $\delta>0$,
such that if $V$ is $\delta$-invariant under $H$, then there is
an $H$-invariant subspace $V'$ with $d(V,V')<\varepsilon$.\end{lem}
\begin{proof}
Let $\mathcal{S}_{H}$ denote the space of $H$-invariant subspaces
of $\mathbb{R}^{d}$. If the statement were false there would be some
$\varepsilon>0$ and a sequence $V_{n}\leq\mathbb{R}^{d}$ of subspaces
such that $V_{n}$ is $1/n$-invariant for $H$, but $d(V_{n},V')\geq\varepsilon$
for every $V'\in\mathcal{S}_{H}$. Using compactness of the space
of subspaces, we can pass to a subsequence $V_{n_{k}}$ converging
to some $V$. Since the linear action is continuous, $d(V,hV)=\lim d(V_{n_{k}},hV_{n_{k}})=0$
for all $h\in H$, so $V\in\mathcal{S}_{H}$. But by hypothesis $d(V_{n_{k}},V)\geq\varepsilon$
for all $k$, a contradiction.
\end{proof}
In fact the choice $\delta=c\cdot\varepsilon^{d+1}$ works for an
appropriate constant $c$ (or $c\cdot\varepsilon^{k+1}$ if one fixes
the dimension $k$ of the subspace in question), but we will not use
this.

Our second tool will be to construct almost-invariant subspaces from
almost-invariant sets of vectors. 
\begin{lem}
\label{lem:almost-invariant-subspace-from-almost-invariant-set}Let
$0<\varepsilon<1$ and write $\varepsilon_{n}=\varepsilon^{n!}$.
Let $H<G_{0}$ be a closed subgroup and let $E\subseteq B_{1}(0)\subseteq\mathbb{R}^{d}$
be a set such that $d(hv,E)<\varepsilon_{n}$ for all $v\in E$ and
$h\in H$. Let $v_{1},\ldots,v_{k}\in E$ be a maximal sequence of
vectors satisfying $d(v_{i},\spn\{v_{1},\ldots,v_{i-1}\})>\varepsilon_{i}$
for $1<i\leq k$, and set $V=\spn\{v_{1},\ldots,v_{k}\}$. Then $V$
is $(H,O(\varepsilon))$-invariant and $E\subseteq V^{(\varepsilon_{k+1})}$.\end{lem}
\begin{proof}
We may assume $k<d$ since otherwise $V=\mathbb{R}^{d}$ and the statements
is trivial. To see that $E\subseteq V^{(\varepsilon_{k+1})}$, note
that if $v\in E\setminus V^{(\varepsilon_{k+1})}$ then the vector
$v_{k+1}=v$ would extend the given sequence of vectors in a way that
contradicts its maximality. For invariance, let $h\in H$ and set
$w_{i}=hv_{i}$ and $W=hV=\spn\{w_{i}\}$. By assumption, for each
$i$ there is a $w'_{i}\in E$ with $d(w_{i},w'_{i})<\varepsilon_{d}\leq\varepsilon_{k+1}$,
and we saw above that $w'_{i}\in V^{(\varepsilon_{k+1})}$, hence
$w_{i}\in V^{(2\varepsilon_{k+1})}$. Also, $h$ is an isometry, so
$d(w_{i},\spn\{w_{1},\ldots,w_{i-1}\})>\varepsilon_{i}\geq\varepsilon_{k}$
for all $1\leq i\leq k$, since the same is true for the $v_{i}$.
Therefore, by Corollary \ref{cor:span-is-close-to-V-if-basis-is-close-to-V},
$\spn\{w_{1},\ldots,w_{k}\}\sqsubseteq V^{(c\cdot\varepsilon_{k+1}/\varepsilon_{k}^{k})}$,
and using the fact that $\dim V=\dim W$, this implies 
\[
d(W,V)=O(\frac{2\varepsilon_{k+1}}{\varepsilon_{k}^{k}})=O(\varepsilon^{(k+1)!-k!\cdot k})=O(\varepsilon^{k!})=O(\varepsilon),
\]
as desired. 
\end{proof}

\subsection{\label{sub:saturation-and-self-similar-measures}Saturation at level
$n$}

We will be interested in the situation where the components of a measure
at some scale typically are highly saturated on a subspace. More precisely,
\begin{defn}
\label{def:saturation-at-level-n}For $V\leq\mathbb{R}^{d}$ , a measure
$\mu\in\mathcal{P}(\mathbb{R}^{d})$ is $(V,\varepsilon,m)$-\emph{saturated
at level }$n$ if 
\[
\mathbb{P}_{i=n}\left(\mu^{x,i}\mbox{ is }(V,\varepsilon,m)\mbox{-saturated}\right)>1-\varepsilon.
\]

\end{defn}
We write 
\[
\sat(\mu,\varepsilon,m,n)=\{V\leq\mathbb{R}^{d}\,:\,\mu\mbox{ is }(V,\varepsilon,m)\mbox{-saturated at level }n\}.
\]
Some technical properties related to this notion are summarized in
the next lemma. In the formulation we write $\sum A$ for $\sum_{a\in A}a$.
\begin{lem}
\label{lem:saturation-at-level-n} Let $\varepsilon,R>0$ and $V\leq\mathbb{R}^{d}$.
Let $\mu\in\mathcal{P}([-R,R]^{d})$ and suppose that $\mu$ is given
as a convex combination of probability measures, $\mu=\sum_{i=1}^{k}\alpha_{i}\mu_{i}$. 
\begin{enumerate}
\item If $\mu$ is $(V,\varepsilon,m)$-saturated, then
\[
\sum\{\alpha_{i}\,:\,\mu_{i}\mbox{ is }(V,\varepsilon',m)\mbox{-saturated}\}>1-\varepsilon',
\]
where $\varepsilon'=O(\sqrt{\varepsilon+(\log kR)/m)})$.
\item \label{enu:saturation-at-level-n-transfers-to-components}For $n$
sufficiently large in a manner depending on $\mu,\alpha_{i},\nu_{i}$,
if $V\in\sat(\mu,\varepsilon,m,n)$ then 
\[
\sum\{\alpha_{i}\,:\,V\in\sat(\mu_{i},\varepsilon',m,n)\}>1-\varepsilon',
\]
where $\varepsilon'=O(\sqrt{\varepsilon})$.
\item \label{enu:saturation-at-level-n-of-convex-combinations}If for some
$n$ we have 
\[
\sum\{\alpha_{i}\,:\,V\in\sat(\mu_{i},\varepsilon,m,n)\}>1-\varepsilon,
\]
 then $V\in\sat(\mu,\varepsilon',m,n)$, where $\varepsilon'=O(\sqrt{\varepsilon})$.
\item \label{enu:saturation-at-level-n-similarity-1}Let $g=2^{-t}U+a\in G$.
If $V\in\sat(\mu,\varepsilon,m,n)$ then $UV\in\sat(g\mu,\varepsilon',m,[n-t])$
where $\varepsilon'\rightarrow0$ as $(\varepsilon,\frac{1}{m})\rightarrow0$.
\item \label{enu:saturation-at-level-n-similarity-2}Under the same assumptions
as in \eqref{enu:saturation-at-level-n-similarity-1}, $UV\in\sat(g\mu,\varepsilon'',m,n)$
where $\varepsilon''\rightarrow0$ as $(\varepsilon,\frac{1}{m})\rightarrow0$.
\end{enumerate}
\end{lem}
\begin{proof}
For (1), by absorbing an $O(1/m)$ error into $\varepsilon$ we can
assume that $\mathcal{D}_{m}=\mathcal{D}_{m}^{V}\lor\mathcal{D}_{m}^{V^{\perp}}$
(Lemma \ref{lem:saturation-under-coordinate-change}). By Lemmas \ref{lem:entropy-combinatorial-properties}
\eqref{enu:entropy-almost-convexity} and the hypothesis, we have
\begin{eqnarray*}
\sum_{i=1}^{k}\alpha_{i}\cdot\frac{1}{m}H(\mu_{i},\mathcal{D}_{m}|\mathcal{D}_{m}^{V^{\perp}}) & \geq & \frac{1}{m}H(\mu,\mathcal{D}_{m}|\mathcal{D}_{m}^{V^{\perp}})-\frac{\log k}{m}\\
 & > & \dim V-(\varepsilon+\frac{\log k}{m}).
\end{eqnarray*}
On the other hand, each $\mu_{i}$ is supported on $[-R,R]^{d}$ so
each term in the average on the left hand side is bounded above by
$\dim V+O(\frac{\log R}{m})$. Now (1) follows by Markov's inequality.

For (2), fix for convenience $\delta=\sqrt{\varepsilon}$. By standard
differentiation theorems, for $\mu_{i}$-a.e. $x$, $\left\Vert \mu_{x,\ell}-(\mu_{i})_{x,\ell}\right\Vert \rightarrow0$
as $\ell\rightarrow\infty$. In particular for large $n$, for a set
of $x$ of $\mu_{i}$-mass at least $1-\delta$, we have $\mu^{x,n}=(1-\delta)\mu_{i}^{x,n}+\delta\theta$
for some $\theta\in\mathcal{P}([0,1]^{d})$ (depending on $x,i$).
For any such $n$ let 
\[
A=\{x\in[0,1]^{d}\,:\,\mu^{x,n}\mbox{ is }(V,m,\varepsilon)\mbox{-saturated}\}.
\]
By hypothesis $\mu(A)>1-\varepsilon$. Since $\mu=\sum\alpha_{i}\mu_{i}$,
by Markov's inequality we have
\begin{equation}
\sum\{\alpha_{i}\,:\,\mu_{i}(A)>1-\sqrt{\varepsilon}\}>1-\sqrt{\varepsilon}.\label{eq:102}
\end{equation}
For $i$ satisfying $\mu_{i}(A)>1-\sqrt{\varepsilon}$, for a set
$x$ of points having $\mu_{i}$-mass $1-\delta-\sqrt{\varepsilon}$
we have that $\mu^{x,n}$ is $(V,\varepsilon,m)$-saturated and $\mu^{x,n}=(1-\delta)\mu_{i}^{x,n}+\delta\theta$
for some $\theta\in\mathcal{P}[0,1]^{d})$. Now we can apply part
(1) of this lemma to $\mu^{x,n}$, which is written as a combination
of two measures ($k=2$) and supported on $[0,1)$ (so $R=1$), and
conclude that $\mu_{i}^{x,n}$ is $(V,O(\sqrt{\varepsilon}),m)$-saturated.
This holds for at least a $(1-\delta-\sqrt{\varepsilon})$-fraction
of the components $\mu_{i}^{x,n}$. Since $\delta=\sqrt{\varepsilon}$
we find that $\mu_{i}$ is $(V,O(\sqrt{\varepsilon}),m,n)$-saturated.
This together with \eqref{eq:102} is what we wanted to prove.

For (3), observe that $\mu^{x,n}$ is a convex combination of components
$\mu_{i}^{x,n}$ (the weights are proportional to $\alpha_{i}\mu_{i}(\mathcal{D}_{n}(x))$).
By Lemmas \ref{lem:saturation-passes-to-convex-combinations} and
\ref{lem:saturation--of-total-variation-perturbations}, we will be
done if we show with $\mu$-probability $>1-\sqrt{2\varepsilon}$
over the choice of $x$, the components $\mu_{i}^{x,n}$ which are
$(V,\varepsilon,m)$-saturated constitute a $(1-\sqrt{2\varepsilon})$-fraction
of the mass of $\mu^{x,n}$.

To show this, let $I=\{1,\ldots,k\}$ and let $\alpha$ be the probability
measure on $I$ arising from the weights $\alpha_{i}$. Consider the
space $I\times\mathbb{R}^{d}$ with the probability measure $\theta$
given by $\theta(\{i\}\times A)=\alpha_{i}\mu_{i}(A)$. Define $f:I\times\mathbb{R}^{d}\rightarrow\mathbb{R}$
by 
\[
f(i,x)=\left\{ \begin{array}{cc}
1 & \mbox{if }\mu_{i}^{x,n}\mbox{ is }(V,\varepsilon,m)\mbox{-saturated}\\
0 & \mbox{otherwise}
\end{array}\right..
\]
Note that $f$ is $2^{I}\times\mathcal{D}_{n}$-measurable. Writing
$I_{0}=\{i\in I\,:\,\mu_{i}\mbox{ is }(V,\varepsilon,m,n)\mbox{-saturated}\}$,
we have 
\begin{eqnarray*}
\int fd\theta & = & \sum_{i\in I}\alpha_{i}\int f(i,x)d\mu_{i}(x)\\
 & \geq & \sum_{i\in I_{0}}\alpha_{i}\int f(i,x)d\mu_{i}(x)\\
 & = & \sum_{i\in I_{0}}\alpha_{i}\mu_{i}(x\,:\,\mu^{x,n}\mbox{ is }(V,\varepsilon,m)\mbox{-saturated})\\
 & > & \sum_{i\in I_{0}}\alpha_{i}(1-\varepsilon)\\
 & > & (1-\varepsilon)^{2}\\
 & > & 1-2\varepsilon
\end{eqnarray*}
(the passage from the third to fourth equation is by the hypothesis).
Let $\mathcal{B}$ be smallest the $\sigma$-algebra that makes the
map $I\times\mathbb{R}^{d}\rightarrow\mathbb{R}^{d}$, $(i,x)\mapsto x$,
measurable. The function $g=\mathbb{E}(f|\mathcal{B})$ also satisfies
$g\leq1$ and $\int gd\theta=\int fd\theta>1-2\varepsilon$, so by
Markov's inequality, 
\[
\theta((i,x)\,:\,g(x)>1-\sqrt{2\varepsilon})>1-\sqrt{2\varepsilon}.
\]
But, writing $D=\mathcal{D}_{n}(x)$, the inequality $g(x)>1-\sqrt{2\varepsilon}$
just means that in the convex combination $\mu^{x,n}=\sum\frac{\alpha_{i}\mu_{i}(D)}{\sum\alpha_{i}\mu_{i}(D)}(\mu_{i})^{x,n}$,
at least $1-\sqrt{2\varepsilon}$ of the mass originates in terms
for which $(\mu_{i})^{x,n}$ is $(V,\varepsilon,m)$-saturated. Since
the distribution on $x$ induced by $\theta$ is equal to $\mu$,
this completes the proof. 

For (4), consider $D\in\mathcal{D}_{n}$ such that $\mu^{D}$ is $(V,\varepsilon,m)$-saturated.
Let $\nu=g(\mu_{D})$. Then $\nu'=S_{[n-t]}\nu$ is the image of $\mu^{D}$
under a similarity that contracts by $O(1)$ and rotates by $U$,
and so by Lemma \ref{lem:saturation-transformed-by-similarity}, $\nu'$
is $(UV,\varepsilon+O(1/m),m)$-saturated. Writing $\nu'=\sum_{D\in\mathcal{D}_{1}}\nu'(D)\cdot\nu'_{D}$
we can apply (1) and conclude that, with $\varepsilon$ small and
$m$ large, most mass in this convex combination comes from terms
that are $(UV,\varepsilon',m)$-saturated. This means precisely that
$\nu'$ is $(UV,\varepsilon',m,n)$-saturated. Now, since $g\mu$
is the convex combination of measures $\nu$ of which a $(1-\varepsilon)$-fraction
are as above, (4) follows from (2).

(5) is proved in the same manner as (4), using $\nu'=S_{n}\nu$ instead
of $S_{[n-t]}\nu$. 
\end{proof}

\subsection{Saturated subspaces of self-similar measures}

From here until the end of the paper we again denote by $\mu$ a self-similar
measure on $\mathbb{R}^{d}$ defined by an IFS $\Phi=\{\varphi_{i}\}_{i\in\Lambda}$
and a positive probability vector $p=(p_{i})_{i\in\Lambda}$. As usual
we write $\varphi_{i}=r_{i}U_{i}+a_{i}$, and for $i\in\Lambda^{k}$
we set $\varphi_{i}=\varphi_{i_{1}}\circ\ldots\circ\varphi_{i_{k}}$,
$p_{i}=p_{i_{1}}\cdot\ldots\cdot p_{i_{k}}$ , and define $r_{i}$,
$U_{i}$, similarly. Denote by $G_{\Phi}\subseteq G_{0}$ the smallest
closed group containing the orthogonal parts $U_{i}$, $i\in\Lambda$,
of the maps $\varphi_{i}\in\Phi$.

In the next few results, all dependences between parameters and implicit
constants depend on $\mu$ and $\Phi$.

The first lemma says that the set of subspaces that are $(\mu,\varepsilon,m)$-saturated
at level $n$ is almost invariant under $G_{\Phi}$:
\begin{lem}
\label{lem:almost-invariant subspaces}For every $\varepsilon>0$
there is a $\delta=\delta(\varepsilon)>0$ such that, if $m>m(\varepsilon)$
and $n>n(\varepsilon,m)$, the following holds. For any $V\in\sat(\mu,\delta,m,n)$
and $g\in G_{\Phi}$ there exists $W\in\sat(\mu,\varepsilon,m,n)$
such that $d(W,gV)<\varepsilon$.\end{lem}
\begin{proof}
Let $\Lambda^{\leq n}=\bigcup_{i=1}^{n}\Lambda^{i}$. Then $S=\{U_{i}\,:\,i\in\bigcup_{n=1}^{\infty}\Lambda^{n}\}$
is a sub-semigroup of $G_{\Phi}$, and $\overline{S}$ is a closed
subgroup of $G_{\Phi}$ (it is a general fact that a closed sub-semigroup
of a compact group is a group). Since $\{U_{i}\}_{i\in\Lambda}\subseteq S$
in fact $\overline{S}=G_{\Phi}$. Since $S$ is dense in $\overline{S}$
and $S$ is the increasing union $S=\bigcup_{n=1}^{\infty}\{U_{i}\,:\,i\in\Lambda^{\leq n}\}$,
we can choose $k_{0}$ large enough that for every $V$, $g\in G_{\Phi}$
there is a $i\in\Lambda^{\leq k_{0}}$ with $d(U_{i}^{-1}V,gV)<\varepsilon$.

Fix $0\leq k\leq k_{0}$. Since $\mu=\sum_{i\in\Lambda^{k}}p_{i}\cdot\varphi_{i}\mu$
we can apply Lemma \ref{lem:saturation-at-level-n} \eqref{enu:saturation-at-level-n-transfers-to-components}
with a small parameter $\delta$. Writing $\varepsilon'=\sqrt{\delta}$,
it follows that if $V\in\sat(\mu,\delta,m,n)$ for some $m$ and $n>n_{0}$,
then $V\in\sat(\varphi_{j}\mu,\varepsilon',m,n)$ for all $j\in\Lambda^{k}$
outside a set $J\subseteq\Lambda^{k}$ with $\sum_{j\in J}p_{j}<\varepsilon'$.
Choose $\delta$ small enough that $p_{j}>\varepsilon'$ for all $j\in\Lambda^{k}$
(this requires $\delta$ small in a manner depending only on $k_{0}$,
and hence only on $\varepsilon$). Thus we have shown that if and
$V\in\sat(\mu,\delta,m,n)$ for some $m$ and $n>n_{0}$, then $V\in\sat(\varphi_{j}\mu,\varepsilon',m,n)$
for all $j\in\Lambda^{k}$. By Lemma \ref{lem:saturation-at-level-n}
\eqref{enu:saturation-at-level-n-similarity-2}, this in turn implies
that is $U_{j}^{-1}V\in\sat(\mu,\varepsilon'',m,n)$, where $\varepsilon''$
can be made $<\varepsilon$ if $\varepsilon'$ and $k_{0}/m$ (and
hence $k/m)$ are small enough. This holds if $\delta$ is small and
$m$ large relative to $\varepsilon$ (and hence $k_{0}$), and the
claim follows from our choice of $k_{0}$.
\end{proof}
The next proposition says, roughly, that there is an essentially maximal
$(\varepsilon,m)$-saturated subspace at each small enough scale $n$,
and that it is $(G_{\Phi},\varepsilon)$-invariant.
\begin{prop}
\label{prop:maximal-level-n-saturated-subspace}For every $0<\varepsilon<\frac{1}{10}$
there exists a $\delta=\delta(\varepsilon)>0$ such that, for $m>m(\varepsilon)$
and $n>n(\varepsilon,m)$ there exists a $(G_{\Phi},\varepsilon)$-invariant
subspace $V_{n}^{*}\in\sat(\mu,\varepsilon,m,n)$ such that every
$W\in\sat(\mu,\delta,m,n)$ satisfies $W\sqsubseteq(V_{n}^{*})^{(\varepsilon)}$.\end{prop}
\begin{proof}
Fix $\varepsilon>0$ and apply the previous lemma to $\varepsilon'=\varepsilon^{d!}/2$
to obtain $\delta'$ and set $\delta=\delta'/2$. Suppose $m$ and
$n$ are large enough to satisfy the conclusion of that lemma. Assume
that $m$ is also large enough that, for a suitable parameter $\varepsilon''<\varepsilon$,
the following holds: if $V_{1},V_{2}\leq\mathbb{R}^{d}$ are subspaces
with $\angle(V_{1},V_{2})>\varepsilon'$ and $\mu$ is $(V_{i},\varepsilon'',m)$-saturated
for $i=1,2$ then $\mu$ is $(V_{1}+V_{2},3\varepsilon',m)$-saturated
(such an $m$ and $\varepsilon''$ exists by Lemma \ref{lem:concentration-saturation-uniformity-under-change-of-subspace}
\eqref{enu:saturation-direct-sum}). Also assume that $m$ is large
enough that if $\mu$ is $(V,\delta,m)$-saturated and $V'\leq V$,
then $\mu$ is $(V,\delta',m)$-saturated (such $m$ exist by Lemma
\ref{lem:concentration-saturation-uniformity-under-change-of-subspace}
\eqref{enu:saturation-subspaces}, using $\delta'=2\delta$). 

By the choice of $\delta'$, if $V\in\sat(\mu,\delta',m,n)$ and $g\in G_{\Phi}$,
then there is a subspace $W\leq\mathbb{R}^{d}$ such that $d(W,gV)<\varepsilon'$
and $W\in\sat(\mu,\varepsilon',m,n)$. Let $\mathcal{W}$ denote the
set of all one-dimensional subspaces $W$ that arise in this way,
and write 
\[
E=\{w\in\mathbb{R}^{d}\,:\,\left\Vert w\right\Vert =1\mbox{ and }\mathbb{R}w\in\mathcal{W}\}.
\]
Observe that if $w\in E$ then $W=\mathbb{R}w\in\mathcal{W}$ and
there exists a $V\in\sat(\mu,\delta',m,n)$ and $g\in G_{\Phi}$ with
$d(W,gV)<\varepsilon'$, hence for every $h\in G_{\Phi}$ we have
$d(hW,hgV)<\varepsilon'$. By definition of $\mathcal{W}$ there is
some $W'\in\mathcal{W}$ such that $d(W',hgV)<\varepsilon'$, so $d(hW,W')<2\varepsilon'=\varepsilon^{d!}$.
Thus there is $w'\in E$ with $W'=\mathbb{R}w'$ and $d(w,w')<\varepsilon^{d!}$. 

It follows that the set $E$ satisfies the hypothesis of Lemma \ref{lem:almost-invariant-subspace-from-almost-invariant-set}
for $\varepsilon$ and the group $G_{\Phi}$. Choosing a maximal sequence
of unit vectors $v_{1},\ldots,v_{k}\in E$ such that $d(v_{i},\spn\{v_{1},\ldots,v_{i-1}\})>\varepsilon^{i!}$
and setting $V=\spn\{v_{1},\ldots,v_{k}\}$, we conclude that $V$
is $(G_{\Phi},O(\varepsilon))$-invariant and $E\subseteq V^{(\varepsilon)}$. 

Since $V=\oplus_{i=1}^{k}\mathbb{R}v_{i}$ and $\angle(v_{i},\spn\{v_{1},\ldots,v_{i-1}\})>\varepsilon^{i!}$,
and $\mathbb{R}v_{i}\in\sat(\mu,\varepsilon^{d!}/2,m,n)$ for all
$i$, repeated application of Lemma \ref{lem:concentration-saturation-uniformity-under-change-of-subspace}
\eqref{enu:saturation-direct-sum}, assuming $m$ large enough relative
to $\varepsilon$ (and hence $\varepsilon'$), gives that $V\in\sat(\mu,O(\varepsilon),m,n)$. 

Finally, if $W\in\sat(\mu,\delta,m,n)$ then we can choose an orthonormal
basis $\{w_{i}\}$ for $W$, so by choice of $m$, $\mathbb{R}w_{i}\in\sat(\mu,\delta',m,n)$,
so $w_{i}\in E$. By Lemma \ref{lem:almost-invariant-subspace-from-almost-invariant-set},
$w_{i}\in V^{(\varepsilon)}$. The $w_{i}$ are orthonormal, so $d(w_{i},\spn\{w_{1},\ldots,w_{j-1}\}=1$.
Hence by Corollary \ref{cor:span-is-close-to-V-if-basis-is-close-to-V},
$W\sqsubseteq V^{(O(\varepsilon))}$.

We have proved the claim for $V_{n}^{*}=V$, up to some constant factors,
to remove them begin with a small multiple of $\varepsilon$ instead
of $\varepsilon$.
\end{proof}
The next proposition allows us to replace a saturated almost-invariant
subspace with a truly invariant one, of some lesser saturation. It
also shows that this new subspace is saturated at many levels, even
though the original subspace a-priori was saturated at a single level.
\begin{prop}
\label{prop:from-almost-to-truely-invariant-subspaces} For every
$\varepsilon>0$, $0<\delta<\delta(\varepsilon)$, $m>m(\varepsilon,\delta)$
and every $n\in\mathbb{N}$, the following holds. If $W\in\sat(\mu,\delta,m,n)$
is $(G_{\Phi},\delta)$-invariant and $\widetilde{W}$ is a $G_{\Phi}$-invariant
subspace with $d(W,\widetilde{W})<\delta$, then for $m'=[\log(2/\delta)]$
and all large enough $n'$ we have $\widetilde{W}\in\sat(\mu,\varepsilon,m',n')$.\end{prop}
\begin{proof}
Fix $0<\delta<\varepsilon$. Also fix $m$ large relative to $\delta$
(we shall see how large later). Let $n\in\mathbb{N}$, $W\leq\mathbb{R}^{d}$
and $m',n'$ be as in the statement, so our assumption is that 
\[
\mathbb{P}_{i=n}\left(\mu^{x,i}\mbox{ is }(W,\delta,m)\mbox{-saturated}\right)>1-\delta.
\]
For each measure $\theta=\mu^{x,n}$ in the event above, writing $\delta_{1}=\sqrt{d\delta+O(\frac{m'}{m})}$,
Lemma \ref{lem:saturation-passes-to-components} implies 
\[
\mathbb{P}_{0\leq j\leq m}\left(\theta^{y,j}\mbox{ is }(W,\delta_{1},m')\mbox{-saturated}\right)>1-\delta_{1}.
\]
Assuming $m$ is large relative to $\delta$ (and hence $m'$), we
can arrange $\delta_{1}<2\sqrt{d\delta}$. Combining the two inequalities
above, we can find a $0\leq k\leq m$ such that
\[
\mathbb{P}_{i=n+k}\left(\mu^{x,i}\mbox{ is }(W,\delta_{1},m')\mbox{-saturated}\right)>1-2\delta_{1}.
\]
Let $\mu^{x,n+k}$ be as in this last event. Since $d(W,\widetilde{W})<\delta<2^{-m'}$,
by Lemma \ref{lem:concentration-saturation-uniformity-under-change-of-subspace}
\eqref{enu:saturation-passes-to-nearby-subspaces}, $\mu^{x,n+k}$
is also $(\widetilde{W},\delta_{2},m')$-saturated, where $\delta_{2}=\delta_{1}+O(1/m')$.
Since this holds for a $1-2\delta_{1}>1-\delta_{2}$ proportion of
components $\mu^{x,n+k}$, (because, if $\delta$ is small, $\delta_{2}\geq2\delta_{1}$),
we have $\widetilde{W}\in\sat(\mu,\delta_{2},m',n+k)$. Note that
$\delta_{2}$ can be made arbitrarily small by choosing $\delta$
small enough.

Finally, let $n'>n+k$. Let $\Lambda(n')\subseteq\bigcup_{j=1}^{\infty}\Lambda^{j}$
denote the set of sequences $i=i_{1}\ldots i_{\ell}$ such that $r_{i_{1}}\cdot\ldots\cdot r_{i_{\ell}}<2^{-(n'-k)}\leq r_{i_{1}}\cdot\ldots\cdot r_{i_{\ell-1}}$.
Then $\sum_{i\in\Lambda(n')}p_{i}=1$ and $\mu=\sum_{i\in\Lambda(n')}p_{i}\varphi_{i}\mu$.
By Lemma \ref{lem:saturation-at-level-n} \eqref{enu:saturation-at-level-n-similarity-1},
$\widetilde{W}=U_{i}\widetilde{W}\in(\varphi_{i}\mu,\delta_{3},m',n')$
for all $i\in\Lambda(n')$, where $\delta_{3}\rightarrow0$ as $\delta\rightarrow0$
and $m'\rightarrow\infty$. Since $\mu$ is a convex combination of
the measures $\varphi_{i}\mu$, $i\in\Lambda(n')$, by Lemma \ref{lem:saturation-at-level-n}
\eqref{enu:saturation-at-level-n-of-convex-combinations} we have
$\widetilde{W}\in(\mu,\delta_{4},m',n')$ for $\delta_{4}$ which
can be made arbitrarily small (and in particular $<\varepsilon$)
if $\delta$ is small and $m'$ large. This completes the proof.
\end{proof}
Finally, we show the existence of a ``maximal'' invariant subspace
which is saturated to all degrees at sufficiently deep levels. Let
us say that a $\mu$ is $V$-saturated if $\mu\in\sat(V,\varepsilon,m,n)$
for all $\varepsilon>0$, $m\geq m(\varepsilon)$ and all $n>n(\varepsilon,m)$.
\begin{prop}
\label{prop:maximal-invariant-saturated-subspace}There exists a unique
subspace $\widetilde{V}\leq\mathbb{R}^{d}$ such that 
\begin{enumerate}
\item $\mu$ is $\widetilde{V}$-saturated.
\item $V\subseteq\widetilde{V}$ whenever $\mu$ is $V$-saturated. 
\item [3.] $\widetilde{V}$ is $G_{\Phi}$-invariant.
\end{enumerate}
\end{prop}
\begin{proof}
A formal consequence of Lemma \ref{lem:concentration-saturation-uniformity-under-change-of-subspace}
\eqref{enu:saturation-direct-sum} is that if $\mu$ is $\widetilde{V}_{1}$-saturated
and $\widetilde{V}_{2}$-saturated then $\mu$ is $\widetilde{V}_{1}+\widetilde{V}_{2}$-saturated.
Thus we can take $\widetilde{V}$ to be the sum of all subspaces $V$
on which $\mu$ is saturated. (1) and (2) are then obvious, and (3)
is a formal consequence of Lemma \ref{lem:invariant-subspace-near-almost-invariant-one},
Proposition \ref{prop:maximal-level-n-saturated-subspace} and Propositions
\ref{prop:from-almost-to-truely-invariant-subspaces}, because taken
together they show that if $\mu$ is $V$-saturated then $\mu$ is
$V'$-saturated for a $G_{\Phi}$-invariant subspace $V'$, and $\dim V'\geq\dim V$
. Applying this to $V=\widetilde{V}$ we conclude $V'\subseteq\widetilde{V}$
and $\dim V'\geq\dim\widetilde{V}$ so $\widetilde{V}=V'$ is $G_{\Phi}$-invariant.
\end{proof}
We now need sufficient conditions for the subspace $\widetilde{V}$
from the last proposition to be of dimension $>1$. To this end, we
have the following.
\begin{prop}
\label{prop:largeness-condition-for-Vtilde}If there exists a sequence
$V_{i}\in\sat(\mu,\varepsilon_{i},m_{i},n_{i})$ with $\varepsilon_{i}\rightarrow0$,
$m_{i}>m(\varepsilon_{i})$ and $n_{i}>n(\varepsilon_{i},m_{i})$,
and if $V_{i}\rightarrow V$, then $V\subseteq\widetilde{V}$, where
$\widetilde{V}$ is as in Proposition \ref{prop:maximal-invariant-saturated-subspace}.\end{prop}
\begin{proof}
In each of the three previous propositions, a $\delta=\delta(\varepsilon)$
was associated to an $\varepsilon$. We can assume that these functions
$\delta$ are increasing (so decreasing $\varepsilon$ leads to no
increase in $\delta(\varepsilon)$).

Let $V_{i},\varepsilon_{i},m_{i},n_{i}$ be given. Assuming that $m_{i},n_{i}$
are large enough relative to $\varepsilon_{i}$, by Proposition \ref{prop:maximal-level-n-saturated-subspace}
there is a sequence $\varepsilon'_{i}=\varepsilon'_{i}(\varepsilon_{i})\rightarrow0$
depending monotonely on $\varepsilon_{i}$, such that for each $i$
there is a $(G_{\Phi},\varepsilon'_{i})$-invariant subspace $V_{i}^{*}\in\sat(\mu,\varepsilon_{i}',m_{i},n_{i})$
with $\angle(V_{i},V_{i}^{*})<\varepsilon'_{i}$ and $\dim V_{i}^{*}$$\geq\dim V_{i}$
(if $\delta(\cdot)$ is the function in that proposition than we choose
$\varepsilon'_{i}=\delta^{-1}(\varepsilon_{i})$). 

We can henceforth assume that $m_{i}$ are large enough relative to
$\varepsilon'_{i}$, and $n_{i}$ relative to $\varepsilon'_{i},m_{i}$
(here we use that $\varepsilon_{i}'$ depends on $\varepsilon_{i}$
in a monotone way, so being large with respect to $\varepsilon_{i}'$
is the same as being large with respect to $\varepsilon_{i}$, which
was assumed). 

Passing to a subsequence we may assume that $V_{i}^{*}$ converge
to some subspace $V^{*}$. Note that $V^{*}$ is $G_{\Phi}$-invariant,
being the limit of $(G_{\Phi},\varepsilon'_{i})$-invariant subspaces.

By increasing $\varepsilon'_{i}$ if needed, we can assume that $m'_{i}=[\log(2/\varepsilon'_{i})]\rightarrow\infty$
more slowly than linearly. 

By Proposition \ref{prop:from-almost-to-truely-invariant-subspaces}
we can choose $\varepsilon''_{i}\rightarrow0$, depending monotonely
on $\varepsilon'_{i}$ such that if $W\in\sat(\mu,\varepsilon'_{i},m_{i},n_{i})$
is a $G_{\Phi}$-invariant subspace, and $d(V_{i}^{*},W)<\varepsilon'_{i}$,
then $W\in\sat(\mu,\varepsilon''_{i},m'_{i},n')$, for all $n'>n_{i}+m_{i}'$
(recall that $m'_{i}=[\log(2/\varepsilon'_{i})]$; if $\delta(\cdot)$
is the function in that proposition, choose $\varepsilon''_{i}=\delta^{-1}(\varepsilon'_{i})$).
Note that since we assumed that $m'_{i}\rightarrow\infty$ more slowly
than linearly, every large enough integer occurs as $m'_{i}$ for
some $i$.

Applying the previous paragraph to $W=V^{*}$, and since we have arranged
that $\{m'_{i}\}$ includes all large enough integers, we see that
$\mu$ if $V^{*}$-saturated. Thus $V^{*}\subseteq\widetilde{V}$.

Finally, combining $\angle(V_{i},V_{i}^{*})\rightarrow0$ with $V_{i}^{*}\rightarrow V^{*}$
and $V_{i}\rightarrow V$, we conclude that  $\angle(V^{*},V)=0$.
Since $\dim V^{*}=\lim\dim V_{i}^{*}\geq\lim\dim V_{i}=\dim V$, we
must have $V\subseteq V^{*}$. Since $V^{*}\subseteq\widetilde{V}$
we get $V\subseteq\widetilde{V}$, as claimed.
\end{proof}

\subsection{\label{sub:Entropy-and-dimension}Entropy and dimension for self-similar
measures}

If $\mu\in\mathcal{P}([0,1]^{d})$ is exact dimensional, as self-similar
measures are, the dimension of $\mu$ is given by the so-called entropy
dimension:

\begin{equation}
\dim\mu=\lim_{n\rightarrow\infty}H_{n}(\mu).\label{eq:103}
\end{equation}

We require a similar expression relating the dimension of conditional
measures on affine subspaces to entropy. We parametrized affine subspaces
as the set of fibers $\pi^{-1}(y)$ where $\pi$ is a linear map $\pi:\mathbb{R}^{d}\rightarrow\mathbb{R}^{k}$
and $y$ ranges over $\mathbb{R}^{k}$. The conditional measure $\mu_{\pi^{-1}(y)}$
of $\mu$ on $\pi^{-1}(y)$ is defined for $\pi\mu$-a.e. $y$ by
the weak-{*} limit 
\[
\mu_{\pi^{-1}(y)}=\lim_{\ell\rightarrow\infty}\mu_{\pi^{-1}D_{\ell}^{k}(y)},
\]
which exists by the measure-valued version of the Martingale convergence
theorem. 
\begin{thm}
\label{thm:local-entropy-averages-for-fibers}Let $\mu\in\mathcal{P}(\mathbb{R}^{d})$
be a self similar measure for the IFS $\Phi$ and let $\pi:\mathbb{R}^{d}\rightarrow\mathbb{R}^{k}$
be a linear map such that $\ker\pi$ is $D\Phi$-invariant. Then the
conditional measure $\mu_{\pi^{-1}(y)}$ is exact dimensional for
$\pi\mu$-a.e. $y$, and the dimension is given by
\begin{eqnarray*}
\dim\mu_{\pi^{-1}(y)} & = & \lim_{p\rightarrow\infty}\left(\liminf_{n\rightarrow\infty}\mathbb{E}_{0\leq i\leq n}\left(\frac{1}{p}H(\mu_{x,i},\mathcal{D}_{i+p}|\pi^{-1}\mathcal{D}_{i+p}^{k})\right)\right)\\
 & = & \lim_{p\rightarrow\infty}\left(\limsup_{n\rightarrow\infty}\mathbb{E}_{0\leq i\leq n}\left(\frac{1}{p}H(\mu_{x,i},\mathcal{D}_{i+p}|\pi^{-1}\mathcal{D}_{i+p}^{k})\right)\right)
\end{eqnarray*}

\end{thm}
We will apply this theorem via the following corollary:
\begin{cor}
If $\mu$ is self-similar and $V$ is a saturated and $G_{\Phi}$-invariant
subspace, then the conditional measures of $\mu$ on translates of
$V$ are a.s. exact dimensional and of dimension $\dim V$. In particular
this holds for the subspace described in Proposition \ref{prop:maximal-invariant-saturated-subspace}.
\end{cor}
The proof we present for Theorem \ref{thm:local-entropy-averages-for-fibers}
has two ingredients. The first is exact dimensionality and dimension
conservation:
\begin{thm}
\label{thm:dimension-conservation}Let $\mu\in\mathcal{P}(\mathbb{R}^{d})$
be a self similar measure for the IFS $\Phi$ and let $\pi:\mathbb{R}^{d}\rightarrow\mathbb{R}^{k}$
be a linear map such that $\ker\pi$ is $D\Phi$-invariant. Then $\pi\mu$
is exact dimensional, $\mu_{\pi^{-1}(y)}$ is exact dimensional for
$\pi\mu$-a.e. $y$, its dimension is $\pi\mu$-a.s. independent of
$y$, and 
\[
\dim\pi\mu+\dim\mu_{\pi^{-1}(y)}=\dim\mu\qquad\mbox{for }\pi\mu\mbox{-a.e. }y.
\]

\end{thm}
This theorem follows from work of Falconer and Jin \cite{FalconerJin2014}
(which in turn relies on methods of Feng and Hu \cite{FengHu09}).
Next, we require an expression for $\dim\pi\mu$ in terms of entropy
of dyadic partitions. A special case of this result appears in \cite{HochmanShmerkin2012}
for the case that $G_{\Phi}$ is the full orthogonal group.
\begin{thm}
\label{thm:local-entropy-averages-for-projections}Let $\mu\in\mathcal{P}(\mathbb{R}^{d})$
be a self similar measure for the IFS $\Phi$ and let $\pi:\mathbb{R}^{d}\rightarrow\mathbb{R}^{k}$
be a linear map such that $\ker\pi$ is $D\Phi$-invariant. Then 
\begin{eqnarray*}
\dim\pi\mu & = & \lim_{p\rightarrow\infty}\left(\liminf_{n\rightarrow\infty}\mathbb{E}_{0\leq i\leq n}\left(\frac{1}{p}H(\mu_{x,i},\pi^{-1}\mathcal{D}_{i+p}^{k}\right)\right)\\
 & = & \lim_{p\rightarrow\infty}\left(\limsup_{n\rightarrow\infty}\mathbb{E}_{0\leq i\leq n}\left(\frac{1}{p}H(\mu_{x,i},\pi^{-1}\mathcal{D}_{i+p}^{k}\right)\right).
\end{eqnarray*}
\end{thm}
\begin{proof}
First, note that 
\[
\mathbb{E}_{0\leq i\leq n}\left(\frac{1}{p}H(\mu_{x,i},\pi^{-1}\mathcal{D}_{i+p}^{k})\right)=\mathbb{E}_{0\leq i\leq n}\left(\frac{1}{p}H(\mu,\pi^{-1}\mathcal{D}_{i+p}^{k}|\mathcal{D}_{i})\right),
\]
As we have seen, changing the dyadic partition to one adapted to a
different coordinate system changes the right hand side of the last
equation by $O(1/p)$, and in the statement of the theorem we consider
the limit as $p\rightarrow\infty$. Thus, the statement is unaffected
by changes to the coordinate system, and we may assume that $\pi$
is a coordinate projection. Therefore we can apply the local entropy
averages lemma for projections \cite{HochmanShmerkin2012}. The lemma
is usually formulated for lower pointwise dimension, but the same
proof exactly, replacing $\liminf$ by $\limsup$, shows that
\[
\limsup_{n\rightarrow\infty}-\frac{1}{n}\log(\mu((\pi^{-1}\mathcal{D}_{n}^{k})(x)))\geq\limsup_{n\rightarrow\infty}\frac{1}{n}\sum_{i=0}^{n-1}\frac{1}{p}H(\mu_{x,i},\pi^{-1}\mathcal{D}_{i+p}^{k})-O(\frac{1}{p})\qquad\mu\mbox{-a.e. }x.
\]
Since $\pi\mu$ is exact dimensional, the left hand side is $\mu$-a.s.
equal to $\dim\pi\mu$, and we have 
\[
\dim\pi\mu\geq\limsup_{n\rightarrow\infty}\frac{1}{n}\sum_{i=0}^{n-1}\frac{1}{p}H(\mu_{x,i},\pi^{-1}\mathcal{D}_{i+p}^{k})-O(\frac{1}{p})\qquad\mu\mbox{-a.e. }x
\]
Integrating this $d\mu$ and using Fatou's lemma, for all $p$, 
\begin{eqnarray}
\dim\pi\mu & \geq & \limsup_{n\rightarrow\infty}\int\left(\frac{1}{n}\sum_{i=0}^{n-1}\frac{1}{p}H(\mu_{x,i},\pi^{-1}\mathcal{D}_{i+p}^{k})\right)d\mu(X)-O(\frac{1}{p})\nonumber \\
 & = & \limsup_{n\rightarrow\infty}\mathbb{E}_{0\leq i\leq n}\left(\frac{1}{p}H(\mu_{x,i},\pi^{-1}\mathcal{D}_{i+p}^{k})\right)-O(\frac{1}{p}).\label{eq:104}
\end{eqnarray}

Equation \eqref{eq:104} is one half of the inequality we are after,
and its proof only used exact dimensionality of $\mu$. For the reverse
inequality we will use self-similarity. Fix $p$, and note that we
have the identity
\[
\mathbb{E}_{0\leq i\leq n}\left(\frac{1}{p}H(\mu_{x,i},\pi^{-1}\mathcal{D}_{i+p}^{k})\right)=\mathbb{E}_{0\leq i\leq n}\left(\frac{1}{p}H(\mu,\pi^{-1}\mathcal{D}_{i+p}^{k}|\mathcal{D}_{i})\right)
\]
where the expectation on the left is over $i$ and $x$, and on the
right only over $i$. Let $r=\min\{r_{i}\,:\,i\in\Lambda\}$, and
for each $i$ let $I_{i}\subseteq\Lambda^{*}$ denote the set of sequences
$j_{1}\ldots j_{k}\in\Lambda^{*}$ such that $r\cdot2^{-i}<r_{j_{1}}\ldots r_{j_{k}}<2^{-i}\leq r_{j_{1}}\ldots r_{j_{k-1}}$.
It is a standard (and easy) fact that $\mu=\sum_{j\in I_{i}}p_{j}\cdot\varphi_{j}\mu$.
By concavity of conditional entropy (Lemma \ref{lem:entropy-combinatorial-properties}\eqref{enu:entropy-concavity}),
for each $i$,
\begin{eqnarray*}
\frac{1}{p}H(\mu,\pi^{-1}\mathcal{D}_{i+p}^{k}|\mathcal{D}_{i}) & \geq & \frac{1}{p}\sum_{i\in I_{i}}p_{i}H(\varphi_{i}\mu,\pi^{-1}\mathcal{D}_{i+p}|\mathcal{D}_{i})+O(\frac{1}{p})\\
 & = & \frac{1}{p}\sum_{i\in I_{i}}p_{i}H(\varphi_{i}\mu,\pi^{-1}\mathcal{D}_{i+p})+O(\frac{1}{p}).
\end{eqnarray*}
where we used the fact that each $\varphi_{i}\mu$, $i\in I_{i}$,
has diameter $O(2^{-i})$, and Lemma \ref{lem:entropy-weak-continuity-properties}\eqref{enu:entropy-combinatorial-distortion}.
Finally, since $\varphi_{i}$ contracts by $2^{-i}$ up to a constant
factor, by changing scale, applying Lemma \ref{lem:entropy-weak-continuity-properties}\eqref{enu:entropy-change-of-scale},
and changing the coordinates system, we have
\[
\frac{1}{p}\sum_{i\in I_{i}}p_{i}H(\varphi_{i}\mu,\pi^{-1}\mathcal{D}_{i+p})=\frac{1}{p}H(\mu,\pi^{-1}\mathcal{D}_{p})+O(\frac{1}{p}).
\]
Note that we used here the fact that $\ker\pi$ is invariant under
the linear part of $\varphi_{i}$.

Putting this all together, we have shown that for every $p$,
\[
\mathbb{E}_{0\leq i\leq n}\left(\frac{1}{p}H(\mu_{x,i},\pi^{-1}\mathcal{D}_{i+p}^{k})\right)\geq\frac{1}{p}H(\pi\mu,\mathcal{D}_{p})+O(\frac{1}{p}).
\]
Taking the $\liminf$ as $n\rightarrow\infty$, we have
\[
\liminf_{n\rightarrow\infty}\mathbb{E}_{0\leq i\leq n}\left(\frac{1}{p}H(\mu_{x,i},\pi^{-1}\mathcal{D}_{i+p}^{k})\right)\geq\frac{1}{p}H(\pi\mu,\mathcal{D}_{p})+O(\frac{1}{p}).
\]
But, since $\pi\mu$ is exact dimensional, as $p\rightarrow\infty$
the right hand side tends to $\dim\pi\mu$. Combined with inequality
\eqref{eq:104}, this proves the statement.
\end{proof}
We can now prove Theorem \ref{thm:local-entropy-averages-for-fibers}.
Begin with the identity
\begin{eqnarray*}
\mathbb{E}_{0\leq i\leq n}\left(\frac{1}{p}H(\mu_{x,i},\mathcal{D}_{i+p}|\pi^{-1}\mathcal{D}_{i+p}^{k})\right) & = & \mathbb{E}_{0\leq i\leq n}\left(\frac{1}{p}H(\mu_{x,i},\mathcal{D}_{i+p})\right)\\
 &  & \;-\quad\mathbb{E}_{0\leq i\leq n}\left(\frac{1}{p}H(\mu_{x,i},\pi^{-1}\mathcal{D}_{i+p}^{k})\right)
\end{eqnarray*}
(this is just Lemma \ref{lem:entropy-combinatorial-properties} \eqref{enu:entropy-conditional-formula}
and linearity of expectation). Taking $n\rightarrow\infty$ and then
$p\rightarrow\infty$, and using \eqref{eq:103} and Theorem \ref{thm:local-entropy-averages-for-projections},
the right hand side becomes $\dim\mu-\dim\pi\mu$, which by Theorem
\ref{thm:dimension-conservation} is the a.s. dimension of fibers.

\subsection{\label{sub:proof-of-main-theorem-Rd}Proof of Theorem \ref{thm:main-individual-RD-entropy}}

Recall from the introduction that $r=\prod_{i\in\Lambda}r_{i}^{p_{i}}$,
$n'=n\log(1/r)$ and $\nu^{(n)}=\sum_{i\in\Lambda^{n}}p_{i}\cdot\delta_{\varphi_{i}}$.
Also recall the definition of the dyadic partition $\mathcal{D}_{n}=\mathcal{D}_{n}^{G}$,
and the partition $\mathcal{E}_{n}=\mathcal{E}_{n}^{G}$ of $G$ according
to the level-$n$ dyadic partition of the translation part of the
maps. In this section we prove the following:
\begin{thm}
\label{thm:main-individual-RD-entropy-1}Let $\Phi=\{\varphi_{i}\}$
be an IFS on $\mathbb{R}^{d}$ that does not preserve a non-trivial
affine subspace, and $\mu$ a self-similar measure for $\Phi$. Then
either
\[
\lim_{n\rightarrow\infty}\frac{1}{n'}H(\nu^{(n)},\mathcal{D}_{qn}^{G}|\mathcal{E}_{n'}^{G})=0\quad\mbox{for all }q>1,
\]
or else there is a $D\Phi$-invariant subspace $V$ such that $\dim\mu_{V+x}=\dim V$
for $\mu$-a.e. $x$.
\end{thm}
This implies Theorem \ref{thm:main-individual-RD-entropy}, see remark
after its statement.

We begin the proof. First, note that $\mu(V)=0$ for every proper
affine subspace $V\subseteq\mathbb{R}^{d}$, since if $\mu(V)>0$
for some $V$ then it is easily shown that $\mu$ is supported on
$V$, and hence $\Phi$ preserves $V$, contrary to hypothesis. 

We now argue by contradiction: suppose that there is a $\delta_{0}>0$
and $q>1$ such that 
\[
\limsup_{n\rightarrow\infty}\frac{1}{qn'}H(\nu^{(n)},\mathcal{D}_{(q+1)n'}|\mathcal{E}_{n'})>\delta_{0}.
\]
Let $\mathcal{F}$ denote the partition of $G$ according to the contraction
ratio. This is an uncountable partition, but the possible contractions
of $\varphi_{i}$, $i\in\Lambda^{n}$, are just all the $n$-fold
products of the contractions $r_{i}$, $i\in\Lambda$. Thus only $O(n^{|\Lambda|+1})$
distinct contraction ratios occur in the support of $\nu^{(n)}$,
so 
\[
\lim_{n\rightarrow\infty}\frac{1}{qn'}H(\nu^{(n)},\mathcal{F})=\lim_{n\rightarrow\infty}\frac{O(\log n)}{n'}=0.
\]
Using the identities $H(\cdot,\mathcal{D}|\mathcal{E}\lor\mathcal{F})=H(\cdot,\mathcal{D}|\mathcal{E})+H(\cdot,\mathcal{F}|\mathcal{E})$
and $H(\cdot,\mathcal{F}|\mathcal{E})\leq H(\cdot,\mathcal{F})$,
the two limits above imply 
\begin{equation}
\limsup_{n\rightarrow\infty}\frac{1}{qn'}H(\nu^{(n)},\mathcal{D}_{(q+1)n'}|\mathcal{E}_{n'}\lor\mathcal{F})>\delta_{0}.\label{eq:38}
\end{equation}

\begin{lem}
$\lim_{n\rightarrow\infty}\int\frac{1}{qn'}H(g\conv\mu,\mathcal{D}_{(q+1)n'}|\mathcal{D}_{n'})\,d\nu^{(n)}(g)=\dim\mu$\end{lem}
\begin{proof}
If $g=2^{-t}U+a$, then $g\conv\mu$ is supported on a set of diameter
$O(2^{-t})$, hence $H(g\conv\mu,\mathcal{D}_{n'})=O(|t-n'|)$. Similarly,
by Lemma \ref{lem:entropy-weak-continuity-properties} \eqref{enu:entropy-change-of-scale},
$H(g\conv\mu,\mathcal{D}_{(q+1)n'})=H(\mu,\mathcal{D}_{qn'})+O(|t-n'|)$. 

If we choose $g=2^{-t}U+a$ randomly according to $\nu^{(n)}$, then
$t$ is distributed as the sum of $n$ independent random variables,
each of which takes value $\log(1/r_{i})$ with probability $p_{i}$
for $i\in\Lambda$, so by the law of large numbers, $t-n'=o(n')$
in probability. We also have a worst-case bound of $t\leq Cn$ (a.s.
for $g\sim\nu^{(n)}$), because $\varphi_{i_{1},\ldots,i_{n}}$ contracts
by at least $(\min_{i\in\Lambda}r_{i})^{n}$, and $\min_{i\in\Lambda}r_{i}<1$.
Hence the bound $t-n'=o(n')$ holds also in the mean sense. It follows
from the first paragraph that 
\begin{eqnarray*}
\frac{1}{qn'}\int H(g\conv\mu,\mathcal{D}_{n'})d\nu^{(n)}(g) & = & o(1)\\
\frac{1}{qn'}\int H(g\conv\mu,\mathcal{D}_{(q+1)n'})d\nu^{(n)}(g) & = & \frac{1}{qn'}H(\mu,\mathcal{D}_{qn'})+o(1)\\
 & = & \dim\mu+o(1)
\end{eqnarray*}
Subtracting the first line from the second proves the claim.\end{proof}
\begin{lem}
$\lim_{n\rightarrow\infty}\frac{1}{qn'}H(\mu,\mathcal{D}_{(q+1)n'}|\mathcal{D}_{n'})=\dim\mu$.\end{lem}
\begin{proof}
Using the conditional entropy formula, 
\begin{eqnarray*}
\frac{1}{(q+1)n'}H(\mu,\mathcal{D}_{(q+1)n'}) & = & \frac{1}{(q+1)}\cdot\frac{1}{n'}H(\mu,\mathcal{D}_{n'})+\\
 &  & \;+\;\frac{q}{(q+1)}\cdot\frac{1}{qn'}H(\mu,\mathcal{D}_{q(n'+1)}|\mathcal{D}_{n'}).
\end{eqnarray*}
The lemma follows by taking $n\rightarrow\infty$ and using the fact
that $\frac{1}{n}H(\mu,\mathcal{D}_{n})\rightarrow\dim\mu$.
\end{proof}
Let $\nu_{I}^{(n)}$ denote, as usual the, conditional measure of
$\nu^{(n)}$ on $I$.
\begin{lem}
$\lim_{n\rightarrow\infty}\left(\sum_{I\in\mathcal{E}_{n'}\lor\mathcal{F}}\nu(I)\cdot\frac{1}{qn'}H(\nu_{I}^{(n)}\conv\mu,\mathcal{D}_{(q+1)n'})\right)=\dim\mu$.\end{lem}
\begin{proof}
Write 
\[
\mu=\sum_{I\in\mathcal{E}_{n'}\lor\mathcal{F}}\nu(I)\cdot\left(\nu_{I}^{(n)}\conv\mu\right).
\]
and note that
\[
\nu_{I}^{(n)}\conv\mu=\int g\conv\mu\;d\nu_{I}^{(n)}(g)
\]
Combining this with concavity of conditional entropy (Lemma \ref{lem:entropy-combinatorial-properties}
\eqref{enu:entropy-concavity}) and the previous two lemmas, 
\begin{eqnarray*}
\dim\mu & = & \lim_{n\rightarrow\infty}\frac{1}{qn'}H(\mu,\mathcal{D}_{(q+1)n'}|\mathcal{D}_{n'})\\
 & \geq & \limsup_{n\rightarrow\infty}\sum_{I\in\mathcal{E}_{n'}\lor\mathcal{F}}\nu^{(n)}(I)\cdot\frac{1}{qn'}H(\nu_{I}^{(n)}\conv\mu,\mathcal{D}_{(q+1)n'}|\mathcal{D}_{n'})\\
 & \geq & \liminf_{n\rightarrow\infty}\sum_{I\in\mathcal{E}_{n'}\lor\mathcal{F}}\nu^{(n)}(I)\cdot\frac{1}{qn'}H(\nu_{I}^{(n)}\conv\mu,\mathcal{D}_{(q+1)n'}|\mathcal{D}_{n'})\\
 & \geq & \liminf_{n\rightarrow\infty}\sum_{I\in\mathcal{E}_{n'}\lor\mathcal{F}}\nu^{(n)}(I)\cdot\int\frac{1}{qn'}H(g\conv\mu,\mathcal{D}_{(q+1)n'}|\mathcal{D}_{n'})\,d\nu_{I}^{(n)}(g)\\
 & = & \lim_{n\rightarrow\infty}\int\frac{1}{qn'}H(g\conv\mu,\mathcal{D}_{(q+1)n'}|\mathcal{D}_{n'})\,d\nu^{(n)}(g)\\
 & = & \dim\mu,
\end{eqnarray*}
as claimed.
\end{proof}
For $I\in\mathcal{E}_{n'}\lor\mathcal{F}$ consisting of similarities
with contraction $2^{-t}$, define 
\[
\widetilde{\nu}_{I}^{(n)}=S_{t}\nu_{I}^{(n)}
\]
This is a measure on the isometry group $G_{0}$. 
\begin{lem}
For every $\delta>0$ and for arbitrarily large $n$ we can find $I\in\mathcal{E}_{n'}\lor\mathcal{F}$
with $\nu(I)>0$ and such that
\[
\frac{1}{qn'}H(\widetilde{\nu}_{I}^{(n)},\mathcal{D}_{qn'})>\delta_{0}
\]
and
\[
\frac{1}{qn'}H(\widetilde{\nu}_{I}^{(n)}\conv\mu,\mathcal{D}_{qn'})<\frac{1}{qn'}H(\mu,\mathcal{D}_{qn'})+\delta.
\]
\end{lem}
\begin{proof}
By \eqref{eq:38}, for infinitely many $n$ we have 
\begin{gather}
\frac{1}{qn'}\sum_{I\in\mathcal{E}_{n'}\lor\mathcal{F}}\nu^{(n)}(I)\cdot H(\nu_{I}^{(n)},\mathcal{D}_{(q+1)n'})\qquad\qquad\label{eq:38a}\\
\begin{aligned} & =\frac{1}{qn'}H(\nu^{(n)},\mathcal{D}_{(q+1)n'}|\mathcal{E}_{n'}\lor\mathcal{F})\\
 & >\delta_{0}
\end{aligned}
\nonumber 
\end{gather}
Suppose $I\in\mathcal{E}_{n'}\lor\mathcal{F}$ contains similitudes
of contraction $t$. Since the action of $S_{t}$ on $G$ is just
ordinary scaling in our coordinates on $G$, we have
\[
\left|H(\nu_{I}^{(n)},\mathcal{D}_{(q+1)n'})-H(\widetilde{\nu}_{I}^{(n)},\mathcal{D}_{qn'})\right|=O(|t-n'|),
\]
Using the fact that for $g=2^{-t}U+a\sim\nu^{(n)}$ we have $t-n'=o(n')$
in probability as $n\rightarrow\infty$, and the pointwise bound $t=O(n')$,
this and \eqref{eq:38a} imply that there are infinitely many $n$
such that
\begin{equation}
\sum_{I\in\mathcal{E}_{n'}\lor\mathcal{F}}\nu^{(n)}(I)\cdot\frac{1}{qn'}H(\widetilde{\nu}_{I}^{(n)},\mathcal{D}_{qn'})>\delta_{0}.\label{eq:39a}
\end{equation}

Similarly, we have $\widetilde{\nu}_{I}^{(n)}\conv\mu=S_{t}(\nu_{I}^{(n)}\conv\mu)$,
so by Lemma \ref{lem:entropy-combinatorial-properties} \eqref{enu:entropy-change-of-scale},
\[
\left|H(\nu_{I}^{(n)}\conv\mu,\mathcal{D}_{(q+1)n'})-H(\widetilde{\nu}_{I}^{(n)}\conv\mu,\mathcal{D}_{qn'})\right|=O(|t-n'|).
\]
Using the previous lemma and again the fact that $|t-n'|=o(n')$ in
probability as $g=2^{-t}U+a\sim\nu^{(n)}$, 
\[
\lim_{n\rightarrow\infty}\left(\frac{1}{qn'}\sum_{I\in\mathcal{E}_{n'}\lor\mathcal{F}}\nu^{(n)}(I)\cdot H(\widetilde{\nu}_{I}^{(n)}\conv\mu,\mathcal{D}_{qn'})\right)=\dim\mu.
\]
On the other hand we know that also
\[
\lim_{n\rightarrow\infty}\left(\frac{1}{qn'}H(\mu,\mathcal{D}_{qn'})\right)=\dim\mu.
\]
Therefore (using boundedness of the normalized entropy), 
\begin{equation}
\lim_{n\rightarrow\infty}\sum_{I\in\mathcal{E}_{n'}\lor\mathcal{F}}\nu^{(n)}(I)\cdot\left|\frac{1}{qn'}H(\widetilde{\nu}_{I}^{(n)}\conv\mu,\mathcal{D}_{qn'})-\frac{1}{qn'}H(\mu,\mathcal{D}_{qn'})\right|=0.\label{eq:39}
\end{equation}
Combining this with \eqref{eq:39a}, for infinitely many $n$ we can
find $I\in\mathcal{E}_{n'}\lor\mathcal{F}$ with the desired properties.
\end{proof}
Now fix a parameter $\varepsilon>0$, and let $\sigma>0$ be such
that $\mu$ is $((\varepsilon/5d)^{2(d+1)},\sigma)$-non-affine (recall
Definition \ref{def:epsilon-sigma-continuity}). Such $\sigma$ exists
because by assumption $\mu$ gives mass $0$ to every proper affine
subspace. Choose large $m\in\mathbb{N}$, and let $\delta>0$ and
$k\in\mathbb{N}$, be as in the conclusion of Theorem \ref{thm:inverse-theorem-for-isometries}.
Apply the theorem to the measures $\widetilde{\nu}_{I}^{(n)}$ for
the set $I\in\mathcal{E}_{n'}\lor\mathcal{F}$ found in the previous
lemma for the parameter $\delta$. We have arrived at the following
conclusion:
\begin{quotation}
For every $\varepsilon>0$, for arbitrarily large $n$, a $(1-\varepsilon)$-fraction
of the level-$k$ components $\theta=\mu_{x,k}$ of $\mu$ have associated
to them a sequence of subspaces $V_{1},\ldots,V_{n}$ of which at
least a $c\delta_{0}$-fraction are of dimension $\geq1$, and which
satisfy 
\begin{equation}
\mathbb{P}_{0\leq i\leq n}\left(\theta^{y,i}\mbox{ is }(V_{i},\varepsilon,m)\mbox{-saturated}\right)>1-\varepsilon.\label{eq:32}
\end{equation}

\end{quotation}
If the last equation held for $\mu$ instead of $\theta$ (possibly
for a different sequence of subspaces), we would be in a position
to apply Proposition \ref{prop:maximal-invariant-saturated-subspace}
(4), which would give the second alternative of the present theorem.
This ``bootstrapping'' from the component $\theta$ to $\mu$ is
accomplished as follows. Let us say that a probability measure $\eta\in\mathcal{P}(\mathbb{R}^{d})$
is fragmented at level $k$ if $\nu(D)>0$ for at least two distinct
$D\in\mathcal{D}_{k}$, otherwise it is unfragmented. We again abbreviate
$\sum A=\sum_{a\in A}A$.
\begin{lem}
Given $k$, if $s\in\mathbb{N}$ is large enough, then 
\[
\sum\{p_{i}\,:\,i\in\Lambda^{s}\mbox{ and }\varphi_{i}\mu\mbox{ is unfragmented at level }k\}>1-\varepsilon.
\]
\end{lem}
\begin{proof}
Let $E=\bigcup\partial D$, where the union is over $D\in\mathcal{D}_{k}$
such that $\supp\mu\cap\overline{D}\neq\emptyset$. Then $E$ is contained
in the union of finitely many proper affine subspaces, so for a small
enough $\rho>0$ we will have $\mu(E^{(\rho)})<\varepsilon$. Let
$s$ be large enough that for $i\in\Lambda^{s}$ the measure $\varphi_{i}\mu$
is supported on a set of diameter $<\rho$. This means that if $\varphi_{i}\mu$
is fragmented then it is supported on $E^{(\rho)}$. Since $\mu=\sum_{i\in\Lambda^{s}}p_{i}\cdot\varphi_{i}\mu$,
we conclude that
\begin{eqnarray*}
\sum\{p_{i}\,:\,i\in\Lambda^{s}\mbox{ and }\varphi_{i}\mu\mbox{ is fragmented at level }k\} & \leq & \mu(E^{(s)})\\
 & < & \varepsilon,
\end{eqnarray*}
as required.
\end{proof}
Let $s$ as in the lemma for the $k$ we found previously. Assuming
$\varepsilon<1/2$, by the lemma and our previous discussion we can
find a level-$k$ component $\theta=\mu_{D}$, $D\in\mathcal{D}_{k}$,
of $\mu$, for which \eqref{eq:32} holds and, furthermore, $1-\varepsilon$
of the mass of $\theta$ comes from components $\varphi_{i}\mu$,
$i\in\Lambda^{s}$, supported entirely on $D$. We can now apply Lemma
\ref{lem:saturation-at-level-n} \eqref{enu:saturation-at-level-n-transfers-to-components}
to conclude that there is an $i\in\Lambda^{s}$ such that for arbitrarily
large $j$ there is a $V_{n}\in\sat(\varphi_{i}\mu,\varepsilon',m,j)$,
where $\varepsilon'\rightarrow0$ as $\varepsilon\rightarrow0$. By
Lemma \ref{lem:saturation-at-level-n} the same is true of $\mu$
for the subspace $U_{i}^{-1}V_{j}$ and some $\varepsilon''$ that
also vanishes as $\varepsilon\rightarrow0$. We can now invoke Proposition
\ref{prop:largeness-condition-for-Vtilde}, which completes the proof.

\subsection{Transversality and the dimension of exceptions}

In this section we prove Theorems \ref{thm:main-parametric-Rd} and
\ref{thm:main-parametric-Rd-2} on the dimension of exceptional parameters
for parametric families of self-similar sets and measures. We adopt
the notation from the introduction, so $\varphi_{i,t}(x)=r_{i}(t)U_{i}(t)x+a_{i}(t)$
are contracting similarities for $t$ in a compact connected set $I\subseteq\mathbb{R}^{m}$,
for $i=i_{1}\ldots i_{n}\in\Lambda^{n}$ we define $\varphi_{i,t}=\varphi_{i_{1},t}\circ\ldots\circ\varphi_{i_{n},t}$
and similarly $r_{i}(t)$ and $U_{i}(t)$. Recall that $\Delta_{i,j}(t)=\varphi_{i,t}(0)-\varphi_{j,t}(0)$
and define
\[
\Delta'_{n}(t)=\min_{i\neq j\in\Lambda^{n}}\left\Vert \Delta_{i,j}(t)\right\Vert .
\]
If, as in the introduction, we write $\Delta_{n}(t)$ for the minimum
of $d(\varphi_{i,t},\varphi_{j,t})$ over distinct $i,j\in\Lambda^{n}$,
then we have $\Delta'_{n}\leq\Delta_{n}$ and hence $(\Delta'_{n})^{-1}((-\varepsilon,\varepsilon)^{d})\supseteq(\Delta_{n})^{-1}((-\varepsilon,\varepsilon)^{d})$.
In particular, in order to prove Theorem \ref{thm:main-parametric-Rd},
one may replace the set $E$ there with the set $E'=\bigcap_{\varepsilon>0}E'_{\varepsilon}$,
where 
\[
E'_{\varepsilon}=\bigcup_{N=1}^{\infty}\,\bigcap_{n>N}\left(\bigcup_{i,j\in\Lambda^{n}}(\Delta_{i,j})^{-1}((-\varepsilon^{n},\varepsilon^{n})^{d})\right).
\]
Thus we wish to show that, under suitable hypotheses, $\pdim E'_{\varepsilon}\rightarrow0$
as $\varepsilon\rightarrow0$.

We begin with the proof of Theorem \ref{thm:main-parametric-Rd-2}.
We require an elementary fact whose proof we include for completeness.
\begin{lem}
Let $V\subseteq\mathbb{R}^{m}$ be open and let $F:V\rightarrow\mathbb{R}^{k}$
be a $C^{2}$ map. Suppose that $K\subseteq V$ is compact and that
$\rank DF\geq r$ everywhere in $K$. Then $K\cap F^{-1}((-\delta,\delta)^{k})$
can be covered by at most $C\cdot1/\delta^{m-r}$ balls of radius
$\delta$, where $C$ depends only on the diameter of $K$ and the
magnitude of the first and second partial derivatives of $F$ on $K$.\end{lem}
\begin{proof}
We first reduce to the case that $k=r$. Assume this case is known.
Consider the general case $k\geq r$. For each $r$-tuple of distinct
coordinates $i=(i_{1},\ldots,i_{r})\in\{1,\ldots,k\}^{r}$, let $\pi_{i}:\mathbb{R}^{k}\rightarrow\mathbb{R}^{r}$
denote the projection to these coordinates. Now, if $\rank DF(x)\geq r$
then $\rank D(\pi_{i}\circ F)(x)\geq r$ for some $r$-tuple $i$,
so we can find an open cover $V=\bigcup V_{i}$ indexed by tuples
as above such that $D(\pi_{i}\circ F)$ has rank $r$ everywhere in
$K\cap V_{i}$. Choose compact sets $K_{i}\subseteq K$ such that
$K_{i}\subseteq V_{i}$ and $K=\bigcup K_{i}$. By our assumption,
for each $i$ the set $K_{i}\cap(\pi_{i}\circ F)^{-1}((-\delta,\delta)^{r})$
can be covered by $O(1/\delta^{m-r})$ balls of radius $\delta$.
If $x\in F^{-1}((-\delta,\delta)^{k})$ then certainly $x\in(\pi_{i}\circ F)^{-1}((-\delta,\delta)^{r})$
for every tuple $i$, so the union of these $\binom{k}{r}$ covers
is a cover of $K\cap F^{-1}((-\delta,\delta)^{k})$ containing at
most $\binom{k}{r}O(1/\delta^{m-r})$ balls of radius $\delta$, as
required (note that restricting the function and composing with a
projection can only decrease its $C^{2}$ norm, so the constant does
not get worse).

Thus we may from the start assume that $k=r$ and that $\rank DF=r$
everywhere in $K$. Let $M$ denote the bound on the first and second
derivatives of $F|_{K}$. Applying the constant rank theorem \cite[Theorem 7.8]{Lee2003},
for each $x\in K$ there is a neighborhood $W_{x}\subseteq\mathbb{R}^{m}$
of $x$ and an open set $W'_{x}\subseteq\mathbb{R}^{r}$ such that
$F|_{W_{x}}:W_{x}\rightarrow W'_{x}$ is a diffeomorphism and is $C^{2}$-conjugate
to the projection $\pi_{1,\ldots,r}:\mathbb{R}^{m}\rightarrow\mathbb{R}^{r}$.
The distortion of the conjugating maps is controlled by $M$. Since
for $\pi_{1,\ldots,r}$ the statement is clear, the conclusion follows
for $F|_{W_{x}}$. Finally, the neighborhoods $W_{x}$ contain balls
centered at $x$ with radius again bounded in terms of $M$. Only
$O((\diam K)^{m})$ of these neighborhoods are needed to cover $K$,
and the statement follows.
\end{proof}
Returning to our parametrized family of IFSs, assume that $D\Delta_{i,j}$
has rank at least $r$ at every point in $I$ and every distinct pair
$i,j\in\Lambda^{\mathbb{N}}$. 
\begin{lem}
For large enough $n$ and all $i,j\in\Lambda^{n}$ the rank of $\Delta_{i,j}$
is at least $r$ everywhere in $I$. \end{lem}
\begin{proof}
It is easy to check that the power series for the functions $\Delta_{i,j}$
converge on a common neighborhood of $I$ in $\mathbb{C}^{m}$ (each
function being defined by its complex power series), and since $\Delta_{i_{1}\ldots i_{n},j_{1}\ldots j_{n}}\rightarrow\Delta_{i,j}$
uniformly on this neighborhood, we find that $D_{v}\Delta_{i_{1},\ldots,i_{n},j_{1},\ldots,j_{n}}\rightarrow D_{v}\Delta_{i,j}$
as $n\rightarrow\infty$ for all $v$. The lemma now follows by a
compactness argument.
\end{proof}
Finally, it is again clear that there is a uniform bound $M$ for
the first and second derivatives of all the functions $\Delta_{i,j}$,
$i,j\in\Lambda^{n}$. The proof of Theorem \ref{thm:main-parametric-Rd-2}
is now concluded as follows. For large enough $n$, for each $i,j\in\Lambda^{n}$,
the set
\[
(\Delta_{i,j})^{-1}((-\varepsilon^{n},\varepsilon^{n})^{d})
\]
can be covered $O_{M}(1/\varepsilon^{n(m-r)})$ balls of radius $\varepsilon^{n}$.
Thus the set
\[
E'_{n,\varepsilon}=\bigcup_{i,j\in\Lambda^{n}}(\Delta_{i,j})^{-1}((-\varepsilon^{n},\varepsilon^{n})^{d})
\]
satisfies
\[
N(E'_{n,\varepsilon},\varepsilon^{n})\leq|\Lambda|^{n}\cdot O_{M}(1/\varepsilon^{n(m-r)})
\]
where $N(X,\delta)$ is the $\delta$-covering number of $X$. Thus
for every $\varepsilon$ and $n$, 
\begin{eqnarray*}
N(\bigcap_{k>n}E'_{k,\varepsilon},\varepsilon^{n}) & \leq & N(E'_{n,\varepsilon}.\varepsilon^{n})\leq\\
 & \leq & |\Lambda|^{n}\cdot O_{M}(1/\varepsilon^{n(m-r)}),
\end{eqnarray*}
hence for each $\varepsilon>0$,
\begin{eqnarray*}
\bdim\left(\bigcap_{k>n}E'_{k,\varepsilon}\right) & = & \limsup_{n\rightarrow\infty}\frac{\log N\left(\bigcap_{k>n}E'_{k,\varepsilon},\varepsilon^{n}\right)}{\log(1/\varepsilon^{n})}\\
 & \leq & m-r+\frac{\log|\Lambda|}{\log(1/\varepsilon)},
\end{eqnarray*}
It follows that $\pdim E'_{\varepsilon}\leq m-r+\log|\Lambda|/\log(1/\varepsilon)$,
and this tends to $m-r$ as $\varepsilon\rightarrow0$, as required.

We now turn to the the proof of Theorem \ref{thm:main-parametric-Rd},
which is very similar to the one-parameter case from \cite{Hochman2014}.
Let $|i|$ denote the length of a sequence $i$ and for sequences
$i,j$ let $i\land j$ denote the longest common initial segment of
sequences (which may be $0$). Let 
\[
\mathcal{B}_{m}=\{e_{1},\ldots,e_{m}\}
\]
denote the standard basis of $\mathbb{R}^{m}$ and let $D_{v}$ denote
the directional derivative operator in direction $v$. Thus for $F=(F_{1},\ldots,F_{d}):I\rightarrow\mathbb{R}^{d}$
we have $D_{v}F=(D_{v}F_{1},\ldots,D_{v}F_{d}):I\rightarrow\mathbb{R}^{d}$.
We also write $D$ for the differentiation operator for functions
$\mathbb{R}^{m}\rightarrow\mathbb{R}^{d}$. It will be convenient
for the rest of this section to use the supremum norm on vectors and
matrices.
\begin{defn}
Let $I\subseteq\mathbb{R}^{m}$ be a connected compact set. A family
$\{\Phi_{t}\}_{t\in I}$ of IFSs is transverse of order $k$ if the
associated functions $r_{i}(\cdot)$, $a_{i}(\cdot)$, $U_{i}(\cdot)$
are $(k+1)$-times continuously differentiable in a neighborhood of
$I$, and there is a constant $c>0$ such that for all $n\in\mathbb{N}$
and all $i,j\in\Lambda^{n}$, 
\[
\forall\,t_{0}\in I\quad\exists\,p\in\{0,\ldots,k\}\quad\exists\,v_{1},\ldots,v_{p}\in\mathcal{B}_{m}\qquad\qquad
\]
\[
\qquad\qquad\mbox{such that }\quad\left\Vert (D_{v_{p}}\ldots D_{v_{1}}\Delta_{i,j})(t_{0})\right\Vert >c\cdot|i\land j|^{-p}\cdot r_{i\land j}(t_{0}).
\]

\end{defn}
A real-analytic function defined $F:I\rightarrow\mathbb{R}^{d}$ can
be extended to a complex-analytic function on an open complex neighborhood
of $I$. Such an $F$ is identically $0$ if and only if at some point
$t_{0}\in U$ we have $D_{v_{1}}\ldots D_{v_{n}}F(t_{0})=0$ for every
$n$ and $v_{1}\ldots v_{n}\in\mathcal{B}_{m}$. For $i,j\in\Lambda^{\mathbb{N}}$
the functions $\Delta_{i,j}$ are real analytic if $r_{i},a_{i},U_{i}$
are, because on the common neighborhood of $I$ in which these functions
are analytic, $\Delta_{i,j}$ is given as an absolutely convergent
powers series in these functions. Thus the $\Delta_{i,j}$ extend
to complex-analytic functions on a common neighborhoods of $I$. 

We have the following analog of \cite[Proposition 5.7]{Hochman2014}: 
\begin{prop}
\label{prop:analytic-implies-transverse-Rd}Let $I\subseteq\mathbb{R}^{m}$
be a connected compact set and $\{\Phi_{t}\}_{t\in I}$ a family of
IFSs on $\mathbb{R}^{d}$ whose associated functions $r_{i}(\cdot),a_{i}(\cdot),U_{i}(\cdot)$
are real analytic on $I$. For $i,j\in\Lambda^{\mathbb{N}}$, suppose
that $\Delta_{i,j}\equiv0$ on $I$ if and only if $i=j$. Then \textup{\emph{$\{\Phi_{t}\}$
is transverse of order $k$ for some $k$.}}\end{prop}
\begin{proof}
For $i,j\in\Lambda^{n}$, let $\ell=|i\land j|$ and let $u,v\in\Lambda^{n-\ell}$
denote the sequences obtained by deleting the first $\ell$ symbols
of $i,j$. Define the function $\widetilde{\Delta}_{i,j}$ by 
\[
\widetilde{\Delta}_{i,j}(t)=\Delta_{u,v}(t).
\]
We find that 
\[
\Delta_{i,j}(t)=r_{i\land j}(t)\cdot U_{i\land j}(t)(\widetilde{\Delta}_{i,j}(t)),
\]
Let $n(u)$ denote the number of times that the symbol $u\in\Lambda$
appears in $i\land j$ and let $U^{T}$ be the transpose of $U$.
Then (since $U_{i\land j}^{T}=U_{i\land j}^{-1}$), 
\[
\widetilde{\Delta}_{i,j}(t)=(\prod_{u\in\Lambda}r_{u}(t)^{n(u)})\cdot U_{i\land j}^{T}(t)\Delta_{i,j}(t).
\]
From here the analysis is entirely analogous to the proof of \cite[Proposition 5.7]{Hochman2014},
bounding iterated directional derivatives rather than the higher derivative
$\widetilde{\Delta}_{i,j}^{(p)}$ from the original proof. We omit
the details.
\end{proof}
Our next task is to show that transversality of order $k$ provides
efficient coverings of pre-images $(\Delta_{i,j})^{-1}((-\varepsilon,\varepsilon)^{d})$.
The argument is again very similar to the one-dimensional case but
with some additional technicalities. The key part of the argument
in dimension $1$ was the fact that if $F:[a,b]\rightarrow\mathbb{R}$
satisfies $|F'|>c$, then $F^{-1}(-\rho,\rho)$ is an interval of
length $\leq2\rho/c$. We now generalize this to higher dimensions.

Let $U\subseteq\mathbb{R}^{m-1}$, let $f:U\rightarrow\mathbb{R}$
be a Lipschitz function with Lipschitz constant $c$, and $E=\{(x,f(x))\in\mathbb{R}^{m}\,:\,x\in U\}$
be its graph. Then we say that $E$ is a $c$-Lipschitz graph in $\mathbb{R}^{m}$
with domain $U$. More generally we apply this name to any isometric
image of $E$ in $\mathbb{R}^{m}$.
\begin{lem}
Let $E\subseteq\mathbb{R}^{m}$ be a $c$-Lipschitz graph with domain
$U=B_{r}(x)\subseteq\mathbb{R}^{m-1}$ and let $0<\varepsilon<r$.
Then the $\varepsilon$-neighborhood of $E$ can be covered by $O((r/\varepsilon)^{m-1})$
balls of radius $\varepsilon$ if $c<1$, and by $O((cr/\varepsilon)^{m-1})$
such balls if $c\geq1$. \end{lem}
\begin{proof}
Assume that $c\leq1$. Let $y=(u,f(u))$ be a point in the graph.
Let $y^{\pm}=(u,f(u)\pm\varepsilon/2)$. Then the union $C=C(u)=B_{\varepsilon}(y^{+})\cup B_{\varepsilon}(y^{-})$
contains the cylinder $B_{\varepsilon/2}(u)\times[-3\varepsilon/2,3\varepsilon/2]$.
Since $f$ is $c$-Lipschitz, this implies that $C$ contains the
$\varepsilon$-neighborhood of the graph over $B_{\varepsilon/2}(u)$.
Now cover $B_{r}(x)$ by $O((r/\varepsilon)^{m-1})$ balls $B_{\varepsilon/2}(u_{i})$.
Then $\bigcup C(u_{i})$ is covered by $O((r/\varepsilon)^{m-1})$
$\varepsilon$-balls, and contains the $\varepsilon$-neighborhood
of the graph.

If $c\geq1$, then $C(u)$ contains an $\varepsilon$-neighborhood
of the graph over $B_{\varepsilon/2c}(u)$, and we obtain the desired
bound by covering $B_{r}(x)$ by $O((cr/\varepsilon)^{m-1})$ balls
of radius $\varepsilon/2c$. \end{proof}
\begin{lem}
Let $I\subseteq\mathbb{R}^{m}$ be a compact set, let $0<\delta<1$
and let $I^{(\delta)}$ denote the $\delta$-neighborhood of $I$,
let $F:I^{(\delta)}\rightarrow\mathbb{R}$ be twice continuously differentiable
with $0<c\leq\left\Vert DF\right\Vert \leq M$ and $\left\Vert D^{2}F\right\Vert \leq M$
on $I^{(\delta)}$. We assume $c\leq1$. Then for $0<\rho<\min\{\delta,c/M\}$,
the set $I\cap F^{-1}(-\rho,\rho)$ can be covered by $O_{M,\vol(I^{(\delta)})}((c/\rho)^{m-1})$
balls of radius $\rho/c$. \end{lem}
\begin{proof}
Let $t\in I$. Under our hypotheses, there is a ball $B_{r}(t)\subseteq I^{(\delta)}$,
with radius $r$ less than $\min\{\delta,c/M\}$ and of this order,
such that $\left\Vert DF(t)-DF(t')\right\Vert <\frac{1}{100}c$ for
$t'\in B_{r}(t)$ (here we use the upper bound on the second derivative
of $F$). It is then an easy fact from calculus, essentially, the
implicit function theorem, that the level set $S=F^{-1}(0)\cap B_{r}(t)$
is the graph of a $1$-Lipschitz function and that in the transverse
direction to $S$ the function $F$ grows at a rate proportional to
$c$ as long as we remain in $B_{r}(t)$. Thus, given $\rho>0$, the
set $F^{-1}((-\rho,\rho))\cap B_{r}(t)$ is contained in the $O(\rho/c)$-neighborhood
of the graph of a $1$-Lipschitz function with domain $B_{r}(t)$
for $r=O_{M}(c)$, and by the previous lemma, if $\rho<\min\{\delta,c/M\}$,
it can be covered by $O_{M}((r^{m-1}/(\rho/c)^{m-1}))$ balls of diameter
$\rho/c$. Also, $I$ can be covered by $O(\vol(I^{(\delta)})/r^{m})$
balls $B_{r}(t)$ as above, so it can be covered by $O_{M}((c/\rho)^{m-1})$
balls of diameter $\rho/c$. \end{proof}
\begin{cor}
For $F:I^{(\delta)}\rightarrow\mathbb{R}^{d}$ and under the same
assumptions as above, the same conclusion holds.\end{cor}
\begin{proof}
We can write $I=I_{1}\cup\ldots\cup I_{d}$ such that on each of the
closed sets $I_{i}$ the assumption of the previous lemma holds for
$F_{i}$ (the $i$-th component of $F$) with some degradation of
$c$. Then $I\cap F^{-1}((-\rho,\rho)^{d})\subseteq\bigcup_{i=1}^{d}I_{i}\cap F_{i}^{-1}(-\rho,\rho)$
and the lemma can be applied to each set in the union to obtain the
desired result. \end{proof}
\begin{prop}
Let $I\subseteq\mathbb{R}^{m}$ be a compact set, $I^{(\delta)}$
the $\delta$-neighborhood of $I$, and $F:I^{(\delta)}\rightarrow\mathbb{R}^{d}$
a $(k+1)$-times differentiable function. Suppose that there are constants
$M>0$ and $0<b<1$  such that 
\begin{enumerate}
\item For every $t\in I$, $0\leq p\leq k+1$ and $v_{1},\ldots,v_{p}$
$\in\mathcal{B}_{m}$ we have $|D_{v_{1}}\ldots D_{v_{p}}F(t)|\leq M$
(for $p=0$ this means $|F(t)|\leq M$). 
\item For every $t\in I$ there exist $p\in\{0,\ldots,k\}$ and $v_{1},\ldots,v_{p}\in\mathcal{B}_{m}$
such that $\left\Vert D_{v_{1}}\ldots D_{v_{p}}F(t)\right\Vert >b$
(for $p=0$ this means $F(t)>b$). 
\end{enumerate}
Then there exists $C=C(b,M,\vol I^{(\delta)})\geq1$ such that for
every $0<\rho<b\cdot b^{2^{k}}$, the set 
\[
Z_{\rho}=I\cap F^{-1}((-\rho,\rho)^{d})
\]
can be covered by $C^{k}(b/\rho)^{(m-1)/2^{k}}$ balls of radius $(\rho/b)^{1/2^{k}}$.\end{prop}
\begin{proof}
Take $C$ large enough to play the role of the constant in the bound
in the previous corollary, and large enough that $mC^{k-1}+C\leq C^{k}$
for $k\geq1$. 

We argue by induction on $k$. The case $k=0$ is trivial (because
$|F(t)|>b$ and $\rho<b$ implies $Z_{\rho}=\emptyset$). 

Now fix $k$ and suppose we have proved the claim for $k-1$. First,
note that we can assume without loss of generality that $I\subseteq\overline{Z_{b}}=\{t\in I\,:\,\left\Vert F(t)\right\Vert \leq b\}$,
since clearly $Z_{\rho}\subseteq\overline{Z_{b}}$ and if we did not
have $I\subseteq\overline{Z_{b}}$ we could simply replace $I$ by
$I\cap\overline{Z_{b}}$, to make it hold.

Since $\left\Vert F(t)\right\Vert \leq b$ on $I$, the hypothesis
(2) necessarily holds at each point with $p\geq1$. Thus we can write
$I$ as a union of closed sets $I_{v}$, $v\in\mathcal{B}_{m}$, on
each of which the induction hypothesis holds for one of the functions
$G_{v}=D_{v}F$. 

Fix $v\in\mathcal{B}_{m}$, take $\rho'=\sqrt{b\rho}$. Note that
$0<b<1$ and $0<\rho<b^{2^{k}}$, so $0<\rho'<b^{2^{k-1}}$. Define
\begin{eqnarray*}
I'_{v} & = & I_{v}\cap G_{v}^{-1}((-\rho',\rho')^{d})\\
I''_{v} & = & I_{v}\setminus I'_{v}
\end{eqnarray*}
We cover $Z_{\rho}$ in each of these sets separately.

First, we actually cover the entire set $I'_{v}$. Indeed, by the
induction hypothesis, it can be covered by 
\[
C^{k-1}(\frac{b}{\rho'})^{(m-1)/2^{k-1}}=C^{k-1}(\frac{b}{\rho})^{(m-1)/2^{k}}
\]
balls of radius $(\rho'/b)^{1/2^{k-1}}=(\rho/b)^{1/2^{k}}$. 

On the other hand, on $I''_{v}$ we have $\left\Vert DF\right\Vert \geq\left\Vert G_{v}\right\Vert \geq\rho'$.
By the previous corollary, $Z_{\rho}\cap I''_{v}=I''_{v}\cap F^{-1}((-\rho,\rho)^{d})$
can be covered by 
\[
C(\rho'/\rho)^{m-1}=C(\frac{b}{\rho})^{(m-1)/2}
\]
balls of diameter $\rho/\rho'=\sqrt{\rho/b}$, hence, since $\sqrt{\rho/b}<(\rho/b)^{1/2^{k}}$,
we can cover $Z_{\rho}\cap I''$ by at most this many balls of radius
$(\rho/b)^{1/2^{k}}$.

Taking the union of the covers we have found for $Z_{\rho}\cap I'_{v}$
and $Z_{\rho}\cap I''_{v}$, we obtain a cover of $Z_{\rho}\cap I_{v}$
by $(C^{k-1}+C)(b/\rho)^{(m-1)/2}$ balls of radius $(\rho/b)^{1/2^{k}}$.
Summing over the $m$ elements $v\in\mathcal{B}_{m}$, we have covered
$Z_{\rho}$ by 
\[
m(C^{k-1}+C)(\frac{1}{\rho})^{(m-1)/2})\leq C^{k}(\frac{b}{\rho})^{(m-1)/2^{k}}
\]
balls of radius $(\rho/b)^{1/2^{k}}$ (using our assumption $mC^{k-1}+C\leq C^{k}$).
This is the desired cover.
\end{proof}
Theorem \ref{thm:main-parametric-Rd} now follows from Proposition
\ref{prop:analytic-implies-transverse-Rd} and the next result:
\begin{thm}
\label{thm:transverse-implies-small-exceptions}If $\{\Phi_{t}\}_{t\in I}$
satisfies  transversality of order $k\geq1$ on the compact set $I\subseteq\mathbb{R}^{m}$,
then\textup{\emph{ the set $E$ of ``exceptional'' parameters in
Theorem \ref{thm:description-of-exceptional-params-Rd} has}} packing
(and hence Hausdorff) dimension at most $m-1$. \end{thm}
\begin{proof}
Let $M$ be a uniform bound for $\left\Vert D_{v_{1}}\ldots D_{v_{k+1}}\Delta_{i,j}(t)\right\Vert $
taken over $v_{i}\in\mathcal{B}_{m}$, $t\in I$ and $i,j\in\Lambda^{*}$.
Such $M$ exists from $k$-fold continuous differentiability of $r_{i}(\cdot),a_{i}(\cdot)$
and the fact that $|r_{i}|$ are bounded away from $1$ on $I$. By
transversality there is a constant $c>0$ such that for all $n\in\mathbb{N}$
and all $i,j\in\Lambda^{n}$, 
\[
\forall\,t_{0}\in I\quad\exists\,p\in\{0,\ldots,k\}\quad\exists\,v_{1},\ldots,v_{p}\in\mathcal{B}_{m}\qquad\qquad
\]
\[
\qquad\qquad\mbox{such that }\quad\left\Vert (D_{v_{p}}\ldots D_{v_{1}}\Delta_{i,j})(t_{0})\right\Vert >c\cdot|i\land j|^{-p}\cdot r_{min}^{|i\land j|}(t_{0}),
\]
where 
\[
r_{min}=\min\{r_{i}(t)\,:\,i\in\Lambda\,,\,t\in I\}.
\]
We may assume that $c<1$ and $k\geq2$. In what follows we suppress
the dependence on $k,M,c$ and $I$ in the $O(\cdot)$ notation: $O(\cdot)=O_{k,M,c,|I|}(\cdot)$. 

Fix $n$ and distinct $i,j\in\Lambda^{n}$. Let $b=b_{n}=cn^{-k}r_{min}^{n}$,
so that the hypothesis of the previous proposition is satisfied for
the function $F=\Delta_{i,j}$ and this $b$. Therefore, for all $0<\rho<n^{2^{k}}$,
the set $\{t\in I\,:\,|\Delta_{i,j}|<\rho\}$ can be covered by at
most $O((b/\rho)^{(m-1)/2^{k}})$ balls of radius $(\rho/b)^{1/2^{k}}$
each. 

Now let $\varepsilon>0$ be such that $\rho=\varepsilon^{n}$ satisfies
$\rho<(b_{n})^{2^{k}}=(cn^{-k}r_{min}^{n})^{2^{k}}$ for all $n$
(this holds for all sufficiently small $\varepsilon>0$). Fixing $n$
again, the discussion above applies to $(\Delta_{i,j})^{-1}(-\varepsilon^{n},\varepsilon^{n})$
for every distinct pair $i,j\in\Lambda^{n}$, so ranging over all
such pairs we find that
\[
E'_{\varepsilon,n}=\bigcup_{i,j\in\Lambda^{n}\,,\,i\neq j}(\Delta_{i,j})^{-1}(-\varepsilon^{n},\varepsilon^{n})
\]
can be covered by $O(|\Lambda|^{n}(b_{n}/\varepsilon^{n})^{(m-1)/2^{k}})$
balls of radius $(\varepsilon^{n}/b_{n})^{1/2^{k}}$. Now, 
\[
E\subseteq E'_{\varepsilon}=\bigcup_{N=1}^{\infty}\bigcap_{n>N}E'_{\varepsilon,n}.
\]
By the above, for each $\varepsilon$ and $N$ we have 
\begin{eqnarray*}
\bdim\left(\bigcap_{n>N}E'_{\varepsilon,n}\right) & \leq & \lim_{n\rightarrow\infty}\frac{\log\left(|\Lambda|^{n}(b_{n}/\varepsilon^{n})^{(m-1)/2^{k}}\right)}{\log\left((b_{n}/\varepsilon^{n})^{1/2^{k}}\right)}\\
 & = & O(\frac{\log(|\Lambda|(r_{min}/\varepsilon)^{(m-1)/2^{k}})}{\log(r_{min}/\varepsilon){}^{1/2^{k}}}).
\end{eqnarray*}
The last expression tends to $m-1$ as $\varepsilon\rightarrow0$,
uniformly in $N$. Thus the same is true of $E'_{\varepsilon}$, and
$E\subseteq E'_{\varepsilon}$ for all $\varepsilon$, so $E$ has
packing (and Hausdorff) dimension $m-1$.
\end{proof}

\subsection{\label{sub:Applications-Rd-proofs}Applications and further comments}
\begin{proof}
[Proof of Theorem \ref{thm:generic-IFS}] Fix $\Lambda$. For $i,j\in\Lambda^{\mathbb{N}}$,
given an IFS $\Phi=\{(\varphi_{i}\}_{i\in\Lambda}=(r_{i}U_{i}+a_{i})_{i\in\Lambda}$,
evidently 
\[
\Delta_{i,j}(\Phi)=\sum_{n=0}^{\infty}\left(r_{i_{1}\ldots i_{n-1}}U_{i_{1}\ldots i_{n-1}}a_{i_{n}}-r_{j_{1}\ldots j_{n-1}}U_{j_{1}\ldots j_{n-1}}a_{j_{n}}\right).
\]
As a function of $(r_{u},U_{u},a_{u})_{u\in\Lambda}\in(\mathbb{R}^{+}\times\mathbb{R}^{d^{2}}\times\mathbb{R}^{d})^{\Lambda}$
this is clearly a non-constant expression. The parametrization is
trivially real-analytic, and the conclusion follows from Theorem \eqref{thm:main-parametric-Rd}.
\end{proof}

\begin{proof}
[Proof of Theorem \ref{thm:generic-translation-IFS}]Fix $\{U_{i}\}_{i\in\Lambda}\in G_{0}^{\Lambda}$
and $\{r_{i}\}_{i\in\Lambda}\in(0,1/2)^{\Lambda}$. Given distinct
$i,j\in\Lambda^{\mathbb{N}}$ let $k=k(i,j)$ be the first index where
they differ. For $a=(a_{u})_{u\in\Lambda}\in(\mathbb{R}^{d})^{\Lambda}$,
let $\Phi_{a}=\{r_{u}U_{u}+a_{u}\}_{u\in\Lambda}$, so 
\[
\Delta_{i,j}(a)=\sum_{n\geq k(i,j)}\left(r_{i_{1}\ldots i_{n-1}}U_{i_{1}\ldots i_{n-1}}a_{i_{n}}-r_{j_{1}\ldots j_{n-1}}U_{j_{1}\ldots j_{n-1}}a_{j_{n}}\right).
\]
This is linear in the $a$ variables. Differentiating by the coordinates
in $a_{i_{k}}=(a_{i_{k}}^{1},\ldots,a_{i_{k}}^{d})$, we obtain a
derivative matrix of the form
\begin{equation}
\left(\frac{\partial\Delta_{i,j}}{\partial a_{i_{k}}}\right)=r_{i_{1}\ldots i_{k-1}}U_{i_{1}\ldots i_{k-1}}+\sum_{n\in I}r_{i_{1}\ldots i_{n}}U_{i_{1}\ldots i_{n}}-\sum_{n\in J}r_{j_{1}\ldots j_{n}}U_{j_{1}\ldots j_{n},}\label{eq:partial-derivative-1}
\end{equation}
where $I=\{n>k\,:\,i_{n}=i_{k}\}$ and $J=\{n>k\,:\,j_{n}=i_{k}\}$.
Similarly, setting $I'=\{n>k\,:\,i_{n}=j_{k}\}$ and $J'=\{n>k\,:\,j_{n}=j_{k}\}$
and differentiating $\Delta_{i,j}$ by the $a_{j_{k}}$ variable (and
using $r_{i_{1}\ldots i_{k-1}}=r_{j_{1}\ldots j_{k-1}}$ and $U_{i_{1}\ldots i_{k-1}}=U_{j_{1}\ldots j_{k-1}}$),
\begin{equation}
\left(\frac{\partial\Delta_{i,j}}{\partial a_{j_{k}}}\right)=r_{j_{1}\ldots j_{k-1}}U_{j_{1}\ldots j_{k-1}}+\sum_{n\in I'}r_{i_{1}\ldots i_{n}}U_{i_{1}\ldots i_{n}}-\sum_{n\in J'}r_{j_{1}\ldots j_{n}}U_{j_{1}\ldots j_{n}}.\label{eq:partial-derivative-2}
\end{equation}
In order for these matrices to be invertible, it is enough that on
the right hand sides of equations \eqref{eq:partial-derivative-1}
and \eqref{eq:partial-derivative-2}, the norm of the sum of the last
two terms is less than the norm of the first term. Let 
\begin{eqnarray*}
R & = & \sum_{n\in I}r_{i_{1}\ldots i_{n}}+\sum_{n\in J}r_{j_{1}\ldots j_{n}}\\
R' & = & \sum_{n\in I'}r_{i_{1}\ldots i_{n}}+\sum_{n\in J'}r_{j_{1}\ldots j_{n}}.
\end{eqnarray*}
These are upper bounds for the norms in question. We have
\begin{eqnarray*}
R+R' & = & \left(\sum_{n\in I}r_{i_{1}\ldots i_{n}}+\sum_{n\in I'}r_{i_{1}\ldots i_{n}}\right)+\left(\sum_{n\in J}r_{j_{1}\ldots j_{n}}+\sum_{n\in J'}r_{j_{1}\ldots j_{n}}\right)\\
 & \leq & r_{i_{1}\ldots i_{k-1}}\prod_{n\in I\cap I'}(r_{i_{n}}+r_{i_{n}})+r_{j_{1}\ldots j_{k-1}}\prod_{n\in J\cap J'}(r_{j_{n}}+r_{j_{n}})\\
 & < & 2r_{i_{1}\ldots i_{k-1}}.
\end{eqnarray*}
(In the first inequality we used $r_{i}<1$. In the second we used
the fact that if $n\in I\cap I'$ then $i_{n}\neq j_{n}$ and hence
$r_{i_{n}}+r_{j_{n}}<1$, and similarly for $n\in J\cap J'$, and
that $r_{i_{1}\ldots i_{k-1}}=r_{j_{1}\ldots j_{k-1}}$ by choice
of $k$). Now, $R+R'<r_{i_{1}\ldots r_{k-1}}$ implies that either
$R<r_{i_{1}\ldots i_{k-1}}$ or $R'<r_{i_{1}\ldots i_{k-1}}$. In
the first case, the first term in \eqref{eq:partial-derivative-1}
is a similarity with contraction $r_{i_{1}\ldots i_{k-1}}$, and the
latter two terms together give a matrix whose norm is at most $R<r_{i_{1}\ldots i_{k-1}}$.
Hence the sum is invertible, and $\rank D\Delta_{i,j}\geq d$. The
same argument applies to \eqref{eq:partial-derivative-2} if $R'<r_{i_{1}\ldots i_{k-1}}$.
The conclusion now follows from Theorem \ref{thm:main-parametric-Rd-2}.
\end{proof}

\begin{proof}
[Proof of Theorem \ref{thm:generic-projection-IFS}] Let $(\varphi_{i})_{i\in\Lambda}$
be given. For $i\in\Lambda^{\mathbb{N}}$ write $\varphi_{i}=\lim_{n\rightarrow\infty}\varphi_{i_{1}\ldots i_{n}}(0)$.
Given distinct $i,j\in\Lambda^{\mathbb{N}}$ and $\pi\in\Pi_{d,k}$,
evidently
\[
\Delta_{i,j}(\pi)=\pi(\varphi_{i})-\pi(\varphi_{j})=\pi(\varphi_{i}-\varphi_{j})
\]
Now, it is easy to verify that for a fixed $0\neq v\in\mathbb{R}^{d}$
the map $\pi\mapsto\pi(v)$, $\Pi_{d,k}\rightarrow\mathbb{R}^{k}$,
has rank $k$ at every point. Taking $v=\varphi_{i}-\varphi_{j}$
this shows that $\Delta_{i,j}$ has rank $k$ at every point. An application
of Theorem \eqref{thm:main-parametric-Rd-2} completes the proof.
\end{proof}

\begin{proof}
[Proof of Theorem \ref{thm:nonuniform-BC}] Writing $\Delta_{i,j}(\beta,\gamma)$
explicitly and noting that it is not constant and real-analytic, Theorem
\ref{thm:nonuniform-BC} is immediate from Theorem \ref{thm:main-parametric-Rd}
(since the IFS in on the line, irreducibility is a non-issue).
\end{proof}

\begin{proof}
[Proof of Theorem \ref{thm:fat-Sierpinski}] We would again like to
apply Theorem \ref{thm:fat-Sierpinski}. Analyticity and non-triviality
of $\Delta_{i,j}$ is again a simple matter, but the usual presentation
of the fat Sierpinski gaskets uses an IFS consisting of homotheties,
which act reducibly. However, the attractor of the fat Sierpinski
gaskets are invariant under rotation by $2\pi/3$ about their center
of mass, and hence they can be presented also as attractors of an
IFS $x\mapsto\lambda U_{i}x+a_{i}$ where $a_{i}$ are the vertices
of a triangle in $\mathbb{R}^{2}$ and the $U_{i}$ are rotations
by $2\pi/3$. Unlike the usual presentation this IFS is irreducible.
Theorem \ref{thm:fat-Sierpinski} now does the job.
\end{proof}
The argument in the last proof relied heavily on the possibility of
presenting the attractor using an irreducible IFS. This is not always
possible. For instance, if we take the fat Sierpinski gasket with
the usual homothetic presentation, and augment it with an additional
homothety, then the symmetry breaks down and there is no irreducible
presentation. In this case Theorem \ref{thm:fat-Sierpinski} no longer
gives information about the set of exceptional parameters, because
the set of reducible parameters is large. Some additional argument
is needed in this case. 

Finally, the proof of Corollary \ref{cor:irr-algebraic-IFS} is based
on the classical fact that polynomials of bounded height in a fixed
set of algebraic numbers either vanishes or is exponentially large
in the degree of the polynomial. For completeness we include a proof,
noting that the version in \cite[Lemma 5.10]{Hochman2014} erroneously
omitted the height assumption:
\begin{lem}
Let $\mathcal{A}\subseteq\mathbb{R}$ be a finite set of algebraic
numbers over $\mathbb{Q}$. If $x$ is a polynomial expression in
the elements of $\mathcal{A}$ with coefficients of magnitude at most
$h$, then either $x=0$ or $|x|>s^{n}$.\end{lem}
\begin{proof}
Let $\mathcal{A}=\{a_{1},\ldots,a_{k}\}$. Let $f(x_{1},..,x_{k})$
be an integer polynomial of degree $n$ and coefficients bounded by
$h$ in absolute value. Assuming $x=f(a_{1},...,a_{k})$ is not zero,
it suffices to show that $|x|>c^{n}/h^{u}$ for some $c,u>0$ depending
only on $\mathcal{A}$. 

Let $\mathbb{F}=\mathbb{Q}(a_{1},\ldots,a_{k})$ be the field over
$\mathbb{Q}$ generated by $\{a_{i}\}$.

We may assume that $a_{i}$ are algebraic integers. This is because
we can choose positive integers $p_{1},...,p_{k}$ such that $b_{i}=p_{i}\cdot a_{i}$
is an algebraic integer. Let $p=p_{1}\cdot\ldots\cdot p_{k}$ (note
that this depends only on the $a_{i}$). Then 
\[
p^{n}\cdot f(a_{1},...,a_{k})=g(b_{1},...,b_{k}),
\]
and $g$ is an integer polynomial of degree $n$ with coefficients
bounded by $h\cdot p^{n}$. So if we have $c=c(b_{1},...,b_{k})>0$
such that $g(b_{1},...,b_{k})>c^{n}/(hp^{n})^{u}$, then $f(a_{1},...,a_{k})>c^{n}/(h^{u}\cdot p^{(u+1)n})$,
which is what we wanted (using the constant $c/p^{u+1}$ instead of
$c$).

Assuming now that $a_{i}$ are algebraic integers, let $\mathbb{F}'$
be the normal closure of $\mathbb{F}=\mathbb{Q}(a_{1},...,a_{k})$
and $\Gamma=\Gal(\mathbb{F}'/\mathbb{Q})$, so the fixed field of
$\Gamma$ is $\mathbb{Q}$. Note that $\mathbb{F}'$, hence $\Gamma$,
depends only on the $a_{i}$, and $\Gamma$ is finite.

Now we do the usual thing: if $f(x_{1},...,a_{k})$ is not zero then
also $\prod_{s\in\Gamma}s(f(x))$ is non-zero, but it is both an algebraic
integer and rational, so its absolute value is at least 1. Hence 
\[
1\le\prod_{s\in\Gamma}|f(sx)|=|f(x)|\cdot\prod_{s\in\Gamma\setminus\{\id\}}|f(sx)|.
\]
The last product has $|\Gamma|-1$ factors $|f(sx)|$, each of size
at most $h\cdot\max\{|\Gamma\mbox{-conjugates of}\ a_{i}|\}^{n}$.
Dividing gives the bound that we want.
\end{proof}
{

}

\bibliographystyle{plain}
\bibliography{bib}

\bigskip{}

\lyxaddress{Email: mhochman@math.huji.ac.il \\
Address: Einstein Institute of Mathematics, Givat Ram, Jerusalem 91904,
Israel}
\end{document}